\documentclass[reqno]{amsart}
\usepackage[active]{srcltx}
\usepackage{amssymb}
\usepackage[active]{srcltx}
\usepackage{amsmath}
\usepackage{times}
\usepackage[T1]{fontenc}
\usepackage{textcomp}
\usepackage{txfonts}
\usepackage{tikz}
\usetikzlibrary{arrows,decorations.markings}
\usetikzlibrary{shapes}
\usepackage{mathrsfs}
\usepackage{multicol}
\usepackage{hyperref}
\usepackage{amsthm}
\usepackage{cases}
\usepackage{leftidx}

\usepackage{titlesec}

\titleformat{\section}{\Large\bfseries}{\thesection.}{0.5em}{}
\titleformat{\subsection}{\large\itshape}{\thesubsection.}{0.5em}{}

\usepackage{pifont}

\usepackage{bm}

\DeclareMathAlphabet{\mathpzc}{OT1}{pzc}{m}{it}

{\catcode`\|=\active
  \gdef\set#1{\mathinner{\lbrace\,{\mathcode`\|"8000%
                                   \let|\midvert #1}\,\rbrace}}
}
\def\midvert{\egroup\mid\bgroup}

\usepackage{aliascnt}
\def\NewTheorem#1{%
  \newaliascnt{#1}{lemma}
  \newtheorem{#1}[#1]{#1}
  \aliascntresetthe{#1}
  \expandafter\def\csname #1autorefname\endcsname{#1}
}

\newcounter{main}
\theoremstyle{plain}

\newtheorem{THEOREM}[main]{Theorem}

\numberwithin{equation}{section}

\swapnumbers

\numberwithin{lemma}{section}
\NewTheorem{Assumption}
\NewTheorem{Proposition}
\NewTheorem{Theorem}
\NewTheorem{Corollary}
\NewTheorem{Conjecture}
\NewTheorem{Lemma}
\theoremstyle{definition}
\NewTheorem{Definition}
\NewTheorem{Remark}
\NewTheorem{Remarks}
\NewTheorem{Example}

\def\@thm#1#2#3{%
  \ifhmode\unskip\unskip\par\fi
  \normalfont
  \trivlist
  \let\thmheadnl\relax
  \let\thm@swap\@gobble
  \let\thm@indent\noindent
  \thm@headfont{\SMStheoremfont}
  \thm@notefont{\fontseries\mddefault\upshape}%
  \thm@headpunct{.}
  \thm@headsep 5\p@ plus\p@ minus\p@\relax
  \thm@space@setup
  #1
  \@topsep \thm@preskip               
  \@topsepadd \thm@postskip           
  \def\@tempa{#2}\ifx\@empty\@tempa
    \def\@tempa{\@oparg{\@begintheorem{#3}{}}[]}%
  \else
    \refstepcounter{#2}%
    \def\@tempa{\@oparg{\@begintheorem{#3}{\csname the#2\endcsname}}[]}%
  \fi
  \@tempa
}

\renewcommand{\P}{\mathscr{P}^\Lambda}
\renewcommand{\H}{\widehat H}
\renewcommand{\O}{{\mathscr O}}
\renewcommand{\L}{{\mathscr L}}
\newcommand{\Z}{\mathbb{Z}}
\newcommand{\N}{\mathbb{N}}
\newcommand{\R}{\mathscr{R}^\Lambda_n}

\newcommand{\K}{\mathcal{K}}

\newcommand{\Sym}{\mathfrak{S}}
\newcommand{\B}{\mathscr B_n(\delta)}

\newcommand{\Bx}{\mathscr B^{\F}_n(x)}
\newcommand{\BOx}{\mathscr B^{\mathscr O}_n(x)}
\newcommand{\F}{\mathbb{F}}
\newcommand{\Tud}{\mathscr T^{ud}}

\def\bi{\mathbf{i}}
\def\bj{\mathbf{j}}
\def\bk{\mathbf{k}}
\def\bl{\mathbf{l}}
\def\bm{\mathbf{m}}

\def\t{\mathsf{t}}
\def\s{\mathsf{s}}
\def\u{\mathsf{u}}
\def\v{\mathsf{v}}
\def\x{\mathsf{x}}
\def\y{\mathsf{y}}
\def\w{\mathsf{w}}
\def\m{\mathfrak{m}}

\def\l{\ell}

\newcommand\f[1]{\frac{f_{#1 #1}}{\gamma_{#1}}}

\usepackage{tikz}  

\colorlet{darkgreen}{green!50!black}
\tikzset{dots/.style={ultra thick,loosely dotted},
         belt/.style={draw,blue,thick,fill=blue!50},
         greendot/.style={fill,circle,color=darkgreen,inner sep=1.5pt,outer sep=0}
}
\newenvironment{braid}{
  \begin{tikzpicture}[baseline=6mm,blue,line width=1pt, xscale = 0.35, yscale=0.4,
                      draw/.append style={rounded corners},
                      every node/.append style={font=\fontsize{5}{5}\selectfont}]%
  }{\end{tikzpicture}
}

\def\Grid(#1,#2){
  \draw[very thin,gray,step=2mm] (0,0)grid(#1,#2);
  \draw[very thin,darkgreen,step=10mm] (0,0)grid(#1,#2);
}

\newcommand\Tableau[2][\relax]{
  \begin{tikzpicture}[scale=0.5,draw/.append style={thick,black},baseline=-3mm]
    \ifx\relax#1\relax%
    \else 
      \foreach\box in {#1} {
        \ifx\box\relax\else
          \filldraw[blue!10]\box+(-.5,-.5) rectangle ++(.5,.5);
        \fi
      }
    \fi;
    \newcount\row\newcount\col
    \row=0
    \foreach \Row in {#2} {
       \col=1
       \foreach \k in \Row {
          \draw (\the\col,\the\row) +(-.5,-.5) rectangle ++(.5,.5);
          \ifnum\k<0
            \draw(\the\col,\the\row) node[fill=gray!20]{-\k};
          \else \draw (\the\col,\the\row) node{\k};
          \fi
          \global\advance\col by 1
       }
       \global\advance\row by -1
    }
  \end{tikzpicture}
}

\newcommand\YoungDiagram[2][\relax]{
  \begin{tikzpicture}[scale=0.5,draw/.append style={thick,black},baseline=-1mm]
    \ifx\relax#1\relax%
    \else 
    \foreach\box in {#1} {
      \filldraw[blue!10]\box rectangle ++(1,1);
    }
    \fi
    \newcount\row
    \row=0
    \foreach \col in {#2} {
       \draw(1,\the\row)grid ++(\col,1);
       \global\advance\row by -1
    }
  \end{tikzpicture}
}

\DeclareMathOperator{\res}{res} 
 \DeclareMathOperator{\node}{node}
 \DeclareMathOperator{\cont}{cont}
 
\DeclareMathOperator{\Std}{Std} 
\DeclareMathOperator{\Shape}{Shape}

\newcommand{\map}[2]{\,{:}\,#1\!\longrightarrow\!#2}
\newcommand{\G}[1]{\mathscr G_{#1}(\delta)}

\usepackage[enableskew,vcentermath]{youngtab}
\def\tab(#1){\mbox{\tiny$\young(#1)$}\,}
\def\ydiag(#1){\mbox{\tiny$\yng(#1)$}\,}

\newpage

\begin{document}    
\bibliographystyle{andrew}

\title{A KLR Grading of the Brauer Algebras}


\author{Ge Li}

\address{School of Mathematics and Statistics\\
		University of Sydney\\
		Sydney, NSW 2006}
		
\email{geli@maths.usyd.edu.au}

\begin{abstract}
We construct a naturally $\Z$-graded algebra $\G{n}$ over $R$ with KLR-like relations and give an explicit isomorphism between $\G{n}$ and $\B$, the Brauer algebras over $R$, when $R$ is a field of characteristic $0$. This isomorphism allows us to exhibit a non-trivial $\Z$-grading on the Brauer algebras over a field of characteristic $0$. As a byproduct of the proof, we also construct an explicit homogeneous cellular basis for $\G{n}$.
\end{abstract}

%

\maketitle          
\pagenumbering{arabic}

\section{Introduction}

Richard Brauer~\cite{BrauerAlg} introduced a class of finite dimensional algebras $\B$ over a field $R$, which are called Brauer algebras, in order to study the $n$-th tensor power of the defining representations of orthogonal groups and symplectic groups. It is well known that the symmetric group algebras $R\Sym_n$ is a subalgebra of $\B$.

Khovanov and Lauda~\cite{KhovLaud:diagI,KhovLaud:diagII} and Rouquier~\cite{Rouq:2KM} have introduced a remarkable new family of algebras $\mathscr R_n$, the quiver Hecke algebras, for each oriented quiver, and they showed that they categorify the positive part of the enveloping algebras of the corresponding quantum groups. The algebras $\mathscr R_n$ are naturally $\Z$-graded. Brundan and Kleshchev \cite{BK:GradedKL} proved that every degenerate and non-degenerate cyclotomic Hecke algebra $H^\Lambda_n$ of type $G(r,1,n)$ over a field is isomorphic to a cyclotomic quiver Hecke algebra $\R$ of type $A$ by constructing an explicit isomorphisms between these two algebras. Hu and Mathas~\cite{HuMathas:SemiQuiver} gave another proof of Brundan and Kleshchev's result using seminormal forms. Moreover, Hu and Mathas~\cite{HuMathas:GradedCellular} defined a homogeneous basis of the cyclotomic quiver algebras $\R$ which showed that $H^\Lambda_n$ is a graded cellular algebra.

Because $R\Sym_n$ is a special case of $H_n^\Lambda$, all above results hold in $R\Sym_n$. It is natural to ask the question that whether Brauer algebras $\B$ are graded cellular algebras. Ehrig and Stroppel~\cite{ES:Graded} proved this result, but they were unable to give a presentation of the graded Brauer algebras similar to the KLR presentation of cyclotomic quiver Hecke algebras.

Let $R$ be a field of characteristic $0$. The main purpose of this paper is to construct a $\Z$-graded algebra over $R$ with a parameter $\delta \in R$ analogues to the cyclotomic quiver Hecke algebras and to prove that this algebra is isomorphic to $\B$. In Section~\ref{sec:algebra} we define a $\Z$-graded algebra, $\G{n}$, by generators and relations. It is generated by elements
\begin{equation} \label{generator:1}
G_n(\delta) = \set{e(\bi)|\bi \in P^n} \cup \set{y_k| 1 \leq k \leq n} \cup \set{\psi_k| 1\leq k \leq n-1} \cup \set{\epsilon_k| 1\leq k \leq n-1},
\end{equation}
and the relations are similar to the KLR relations for the cyclotomic quiver Hecke algebras of type $A$. In Section~\ref{sec:ir} we construct a set of homogeneous elements in $\G{n}$
\begin{equation} \label{basis:cellular}
\set{\psi_{\s\t} | (\lambda,f) \in \widehat B_n, \s,\t \in \Tud_n(\lambda)},
\end{equation}
with $\deg \psi_{\s\t} = \deg \s + \deg \t$, which are the Brauer-algebra-analogue of the graded cellular basis of the cyclotomic quiver Hecke algebras~\cite{HuMathas:GradedCellular}. In Section~\ref{sec:span} we prove:
\begin{THEOREM} \label{main:A}
The algebra $\G{n}$ is spanned by the elements $\set{\psi_{\s\t} | (\lambda,f) \in \widehat B_n, \s,\t \in \Tud_n(\lambda)}$.
\end{THEOREM}

We prove Theorem~\ref{main:A} via showing that
\begin{enumerate}
\item $1_R$ is a linear combination of~\eqref{basis:cellular} (cf.~\autoref{I:end}).

\item For any $(\lambda,f) \in \widehat B_n$, $\s,\t \in \Tud_n(\lambda)$ and $a \in \G{n}$, we have
\begin{equation} \label{cellular:property}
\psi_{\s\t} a = \sum_{\v \in \Tud_n(\lambda)} c_{\v} \psi_{\s\v} + \sum_{\substack{(\mu,\l) > (\lambda,f) \\ \u,\v \in \Tud_n(\mu)}} c_{\u\v} \psi_{\u\v},
\end{equation}
where $c_\v, c_{\u\v} \in R$ and $>$ is a total ordering on $\widehat B_n$ (cf.~\autoref{basis:main}).
\end{enumerate}

Theorem~\ref{main:A} shows that $\dim \G{n} \leq (2n-1)!! = \dim \B$ and as a byproduct,~\eqref{cellular:property} shows that~\eqref{basis:cellular} has a cellular-like property.

In Section~\ref{sec:gene} we construct a new set of generators of the Brauer algebra $\B$
\begin{equation} \label{generator:2}
\set{e(\bi)|\bi \in I^n} \cup \set{y_k| 1 \leq k \leq n} \cup \set{\psi_k| 1\leq k \leq n-1} \cup \set{\epsilon_k| 1\leq k \leq n-1},
\end{equation}
and in Section~\ref{sec:grade} we prove that the map $\G{n} \longrightarrow \B$ given by sending the generators in (\ref{generator:1}) to those in (\ref{generator:2}) is a surjective algebra homomorphism. To show that the generators in (\ref{generator:2}) satisfy the relations of $\G{n}$ we make extensive use of seminormal forms of $\B$, following the Hu-Mathas~\cite{HuMathas:SemiQuiver} approach in type A. In turn this relies Naz~\cite{Nazarov:brauer}, or Rui-Si~\cite{RuiSi:BrauerDet}. By construction our map $\G{n}\longrightarrow\B$ is surjective and it is injective by Theorem~\ref{main:A}, we obtain the main result of this paper:

\begin{THEOREM} \label{main:B}
Suppose $R$ is a field of characteristic $p = 0$ and $\delta \in R$. Then $\B \cong \G{n}$.
\end{THEOREM}

By given $\G{n} \cong \B$, (\ref{basis:cellular}) is a basis of $\G{n}$ because it spans $\G{n}$ by Theorem~\ref{main:A} and it has the right number of elements($ = (2n-1)!!$). Because we proved the cellularity of (\ref{basis:cellular}) by~\eqref{cellular:property}, it is a graded cellular basis of $\G{n}$:

\begin{THEOREM} \label{main:C}
The algebra $\G{n}$ is a graded cellular algebra with a graded cellular basis~\eqref{basis:cellular}.
\end{THEOREM}

Because $\G{n} \cong \B$ and $\R \cong R\Sym_n$, similar to the Brauer algebras, by removing the elements $\epsilon_k$ for $1 \leq k \leq n$, the quotient of $\G{n}$ is isomorphic to $\R$ with weight $\Lambda_k$ for any $k \in \Z$.


Finally we remark that the strategy that we use in this paper can be extended to the Brauer algebra over fields of positive characteristic, the degenerate cyclotomic Nazarov-Wenzl algebras and partition algebras. As this paper is already long enough, the details will appear in subsequent papers.


\section{Preliminaries}


\subsection{The symmetric groups and Brauer algebras} \label{sec:basis:def}

Let $\Sym_n$ be the symmetric group acting on the integers $\{1,2,\ldots,n\}$. For $1 \leq k \leq n-1$, define $s_k = (k,k+1)$ as the elementary transpositions in $\Sym_n$. Hence $\Sym_n$ is generated by
$$
\set{s_k | 1 \leq k \leq n-1}
$$
subject to the relations:
\begin{align*}
& s_k^2 = 1, && \text{for $1 \leq k \leq n-1$,}\\
& s_k s_r = s_r s_k, && \text{for $1 \leq k,r \leq n-1$ and $|k - r| > 1$,}\\
& s_k s_{k+1} s_k = s_{k+1} s_k s_{k+1}, && \text{for $1 \leq k \leq n-2$.}
\end{align*}

An expression $w = s_{i_1} s_{i_2} \ldots s_{i_m}$ for $w$ in terms of elementary transpositions is \textit{reduced} if $w$ cannot be expressed as a proper sub-expression of $s_{i_1} s_{i_2} \ldots s_{i_m}$. For example, $w = s_2 s_3 s_2$ is a reduced expression and $w = s_2 s_3 s_2 s_3$ is not a reduced expression, because $w = s_2 s_3 s_2 s_3 = s_3 s_2 s_3 s_3 = s_3 s_2$. Notice that generally there are more than one reduced expressions for an element of $\Sym_n$. For example, $s_2 s_3 s_2 = s_3 s_2 s_3$ and both of expressions are reduced.

Let $R$ be a commutative ring with identity $1$ and $\delta \in R$. The Brauer algebra $\B$ is a unital associative $R$-algebra with generators
$$
\{s_1,s_2,\ldots,s_{n-1}\} \cup \{e_1,e_2,\ldots,e_{n-1}\},
$$
associated with relations
\begin{enumerate}
\item (Inverses) $s_k^2 = 1$.

\item (Essential idempotent relation) $e_k^2 = \delta e_k$.

\item (Braid relations) $s_k s_{k+1} s_k = s_{k+1} s_k s_{k+1}$ and $s_k s_r = s_r s_k$ if $|k - r| > 1$.

\item (Commutation relations) $s_k e_l = e_l s_k$ and $e_k e_r = e_r e_k$ if $|k - r| > 1$.

\item (Tangle relations) $e_k e_{k+1} e_k = e_k$, $e_{k+1} e_k e_{k+1} = e_{k+1}$, $s_k e_{k+1} e_k = s_{k+1} e_k$ and $e_k e_{k+1} s_k = e_k s_{k+1}$.

\item (Untwisting relations) $s_k e_k = e_k s_k = e_k$.
\end{enumerate}

The symmetric group algebra $R\Sym_n$ can be considered as a subalgebra of the Brauer algebra $\B$ for any $\delta$. The Brauer algebra $\B$ (cf.~\cite{BrauerAlg},~\cite{Wenzl:ssbrauer}) has $R$-basis consisting of Brauer diagrams $D$, which consist of two rows of $n$ dots, labelled by $\{1,2,\ldots, n\}$, with each dot joined to one other dot. See the following diagram as an example:
$$
D = 
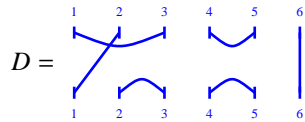
\begin{tikzpicture} [baseline=3mm,blue,line width=1pt, xscale = 0.6, yscale=0.8,
                      draw/.append style={rounded corners},
                      every node/.append style={font=\fontsize{5}{5}\selectfont}]
\draw (0, -0.1) node[below] {$1$} -- (0,0.1);
\draw (1, -0.1) node[below] {$2$} -- (1,0.1);
\draw (2, -0.1) node[below] {$3$} -- (2,0.1);
\draw (3, -0.1) node[below] {$4$} -- (3,0.1);
\draw (4, -0.1) node[below] {$5$} -- (4,0.1);
\draw (5, -0.1) node[below] {$6$} -- (5,0.1);

\draw (0, 0.9) -- (0,1.1) node[above] {$1$};
\draw (1, 0.9) -- (1,1.1) node[above] {$2$};
\draw (2, 0.9) -- (2,1.1) node[above] {$3$};
\draw (3, 0.9) -- (3,1.1) node[above] {$4$};
\draw (4, 0.9) -- (4,1.1) node[above] {$5$};
\draw (5, 0.9) -- (5,1.1) node[above] {$6$};

\draw (0,1) .. controls (1,0.7) .. (2,1);
\draw (3,1) .. controls (3.5,0.7) .. (4,1);
\draw (1,0) .. controls (1.5,0.3) .. (2,0);
\draw (3,0) .. controls (3.5,0.3) .. (4,0);
\draw (0,0) -- (1,1);
\draw (5,0) -- (5,1);
\end{tikzpicture}
$$

Two diagrams $D_1$ and $D_2$ can be composed to get $D_1 \circ D_2$ by placing $D_1$ above $D_2$ and joining corresponding points and deleting all the interior loops. The multiplication of $\B$ is defined by
$$
D_1 {\cdot} D_2 = \delta^{n(D_1, D_2)} D_1 \circ D_2,
$$
where $n(D_1, D_2)$ is the number of deleted loops. For example:
$$
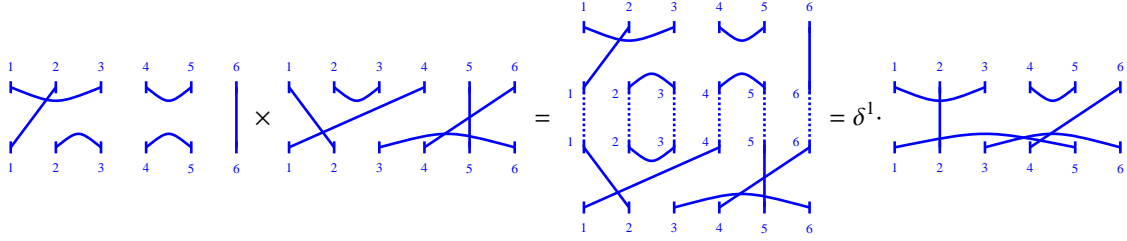
\begin{tikzpicture} [baseline=3mm,blue,line width=1pt, xscale = 0.6, yscale=0.8,
                      draw/.append style={rounded corners},
                      every node/.append style={font=\fontsize{5}{5}\selectfont}]
\draw (0, -0.1) node[below] {$1$} -- (0,0.1);
\draw (1, -0.1) node[below] {$2$} -- (1,0.1);
\draw (2, -0.1) node[below] {$3$} -- (2,0.1);
\draw (3, -0.1) node[below] {$4$} -- (3,0.1);
\draw (4, -0.1) node[below] {$5$} -- (4,0.1);
\draw (5, -0.1) node[below] {$6$} -- (5,0.1);

\draw (0, 0.9) -- (0,1.1) node[above] {$1$};
\draw (1, 0.9) -- (1,1.1) node[above] {$2$};
\draw (2, 0.9) -- (2,1.1) node[above] {$3$};
\draw (3, 0.9) -- (3,1.1) node[above] {$4$};
\draw (4, 0.9) -- (4,1.1) node[above] {$5$};
\draw (5, 0.9) -- (5,1.1) node[above] {$6$};

\draw (0,1) .. controls (1,0.7) .. (2,1);
\draw (3,1) .. controls (3.5,0.7) .. (4,1);
\draw (1,0) .. controls (1.5,0.3) .. (2,0);
\draw (3,0) .. controls (3.5,0.3) .. (4,0);
\draw (0,0) -- (1,1);
\draw (5,0) -- (5,1);
\end{tikzpicture}
\times
\begin{tikzpicture} [baseline=3mm,blue,line width=1pt, xscale = 0.6, yscale=0.8,
                      draw/.append style={rounded corners},
                      every node/.append style={font=\fontsize{5}{5}\selectfont}]
\draw (0, -0.1) node[below] {$1$} -- (0,0.1);
\draw (1, -0.1) node[below] {$2$} -- (1,0.1);
\draw (2, -0.1) node[below] {$3$} -- (2,0.1);
\draw (3, -0.1) node[below] {$4$} -- (3,0.1);
\draw (4, -0.1) node[below] {$5$} -- (4,0.1);
\draw (5, -0.1) node[below] {$6$} -- (5,0.1);

\draw (0, 0.9) -- (0,1.1) node[above] {$1$};
\draw (1, 0.9) -- (1,1.1) node[above] {$2$};
\draw (2, 0.9) -- (2,1.1) node[above] {$3$};
\draw (3, 0.9) -- (3,1.1) node[above] {$4$};
\draw (4, 0.9) -- (4,1.1) node[above] {$5$};
\draw (5, 0.9) -- (5,1.1) node[above] {$6$};

\draw (1,1) .. controls (1.5,0.7) .. (2,1);
\draw (2,0) .. controls (3.5,0.3) .. (5,0);
\draw (0,0) -- (3,1);
\draw (1,0) -- (0,1);
\draw (3,0) -- (5,1);
\draw (4,0) -- (4,1);
\end{tikzpicture}
=
\begin{tikzpicture} [baseline=11mm,blue,line width=1pt, xscale = 0.6, yscale=0.8,
                      draw/.append style={rounded corners},
                      every node/.append style={font=\fontsize{5}{5}\selectfont}]
\draw (0, 1.9) node[below,left] {$1$} -- (0,2.1);
\draw (1, 1.9) node[below,left] {$2$} -- (1,2.1);
\draw (2, 1.9) node[below,left] {$3$} -- (2,2.1);
\draw (3, 1.9) node[below,left] {$4$} -- (3,2.1);
\draw (4, 1.9) node[below,left] {$5$} -- (4,2.1);
\draw (5, 1.9) node[below,left] {$6$} -- (5,2.1);

\draw (0, 2.9) -- (0,3.1) node[above] {$1$};
\draw (1, 2.9) -- (1,3.1) node[above] {$2$};
\draw (2, 2.9) -- (2,3.1) node[above] {$3$};
\draw (3, 2.9) -- (3,3.1) node[above] {$4$};
\draw (4, 2.9) -- (4,3.1) node[above] {$5$};
\draw (5, 2.9) -- (5,3.1) node[above] {$6$};

\draw (0,3) .. controls (1,2.7) .. (2,3);
\draw (3,3) .. controls (3.5,2.7) .. (4,3);
\draw (1,2) .. controls (1.5,2.3) .. (2,2);
\draw (3,2) .. controls (3.5,2.3) .. (4,2);
\draw (0,2) -- (1,3);
\draw (5,2) -- (5,3);

\draw (0, -0.1) node[below] {$1$} -- (0,0.1);
\draw (1, -0.1) node[below] {$2$} -- (1,0.1);
\draw (2, -0.1) node[below] {$3$} -- (2,0.1);
\draw (3, -0.1) node[below] {$4$} -- (3,0.1);
\draw (4, -0.1) node[below] {$5$} -- (4,0.1);
\draw (5, -0.1) node[below] {$6$} -- (5,0.1);

\draw (0, 0.9) -- (0,1.1) node[above,left] {$1$};
\draw (1, 0.9) -- (1,1.1) node[above,left] {$2$};
\draw (2, 0.9) -- (2,1.1) node[above,left] {$3$};
\draw (3, 0.9) -- (3,1.1) node[above,left] {$4$};
\draw (4, 0.9) -- (4,1.1) node[above,left] {$5$};
\draw (5, 0.9) -- (5,1.1) node[above,left] {$6$};

\draw (1,1) .. controls (1.5,0.7) .. (2,1);
\draw (2,0) .. controls (3.5,0.3) .. (5,0);
\draw (0,0) -- (3,1);
\draw (1,0) -- (0,1);
\draw (3,0) -- (5,1);
\draw (4,0) -- (4,1);

\draw[densely dotted] (0,1) -- (0,2);
\draw[densely dotted] (1,1) -- (1,2);
\draw[densely dotted] (2,1) -- (2,2);
\draw[densely dotted] (3,1) -- (3,2);
\draw[densely dotted] (4,1) -- (4,2);
\draw[densely dotted] (5,1) -- (5,2);
\end{tikzpicture}
=
\delta^1{\cdot}
\begin{tikzpicture} [baseline=3mm,blue,line width=1pt, xscale = 0.6, yscale=0.8,
                      draw/.append style={rounded corners},
                      every node/.append style={font=\fontsize{5}{5}\selectfont}]
\draw (0, -0.1) node[below] {$1$} -- (0,0.1);
\draw (1, -0.1) node[below] {$2$} -- (1,0.1);
\draw (2, -0.1) node[below] {$3$} -- (2,0.1);
\draw (3, -0.1) node[below] {$4$} -- (3,0.1);
\draw (4, -0.1) node[below] {$5$} -- (4,0.1);
\draw (5, -0.1) node[below] {$6$} -- (5,0.1);

\draw (0, 0.9) -- (0,1.1) node[above] {$1$};
\draw (1, 0.9) -- (1,1.1) node[above] {$2$};
\draw (2, 0.9) -- (2,1.1) node[above] {$3$};
\draw (3, 0.9) -- (3,1.1) node[above] {$4$};
\draw (4, 0.9) -- (4,1.1) node[above] {$5$};
\draw (5, 0.9) -- (5,1.1) node[above] {$6$};

\draw (0,1) .. controls (1,0.7) .. (2,1);
\draw (3,1) .. controls (3.5,0.7) .. (4,1);
\draw (2,0) .. controls (3.5,0.3) .. (5,0);
\draw (0,0) .. controls (2,0.3) .. (4,0);
\draw (1,0) -- (1,1);
\draw (3,0) -- (5,1);
\end{tikzpicture}
$$

It is easy to see that we have $2n-1$ possibilities to join the first dot with another one, then $2n-3$ possibilities for the next dot and so on. So there are $(2n-1)!! = (2n-1){\cdot}(2n-3){\cdot} \ldots 3 {\cdot}1$ number of Brauer diagrams, which implies the dimension of $\B$ is $(2n-1)!!$.

\subsection{(Graded) cellular algebras} \label{sec:cellular}

Following~\cite{GL}, we now introduce the graded cellular algebras. Reader may also refer to~\cite{HuMathas:GradedCellular}. Let $R$ be a commutative ring with $1$ and let $A$ be a unital $R$-algebra.

\begin{Definition} \label{def:GradedCellular}
A \textit{graded cell datum} for $A$ is a triple $(\Lambda,T,C,\deg)$ where $\Lambda = (\Lambda,>)$ is a poset, either finite or infinite, and $T(\lambda)$ is a finite set for each $\lambda \in \Lambda$, $\deg$ is a function from $\coprod_\lambda T(\lambda)$ to $\Z$, and
$$
C\map{\prod_{\lambda\in\Lambda} T(\lambda)\times T(\lambda)}{A}
$$
is an injective map which sends $(s,t)$ to $a_{st}^\lambda$ such that:
\begin{enumerate}
\item $\set{ a_{\s\t}^\lambda|\lambda\in\Lambda,s,t\in T(\lambda)}$ is an $R$-free basis of $A$;

\item for any $r\in A$ and $\t\in T(\lambda)$, there exists scalars $c_{\t}^v(r)$ such that, for any $\s\in T(\lambda)$,
$$
a_{\s\t}^\lambda{\cdot}r \equiv \sum_{\v\in T(\lambda)} c_{\t}^{\v}(r) a_{\s\v}^\lambda \mod{A^{>\lambda}}
$$
where $A^{>\lambda}$ is the $R$-submodule of $A$ spanned by $\set{ a_{\x\y}^\mu|\mu>\lambda, \x,\y\in T(\mu)}$;

\item the $R$-linear map $*\map{A}{A}$ which sends $a_{\s\t}^\lambda$ to $a_{\t\s}^\lambda$, for all $\lambda\in \Lambda$ and $\s,\t\in T(\lambda)$, is an anti-isomorphism of~$A$.

\item each basis element $a_{\s\t}^\lambda$ is homogeneous of degree $\deg a_{\s\t}^\lambda = \deg \s + \deg \t$, for $\lambda\in \Lambda$ and all $\s,\t\in T(\lambda)$.
\end{enumerate}
\end{Definition}

If a graded cell datum exists for $A$ then $A$ is a \textit{graded cellular algebra}. Similarly, by forgetting the grading we can define a \textit{cell datum} and hence a \textit{cellular algebra}.

Suppose $A$ is a graded cellular algebra with graded cell datum $(\Lambda,T,C,\deg)$. For any $\lambda \in \Lambda$, define $A^{\geq \lambda}$ to be the $R$-submodule of $A$ spanned by
$$
\set{c_{\s\t}^\mu|\mu \geq \lambda, \s,\t\in T(\mu)}.
$$

Then $A^{>\lambda}$ is an ideal of $A^{\geq\lambda}$ and hence $A^{\geq\lambda}/A^{>\lambda}$ is a $A$-module. For any $\s\in T(\lambda)$ we define $C_\s^\lambda$ to be the $A$-submodule of $A^{\geq\lambda}/A^{>\lambda}$ with basis $\set{a^\lambda_{\s\t} + A^{>\lambda}|\t\in T(\lambda)}$. By the cellularity of $A$ we have $C_\s^\lambda \cong C_\t^\lambda$ for any $\s,\t\in T(\lambda)$.

\begin{Definition} \label{def:cell}
Suppose $\lambda \in \P_n$. Define the \textit{cell module} of $A$ to be $C^\lambda = C_\s^\lambda$ for any $\s\in T(\lambda)$, which has basis $\set{a^\lambda_\t|\t\in T(\lambda)}$ and for any $r \in A$,
$$
a^\lambda_\t{\cdot}r= \sum_{\u\in T(\lambda)} c_\u^r a_\u^\lambda
$$
where $c_\u^r$ are determined by
$$
a^\lambda_{\s\t}{\cdot}r = \sum_{\u\in T(\lambda)} c_\u^r a^\lambda_{\s\u} + A^{>\lambda}.
$$
\end{Definition}

We can define a bilinear map $\langle\cdot,\cdot\rangle\map{C^\lambda\times C^\lambda}{\Z}$ such that
$$
\langle a^\lambda_\s,a^\lambda_\t\rangle a^\lambda_{\u\v} = a^\lambda_{\u\s}a^\lambda_{\t\v} + A^{>\lambda}
$$
and let $\text{rad }C^\lambda = \set{\s\in C^\lambda|\langle \s,\t\rangle = 0\text{ for all }\t\in C^\lambda}$. As $\langle\cdot,\cdot \rangle$ is homogeneous of degree $0$, $\text{rad }C^\lambda$ is a graded $A$-submodule of $C^\lambda$.

\begin{Definition}\label{def:simple}
Suppose $\lambda \in \P_n$. Let $D^\lambda = C^\lambda/\text{rad }C^\lambda$ as a graded $A$-module.
\end{Definition}

Exactly as in the ungraded case~\cite[Theorem 3.4]{GL}, we obtain the following:

\begin{Theorem} [\protect{Hu-Mathas~\cite[Theorem 2.10]{HuMathas:GradedCellular}}] ~\label{simpleMod}
The set $\set{D^\lambda\langle k\rangle|\lambda\in \Lambda, D^\lambda\neq 0, k\in\Z}$ is a complete set of pairwise non-isomorphic graded simple $A$-modules.
\end{Theorem}

In particular, the symmetric group algebra $R\Sym_n$ is a graded cellular algebra. The details of the graded cellular structure of $R\Sym_n$ will be introduced in Section 2.5. It is well-known that the Brauer algebras $\B$ are (ungraded) cellular algebras~\cite[Theorem 4.10]{GL}. In the following two subsections we will construct the cellular basis of the Brauer algebras.

\subsection{Combinatorics} \label{sec:combinatorics}

In this subsection we introduce the combinatorics of up-down tableaux, which will be used to index the cellular basis of the Brauer algebras. Throughout the rest of this paper, we fix $R$ to be a field with characteristic $0$ and $\delta \in R$.

Recall that a \textit{partition} of $n$ is a weakly decreased sequence of nonnegative integers $\lambda = (\lambda_1, \lambda_2, \ldots)$ such that $|\lambda| := \lambda_1 + \lambda_2 + \ldots = n$. In such case we denote $\lambda \vdash n$. Because $|\lambda| < \infty$, there are finite many nonzero $\lambda_i$ for $i \geq 1$. Because $\lambda_1 \geq \lambda_2 \geq \ldots$, there exists $k \geq 1$ such that $\lambda_k \geq 0$ and $\lambda_{k+1} = \lambda_{k+2} = \ldots = 0$. Usually we will write $\lambda = (\lambda_1, \ldots, \lambda_k)$ instead of an infinite sequence.

Let $\H_n$ be the set of all partitions of $n$. We can define a partial ordering $\unlhd$ on $\H_n$, which is called the \textit{dominance ordering}. Given $\lambda,\mu \in \H_n$, we say $\lambda \unlhd \mu$ if for any $k \geq 1$, we have $|\sum_{i=1}^k \lambda_i| \leq |\sum_{i=1}^k \mu_i|$. Write $\lambda \lhd \mu$ if $\lambda \unlhd \mu$ and $\lambda \neq \mu$. The dominance ordering can be extended to a total ordering $\leq$, the \textit{lexicographic ordering}. We write $\lambda < \mu$ if there exists $k$ such that $\lambda_i = \mu_i$ for all $i < k$ and $\lambda_k < \mu_k$. Define $\lambda \leq \mu$ if $\lambda < \mu$ or $\lambda = \mu$. Then $\lambda \unlhd \mu$ implies $\lambda \leq \mu$.

The \textit{Young diagram} of a partition $\lambda$ is the set
$$
[\lambda] := \set{(r,l) | 1 \leq l \leq \lambda_r}.
$$

For example, $\lambda = (3,2,2)$ is a partition of $7$, and the Young diagram of $\lambda$ is
$$
[\lambda] = \ydiag(3,2,2).
$$

Define $\alpha = \pm (r,l)$ to be a \textit{node} for positive integers $r$ and $l$, and denote $\alpha > 0$ if $\alpha = (r,l)$ and $\alpha < 0$ if $\alpha = - (r,l)$. In this paper we allow to work with linear combination of nodes. In more details, suppose $\alpha$ and $\beta$ are two nodes, we write $\alpha + \beta = 0$ if $\alpha = (r,l)$ and $\beta = -(r,l)$, or vise versa. Similarly, we write $\alpha = - \beta$ if $\alpha + \beta = 0$.

Suppose $\alpha = (r,l)$. We say $\alpha$ is a node of $\lambda$ in row $r$ and column $l$ if $\alpha \in [\lambda]$. For example, let $\alpha = (2,2)$, $\beta = (2,3)$ and $\lambda = (3,2,2)$. Then $\alpha$ is a node of $\lambda$ in row $2$ and column $2$, and $\beta$ is not a node of $\lambda$.

Suppose $\lambda$ is a partition. A node $\alpha > 0$ is \textit{addable} if $\lambda \cup \{\alpha\}$ is still a partition, and it is \textit{removable} if $\lambda \backslash \{\alpha\}$ is still a partition. Let $\mathscr A(\lambda)$ and $\mathscr R(\lambda)$ be the sets of addable and removable nodes of $\lambda$, respectively, and set $\mathscr{AR}(\lambda) = \mathscr A(\lambda) \cup \mathscr R(\lambda)$.

A \textit{$\lambda$-tableau} is any bijection $\t\map{\{1,2,\ldots,n\}}{[\lambda]}$. We identify a $\lambda$-tableau $\t$ with a labeling of the diagram of $\lambda$. That is, we label the node $(r,l)\in[\lambda]$ with the integer $\t^{-1}(r,l)$. We say a tableau $\t$ has shape $\lambda$ if it is a $\lambda$-tableau. A tableau $\t$ is \textit{standard} if the entries of each row and each column of $\t$ increase. Suppose $\lambda \vdash n$. Denote $\Std(\lambda)$ to be the set of all standard tableau of shape $\lambda$.

Define $\widehat B_n := \set{(\lambda,f)|\lambda \in \widehat H_{n - 2f}\text{ and }0 \leq f \leq \lfloor \frac{n}{2} \rfloor}$ and $\widehat B$ to be the graph with
\begin{enumerate}
\item vertices at level $n$: $\widehat B_n$, and

\item an edge $(\lambda,f) \rightarrow (\mu,m)$, $(\lambda,f) \in \widehat B_{n-1}$ and $(\mu,m) \in \widehat B_n$, if either $\mu$ is obtained by adding a node to $\lambda$, or by deleting a node from $\lambda$.
\end{enumerate}

We can extend the dominance ordering and lexicographic ordering of partitions to $\widehat B_n$ by defining $(\lambda,f) \unlhd (\mu,m)$ if $f < m$, or $f = m$ and $\lambda\unlhd \mu$; and $(\lambda,f) \leq (\mu,m)$ if $f < m$, or $f = m$ and $\lambda \leq \mu$. We define $\lhd$ and $<$ similarly.

\begin{Definition}
Let $(\lambda,f)\in \widehat B_n$. An \textit{up-down tableau} of shape $(\lambda,f)$ is a sequence
\begin{equation} \label{def:udtab:1}
\t = ((\lambda^{(0)},f_0),(\lambda^{(1)},f_1),\ldots,(\lambda^{(n)},f_n)),
\end{equation}
where $(\lambda^{(0)},f_0) = (\emptyset, 0)$, $(\lambda^{(n)},f_n) = (\lambda,f)$ and $(\lambda^{(k-1)},f_{k-1}) \rightarrow (\lambda^{(k)},f_k)$ is an edge in $\widehat B$, for $k = 1,\ldots,n$. We write $\Shape(\t) = (\lambda,f)$. If $k = 0,1,\ldots,n$, we denote $\t_k = \lambda^{(k)}$ and define the truncation of $\t$ to level $k$ to be the up-down tableau
$$
\t|_k = ((\lambda^{(0)},f_0),(\lambda^{(1)},f_1),\ldots,(\lambda^{(k)},f_k)).
$$

For any $0 \leq f \leq \lfloor \frac{n}{2} \rfloor$ and $\lambda \vdash n - 2f$, define
$$
\Tud_n(\lambda) := \set{\t|\text{$\t$ is an up-down tableau of shape $(\lambda,f) \in \widehat B_n$}}.
$$
\end{Definition}

Suppose $\s,\t \in \Tud_n(\lambda)$. We define the dominance ordering $\s \unlhd \t$ if $\s_k \unlhd \t_k$ for any $k$ with $1 \leq k \leq n$ and $\s \lhd \t$ if $\s \unlhd \t$ and $\s \neq \t$.

An up-down tableau $\t = ((\lambda^{(0)},f_0),(\lambda^{(1)},f_1),\ldots,(\lambda^{(n)},f_n))$ can be identified with a $n$-tuple of nodes:
\begin{equation} \label{def:udtab:2}
\t = (\alpha_1, \alpha_2, \ldots, \alpha_n),
\end{equation}
where $\alpha_k = (r,l)$ if $\lambda^{(k)} = \lambda^{(k-1)} \cup \{(r,l)\}$ and $\alpha_k = -(r,l)$ if $\lambda^{(k)} = \lambda^{(k-1)} \backslash \{(r,l)\}$ for $1 \leq k \leq n$. Therefore an up-down tableau $\t$ can be identified as a map $\t\map{\{1,2,\ldots,n\}}{\set{\pm(r,l) | r, l \in \Z_{>0}}}$. Note that the range of $\t$ is all (positive and negative) nodes and $\t$ is not necessary injective. We have $\t(k) = \alpha_k$ for $1 \leq k \leq n$.


We define a right action of $\Sym_n$ on the up-down tableaux of $n$. Suppose $\t = (\alpha_1, \ldots, \alpha_n)$ and $1 \leq k \leq n-1$. Define $\t{\cdot}s_k = (\alpha_1, \ldots, \alpha_{k-2}, \alpha_k, \alpha_{k-1}, \alpha_{k+1}, \ldots, \alpha_n)$. We note that $\t{\cdot}s_k$ is not necessarily an up-down tableau, and when $\t \in \Tud_n(\lambda)$, then $\t{\cdot}s_k \in \Tud_n(\lambda)$ if $\t{\cdot}s_k$ is an up-down tableau.

Two nodes $\alpha = (i,j) > 0$ and $\beta = (r,l) > 0$ are \textit{adjacent} if $i = r \pm 1$ and $j = l$, or $i = r$ and $j = l \pm 1$. The next Lemma can be verified directly by the construction of up-down tableaux. It gives conditions for $\t{\cdot}s_k$ to be an up-down tableau.

\begin{Lemma} \label{y:h2:8}
Suppose $(\lambda,f) \in \widehat B_n$ and $\t \in \Tud_n(\lambda)$. For $1 \leq k \leq n-1$, $\t{\cdot}s_k$ is an up-down tableau if and only if one of the following conditions hold:
\begin{enumerate}
\item $\t(k) > 0$, $\t(k+1) > 0$, and $\t(k)$ and $\t(k+1)$ are not adjacent.

\item $\t(k) < 0$, $\t(k+1) < 0$, and $-\t(k)$ and $-\t(k+1)$ are not adjacent.

\item $\t(k) > 0$, $\t(k+1) < 0$, and $\t(k) + \t(k+1) \neq 0$.

\item $\t(k) < 0$, $\t(k+1) > 0$, and $\t(k) + \t(k+1) \neq 0$.
\end{enumerate}
\end{Lemma}

Recall $\delta \in R$ and let $x$ be an indeterminate. Suppose $\alpha = (r,l)$ is a positive node. The \textit{content} of $\alpha$ is $\cont(\alpha) = \frac{x-1}{2} + l - r$ and the \textit{residue} of $\alpha$, $\res(\alpha)$, is the evaluation of the content at $x = \delta$. Set $\cont(-\alpha) = -\cont(\alpha)$ and $\res(-\alpha) = -\res(\alpha)$.

Suppose $\t = (\alpha_1, \alpha_2, \ldots, \alpha_n)$ is an up-down tableau. Define $c_\t(k) = \cont(\alpha_k)$ and $r_\t(k) = \res(\alpha_k)$ for $1 \leq k \leq n$. We define the \textit{residue sequence} of $\t$ to be $\bi_\t = (i_1, i_2, \ldots, i_n)$ such that $i_k = r_\t(k)$ for any $1 \leq k \leq n$. Let $P = \frac{\delta - 1}{2} + \Z$. One can see that $\res(\alpha) \in P$ for any node $\alpha$. Therefore one can see that the residue sequence $\bi_\t \in P^n$. Suppose $\bi \in P^n$. Let $\Tud_n(\bi)$ be the set containing all the up-down tableaux with residue sequence $\bi$.

Here we give the notations for subtraction and concatenation of $n$-tuples in $P^n$. For $\bi = (i_1,\ldots,i_n) \in P^n$, denote $\bi|_k = (i_1,\ldots,i_j) \in P^k$ for $1 \leq k \leq n$; and denote $\bi \vee i = (i_1,\ldots,i_n,i) \in P^{n+1}$ for $i \in P$.

We now introduce the degree function of the set of up-down tableaux. Ultimately this degree function will describe the grading on $\B$.

Suppose we have $(\lambda,f) \rightarrow (\mu,m)$. Write $\lambda \ominus \mu = \alpha$ if $\lambda = \mu\cup \{\alpha\}$ or $\mu = \lambda \cup \{\alpha\}$. For any up-down tableau $\t$ and an integer $k$, with $1 \leq k \leq n$, let $\lambda = \t_{k-1}$, $\mu = \t_k$ and $\alpha = (r,l) = \lambda \ominus \mu$. Define
$$
\begin{array}{rcll}
\mathscr A_\t(k) & = & \set{ \beta = (k,c) \in \mathscr A(\lambda) | \text{$\res(\beta) = \res(\alpha)$ and $k > r$}}, &\text{if $\mu = \lambda \cup \{\alpha\}$,}\\
\mathscr {\widehat A}_\t(k) & = & \set{ \beta = (k,c) \in \mathscr A(\mu) | \text{$\res(\beta) = -\res(\alpha)$ and $k \neq r$}}, & \text{if $\mu = \lambda \backslash \{\alpha\}$;}\\
\mathscr R_\t(k) & = & \set{ \beta = (k,c) \in \mathscr R(\lambda) | \text{$\res(\beta) = \res(\alpha)$ and $k > r$}}, & \text{if $\mu = \lambda \cup \{\alpha\}$,}\\
\mathscr {\widehat R}_\t(k) & = & \set{ \beta = (k,c) \in \mathscr R(\mu) | \text{$\res(\beta) = -\res(\alpha)$}}, & \text{if $\mu = \lambda \backslash \{\alpha\}$.}
\end{array}
$$

\begin{Definition}
Suppose $\t$ is an up-down tableau of size $n$. For integer $k$ with $0 \leq k \leq n-1$, write $\lambda = \t_{k-1}$, $\mu = \t_k$ and $\alpha = \lambda \ominus \mu$. Define
$$
\deg (\t|_{k-1} \Rightarrow \t|_k) :=
\begin{cases}
|\mathscr A_\t(k)| - |\mathscr R_\t(k)|, & \text{if $\mu = \lambda \cup \{\alpha\}$,}\\
|\mathscr {\widehat A}_\t(k)| - |\mathscr {\widehat R}_\t(k)| + \delta_{\res(\alpha),-\frac{1}{2}}, & \text{if $\mu = \lambda\backslash\{\alpha\}$,}
\end{cases}
$$
and the \textit{degree} of $\t$ is
$$
\deg \t := \sum_{k = 1}^n \deg (\t|_{k-1} \Rightarrow \t|_k).
$$
\end{Definition}

\begin{Remark} \label{remark:deg:1}
We note that when the characteristic of the field is $0$, we have $|\mathscr A_\t(k)| = |\mathscr R_\t(k)| = 0$ for any $1 \leq k \leq n$. Therefore, we always have $\deg (\t|_{k-1} \Rightarrow \t|_k) = 0$ when $\mu = \lambda \cup \{\alpha\}$.
\end{Remark}

\begin{Example} \label{ex:deg:1}
Let $n = 6$, $\lambda = (1,1)$, $\delta = 1$ and $\t = \left( \emptyset, \ydiag(1), \ydiag(2), \ydiag(2,1), \ydiag(2,2), \ydiag(2,1), \ydiag(1,1) \right) \in \Tud_n(\lambda)$. By~\autoref{remark:deg:1}, we have $\deg \t = \sum_k \deg (\t|_{k-1} \Rightarrow \t|_k)$, where $k$ take values such that $\t|_k$ is obtained by removing a node from $\t|_{k-1}$. Therefore, we have $\deg \t = \deg (\t|_4 \Rightarrow \t|_5) + \deg (\t|_5 \Rightarrow \t|_6)$.

By the definitions, we have $\mathscr {\widehat A}_\t(5) = \mathscr {\widehat R}_\t(5) = \emptyset$, $\mathscr {\widehat A}_\t(6) = \emptyset$ and $\mathscr {\widehat R}_\t(6) = \{(2,1)\}$. Because $\delta = 1$, for any node $\alpha$, we have $\res(\alpha) \in \Z$, which implies $\delta_{\res(\alpha), -\frac{1}{2}} = 0$. Hence, the degree of $\t$ is
$$
\deg\t = \deg (\t|_4 \Rightarrow \t|_5) + \deg (\t|_5 \Rightarrow \t|_6) = 0 - 1 = -1.
$$
\end{Example}

\begin{Example} \label{ex:deg:2}
Let $n = 6$, $\lambda = (1,1)$, $\delta = 0$ and $\t = \left( \emptyset, \ydiag(1), \ydiag(2), \ydiag(2,1), \ydiag(2,2), \ydiag(2,1), \ydiag(1,1) \right) \in \Tud_n(\lambda)$. Following the same argument as in~\autoref{ex:deg:1}, we have $\deg \t = \deg (\t|_4 \Rightarrow \t|_5) + \deg (\t|_5 \Rightarrow \t|_6)$.

By the definitions, we have $\mathscr {\widehat A}_\t(5) = \emptyset$ and $\mathscr {\widehat R}_\t(5) = \{(1,2)\}$. If we set $\lambda = \t_4$ and $\mu = \t_5$, we have $\mu = \lambda\backslash\{\alpha\}$ where $\alpha = (2,2)$. Because $\res(\alpha) = -\frac{1}{2}$, we have
$$
\deg (\t|_4 \Rightarrow \t|_5) = |\mathscr {\widehat A}_\t(5)| - |\mathscr {\widehat R}_\t(5)| + \delta_{\res(\alpha),-\frac{1}{2}} = 0 - 1 + 1 = 0.
$$

Similarly, we have $\mathscr {\widehat A}_\t(6) = \mathscr {\widehat R}_\t(6) = \emptyset$. If we set $\lambda = \t_5$ and $\mu = \t_6$, we have $\mu = \lambda\backslash\{\alpha\}$ where $\alpha = (1,2)$. Because $\res(\alpha) = \frac{1}{2}$, we have
$$
\deg (\t|_5 \Rightarrow \t|_6) = |\mathscr {\widehat A}_\t(6)| - |\mathscr {\widehat R}_\t(6)| + \delta_{\res(\alpha),-\frac{1}{2}} = 0.
$$

Hence, the degree of $\t$ is $\deg\t = \deg (\t|_4 \Rightarrow \t|_5) + \deg (\t|_5 \Rightarrow \t|_6) = 0 + 0 = 0$.
\end{Example}

\subsection{Jucys-Murphy elements and Cellularity of Brauer algebras} \label{sec:cellular:Brauer}

In the Brauer algebra $\B$, Nazarov~\cite{Nazarov:brauer} defined \textit{Jucys-Murphy elements} $L_k$ for $1 \leq k \leq n$ by $L_1 = \frac{\delta-1}{2}$ and
$$
L_{k+1} = s_k - e_k + s_k L_k s_k,\qquad \text{for $1 \leq k \leq n-1$.}
$$

\begin{Lemma} [Nazarov~\cite{Nazarov:brauer}]
The following relations hold in the algebra $\B$:
$$
\begin{array}{c}
s_k L_r = L_r s_k, \quad e_k L_r = L_r e_k;\quad r\neq k, k+1;\\
s_k L_k - L_{k+1} s_k = e_k - 1, \quad L_k s_k - s_k L_{k+1} = e_k - 1;\\
e_k(L_k + L_{k+1}) = 0, \quad (L_k + L_{k+1})e_k = 0.
\end{array}
$$
\end{Lemma}

Graham and Lehrer~\cite{GL} proved that $\B$ is a cellular algebra over any commutative ring $R$. Enyang~\cite{Enyang, Enyang:ss} constructed another cellular basis indexed by pairs $(\s,\t)$, where $\s,\t \in \Tud_n(\lambda)$ and $(\lambda,f) \in \widehat B_n$. We will only state the Theorem here.

\begin{Theorem} [Enyang~\cite{Enyang:ss}] \label{B:basis}
Let $\B$ be a Brauer algebra over a commutative ring $R$ and $*\map{\B}{\B}$ be the $R$-linear involution which fixes $s_k$ and $e_k$ for $1 \leq k \leq n-1$. Then $\B$ has a cellular basis
$$
\set{m_{\s\t} | \s,\t \in \Tud_n(\lambda), (\lambda,f) \in \widehat B_n}
$$
such that
$$
m_{\s\t} L_k = r_\t(k) m_{\s\t} + \sum_{\substack{\v \in \Tud_n(\lambda) \\ \v \rhd \t}} c_\v m_{\s\v} + \sum_{\substack{\u,\v \in \Tud_n(\mu) \\ (\mu,m) \in \widehat B_n \\(\mu,m) \rhd (\lambda,f)}} c_{\u\v} m_{\u\v}.
$$
\end{Theorem}

%
%
%

\subsection{Seminormal forms and idempotents} \label{sec:semi}

In this subsection we develop the theory of seminormal forms for Brauer algebras, summarizing results that are in the literature, such as~\cite{M:seminormal,RuiSi:BrauerDet}.

Recall that $\B$ is a $R$-algebra, where $R$ is a field of characteristic $0$. Define $\F = R(x)$ to be the rational field with indeterminate $x$ and $\O = R[x]_{(x - \delta)} = R[[x - \delta]]$. Let $\m = (x - \delta)\O \subset \O$. Then $\m$ is a maximal ideal of $\O$ and $R \cong \O/\m$.

Let $\Bx$ and $\BOx$ be the Brauer algebras over $\F$ and $\O$, respectively. Then $\Bx = \BOx \otimes_\O \F$ and $\B \cong \BOx \otimes_\O R \cong \BOx/(x - \delta)\BOx$. In order to avoid confusion we will write the generators of $\BOx$ and $\Bx$ as $s_k^\O$ and $e_k^\O$ and generators of $\B$ as $s_k$ and $e_k$. Hence for any element $w \in \B$, we write $w^\O = w \otimes_R 1_\O \in \BOx$, so that $w = w^\O \otimes_\O 1_R$.

Because $\B \cong \BOx\otimes_\O R \cong \BOx/(x-\delta)\BOx$, if $x,y \in \BOx$ and we have $x \equiv y \pmod{(x-\delta)\BOx}$, then $x\otimes_\O 1_R = y\otimes_\O 1_R$ as elements of $\B$. This observation will give us a way to extend the results of $\BOx$ to $\B$.

The next Lemma says that an up-down tableau $\t$ is completely determined by its contents $c_\t(k)$ for $1 \leq k \leq n$.

\begin{Lemma}
Suppose $\s,\t$ are up-down tableaux of size $n$. Then $\s = \t$ if and only if $c_\s(k) = c_\t(k)$ for all $1 \leq k \leq n$.
\end{Lemma}

Hence $\Bx$ and $\F$ satisfies the separation condition in the sense of Mathas~\cite[Definition 2.8]{M:seminormal} or Rui-Si~\cite[Assumption 3.1]{RuiSi:BrauerDet}. The results we have in the rest of this subsection are already included in Mathas~\cite[Section 3, 4]{M:seminormal} and Rui-Si~\cite[Section 3]{RuiSi:BrauerDet}.

\begin{Definition}
Suppose $\Bx$ is the Brauer algebra over $\F$ and $\set{m_{\s\t} | \Tud_n(\lambda), (\lambda,f) \in \widehat B_n}$ is the basis of $\Bx$.
\begin{enumerate}
\item For any $1 \leq k \leq n$, define $\mathscr C(k) = \set{c_\t(k) | \t \in \Tud_n(\lambda), (\lambda,f) \in \widehat B_n}$.

\item $F_\t = \prod_{k=1}^n \prod_{\substack{c \in \mathscr C(k) \\ c_\t(k) \neq c}} \frac{L_k^\O - c}{c_\t(k) - c}$.

\item $f_{\s\t} = F_\s m_{\s\t} F_\t$,
\end{enumerate}
where $\s,\t \in \Tud_n(\lambda)$ and $(\lambda,f) \in \widehat B_n$.
\end{Definition}

By~\autoref{B:basis}, $f_{\s\t} = m_{\s\t} + \sum_{\u \rhd \s, \v \rhd \t} r_{\u\v} m_{\u\v}$, for some $r_{\u\v} \in \F$. Therefore
$$
\set{f_{\s\t} | \s,\t \in \Tud_n(\lambda), (\lambda,f) \in \widehat B_n}
$$
is a basis of $\Bx$. This basis is called the \textit{seminormal basis} of $\B$; see~\cite[Theorem 3.7]{M:seminormal}.

\begin{Lemma} \label{semi:L}
Suppose that $\Bx$ is the Brauer algebra over $\F$. Then we have
$$
f_{\s\t} L_k^\O = c_\t(k) f_{\s\t}, \qquad L_k^\O f_{\s\t} = c_\s(k) f_{\s\t}, \qquad \text{and} \qquad f_{\s\t} f_{\u\v} = \delta_{\t,\u} f_{\s\v},
$$
for $1 \leq k \leq n$ and $\s,\t, \u,\v \in \Tud_n(\lambda)$.
\end{Lemma}

Nazarov~\cite{Nazarov:brauer} gave the actions of $s_k^\O$ and $e_k^\O$ on the $f_{\s\t}$'s. Readers may also check Rui-Si~\cite{RuiSi:BrauerDet}. Suppose $\s,\t \in \Tud_n(\lambda)$ and $(\lambda,f) \in \widehat B_n$. Define
$$
f_{\s\t}e_k^\O = \sum_{\u \in \Tud_n(\lambda)}e_k(\t,\u)f_{\s\u} \qquad \text{and} \qquad f_{\s\t} s_k^\O = \sum_{\u \in \Tud_n(\lambda)} s_k(\t,\u) f_{\s\u}.
$$

\begin{Definition}
Suppose $1 \leq k \leq n-1$ and $(\lambda,f) \in \widehat B_n$. For $\t \in \Tud_n(\lambda)$ with $\t_{k-1} = \t_{k+1}$, define an equivalence relation $\overset{k}\sim$ by declaring that $\t \overset{k}\sim \s$ if $\t_r = \s_r$ whenever $1 \leq r \leq n$ and $r \neq k$, for $\s \in \Tud_n(\lambda)$.
\end{Definition}

Suppose $\bi,\bj \in P^n$ with $i_k + i_{k+1} = j_k + j_{k+1} = 0$ for some $1 \leq k \leq n-1$. We define an equivalence relation $\overset{k}\sim$ on $P^n$ by declaring that $\bi \overset{k}\sim \bj$ if $i_r = j_r$ whenever $1 \leq r \leq n$ and $r \neq k,k+1$. It is easy to see that $\t \overset{k}\sim \s$ only if $\s_{k-1} = \s_{k+1} = \t_{k-1} = \t_{k+1}$ and $\bi_\t \overset{k}\sim \bi_\s$. The next result is a special case of~\cite[4.2]{AMR}.

\begin{Lemma} \label{e_k:h1}
Suppose $\t \in \Tud_n(\lambda)$ with $\t_{k-1} = \t_{k+1} = \mu$. Then there is a bijection between $\mathscr{AR}(\mu)$ and the set $\set{\s \in \Tud_n(\lambda) | \s \overset{k}\sim \t}$.
\end{Lemma}

Suppose $\t \in \Tud_n(\lambda)$ with $\t_{k-1} = \t_{k+1}$ for some $k$ with $1 \leq k \leq n-1$. Define
$$
e_k(\t,\t) := (2c_\t(k) + 1) \prod_{\substack{\u \overset{k}\sim \t\\ \u\neq \t}}\frac{c_\t(k) + c_\u(k)}{c_\t(k) - c_\u(k)} \in \F.
$$
%
%
%
%
%
%
%

\begin{Theorem} \label{semi:B}
Suppose $(\lambda,f) \in \widehat B_n$ and $\t \in \Tud_n(\lambda)$. For any $1 \leq k \leq n-1$ and $\s \in \Tud_n(\lambda)$, we have:
\begin{enumerate}
\item If $\t_{k-1} \neq \t_{k+1}$ and $\t{\cdot}s_k$ does not exist, then
$$
f_{\s\t} s_k^\O = \frac{1}{c_\t(k+1) - c_\t(k)} f_{\s\t}.
$$

\item If $\t_{k-1} \neq \t_{k+1}$ and $\u = \t{\cdot}s_k \in \Tud_n(\lambda)$, then
$$
f_{\s\t} s_k^\O = \begin{cases}
 \frac{1}{c_\t(k+1) - c_\t(k)} f_{\s\t} + f_{\s\u}, & \text{if $\t \rhd \u$,}\\
 \frac{1}{c_\t(k+1) - c_\t(k)} f_{\s\t} + (1 - \frac{1}{(c_\s(k+1) - c_\s(k))^2})f_{\s\u}, & \text{if $\u \rhd \t$.}
\end{cases}
$$

\item If $\t_{k-1} \neq \t_{k+1}$, then $f_{\s\t} e_k^\O = 0$.

\item If $\t_{k-1} = \t_{k+1}$, then
$$
f_{\s\t} s_k^\O = \sum_{\u\overset{k}\sim \t} s_k(\t,\u) f_{\s\u} = \sum_{\u \overset{k}\sim \t} \frac{e_k(\t,\u) - \delta_{\t\u}}{c_\t(k) + c_\u(k)} f_{\s\u}.
$$

\item If $\t_{k-1} = \t_{k+1}$, then
$$
f_{\s\t} e_k^\O = \sum_{\u\overset{k}\sim \t} e_k(\t,\u) f_{\s\u}.
$$

\item If $\t_{k-1} = \t_{k+1}$ and $\u \overset{k}\sim \t \overset{k}\sim \v$, then
$$
e_k(\u,\t) e_k(\t,\v) = e_k(\u,\v) e_k(\t,\t).
$$
\end{enumerate}
\end{Theorem}

The following result gives an explicit construction on (central) primitive idempotents of $\Bx$. Such result has been proved by Mathas~\cite[Theorem 3.16]{M:seminormal} for general cellular algebras under separation condition.

For $\t \in \Tud_n(\lambda)$, define $\gamma_\t \in \F$ such that $f_{\t\t} f_{\t\t} = \gamma_\t f_{\t\t}$. By the cellularity of $\{f_{\s\t}\}$, we have $f_{\s\t} f_{\t\u} = \gamma_\t f_{\s\u}$ for any $\s,\u \in \Tud_n(\lambda)$. Note that $\gamma_\t$ can be computed recursively by Rui-Si~\cite[Proposition 4.9]{RuiSi:BrauerDet}.

\begin{Proposition} \label{idem:semi:1}
We have the following results:
\begin{enumerate}
\item Suppose $\t \in \Tud_n(\lambda)$ and $(\lambda,f) \in \widehat B_n$. Then $\f{\t}$ is a primitive idempotent of $\Bx$.

\item $\sum_{\t \in \Tud_n(\lambda)} \f{\t}$ is a central primitive idempotent. Moreover,
$$
\sum_{(\lambda,f) \in \widehat B_n} \sum_{\t \in \Tud_n(\lambda)} \f{\t} = 1.
$$
\end{enumerate}
\end{Proposition}

\subsection{Graded cellular algebras $R\Sym_n$} \label{sec:KLR}

In this subsection we will introduce a graded cellular structure on $R\Sym_n$ by giving a graded cellular basis of $R\Sym_n$.

Khovanov and Lauda~\cite{KhovLaud:diagI, KhovLaud:diagII} and Rouquier~\cite{Rouq:2KM} have introduced a naturally $\Z$-graded algebra $\mathscr R_n$. Fix an integer $e\in \{0,2,3,4\ldots\}$. Define $P' = \Z \cup \Z/2$ when $e = 0$ and $P' = \Z/e\Z$ when $e > 0$. Let $\Gamma_e$ be the oriented quiver with vertex set $\Z/e\Z$ and directed edges $i\rightarrow i+1$, for $i\in \Z/e\Z$. Thus, $\Gamma_e$ is the quiver of type $A_{\infty}$ if $e = 0$ and if $e \geq 2$ then it is a cyclic quiver of type $A_e^{(1)}$:
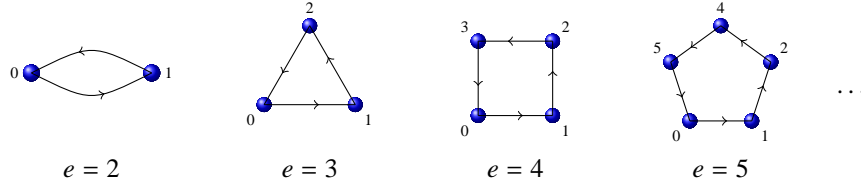
\begin{center}
\begin{tabular}{*5c}
  \begin{tikzpicture}[scale=0.8,decoration={curveto, markings,
            mark=at position 0.6 with {\arrow{>}}
    }]
    \useasboundingbox (-1.7,-0.7) rectangle (1.7,0.7);
    \foreach \x in {0, 180} {
      \shade[ball color=blue] (\x:1cm) circle(4pt);
    }
    \tikzstyle{every node}=[font=\tiny]
    \draw[postaction={decorate}] (0:1cm) .. controls (90:5mm) .. (180:1cm)
           node[left,xshift=-0.5mm]{$0$};
    \draw[postaction={decorate}] (180:1cm) .. controls (270:5mm) .. (0:1cm)
           node[right,xshift=0.5mm]{$1$};
  \end{tikzpicture}
& 
  \begin{tikzpicture}[scale=0.7,decoration={ markings,
            mark=at position 0.6 with {\arrow{>}}
    }]
    \useasboundingbox (-1.7,-0.7) rectangle (1.7,1.4);
    \foreach \x in {90,210,330} {
      \shade[ball color=blue] (\x:1cm) circle(4pt);
    }
    \tikzstyle{every node}=[font=\tiny]
    \draw[postaction={decorate}]( 90:1cm)--(210:1cm) node[below left]{$0$};
    \draw[postaction={decorate}](210:1cm)--(330:1cm) node[below right]{$1$};
    \draw[postaction={decorate}](330:1cm)--( 90:1cm)
           node[above,yshift=.5mm]{$2$};
  \end{tikzpicture}
& 
  \begin{tikzpicture}[scale=0.7,decoration={ markings,
            mark=at position 0.6 with {\arrow{>}}
    }]
    \useasboundingbox (-1.7,-0.7) rectangle (1.7,0.7);
    \foreach \x in {45, 135, 225, 315} {
      \shade[ball color=blue] (\x:1cm) circle(4pt);
    }
    \tikzstyle{every node}=[font=\tiny]
    \draw[postaction={decorate}] (135:1cm) -- (225:1cm) node[below left]{$0$};
    \draw[postaction={decorate}] (225:1cm) -- (315:1cm) node[below right]{$1$};
    \draw[postaction={decorate}] (315:1cm) -- (45:1cm) node[above right]{$2$};
    \draw[postaction={decorate}] (45:1cm) -- (135:1cm) node[above left]{$3$};
  \end{tikzpicture}
& 
  \begin{tikzpicture}[scale=0.7,decoration={ markings,
            mark=at position 0.6 with {\arrow{>}}
    }]
    \useasboundingbox (-1.7,-0.7) rectangle (1.7,0.7);
    \foreach \x in {18,90,162,234,306} {
      \shade[ball color=blue] (\x:1cm) circle(4pt);
    }
    \tikzstyle{every node}=[font=\tiny]
    \draw[postaction={decorate}] (162:1cm) -- (234:1cm)node[below left]{$0$};
    \draw[postaction={decorate}] (234:1cm) -- (306:1cm) node[below right]{$1$};
    \draw[postaction={decorate}] (306:1cm) -- (18:1cm) node[above right]{$2$};
    \draw[postaction={decorate}] (18:1cm) -- (90:1cm)
           node[above,yshift=.5mm]{$4$};
    \draw[postaction={decorate}] (90:1cm) -- (162:1cm) node[above left]{$5$};
  \end{tikzpicture}
&\raisebox{3mm}{$\dots$}
\\[4mm]
  $e=2$&$e=3$&$e=4$&$e=5$&
\end{tabular}
\end{center}

Let $(a_{i,j})_{i,j\in \Z/e\Z}$ be the symmetric Cartan matrix associated with $\Gamma_e$, so that
\begin{equation*} \label{notation: cartan}
a_{i,j} = \begin{cases}
2, & \text{ if $i = j$,}\\
0, & \text{ if $i\neq j \pm 1$,}\\
-1, & \text{ if $e\neq 2$ and $i = j \pm 1$,}\\
-2, & \text{ if $e = 2$ and $i = j + 1$.}
\end{cases}
\end{equation*}

To the quiver $\Gamma_e$ attach the standard Lie theoretic data of a Cartan matrix $(a_{ij})_{i,j\in \Z/e\Z}$, fundamental weights $\{\Lambda_i| i\in \Z/e\Z\}$, positive weights $\sum_{i\in \Z/e\Z}\N \Lambda_i$, positive roots $\bigoplus_{i\in \Z/e\Z}\N\alpha_i$ and let $(\cdot,\cdot)$ be the bilinear form determined by
$$
(\alpha_i,\alpha_j)=a_{ij}\qquad\text{and}\qquad
          (\Lambda_i,\alpha_j)=\delta_{ij},\qquad\text{for }i,j\in \Z/e\Z.
$$
Fix a \textit{weight} $\Lambda = \sum_{i\in \Z/e\Z} a_i \Lambda_i \in \sum_{i\in \Z/e\Z}\N \Lambda_i$. Then $\Lambda$ is a weight of
\textit{level} $l(\Lambda) = \l = \sum_{i\in \Z/e\Z} a_i$. A \textit{multicharge} for $\Lambda$ is a sequence
$\kappa_\Lambda = (\kappa_1,\ldots,\kappa_\l) \in (\Z/e\Z)^\l$ such that
$$
(\Lambda,\alpha_i) = a_i = \#\set{ 1\leq s\leq \l|\kappa_s \equiv i\pmod{e}}
$$
for any $i\in \Z/e\Z$.

The following algebras were introduced by Khovanov and Lauda and Rouquier who defined KLR algebras for arbitrary oriented quivers.


\begin{Definition}[Khovanov and Lauda~\cite{KhovLaud:diagI,KhovLaud:diagII} and Rouquier~\cite{Rouq:2KM}] \label{BK:def}
Suppose $\K$ is an integral ring and $n$ is a positive integer. The \textit{Khovanov-Lauda--Rouquier algebra}, $\mathscr{R}_n(\K)$ of type $\Gamma_e$ is the unital associative $\K$-algebra with generators
  $$\{ \psi_1,\dots, \psi_{n-1}\} \cup
         \{  y_1,\dots, y_n \} \cup \set{ e(\bi)|\bi\in (\Z/e\Z)^n}$$
  and relations
\begin{align*}
 e(\bi)  e(\bj) &= \delta_{\bi\bj} e(\bi),
&{\textstyle\sum_{\bi \in (\Z/e\Z)^n}}  e(\bi)&= 1,\\
 y_r  e(\bi) &=  e(\bi) y_r,
& \psi_r  e(\bi)&=  e(s_r{\cdot}\bi)  \psi_r,
& y_r  y_s &=  y_s  y_r,
\end{align*}
\vskip-21pt
\begin{align*}
 \psi_r  y_s  &=  y_s  \psi_r,&\text{if }s \neq r,r+1,\\
 \psi_r  \psi_s &=  \psi_s  \psi_r,&\text{if }|r-s|>1,
\end{align*}
\vskip-21pt
\begin{align*}
   \psi_r  y_{r+1}  e(\bi) &= \begin{cases}
      ( y_r \psi_r+1) e(\bi),\hspace*{18mm} &\text{if $i_r=i_{r+1}$},\\
     y_r \psi_r  e(\bi),&\text{if $i_r\neq i_{r+1}$}
  \end{cases}\\
   y_{r+1}  \psi_r  e(\bi) &= \begin{cases}
      ( \psi_r  y_r+1)  e(\bi),\hspace*{18mm} &\text{if $i_r=i_{r+1}$},\\
     \psi_r  y_r  e(\bi), &\text{if $i_r\neq i_{r+1}$}
  \end{cases}\\
   \psi_r^2  e(\bi) &= \begin{cases}
       0,&\text{if $i_r = i_{r+1}$},\\
       e(\bi),&\text{if $i_r \ne i_{r+1}\pm1$},\\
      ( y_{r+1}- y_r) e(\bi),&\text{if  $e\ne2$ and $i_{r+1}=i_r+1$},\\
       ( y_r -  y_{r+1}) e(\bi),&\text{if $e\ne2$ and  $i_{r+1}=i_r-1$},\\
      ( y_{r+1} -  y_{r})( y_{r}- y_{r+1})  e(\bi),&\text{if $e=2$ and
$i_{r+1}=i_r+1$}
\end{cases}\\
 \psi_{r} \psi_{r+1}  \psi_{r}  e(\bi) &= \begin{cases}
    ( \psi_{r+1}  \psi_{r}  \psi_{r+1} +1) e(\bi),\hspace*{7mm}
       &\text{if $e\ne2$ and $i_{r+2}=i_r=i_{r+1}-1$},\\
  ( \psi_{r+1}  \psi_{r}  \psi_{r+1} -1) e(\bi),
       &\text{if $e\ne2$ and $i_{r+2}=i_r=i_{r+1}+1$},\\
  \big( \psi_{r+1}  \psi_{r}  \psi_{r+1} + y_r\\
  \qquad -2 y_{r+1}+ y_{r+2}\big) e(\bi),
    &\text{if $e=2$ and $i_{r+2}=i_r=i_{r+1}+1$},\\
   \psi_{r+1}  \psi_{r}  \psi_{r+1}  e(\bi),&\text{otherwise.}
\end{cases}
\end{align*}
for $\bi,\bj\in (\Z/e\Z)^n$ and all admissible $r$ and $s$. Moreover, $\mathscr{R}_n(\mathscr O)$ is naturally
$\Z$-graded with degree function determined by
$$\deg  e(\bi)=0,\qquad \deg  y_r=2\qquad\text{and}\qquad \deg
   \psi_s  e(\bi)=-a_{i_s,i_{s+1}},$$
for $1\le r\le n$, $1\le s<n$ and $\bi\in (\Z/e\Z)^n$.
\end{Definition}

Fix a weight $\Lambda = \sum_{i\in \Z/e\Z} a_i \Lambda_i$ with $a_i \in \N$. Let $N_n^\Lambda(\K)$ be the two-sided ideal of $\mathscr{R}_n$ generated by the elements $e(\bi)y_1^{(\Lambda,\alpha_{i_1})}$, for $\bi \in (\Z/e\Z)^n$. We define the cyclotomic Khovanov-Lauda-Rouquier algebras, which were introduced by Khovanov and Lauda~\cite[Section 3.4]{KhovLaud:diagI}.

\begin{Definition}
The \textit{cyclotomic Khovanov-Lauda-Rouquier algebras} of weight $\Lambda$ and type $\Gamma_e$ is the algebra $\R(\K) = \mathscr{R}_n(\K)/ N_n^\Lambda(\K)$.
\end{Definition}

Brundan and Kleshchev~\cite{BK:GradedKL} proved the remarkable result that when $\K$ is a field of characteristic $p$, $\K\Sym_n \cong \R(\K)$ when we set $e = p$ and $\Lambda = \Lambda_k$ for any $k \in \Z/e\Z$.

\begin{Theorem} [Brundan-Kleshchev~\cite{BK:GradedKL}] \label{BK:iso}
Suppose $\K$ is a field of characteristic $p$ and $\R(\K)$ is the cyclotomic Khovanov-Lauda-Rouquier algebra over $\K$ with $e = p$ and $\Lambda = \Lambda_k$ for any $k \in \Z/e\Z$. Then $\K\Sym_n \cong \R(\K)$.
\end{Theorem}


Murphy~\cite{Murphy:basis} constructed the first cellular basis for $\K\Sym_n$ which shows that $\K\Sym_n$ is a cellular algebra. Hu-Mathas~\cite{HuMathas:GradedCellular} gave a graded cellular basis of $\K\Sym_n$ and prove that $\K\Sym_n$ is a graded cellular algebra. Next we introduce a graded cellular basis of $\K\Sym_n$.

Suppose $\Lambda = \Lambda_k$ for some $k \in \Z/e\Z$ and $\lambda \vdash n$ is a partition of $n$. For $\t \in \Std(\lambda)$, write $\t(\l) = (r_\l, c_\l)$ for $1 \leq \l \leq n$. We define the \textit{residue sequence} of $\t$ to be $\bi = (i_1, \ldots, i_n) \in (\Z/e\Z)^n$ where $i_\l \equiv k + c_\l - r_\l \pmod{e}$ for $1 \leq \l \leq n$.

\begin{Remark}
Note we have defined residue sequence for up-down tableau. A standard $\lambda$-tableau can be considered as a special case of up-down tableau with shape $(\lambda,0)$. For $\t \in \Std(\lambda)$, we have two residue sequence -- $\bi = (i_1, \ldots, i_n) \in (\Z/e\Z)^n$ by considering $\t$ as a standard tableau, and $\bj = (j_1, \ldots, j_n) \in P^n$ by considering $\t$ as an up-down tableau. Essentially these two definitions are equivalent when $e = p = 0$. One can see that $\bi$ is a "shift" of $\bj$. In more details, for $1 \leq \l \leq n$, we have $i_\l = j_\l - (\frac{\delta-1}{2} - k)$.
\end{Remark}

Define $\t^\lambda$ to be the unique standard $\lambda$-tableau such that $\t^\lambda \unrhd \t$ for all standard $\lambda$-tableau $\t$ and let $\bi_\lambda = (i_1, \ldots, i_n)$ to be the residue sequence of $\t^\lambda$ by considering it as a standard tableau. We define $e_\lambda = e(\bi_\lambda)$.

Suppose $w \in \Sym_n$ with reduced expression $w = s_{i_1}s_{i_2} \ldots s_{i_m}$. Define
$$
\psi_{w} = \psi_{i_1} \psi_{i_2} \ldots \psi_{i_\l}\in \R(\K) \hspace*{5mm} \text{and} \hspace*{5mm}\psi_{w}^* = \psi_{i_\l}\psi_{i_{\l-1}}\ldots \psi_{i_2}\psi_{i_1}\in \R(\K).
$$

For any standard tableau $\t$ with shape $\lambda$, we define $d(\t) \in \Sym_n$ such that $\t = \t^\lambda d(\t)$.

\begin{Definition} \label{def:basis of Sym}
Suppose $\lambda \vdash n$ and $\s,\t \in \Std(\lambda)$. Define
$$
\psi_{\s\t} = \psi_{d(\s)}^* e_\lambda \psi_{d(\t)} \in \R(\K).
$$
\end{Definition}

\begin{Theorem}[\protect{Hu-Mathas~\cite[Theorem 5.14]{HuMathas:GradedCellular}}] \label{HM:basis of Sym}
Suppose $\K$ is a field. Then
$$
\set{ \psi_{\s\t} | \s,\t \in \Std(\lambda) \text{ for $\lambda \vdash n$}}
$$
is a graded cellular basis of $\R(\K)$.
\end{Theorem}

\begin{Remark}
We note that all the results of this subsection were originally proved in the cyclotomic Hecke algebras of type $A$, $\mathscr H_n^\Lambda(\K)$, with weight $\Lambda = \sum_{i\in \Z/e\Z} a_i \Lambda_i \in \sum_{i\in \Z/e\Z}\N \Lambda_i$ instead of the symmetric group algebra, $\K\Sym_n$. In this paper we only need the results for $\K\Sym_n$. So we restrict all the results to $\K\Sym_n$.
\end{Remark}

\section{The graded algebras $\G{n}$} \label{sec:algebra}

Let $x$ be an invariant and recall $\O = R[x]_{(x - \delta)} = R[[x-\delta]]$, and $\F = R(x)$. In this section we define a new algebra $\G{n}$ over $R$ associated with KLR-like relations, which is naturally $\Z$-graded.

\subsection{A categorification of $n$-tuple $\bi \in P^n$} \label{sec:bi}

Recall $P = \frac{\delta-1}{2} + \Z$. In this subsection we define a mapping $h_k\map{P^n}{\Z}$ which separates $P^n$ into three mutually exclusive subsets. Such categorification will be used to determine the relations and the degree of generators of the graded algebra $\G{n}$.

First we give a proper definition of $h_k$.

\begin{Definition}
Suppose $\bi = (i_1, i_2, \ldots, i_n) \in P^n$ and $k$ is an integer with $1 \leq k \leq n$. We define
\begin{eqnarray*}
h_k(\bi) & := & \delta_{i_k, -\frac{\delta-1}{2}} + \#\set{1 \leq r \leq k-1 | i_r = -i_k \pm 1} + 2\#\set{1 \leq r \leq k-1 | i_r = i_k}\\
&& - \delta_{i_k, \frac{\delta-1}{2}} - \#\set{1 \leq r \leq k-1 | i_r = i_k \pm 1} - 2\#\set{1 \leq r \leq k-1 | i_r = -i_k}.
\end{eqnarray*}
\end{Definition}

\begin{Remark}
Suppose $\bi = (i_1, \ldots, i_n)$ and $1 \leq k \leq n-1$. If $i_{k+1} = i_k$, we have
\begin{equation} \label{remark:h:eq1}
h_{k+1}(\bi) =
\begin{cases}
h_k(\bi), & \text{if $i_k = 0$,}\\
h_k(\bi) + 3, & \text{if $i_k = \pm \frac{1}{2}$,}\\
h_k(\bi) + 2, & \text{otherwise;}
\end{cases}
\end{equation}
and if $i_{k+1} = -i_k$, we have
\begin{equation} \label{remark:h:eq2}
h_{k+1}(\bi) =
\begin{cases}
- h_k(\bi), & \text{if $i_k = 0$,}\\
- h_k(\bi) - 3, & \text{if $i_k = \pm \frac{1}{2}$,}\\
- h_k(\bi) - 2, & \text{otherwise.}
\end{cases}
\end{equation}

We will use this result frequently in the rest of this paper.
\end{Remark}

Given $(\lambda,f) \in \widehat B_{k-1}$, the key point of $h_k$ is that it gives us a way to understand the structure of $\mathscr{AR}(\lambda)$. In order to connect $h_k$ and $\mathscr{AR}(\lambda)$, we introduce a Lemma which is first proved by Nazarov~\cite{Nazarov:brauer}.


\begin{Lemma} [\protect{Nazarov~\cite[Lemma 3.8]{Nazarov:brauer}}] \label{Naz}
Suppose $u$ is a unknown and $\t$ is an up-down tableau of size $n$. If for $1 \leq k \leq n-1$ we have $\t(k) + \t(k+1) = 0$, then
\begin{equation} \label{eq:Naz}
\frac{u + (x - 1)/2}{u - (x - 1)/2} \prod_{r = 1}^{k-1} \frac{(u + c_\t(r))^2 - 1}{(u - c_\t(r))^2 - 1} \frac{(u - c_\t(r))^2}{(u + c_\t(r))^2} = \sum_{\s \overset{k}\sim \t} \frac{u + c_\s(k)}{u - c_\s(k)}.
\end{equation}
\end{Lemma}

Suppose $\lambda$ is a partition and $\alpha \in \mathscr{AR}(\lambda)$. Define
$$
\res_\lambda(\alpha) =
\begin{cases}
\res(\alpha), & \text{if $\alpha \in \mathscr A(\lambda)$,}\\
-\res(\alpha), & \text{if $\alpha \in \mathscr R(\lambda)$,}
\end{cases}
$$
and for $i \in P$, we denote $\mathscr{AR}_\lambda(i) = \set{ \alpha \in \mathscr{AR}(\lambda) | \res_\lambda(\alpha) = i}$.

The next Lemma gives the most important property of $h_k$.

\begin{Lemma} \label{deg:h1:1}
For any $\bi \in P^n$ such that $\bi|_{k-1}$ is the residue sequence of some up-down tableaux with shape $(\lambda,f)$, we have $h_k(\bi) = |\mathscr{AR}_\lambda(-i_k)| - |\mathscr{AR}_\lambda(i_k)|$.
\end{Lemma}

\begin{proof}
Suppose $\u \in \Tud_{k-1}(\bi|_{k-1})$ with shape $(\lambda,f)$. Choose any up-down tableau $\t$ of size $k+1$ with $\t|_{k-1} = \u$ and $\t(k) + \t(k+1) = 0$. By the construction of up-down tableaux, there exists such $\t$ as long as $|\mathscr{AR}(\lambda)| > 0$, which is always true.

Substitute $u = i$ into (\ref{eq:Naz}). Then we have
$$
\frac{i + (x - 1)/2}{i - (x - 1)/2} \prod_{r = 1}^{k-1} \frac{(i + c_\t(r))^2 - 1}{(i - c_\t(r))^2 - 1} \frac{(i - c_\t(r))^2}{(i + c_\t(r))^2} = \prod_{\s \overset{k}\sim \t} \frac{i + c_\s(k)}{i - c_\s(k)} \in \F.
$$

For convenience, write $d = \frac{x - \delta}{2}$. Then we have
$$
\prod_{\s \overset{k} \sim \t} \frac{i + c_\s(k)}{i - c_\s(k)} = \frac{d^{|\mathscr{AR}_\lambda(-i)|}}{d^{|\mathscr{AR}_\lambda(i)|}} v_1,
$$
for some $v_1$ invertible in $\O$ and
\begin{eqnarray*}
&& \frac{i + (x - 1)/2}{i - (x - 1)/2} \prod_{r = 1}^{k-1} \frac{(i + c_\t(r))^2 - 1}{(i - c_\t(r))^2 - 1} \frac{(i - c_\t(r))^2}{(i + c_\t(r))^2} = \frac{d^{\delta_{i, -\frac{\delta-1}{2}}}}{d^{\delta_{i, \frac{\delta-1}{2}}}} \frac{d^{\#\set{1 \leq r \leq k-1 | i_r = -i \pm 1}}}{d^{\#\set{1 \leq r \leq k-1 | i_r = i \pm 1}}} \frac{d^{2\#\set{1 \leq r \leq k-1 | i_r = i}}}{d^{2\#\set{1 \leq r \leq k-1 | i_r = -i}}} v_2,
\end{eqnarray*}
for some $v_2$ invertible in $\O$.

Hence,
$$
\frac{d^{|\mathscr{AR}_\lambda(-i)|}}{d^{|\mathscr{AR}_\lambda(i)|}} v_1 = \frac{d^{\delta_{i, -\frac{\delta-1}{2}}}}{d^{\delta_{i, \frac{\delta-1}{2}}}} \frac{d^{\#\set{1 \leq r \leq k-1 | i_r = -i \pm 1}}}{d^{\#\set{1 \leq r \leq k-1 | i_r = i \pm 1}}} \frac{d^{2\#\set{1 \leq r \leq k-1 | i_r = i}}}{d^{2\#\set{1 \leq r \leq k-1 | i_r = -i}}} v_2
$$
where $v_1,v_2$ invertible in $\O$. Because $d$ is not invertible in $\O$, we have
\begin{eqnarray*}
|\mathscr{AR}_\lambda(-i)| - |\mathscr{AR}_\lambda(i)| & = & \delta_{i, -\frac{\delta - 1}{2}} + \#\set{1 \leq r \leq k-1 | i_r = -i \pm 1} + 2\#\set{1 \leq r \leq k-1 | i_r = i}\\
&& - \ \delta_{i, \frac{\delta - 1}{2}} - \#\set{1 \leq r \leq k-1 | i_r = i \pm 1} - 2\#\set{1 \leq r \leq k-1 | i_r = -i} = h_k(\bi),
\end{eqnarray*}
which completes the proof.
\end{proof}


The next Corollary is a special case of~\autoref{deg:h1:1}, which shows the connection between $h_k(\bi_\t)$ and $\t$.

\begin{Corollary} \label{deg:h1}
Suppose $\t$ is an up-down tableau of size $n$ and $\bi_\t$ is the residue sequence of $\t$. For $1 \leq k \leq n$, let $\t_{k-1} = \lambda$. Then we have $h_k(\bi_\t) = |\mathscr{AR}_\lambda(-i_k)| - |\mathscr{AR}_\lambda(i_k)|$.
\end{Corollary}

The first application of~\autoref{deg:h1:1} is that when $\bi$ is a residue sequence of some up-down tableaux, the value of $h_k(\bi)$ is bounded.

\begin{Lemma} \label{deg:h2}
Suppose $\bi \in P^n$ and $1 \leq k \leq n$. If $\bi$ is the residue sequence of an up-down tableau, we have $h_k(\bi) \in \{-2, -1, 0\}$.
\end{Lemma}

\begin{proof}
Suppose $\t$ is an up-down tableau with residue sequence $\bi$. Write $\lambda = \t_{k-1}$. By~\autoref{deg:h1:1} we have$|\mathscr{AR}_\lambda(-i_k)| - |\mathscr{AR}_\lambda(i_k)| = h_k(\bi)$.

The existence of $\t$ implies $|\mathscr{AR}_\lambda(i_k)| \geq 1$. By the construction of partitions, we have
\begin{equation} \label{deg:h2:ineq}
 0 \leq |\mathscr{AR}_\lambda(-i_k)| \leq 2, \qquad 0 \leq |\mathscr{AR}_\lambda(i_k)| \leq 2, \qquad\text{and} \qquad 0 \leq |\mathscr{AR}_\lambda(-i_k)| + |\mathscr{AR}_\lambda(i_k)| \leq 2.
\end{equation}

The Lemma follows easily by direct calculations.
\end{proof}

\begin{Lemma} \label{deg:h3}
Suppose $\bi \in P^n$ and $1 \leq k \leq n$. If $\bi$ is the residue sequence of an up-down tableau, we have $h_k(\bi) = 0$ if $i_k = 0$ and $h_k(\bi) \in \{-1, -2\}$ if $i_k = \pm\frac{1}{2}$.
\end{Lemma}

\begin{proof}
Suppose $i_k = 0$. As $i_k = -i_k$, by the definition of $h_k(\bi)$, we have $h_k(\bi) = -h_k(\bi) = 0$. Suppose $i_k = -\frac{1}{2}$ and set $\lambda = \t_{k-1}$. By the construction of $\lambda$, we have $|\mathscr{AR}_\lambda(-i_k)| = 0$ and $|\mathscr{AR}_\lambda(i_k)| \geq 1$, which implies that $h_k(\bi) \leq -1$ by~\autoref{deg:h1:1}. Hence $h_k(\bi) \in \{-1,-2\}$. For $i_k = \frac{1}{2}$ we have the same result following the same argument.
\end{proof}

\autoref{deg:h2} and~\autoref{deg:h3} give us an easy way to test whether $\bi \in P^n$ is the residue sequence of an up-down tableau. If we have $h_k(\bi) \not\in \{-2,-1,0\}$ for some $1 \leq k \leq n$, then $\bi$ is not the residue sequence of an up-down tableau. But the reverse is not always valid. We will discuss this problem further in Section~\ref{sec:residue}.

The next important application of $h_k$ is by giving $\t \in \Tud_n(\bi)$ and $\lambda = \t_{k-1}$ for $1 \leq k \leq n$, we know the exact values of $|\mathscr{AR}_\lambda(-i_k)|$ and $|\mathscr{AR}_\lambda(i_k)|$ by knowing $h_k(\bi)$. Because $\bi$ is the residue sequence of $\t$, we have $|\mathscr{AR}_\lambda(i_k)| \geq 1$. Then by~\autoref{deg:h1:1} and (\ref{deg:h2:ineq}), the following results are straightforward:
\begin{align}
& |\mathscr{AR}_\lambda(-i_k)| = 0, \quad |\mathscr{AR}_\lambda(i_k)| = 2 & &\text{if $h_k(\bi) = -2$,}\label{deg:h4:eq1}\\
& |\mathscr{AR}_\lambda(-i_k)| = 0, \quad |\mathscr{AR}_\lambda(i_k)| = 1 & &\text{if $h_k(\bi) = -1$,}\label{deg:h4:eq2}\\
& |\mathscr{AR}_\lambda(-i_k)| = 1, \quad |\mathscr{AR}_\lambda(i_k)| = 1 & &\text{if $h_k(\bi) = 0$,}\label{deg:h4:eq3}
\end{align}

The following results are implied by (\ref{deg:h4:eq1}) - (\ref{deg:h4:eq3}), which can be used to determine the structure of $\t$. These results will be used frequently in the rest of this paper.

\begin{Lemma} \label{deg:h4}
Suppose $\t \in \Tud_n(\bi)$. For $1 \leq k \leq n$, write $\lambda = \t_{k-1}$. Then we have the following properties:
\begin{enumerate}
\item When $h_k(\bi) = -2$, then $\mathscr{AR}_\lambda(i_k) = \{\alpha, \beta\}$ where $\alpha \in \mathscr A(\lambda)$ and $\beta \in \mathscr R(\lambda)$.

\item When $h_k(\bi) = 0$, then $\mathscr{AR}_\lambda(i_k) = \{\alpha\}$ and $\mathscr{AR}_\lambda(-i_k) = \{\beta\}$, where either $\alpha, \beta \in \mathscr A(\lambda)$ or $\alpha,\beta \in \mathscr R(\lambda)$.
\end{enumerate}
\end{Lemma}

\proof
(1). When $h_k(\bi) = -2$, by (\ref{deg:h4:eq1}) we have $\mathscr{AR}_\lambda(i_k) = \{\alpha, \beta\}$. Suppose $\alpha \in \mathscr A(\lambda)$. If $\beta \in \mathscr A(\lambda)$, then we have $\alpha, \beta \in \mathscr A(\lambda)$ such that $\res(\alpha) = \res(\beta) = i_k$. But for any $\lambda$, there exists at most one addable node with residue $i_k$. Hence we must have $\beta \in \mathscr R(\lambda)$. Suppose $\alpha \in \mathscr R(\lambda)$. If $\beta \in \mathscr R(\lambda)$, then we have $\alpha, \beta \in \mathscr R(\lambda)$ such that $\res(\alpha) = \res(\beta) = -i_k$. But for any $\lambda$, there exists at most one removable node with residue $-i_k$. Hence we must have $\beta \in \mathscr A(\lambda)$. Therefore part (1) follows.

(2). When $h_k(\bi) = 0$, by (\ref{deg:h4:eq3}) we have $\mathscr{AR}_\lambda(i_k) = \{\alpha\}$ and $\mathscr{AR}_\lambda(-i_k) = \{\beta\}$. Suppose $\alpha \in \mathscr A(\lambda)$. If $\beta \in \mathscr R(\lambda)$, we have $\res(\alpha) = \res(\beta) = i_k$. But for any $\lambda$, if there exists an addable node with residue $i_k$, there does not exist a removable node with residue $i_k$. Hence we must have $\beta \in \mathscr A(\lambda)$. Suppose $\alpha \in \mathscr R(\lambda)$. If $\beta \in \mathscr A(\lambda)$, we have $\res(\alpha) = \res(\beta) = -i_k$. But for any $\lambda$, if there exists an addable node with residue $-i_k$, there does not exist a removable node with residue $-i_k$. Hence we must have $\beta \in \mathscr R(\lambda)$. Therefore part (2) follows. \endproof

\begin{Lemma} \label{deg:h4:h}
Suppose $\bi \in P^n$ with $i_k = i_{k+1}$ for $1 \leq k \leq n-1$. For $\t \in \Tud_n(\bi)$, we have $\t(k) > 0$, $\t(k+1) < 0$ or $\t(k) < 0$, $\t(k+1) > 0$. Moreover, we have $\t(k) + \t(k+1) \neq 0$ if and only if $i_k = 0$.
\end{Lemma}

\proof Because of \eqref{remark:h:eq1} and~\autoref{deg:h2}, it forces $i_k \neq \pm\frac{1}{2}$, and
$$
h_k(\bi) =
\begin{cases}
0, & \text{if $i_k = 0$,}\\
-2, & \text{if $i_k \neq 0$,}
\end{cases}
\qquad
\text{and $h_{k+1}(\bi) = 0$.}
$$

Write $\lambda = \t_{k-1}$. Assume $i_k = 0$. By~\eqref{deg:h4:eq3} we have $|\mathscr{AR}_\lambda(i_k)| = 1$. Hence by the construction of up-down tableaux, we require $\t(k) + \t(k+1) = 0$, which also implies that $\t(k) > 0$, $\t(k+1) < 0$ or $\t(k) < 0$, $\t(k+1) > 0$.

Assume $i_k \neq 0$. By~\eqref{deg:h4:eq1} we have $|\mathscr{AR}_\lambda(i_k)| = 2$. Let $\alpha, \beta \in \mathscr{AR}_\lambda(i_k)$ be distinct nodes. By~\autoref{deg:h4}, we set $\alpha \in \mathscr A(\lambda)$ and $\beta \in \mathscr R(\lambda)$. By the construction of up-down tableaux, we require $\t(k) = \alpha$, $\t(k+1) = -\beta$ or $\t(k) = -\beta$, $\t(k+1) = \alpha$. Henceforth, we have $\t(k) > 0$, $\t(k+1) < 0$ or $\t(k) < 0$, $\t(k+1) > 0$. Because $\alpha$ and $\beta$ are distinct, we have $\t(k) + \t(k+1) \neq 0$, which completes the proof. \endproof

\begin{Lemma} \label{deg:h4:1}
Suppose $\t \in \Tud_n(\bi)$. If $\t_{k-1} = \t_{k+1} = \lambda$ for $1 \leq k \leq n-1$, we have the following properties:
\begin{enumerate}
\item $h_k(\bi) = -2$ if and only if there exists an unique $\s \neq \t$ such that $\s \overset{k}\sim \t$ and $\s \in \Tud_n(\bi)$. Moreover, we have $c_\t(k) - i_k = -(c_\s(k) - i_{k+1})$.

\item $h_k(\bi) = 0$ if and only if there exists an unique $\s$ such that $\s \overset{k}\sim \t$ and $\s \in \Tud_n(\bi{\cdot}s_k)$. Moreover, we have $c_\t(k) - i_k = c_\s(k) - i_k$.

\item $h_k(\bi_\t) = -1$ if and only if $\s \in \Tud_n(\bi) \cup \Tud_n(\bi{\cdot}s_k)$ and $\s \overset{k}\sim \t$ implies $\s = \t$.
\end{enumerate}
\end{Lemma}

\begin{proof}
In (1), by (\ref{deg:h4:eq1}) we have $|\mathscr{AR}_\lambda(-i_k)| = 0$ and $|\mathscr{AR}_\lambda(i_k)| = 2$ if and only if $h_k(\bi) = -2$. Hence by~\autoref{e_k:h1}, $|\mathscr{AR}_\lambda(i_k)| = 2$ if and only if there exist exactly two distinct up-down tableau $\u,\v \in \Tud_n(\bi_\t)$ such that $\u \overset{k}\sim \t$ and $\v \overset{k}\sim \t$. It is obvious that one of $\u$ and $\v$ is $\t$. It is Without loss of generality, we set $\u = \t$. Hence by setting $\s = \v$, the uniqueness and existence of $\s$ follows.

Moreover, by~\autoref{deg:h4} we have $\mathscr{AR}_\lambda(i_k) = \{\alpha, \beta\}$ where $\alpha \in \mathscr A(\lambda)$ and $\beta \in \mathscr R(\lambda)$. Then we have $\t(k) = \alpha$ and $\s(k) = -\beta$, or vice versa. In both cases, we have $c_\t(k) - i_k = -(c_\s(k) - i_{k+1})$, which proves part (1).

Using similar arguments, (2) can be implied by (\ref{deg:h4:eq3}) and~\autoref{deg:h4}; and (3) can be implied by (\ref{deg:h4:eq2}).
\end{proof}

\begin{Lemma} \label{deg:h4:2}
Suppose $\t \in \Tud_n(\bi)$. If $i_k + i_{k+1} = 0$ for $1 \leq k \leq n-1$, we have the following properties:
\begin{enumerate}
\item When $h_{k+1}(\bi) = 0$ or $-1$, then $\t(k) + \t(k+1) = 0$.

\item When $h_{k+1}(\bi) = -2$, then either $\t(k) + \t(k+1) = 0$, or $c_\t(k) - i_k = c_\t(k+1) - i_{k+1}$.
\end{enumerate}
\end{Lemma}

\proof Suppose $\t_k = \lambda$ and $\t(k) = \alpha$. Without loss of generality, we assume $\alpha > 0$. When $\alpha < 0$ the Lemma follows by the same argument.

(1). When $h_{k+1}(\bi) = 0$ or $-1$, by (\ref{deg:h4:eq2}) and (\ref{deg:h4:eq3}) we have $|\mathscr{AR}_\lambda(i_{k+1})| = 1$. As $\alpha \in \mathscr R(\lambda)$ and $\res_\lambda(\alpha) = -\res(\alpha) = -i_k = i_{k+1}$, we have $\mathscr{AR}_\lambda(i_{k+1}) = \{\alpha\}$. Hence, it forces $\t(k+1) = -\alpha = -\t(k)$.

(2). When $h_{k+1}(\bi) = -2$, by (\ref{deg:h4:eq1}) we have $|\mathscr{AR}_\lambda(i_{k+1})| = 2$. For the same reason as above, we have $\alpha \in \mathscr{AR}_\lambda(i_{k+1})$. Hence we have $\mathscr{AR}_\lambda(i_{k+1}) = \{\alpha,\beta\}$. By~\autoref{deg:h4}, we have $\alpha \in \mathscr R(\lambda)$ and $\beta \in \mathscr A(\lambda)$. Therefore $\t(k+1) = -\alpha$ or $\beta$. If $\t(k+1) = -\alpha = -\t(k)$, we have $\t(k) + \t(k+1) = 0$; and if $\t(k+1) = \beta$, we have $\t(k) + \t(k+1) \neq 0$ and $c_\t(k) - i_k = \frac{x-\delta}{2} = c_\t(k+1) - i_{k+1}$. \endproof

%
%
%
%
%
%


We now categorize $P^n$ using $h_k$. For $1 \leq k \leq n$, define $P_{k,+}^n$, $P_{k,-}^n$ and $P_{k,0}^n$ as subsets of $P^n$ by
\begin{align*}
P_{k,+}^n & := \set{\bi \in P^n | \text{$i_k \neq 0, -\frac{1}{2}$ and $h_k(\bi) = 0$, or $i_k = -\frac{1}{2}$ and $h_k(\bi) = -1$}}.\\
P_{k,-}^n & := \set{\bi \in P^n | \text{$i_k \neq 0, -\frac{1}{2}$ and $h_k(\bi) = -2$, or $i_k = -\frac{1}{2}$ and $h_k(\bi) = -3$}},\\
P_{k,0}^n & := P^n \backslash (P_{k,+}^n \cup P_{k,-}^n).
\end{align*}

We split $P^n$ into three mutually exclusive subsets, i.e. $P^n = P_{k,+}^n \sqcup P_{k,-}^n \sqcup P_{k,0}^n$. Let $I^n$ be the set containing all the residue sequences of up-down tableaux of size $n$. For $1 \leq k \leq n$, define
$$
I_{k,a}^n := \set{\bi \in P_{k,a}^n | \text{$\bi$ is a residue sequence of some up-down tableaux}}
$$
where $a \in \{+,-,0\}$. It is easy to see that $I^n = I_{k,+}^n \sqcup I_{k,-}^n \sqcup I_{k,0}^n$. By the definitions of $P_{k,+}^n$, $P_{k,-}^n$ and $P_{k,0}^n$,~\autoref{deg:h2} and~\autoref{deg:h3} imply that
\begin{align*}
I_{k,+}^n & = \set{\bi \in I^n | \text{$i_k \neq 0, -\frac{1}{2}$ and $h_k(\bi) = 0$, or $i_k = -\frac{1}{2}$ and $h_k(\bi) = -1$}}.\\
I_{k,-}^n & = \set{\bi \in I^n | \text{$i_k \neq 0, -\frac{1}{2}$ and $h_k(\bi) = -2$}}.\\
I_{k,0}^n & = \set{\bi \in I^n | \text{$i_k \neq 0, -\frac{1}{2}$ and $h_k(\bi) = -1$, or $i_k = -\frac{1}{2}$ and $h_k(\bi) = -2$, or $i_k = 0$}}.
\end{align*}

%

We will give a further explanation of this categorification after we construct our graded algebra $\G{n}$.


\subsection{Graded algebras $\G{n}$} \label{sec:A}

In this subsection we construct a naturally $\Z$-graded algebra $\G{n}$ over $R$ and introduce some of its properties. For $\bi \in P^n$ and $1 \leq k \leq n-1$, define $a_k(\bi) \in \Z$ and $A_{k,1}^\bi, A_{k,2}^\bi, A_{k,3}^\bi, A_{k,4}^\bi \in \{1,2,\ldots,k-1\}$ by
$$
a_k(\bi) =
\begin{cases}
\#\set{1 \leq r \leq k-1 | i_r \in \{-1,1\} } + 1 + \delta_{\frac{i_k - i_{k+1}}{2},\frac{\delta - 1}{2}} , & \text{if $\frac{i_k - i_{k+1}}{2} = 0$,}\\
\#\set{1 \leq r \leq k-1 | i_r \in \{-1,1\} } + \delta_{\frac{i_k - i_{k+1}}{2},\frac{\delta - 1}{2}}, & \text{if $\frac{i_k - i_{k+1}}{2} = 1$,}\\
\delta_{\frac{i_k - i_{k+1}}{2},\frac{\delta - 1}{2}}, & \text{if $\frac{i_k - i_{k+1}}{2} = 1/2$,}\\
\#\set{1 \leq r \leq k-1 | i_r \in \{\frac{i_k - i_{k+1}}{2}, \frac{i_k - i_{k+1}}{2} - 1, -\frac{i_k - i_{k+1}}{2}, -\frac{i_k - i_{k+1}}{2} + 1\} } + \delta_{\frac{i_k - i_{k+1}}{2},\frac{\delta - 1}{2}}, & \text{otherwise,}
\end{cases}
$$
and
\begin{align*}
A_{k,1}^\bi & := \set{1 \leq r \leq k-1 | i_r = -i_k \pm 1}, &
A_{k,2}^\bi & := \set{1 \leq r \leq k-1 | i_r = i_k},\\
A_{k,3}^\bi & := \set{1 \leq r \leq k-1 | i_r = i_k \pm 1}, &
A_{k,4}^\bi & := \set{1 \leq r \leq k-1 | i_r = -i_k};
\end{align*}
and for $\bi \in P_{k,0}^n$ and $1 \leq k \leq n-1$, define $z_k(\bi) \in \Z$ by
$$
z_k(\bi) =
\begin{cases}
0, & \text{if $h_k(\bi) < -2$, or $h_k(\bi) \geq 0$ and $i_k \neq 0$,}\\
(-1)^{a_k(\bi)} (1 + \delta_{i_k, -\frac{1}{2}}), & \text{if $-2 \leq h_k(\bi) < 0$,}\\
\frac{1 + (-1)^{a_k(\bi)}}{2}, & \text{if $i_k = 0$.}
\end{cases}
$$

Let $\G{n}$ be an unital associate $R$-algebra with generators
$$
G_n(\delta) = \set{e(\bi)|\bi \in P^n} \cup \set{y_k| 1 \leq k \leq n} \cup \set{\psi_k| 1\leq k \leq n-1} \cup \set{\epsilon_k| 1\leq k \leq n-1}
$$
associated with the following relations:

\begin{enumerate}
\item Idempotent relations:
Let $\bi, \bj \in P^n$ and $1 \leq k \leq n-1$. Then
\begin{equation} \label{rela:1}
y_1^{\delta_{i_1, \frac{\delta-1}{2}}} e(\bi) = 0, \qquad \sum_{\bi \in P^n} e(\bi) = 1,\qquad e(\bi) e(\bj) = \delta_{\bi,\bj} e(\bi), \qquad e(\bi) \epsilon_k = 0 \text{  if $i_k + i_{k+1} \neq 0$;}
\end{equation}

\item Commutation relations:
Let $\bi \in P^n$. Then
\begin{align}
y_k e(\bi) = e(\bi) y_k, &\hspace*{10mm} \ \psi_k e(\bi) = e(\bi{\cdot}s_k) \psi_k &&&&\text{and} \label{rela:2:1} \\
y_k y_r = y_r y_k, &\hspace*{10mm} \ y_k \psi_r = \psi_r y_k, && \hspace*{-5mm} y_k \epsilon_r = \epsilon_r y_k, &&& \label{rela:2:2}\\
\psi_k \psi_r = \psi_r \psi_k, &\hspace*{10mm} \ \psi_k \epsilon_r = \epsilon_r \psi_k, && \hspace*{-5mm} \epsilon_k \epsilon_r = \epsilon_r \epsilon_k &&\text{if $|k - r| > 1$;} \label{rela:2:3}
\end{align}

\item Essential commutation relations:
Let $\bi \in P^n$ and $1 \leq k \leq n-1$. Then
\begin{align}
e(\bi) y_k \psi_k & = e(\bi) \psi_k y_{k+1} + e(\bi) \epsilon_k e(\bi{\cdot}s_k) - \delta_{i_k, i_{k+1}} e(\bi), \label{rela:3:1}\\
\text{and} \qquad e(\bi) \psi_k y_k & = e(\bi) y_{k+1} \psi_k + e(\bi) \epsilon_k e(\bi{\cdot}s_k) - \delta_{i_k, i_{k+1}} e(\bi). \label{rela:3:2}
\end{align}

\item Inverse relations:
Let $\bi \in P^n$ and $1 \leq k \leq n-1$. Then
\begin{equation} \label{rela:4}
e(\bi) \psi_k^2 =
\begin{cases}
0, & \text{if $i_k = i_{k+1}$ or $i_k + i_{k+1} = 0$ and $h_k(\bi) \neq 0$,}\\
(y_k - y_{k+1}) e(\bi), & \text{if $i_k = i_{k+1} + 1$ and $i_k + i_{k+1} \neq 0$,}\\
(y_{k+1} - y_k) e(\bi), & \text{if $i_k = i_{k+1} - 1$ and $i_k + i_{k+1} \neq 0$,}\\
e(\bi), & \text{otherwise;}
\end{cases}
\end{equation}

\item Essential idempotent relations:
Let $\bi, \bj, \bk \in P^n$ and $1 \leq k \leq n-1$. Then
\begin{align}
e(\bi) \epsilon_k e(\bi) & =
\begin{cases}
(-1)^{a_k(\bi)} e(\bi), & \text{if $\bi \in P_{k,0}^n$ and $i_k = -i_{k+1} \neq \pm\frac{1}{2}$,} \\
(-1)^{a_k(\bi) + 1} (y_{k+1} - y_k) e(\bi), & \text{if $\bi \in P_{k,+}^n$;}
\end{cases} \label{rela:5:1}\\
y_{k+1} e(\bi) & = y_k e(\bi) - 2(-1)^{a_k(\bi)}y_k e(\bi) \epsilon_k e(\bi)\\
& = y_k e(\bi) - 2(-1)^{a_k(\bi)} e(\bi) \epsilon_k e(\bi) y_k, \qquad\qquad\qquad\qquad\qquad \text{if $\bi \in P_{k,0}^n$ and $i_k = -i_{k+1} = \frac{1}{2}$,} \label{rela:5:2}\\
e(\bi) & = (-1)^{a_k(\bi)}e(\bi) \epsilon_k e(\bi) - 2(-1)^{a_{k-1}(\bi)} e(\bi) \epsilon_{k-1} e(\bi) \notag \\
&\hspace*{1cm}+ e(\bi) \epsilon_{k-1} \epsilon_k e(\bi) + e(\bi) \epsilon_k \epsilon_{k-1} e(\bi), \qquad\qquad \text{if $\bi \in P_{k,0}^n$ and $-i_{k-1} = i_k = -i_{k+1} = -\frac{1}{2}$,} \label{rela:5:3}\\
e(\bi) & = (-1)^{a_k(\bi)} e(\bi) (\epsilon_k y_k + y_k \epsilon_k) e(\bi), \qquad\qquad\qquad \ \ \ \ \  \text{if $\bi \in P_{k,-}^n$ and $i_k = -i_{k+1}$,} \label{rela:5:4}
\end{align}
\begin{align}
e(\bj) \epsilon_k e(\bi) \epsilon_k e(\bk) & =
\begin{cases}
z_k(\bi) e(\bj) \epsilon_k e(\bk), & \text{if $\bi \in P_{k,0}^n$,}\\
0, & \text{if $\bi \in P_{k,-}^n$,}\\
(-1)^{a_k(\bi)} (1 + \delta_{i_k, -\frac{1}{2}})(\sum_{r \in A_{k,1}^\bi} y_r - 2 \sum_{r \in A_{k,2}^\bi} y_r, &\\
\hspace*{2cm} + \sum_{r \in A_{k,3}^\bi} y_r - 2\sum_{r \in A_{k,4}^\bi} y_r) e(\bj) \epsilon_k e(\bk), & \text{if $\bi \in P_{k,+}^n$;}
\end{cases} \label{rela:5:5}
\end{align}

\item Untwist relations:
Let $\bi, \bj \in P^n$ and $1 \leq k \leq n-1$. Then
\begin{align}
e(\bi) \psi_k \epsilon_k e(\bj) & =
\begin{cases}
(-1)^{a_k(\bi)} e(\bi) \epsilon_k e(\bj), & \text{if $\bi \in P_{k,+}^n$ and $i_k \neq 0, -\frac{1}{2}$,}\\
0, & \text{otherwise;}
\end{cases} \label{rela:6:1}\\
e(\bj) \epsilon_k \psi_k e(\bi) & = \begin{cases}
(-1)^{a_k(\bi)} e(\bj) \epsilon_k e(\bi), & \text{if $\bi \in P_{k,+}^n$ and $i_k \neq 0, -\frac{1}{2}$,}\\
0, & \text{otherwise;}
\end{cases} \label{rela:6:2}
\end{align}

\item Tangle relations:
Let $\bi, \bj \in P^n$ and $1 < k < n$. Then
\begin{align}
& e(\bj) \epsilon_k \epsilon_{k-1} \psi_k e(\bi) = e(\bj) \epsilon_k \psi_{k-1} e(\bi), &&
e(\bi) \psi_k \epsilon_{k-1} \epsilon_k e(\bj) = e(\bi) \psi_{k-1} \epsilon_k e(\bj),& \label{rela:7:1}\\
& e(\bi) \epsilon_k \epsilon_{k-1} \epsilon_k e(\bj) = e(\bi) \epsilon_k e(\bj); &&
e(\bi) \epsilon_{k-1} \epsilon_k \epsilon_{k-1} e(\bj) = e(\bi) \epsilon_{k-1} e(\bj); &
e(\bi) \epsilon_k e(\bj) (y_k + y_{k+1}) = 0; \label{rela:7:2}
\end{align}

\item Braid relations:
Let $\mathcal B_k = \psi_k \psi_{k-1} \psi_k - \psi_{k-1} \psi_k \psi_{k-1}$, $\bi \in P^n$ and $1 < k < n$. Then
\begin{numcases}{e(\bi) \mathcal B_k =}
  e(\bi) \epsilon_k \epsilon_{k-1} e(\bi{\cdot}s_k s_{k-1} s_k), &if $i_k + i_{k+1} = 0$ and $i_{k-1} = \pm (i_k - 1)$,\label{rela:8:1}\\
- e(\bi) \epsilon_k \epsilon_{k-1} e(\bi{\cdot}s_k s_{k-1} s_k), &if $i_k + i_{k+1} = 0$ and $i_{k-1} = \pm (i_k + 1)$,\label{rela:8:2}\\
  e(\bi) \epsilon_{k-1} \epsilon_k e(\bi{\cdot}s_k s_{k-1} s_k), &if $i_{k-1} + i_k = 0$ and $i_{k+1} = \pm (i_k - 1)$,\label{rela:8:3}\\
- e(\bi) \epsilon_{k-1} \epsilon_k e(\bi{\cdot}s_k s_{k-1} s_k), &if $i_{k-1} + i_k = 0$ and $i_{k+1} = \pm (i_k + 1)$,\label{rela:8:4}\\
- (-1)^{a_{k-1}(\bi)} e(\bi) \epsilon_{k-1} e(\bi{\cdot}s_k s_{k-1} s_k), &if $i_{k-1} = -i_k = i_{k+1} \neq 0, \pm\frac{1}{2}$ and $h_k(\bi) = 0$,\label{rela:8:5}\\
  (-1)^{a_k(\bi)} e(\bi) \epsilon_k e(\bi{\cdot}s_k s_{k-1} s_k), &if $i_{k-1} = -i_k = i_{k+1} \neq 0, \pm\frac{1}{2}$ and $h_{k-1}(\bi) = 0$,\label{rela:8:6}\\
  e(\bi), &if $i_{k-1} + i_k, i_{k-1} + i_{k+1}, i_k + i_{k+1} \neq 0$ \notag\\
  & \hspace*{2cm} and $i_{k-1} = i_{k+1} = i_k - 1$,\label{rela:8:7}\\
- e(\bi), &if $i_{k-1} + i_k, i_{k-1} + i_{k+1}, i_k + i_{k+1} \neq 0$ \notag\\
  & \hspace*{2cm} and $i_{k-1} = i_{k+1} = i_k + 1$,\label{rela:8:8}\\
  0, &otherwise.\label{rela:8:9}
\end{numcases}
\end{enumerate}

The algebra is self-graded, where the degree of $e(\bi)$ is $0$, $y_k$ is $2$ and
$$
\deg e(\bi) \psi_k =
\begin{cases}
1, & \text{if $i_k = i_{k+1} \pm 1$,}\\
-2, & \text{if $i_k = i_{k+1}$,}\\
0, & \text{otherwise;}
\end{cases}
$$
and $\deg e(\bi) \epsilon_k e(\bj) = \deg_k(\bi) + \deg_k(\bj)$, where
$$
\deg_k(\bi) =
\begin{cases}
1, & \text{if $\bi \in P_{k,+}^n$,}\\
-1, & \text{if $\bi \in P_{k,-}^n$,}\\
0, & \text{if $\bi \in P_{k,0}^n$.}
\end{cases}
$$

\begin{Remark}
By the definition of $\G{n}$, the categorification of $P^n$ is used to determine the degree of $e(\bi) \epsilon_k e(\bj)$. By (\ref{rela:1}), we have $e(\bi) \epsilon_k e(\bj) = 0$ only if $i_k + i_{k+1} = 0$ and $j_k + j_{k+1} = 0$. We will see that $e(\bi) = 0$ only if $\bi \in I^n$ at the end of Section~\ref{sec:span} (see~\autoref{idem:nonzero}). Hence rather than categorizing $P^n$, we categorize the set $\set{\bi \in I^n | i_k + i_{k+1} = 0}$ for $1 \leq k \leq n-1$.
\end{Remark}

It is easy to verify that there exists an involution $*$ on $\G{n}$ such that $e(\bi)^* = e(\bi)$, $y_k^* = y_k$, $\psi_r^* = \psi_r$ and $\epsilon_r^* = \epsilon_r$ for $1 \leq k \leq n$ and $1 \leq r \leq n-1$.

We have an diagrammatic representation of $\G{n}$. To do this, we associate to each generator of $\G{n}$ an $P$-labelled decorated planar diagram on $2n$ dots in the following way:
\begin{align*}
e(\bi) &=
\begin{tikzpicture} [baseline=3mm,blue,line width=1pt, xscale = 0.6, yscale=0.8,
                      draw/.append style={rounded corners},
                      every node/.append style={font=\fontsize{5}{5}\selectfont}]
  \draw (0,1)node[above]{$i_1$}--(0,0);
  \draw (1,1)node[above]{$i_2$}--(1,0);
  \draw[dots] (1.2,1)--(2.8,1);
  \draw[dots] (1.2,0)--(2.8,0);
  \draw (3,1)node[above]{$i_{n}$}--(3,0);
\end{tikzpicture},
&
y_s e(\bi) & =
\begin{tikzpicture} [baseline=3mm,blue,line width=1pt, xscale = 0.6, yscale=0.8,
                      draw/.append style={rounded corners},
                      every node/.append style={font=\fontsize{5}{5}\selectfont}]
  \draw (0,1)node[above]{$i_1$}--(0,0);
  \draw[dots] (0.2,1)--(1.8,1);
  \draw[dots] (0.2,0)--(1.8,0);
  \draw (2,1)node[above]{$i_{s-1}$}--(2,0);
  \draw (3,1)node[above]{$i_s$}--(3,0);
  \draw (4,1)node[above]{$i_{s+1}$}--(4,0);
  \draw[dots] (4.2,1)--(5.8,1);
  \draw[dots] (4.2,0)--(5.8,0);
  \draw (6,1)node[above]{$i_{n}$}--(6,0);
  \node[greendot] at (3,0.5){};
\end{tikzpicture}, \\
\psi_re(\bi) & =
\begin{tikzpicture} [baseline=3mm,blue,line width=1pt, xscale = 0.6, yscale=0.8,
                      draw/.append style={rounded corners},
                      every node/.append style={font=\fontsize{5}{5}\selectfont}]
  \draw (0,1)node[above]{$i_1$}--(0,0);
  \draw[dots] (0.2,1)--(1.8,1);
  \draw[dots] (0.2,0)--(1.8,0);
  \draw (2,1)node[above]{$i_{r-1}$}--(2,0);
  \draw (3,1)node[above]{$i_r$}--(4,0);
  \draw (4,1)node[above]{$i_{r+1}$}--(3,0);
  \draw (5,1)node[above]{$i_{r+2}$}--(5,0);
  \draw[dots] (5.2,1)--(6.8,1);
  \draw[dots] (5.2,0)--(6.8,0);
  \draw (7,1)node[above]{$i_{n}$}--(7,0);
\end{tikzpicture},
& \epsilon_re(\bi) & =
\begin{tikzpicture} [baseline=3mm,blue,line width=1pt, xscale = 0.6, yscale=0.8,
                      draw/.append style={rounded corners},
                      every node/.append style={font=\fontsize{5}{5}\selectfont}]
  \draw (0,1)node[above]{$i_1$}--(0,0);
  \draw[dots] (0.2,1)--(1.8,1);
  \draw[dots] (0.2,0)--(1.8,0);
  \draw (2,1)node[above]{$i_{r-1}$}--(2,0);
  \draw (3,1)node[above]{$i_r$} .. controls (3.5,0.7) .. (4,1)node[above]{$i_{r+1}$};
  \draw (3,0) .. controls (3.5,0.3) .. (4,0);
  \draw (5,1)node[above]{$i_{r+2}$}--(5,0);
  \draw[dots] (5.2,1)--(6.8,1);
  \draw[dots] (5.2,0)--(6.8,0);
  \draw (7,1)node[above]{$i_{n}$}--(7,0);
\end{tikzpicture},
\end{align*}
for $\bi = (i_1, \ldots, i_n) \in P^n$, $1 \leq s \leq n$ and $1 \leq r \leq n-1$. The $r$-th string of the diagram is the string labelled with $i_r$.

Diagrams are considered up to isotopy, and multiplication of diagrams is given by concatenation, subject to the relations (\ref{rela:1}) - (\ref{rela:8:9}). We are not going to use the diagrammatic notations in this paper, but the reader may be aware that there exists such representation of $\G{n}$.


We give some identities in $\G{n}$ which will be used later for computational purposes. They can be easily verified using the relations of $\G{n}$. When $1 \leq k < n-1$ and $\bi,\bj \in P^n$, we have
\begin{align}
\epsilon_k \epsilon_{k+1} \psi_k & = \epsilon_k \epsilon_{k+1} \epsilon_k \psi_{k+1} = \epsilon_k \psi_{k+1}, \label{rela:9} \\
\psi_k \epsilon_{k+1} \epsilon_k & = \psi_{k+1} \epsilon_k \epsilon_{k+1} \epsilon_k = \psi_{k+1} \epsilon_k, \label{rela:10}\\
\psi_{k+1} \epsilon_k \psi_{k+1} & = \psi_k \epsilon_{k+1} \epsilon_k \epsilon_{k+1} \psi_k = \psi_k \epsilon_{k+1} \psi_k, \label{rela:11}\\
e(\bi) \epsilon_k \psi_{k+1} \epsilon_k e(\bj) & = e(\bi) \epsilon_k \epsilon_{k+1} \psi_k \epsilon_k e(\bj) = \pm e(\bi) \epsilon_k e(\bj); \label{rela:15}
\end{align}
and when $1 < k < n-1$, we have
\begin{align}
\psi_{k+1} \epsilon_{k-1} \epsilon_k \epsilon_{k+1} & = \epsilon_{k-1} \psi_k \epsilon_{k+1} = \epsilon_{k-1} \epsilon_k \epsilon_{k+1} \psi_{k-1}, \label{rela:12} \\
\epsilon_{k+1} \epsilon_{k-1} \epsilon_k \epsilon_{k+1} & = \epsilon_{k-1} \epsilon_{k+1} = \epsilon_{k-1} \epsilon_k \epsilon_{k+1} \epsilon_{k-1}, \label{rela:13}\\
\psi_{k+1} \epsilon_{k-1} \epsilon_k \psi_{k+1} & =  \epsilon_{k-1} \psi_k \epsilon_{k+1} \psi_k = \epsilon_{k-1} \epsilon_k \epsilon_{k+1} \psi_{k-1} \psi_k, \label{rela:14}\\
y_k \epsilon_{k+1} \epsilon_k & = - \epsilon_{k+1} y_{k+1} \epsilon_k = \epsilon_{k+1} \epsilon_k y_{k+2}. \label{rela:16}
\end{align}

\subsection{A method to determine the residue sequence} \label{sec:residue}

At the end of Section~\autoref{sec:span}, we will prove $e(\bi) = 0$ if $\bi \not\in I^n$ (See~\autoref{idem:nonzero}). Hence, for $\bi \in P^n$, it is important to determine whether $\bi$ is the residue sequence of some up-down tableaux. By~\autoref{deg:h2}, $\bi$ is not a residue sequence if for any $1 \leq k \leq n$ we have $h_k(\bi) \not\in \{-2,-1,0\}$, but the reverse is not true. In this subsection we will discuss the reverse part. The results given in this subsection will be used in Section~\ref{sec:grade} as well.

Suppose $\bi \in P^n$ and for $1 \leq k \leq n$, $\bi|_{k-1}$ is a residue sequence of some up-down tableaux. Choose arbitrary $\u \in \Tud_{k-1}(\bi|_{k-1})$ and let $\Shape(\u) = (\lambda,f)$. If $h_k(\bi) \in \{-2,-1\}$, by~\autoref{deg:h1:1} and $|\mathscr{AR}_\lambda(-i_k)| \geq 0$, it forces $|\mathscr{AR}_\lambda(i_k)| > 0$. Therefore the next Lemma follows.

\begin{Lemma} \label{deg:h9}
Suppose $\bi \in P^n$. If we have $h_k(\bi) \in \{-2,-1\}$ for all $1 \leq k \leq n$, then $\bi$ is the residue sequence of some up-down tableaux.
\end{Lemma}

\begin{proof}
We prove the Lemma by induction. When $n = 1$, the Lemma follows trivially. Assume when $n' < n$ the Lemma holds. When $n' = n$, set $\bj = (i_1, \ldots, i_{n-1})$. By induction, $\bj$ is the residue sequence of some up-down tableau. Suppose $\u$ is the up-down tableau with residue sequence $\bj$ and $(\lambda,f) = \Shape(\u)$. As $h_n(\bi) \in \{-2,-1\}$, we have $|\mathscr{AR}_\lambda(i_n)| > 0$ by~\autoref{deg:h1:1}. Let $\alpha \in \mathscr{AR}_\lambda(i_n)$. Without loss of generality we assume $\alpha \in \mathscr A(\lambda)$. Therefore if we write $\u = (\alpha_1, \ldots, \alpha_{n-1})$, then $\s = (\alpha_1, \ldots, \alpha_{n-1}, \alpha)$ is an up-down tableau with residue sequence $\bi$.
\end{proof}

If we have $h_k(\bi) = 0$ for some $1 \leq k \leq n$, by~\autoref{deg:h1:1} and~\eqref{deg:h2:ineq} we have $|\mathscr{AR}_\lambda(-i_k)| = |\mathscr{AR}_\lambda(i_k)| \in \{0,1\}$. Hence we cannot decide whether $\bi$ is the residue sequence of some up-down tableau. Note that when $i_k = 0$, we have $h_k(\bi) = 0$. In the rest of this subsection we extend~\autoref{deg:h9} and include the case when $i_k = 0$.

\begin{Lemma} \label{deg:h4:6}
Suppose $\bi \in P^n$ and $\t \in \Tud_n(\lambda)$ is an up-down tableau with residue sequence $\bi$. Then we have
\begin{eqnarray*}
&& \#\set{\alpha \in [\lambda] | \res(\alpha) = 1} - \#\set{\alpha \in [\lambda] | \res(\alpha) = -1}\\
& = & \#\set{1 \leq k \leq n | i_k = 1} - \#\set{1 \leq k \leq n | i_k = -1}.
\end{eqnarray*}
\end{Lemma}

\begin{proof}
Apply induction on $n$. The base case $n = 1$ follows trivially. For the induction step, we assume that the Lemma holds for all $\bj \in P^{n-1}$ and prove the Lemma holds for $\bi \in P^n$. For convenience we denote
\begin{align*}
g_1(\lambda) & = \#\set{\alpha \in [\lambda] | \res(\alpha) = 1}, & g_{-1}(\lambda) & = \#\set{\alpha \in [\lambda] | \res(\alpha) = -1}\\
n_1(\bi) & = \#\set{1 \leq k \leq n | i_k = 1}, & n_{-1}(\bi) & = \#\set{1 \leq k \leq n | i_k = -1},
\end{align*}
and we need to prove that
\begin{equation} \label{deg:h4:6:eq1}
g_1(\lambda) - g_{-1}(\lambda) = n_1(\bi) - n_{-1}(\bi).
\end{equation}

Let $\bj = (i_1, \ldots, i_{n-1})$ and $\s = \t|_{n-1} \in \Tud_{n-1}(\mu)$ for some $\mu$. We note that if $\t(n) = \alpha > 0$, we have $\lambda = \mu\cup\{\alpha\}$; and if $\t(n) = -\alpha < 0$, we have $\mu = \lambda\cup\{\alpha\}$.

By the construction, $\s$ is an up-down tableau with residue sequence $\bj$. By induction we have
\begin{equation} \label{deg:h4:6:eq2}
g_1(\mu) - g_{-1}(\mu) = n_1(\bj) - n_{-1}(\bj).
\end{equation}

If $i_n \neq \pm 1$, we have $g_1(\lambda) = g_1(\mu)$, $g_{-1}(\lambda) = g_{-1}(\mu)$, $n_1(\bi) = n_1(\bj)$ and $n_{-1}(\bi) = n_{-1}(\bj)$. Hence (\ref{deg:h4:6:eq1}) holds by (\ref{deg:h4:6:eq2}).

If $i_n = 1$ and $\t(n) > 0$, we have $g_1(\lambda) = g_1(\mu) + 1$, $g_{-1}(\lambda) = g_{-1}(\mu)$, $n_1(\bi) = n_1(\bj) + 1$ and $n_{-1}(\bi) = n_{-1}(\bj)$; and if $i_n = 1$ and $\t(n) < 0$, we have $g_1(\lambda) = g_1(\mu)$, $g_{-1}(\lambda) = g_{-1}(\mu) - 1$, $n_1(\bi) = n_1(\bj) + 1$ and $n_{-1}(\bi) = n_{-1}(\bj)$. Hence (\ref{deg:h4:6:eq1}) holds by (\ref{deg:h4:6:eq2}).

If $i_n = -1$ and $\t(n) > 0$, we have $g_1(\lambda) = g_1(\mu)$, $g_{-1}(\lambda) = g_{-1}(\mu) + 1$, $n_1(\bi) = n_1(\bj)$ and $n_{-1}(\bi) = n_{-1}(\bj) + 1$; and if $i_n = -1$ and $\t(n) < 0$, we have $g_1(\lambda) = g_1(\mu) - 1$, $g_{-1}(\lambda) = g_{-1}(\mu)$, $n_1(\bi) = n_1(\bj)$ and $n_{-1}(\bi) = n_{-1}(\bj) + 1$. Hence (\ref{deg:h4:6:eq1}) holds by (\ref{deg:h4:6:eq2}).
\end{proof}

One can see that by knowing the values of $\delta$ and $\#\set{\alpha \in [\lambda] | \res(\alpha) = 1} - \#\set{\alpha \in [\lambda] | \res(\alpha) = -1}$, we can determine the value of $|\mathscr{AR}_\lambda(0)|$.

\begin{Lemma} \label{deg:h4:3}
Suppose $(\lambda, f) \in \widehat B_n$. Then $|\mathscr{AR}_\lambda(0)| = 1$ if and only if $\delta$ is odd and one of the following conditions holds:

(1) $\delta = 1$ and $\#\set{\alpha \in [\lambda] | \res(\alpha) = 1} - \#\set{\alpha \in [\lambda] | \res(\alpha) = -1} = 0$;

(2) $\delta < 1$ and $\#\set{\alpha \in [\lambda] | \res(\alpha) = 1} - \#\set{\alpha \in [\lambda] | \res(\alpha) = -1} = -1$;

(3) $\delta > 1$ and $\#\set{\alpha \in [\lambda] | \res(\alpha) = 1} - \#\set{\alpha \in [\lambda] | \res(\alpha) = -1} = 1$.
\end{Lemma}

\begin{proof}
The Lemma follows directly by the construction of $[\lambda]$.
\end{proof}

The next result is directly implied by~\autoref{deg:h4:6} and~\autoref{deg:h4:3}. Define
$$
a_k^*(\bi) = \#\set{1 \leq r \leq k-1 | i_r \in \{-1,1\} } + \delta_{0,\frac{\delta - 1}{2}}
$$
for $1 \leq k \leq n$ and $\bi \in P^n$ with $i_k = 0$. It is easy to see that if $1 \leq k \leq n-1$ and $i_k = i_{k+1} = 0$, we have $a_k(\bi) = a_k^*(\bi) + 1$.


\begin{Corollary} \label{deg:h4:5}
Suppose $\bi \in I^{n-1}$ and $\bj = \bi \vee 0 \in P^n$. Let $\t \in \Tud_{n-1}(\bi)$ with shape $(\lambda,f)$. Then $|\mathscr{AR}_\lambda(0)| = 1$ if and only if $a_{n-1}(\bj)$ is even.
\end{Corollary}

\proof Suppose $|\mathscr{AR}_\lambda(0)| = 1$. By~\autoref{deg:h4:6}, the parities of $\#\set{1 \leq r \leq n-1 | i_r = 1} + \#\set{1 \leq r \leq n-1 | i_r = -1}$ and $\#\set{\alpha \in [\lambda] | \res(\alpha) = 1} + \#\set{\alpha \in [\lambda] | \res(\alpha) = -1}$ are the same. As $i_1 = \frac{\delta-1}{2}$, we have $a_{n-1}^*(\bj)$ is odd by~\autoref{deg:h4:3} because $|\mathscr{AR}_\lambda(0)| = 1$. This proves the only if part.

Suppose that $a_{n-1}^*(\bj)$ is odd. It forces $\delta$ to be odd. When $\delta = 1$, by the construction of young diagrams, we have
$$
-1 \leq \#\set{\alpha \in [\lambda] | \res(\alpha) = 1} - \#\set{\alpha \in [\lambda] | \res(\alpha) = -1} \leq 1,
$$
which implies that $\delta_k(\bj)$ is odd if and only if $\#\set{\alpha \in [\lambda] | \res(\alpha) = 1} - \#\set{\alpha \in [\lambda] | \res(\alpha) = -1} = 0$; and when $\delta < 1$, by the construction of young diagrams, we have
$$
-2 \leq \#\set{\alpha \in [\lambda] | \res(\alpha) = 1} - \#\set{\alpha \in [\lambda] | \res(\alpha) = -1} \leq 0,
$$
which implies that $\delta_k(\bj)$ is odd if and only if $\#\set{\alpha \in [\lambda] | \res(\alpha) = 1} - \#\set{\alpha \in [\lambda] | \res(\alpha) = -1} = -1$; and when $\delta > 1$, by the construction of young diagrams, we have
$$
0 \leq \#\set{\alpha \in [\lambda] | \res(\alpha) = 1} - \#\set{\alpha \in [\lambda] | \res(\alpha) = -1} \leq 2,
$$
which implies that $\delta_k(\bj)$ is odd if and only if $\#\set{\alpha \in [\lambda] | \res(\alpha) = 1} - \#\set{\alpha \in [\lambda] | \res(\alpha) = -1} = 1$. Therefore, we have $a_{n-1}^*(\bj)$ is odd if and only if $|\mathscr{AR}_\lambda(0)| = 1$ by~\autoref{deg:h4:3}. This proves the if part. \endproof

The next Corollary is easy to be verified.

\begin{Corollary} \label{deg:h4:4}
Suppose $\bi \in I^{n-1}$ and $\bj = \bi \vee 0 \in P^n$. Then $\bj \in I^n$ if and only if $a_{n-1}(\bj)$ is even.
\end{Corollary}

Now we can extend~\autoref{deg:h9}.

\begin{Lemma} \label{deg:h8}
Suppose $\bi \in P^n$. If for any $1 \leq k \leq n$, we have either $h_k(\bi) \in \{-2,-1\}$, or $i_k = 0$ and $a_k^*(\bi)$ is odd. Then $\bi$ is the residue sequence of some up-down tableau.
\end{Lemma}

\begin{proof}
We prove the Lemma by induction. When $n = 1$, the Lemma follows trivially. Assume when $n' < n$ the Lemma holds. When $n' = n$, set $\bj = (i_1, \ldots, i_{n-1})$. By induction, $\bj$ is the residue sequence of some up-down tableau. If $h_n(\bi) \in \{-2,-1\}$, then $\bi \in I^n$ by~\autoref{deg:h9}; and if $i_k = 0$ and $a_k^*(\bi)$ is odd, then $\bi \in I^n$ by~\autoref{deg:h4:4}.
\end{proof}


\begin{Remark}
\autoref{deg:h8} is not sufficient to determine whether $\bi \in P^n$ is the residue sequence of some up-down tableaux. For instance, we cannot determine whether $\bi \in I^n$ if $i_k \neq 0$ and $h_k(\bi) = 0$ for some $1 \leq k \leq n$. We also remark the whole argument only works when $R$ is a field of characteristic $0$. When $R$ is a field of positive characteristic, there exists residue sequence $\bi \in I^n$ such that $h_k(\bi) \not\in \{-2,-1,0\}$ for some $1 \leq k \leq n$.
\end{Remark}

We close this section by giving two more results which can be used to determine whether $\bi{\cdot}s_k$ is a residue sequence of some up-down tableaux by giving $\bi$ is a residue sequence.

\begin{Lemma} \label{A:9}
Suppose $1 \leq k \leq n-1$ and $\bi = (i_1, \ldots, i_n)$ is a residue sequence with $i_k + i_{k+1} = 0$. Then $\bi{\cdot}s_k$ is a residue sequence if and only if $h_k(\bi) = 0$.
\end{Lemma}

\proof Because $i_k + i_{k+1} = 0$, we have $i_k = -i_{k+1}$. Therefore, by the definition of $h_k$, we have $h_k(\bi{\cdot}s_k) = - h_k(\bi)$. By~\autoref{deg:h2}, we have $-2 \leq h_k(\bi) \leq 0$ as $\bi \in I^n$, which implies $0 \leq h_k(\bi{\cdot}s_k) \leq 2$. By~\autoref{deg:h2}, we have $\bi{\cdot}s_k \in I^n$ only if $h_k(\bi{\cdot}s_k) = 0$, which implies $h_k(\bi) = 0$. This proves the only if part.

Assume $h_k(\bi) = 0$. When $i_k = i_{k+1} = 0$, we have $\bi{\cdot}s_k = \bi \in I^n$; and by~\autoref{deg:h3}, we have $i_k \neq \pm \frac{1}{2}$ as $h_k(\bi) = 0$. Hence, it only left us to prove $\bi{\cdot}s_k$ is a residue sequence when $|i_k - i_{k+1}| > 1$. In this case, choose arbitrary $\t \in \Tud_n(\bi)$.

Suppose $\t(k) + \t(k+1) \neq 0$. When $\t(k) > 0$, $\t(k+1) < 0$ or $\t(k) < 0$, $\t(k+1) > 0$, we have that $\t{\cdot}s_k$ is an up-down tableau by~\autoref{y:h2:8}; and when $\t(k), \t(k+1) > 0$ or $\t(k), \t(k+1) < 0$, because $|i_k - i_{k+1}| > 1$, the nodes $\t(k)$ and $\t(k+1)$ are not adjacent, which implies $\t{\cdot}s_k$ is an up-down tableau by~\autoref{y:h2:8}. This proves $\bi{\cdot}s_k$ is a residue sequence.

Suppose $\t(k) + \t(k+1) = 0$. By~\autoref{deg:h4:1}, there exists a unique up-down tableau $\s$ with residue sequence $\bi{\cdot}s_k$ such that $\s \overset{k}\sim \t$, which implies $\bi{\cdot}s_k$ is a residue sequence. This proves the if part. \endproof

\begin{Lemma} \label{A:8}
Suppose $1 \leq k \leq n-1$ and $\bi = (i_1, \ldots, i_n) \in P^n$ with $|i_k - i_{k+1}| > 1$ and $i_k + i_{k+1} \neq 0$. If we have $\t \in \Tud_n(\bi)$, then $\t{\cdot}s_k$ is an up-down tableau. In another word, we have $\bi{\cdot}s_k \in I^n$ if and only if $\bi \in I^n$.
\end{Lemma}

\proof Suppose $\bi \in I^n$, there exists $\t \in \Tud_n(\bi)$. Because $i_k + i_{k+1} \neq 0$ and $|i_k - i_{k+1}| > 1$, $\t(k)$ and $\t(k+1)$ satisfy the conditions of~\autoref{y:h2:8}. Hence, $\s = \t{\cdot}s_k$ is an up-down tableau. Because $\s \in \Tud_n(\bi{\cdot}s_k)$, we have $\bi{\cdot}s_k \in I^n$. Following the same argument, we have $\bi \in I^n$ if $\bi{\cdot}s_k \in I^n$. \endproof

\section{Induction and restriction of $\G{n}$} \label{sec:ir}

In this section we discuss the induction and restriction properties of $\G{n}$. Instead of working on $\G{n}$ directly, first we construct a set of homogeneous elements
$$
\set{\psi_{\s\t} | (\lambda, f) \in \widehat B_n, \s,\t \in \Tud_n(\lambda)} \subset \G{n}
$$
analogue to the $\psi$-basis of $\K\Sym_n$ given in~\autoref{def:basis of Sym} and define $R_n(\delta)$ to be the $R$-span of $\{\psi_{\s\t}\}$. Then we prove the induction and restriction properties of $R_n(\delta)$. In Section~\ref{sec:span}, we will prove that $R_n(\delta) = \G{n}$ and all the results of this section will directly apply to $\G{n}$.

\subsection{A set of homogeneous elements of $\G{n}$} \label{sec:basis}

In this subsection we construct the set $\set{\psi_{\s\t} | (\lambda, f) \in \widehat B_n, \s,\t \in \Tud_n(\lambda)}$, which is a set of homogeneous elements of $\G{n}$. Then we calculate the degree of $\psi_{\s\t}$ and show that $\deg \psi_{\s\t}$ is determined by $\deg \s$ and $\deg \t$.

Suppose $(\lambda,f) \in \widehat B_n$ and $\t \in \Tud_n(\lambda)$ is an up-down tableau. Recall that we can write $\t = (\alpha_1,\alpha_2,\ldots,\alpha_n)$. For convenience we denote $\alpha_0 = (1,1)$. Let $\alpha_{k_1}, \alpha_{k_2},\ldots,\alpha_{k_f}$ be all negative nodes of $\t$ and, for $1 \leq i \leq f$, let $\alpha_{j_i}$ be the first node on the left of $\alpha_{k_i}$ in $\t$ with $\alpha_{j_i} = - \alpha_{k_i}$. So we can choose $f$ pairs of nodes $(\alpha_{j_1},\alpha_{k_1}),\ldots,(\alpha_{j_f},\alpha_{k_f})$ by re-ordering these nodes such that $k_1 < k_2 < \ldots < k_f$. Define the sequence to be the \textit{remove pairs} of the up-down tableau $\t$. Suppose $(\alpha_1,\ldots,\hat\alpha_{j_i}, \ldots,\hat\alpha_{k_i},\ldots,\alpha_n)$ is an up-down tableau, we say $(\alpha_{j_i},\alpha_{k_i})$ a \textit{removable pair} of $\t$ and $i$ the \textit{removable index} of $\t$. If $(\alpha_{j_i},\alpha_{k_i}) \neq (\alpha_0, -\alpha_0)$, we say $(\alpha_{j_i},\alpha_{k_i})$ is a \textit{nontrivial removable pair} of $\t$ and $i$ the \textit{nontrivial removable index} of $\t$.

One can see that there exists an integer $h$ such that for any $i \leq h$, $j_i = 2i-1$, $k_i = 2i$, $\alpha_{j_i} = \alpha_0$ and $\alpha_{k_i} = -\alpha_0$, and $k_{h+1} \neq 2(h+1)$. We define such integer $h$ to be the \textit{head} of the up-down tableau $\t$ and denote $head(\t) = h$.

\begin{Example} \label{udtab:ex1}
Suppose $n = 9$ and $\lambda = (1)$. We have $(\lambda,4) \in \widehat B_9$. Define
$$
\t = \left( \emptyset, \ydiag(1), \emptyset, \ydiag(1), \ydiag(1,1), \ydiag(2,1), \ydiag(2), \ydiag(3), \ydiag(2), \ydiag(1) \right) \in \Tud_9(\lambda).
$$

We can write $\t$ as
$$
\t = \left( \alpha_1, \ldots, \alpha_9 \right) = \left( (1,1), -(1,1), (1,1), (2,1), (1,2), -(2,1), (1,3), -(1,3), -(1,2) \right),
$$
where the negative nodes are $\alpha_2, \alpha_6, \alpha_8$ and $\alpha_9$ and the head $h = 1$. Moreover, we have a sequence of remove pairs $\left((\alpha_1, \alpha_2), (\alpha_4, \alpha_6), (\alpha_7, \alpha_8), (\alpha_5, \alpha_9) \right)$, but not all of them are removable pairs. In more details, $(\alpha_5, \alpha_9)$ is not removable and all the other pairs are removable.
\end{Example}

\begin{Lemma} \label{udtab:h1}
Suppose $(\lambda,f) \in \widehat B_n$ and $\t \in \Tud_n(\lambda)$ with head $h < f$. Then $h+1$ is a removable index, i.e. $(\alpha_{j_{h+1}}, \alpha_{k_{h+1}})$ is a removable pair of $\t$.
\end{Lemma}

\begin{proof}
One can see that for any $i$ with $h < i \leq f$, $(\alpha_{j_i}, \alpha_{k_i})$ is not removable only if there exists $h < m < i$ such that $j_i < j_m < k_m < k_i$. When $i = h+1$, there is not such $m$ exists, which proves the Lemma.
\end{proof}

\begin{Definition}
Suppose $(\lambda,f) \in \widehat B_n$ and $\t = (\alpha_1,\ldots,\alpha_n) \in \Tud_n(\lambda)$ with head $h < f$ and remove pairs $(\alpha_{j_1},\alpha_{k_1}), \ldots,$ $(\alpha_{j_f},\alpha_{k_f})$. Let $i$ be a non-trivial removable index of $\t$ and $\s = (\alpha_0, -\alpha_0, \alpha_1, \ldots, \hat\alpha_{j_i}, \ldots, \hat\alpha_{k_i},\ldots, \alpha_n)$. We denote $\s \rightarrow \t$, and $\rho(\s,\t) = (j_i, k_i)$.
\end{Definition}

\begin{Example} \label{udtab:ex2}
Suppose $\t$ is defined as in~\autoref{udtab:ex1}. We have $2$ nontrivial removable pairs of $\t$: $(\alpha_4, \alpha_6)$ and $(\alpha_7, \alpha_8)$. Define
\begin{align*}
\s_1 & = \left( (1,1), -(1,1), (1,1), -(1,1), (1,1), (1,2), (1,3), -(1,3), -(1,2) \right)\\
& = \left( \emptyset, \ydiag(1), \emptyset, \ydiag(1), \emptyset, \ydiag(1), \ydiag(2), \ydiag(3), \ydiag(2), \ydiag(1) \right) \in \Tud_n(\lambda); \\
\s_2 & = \left( (1,1), -(1,1), (1,1), -(1,1), (1,1), (2,1), (1,2), -(2,1), -(1,2) \right)\\
& = \left( \emptyset, \ydiag(1), \emptyset, \ydiag(1), \emptyset, \ydiag(1), \ydiag(1,1), \ydiag(2,1), \ydiag(2), \ydiag(1) \right) \in \Tud_n(\lambda).
\end{align*}

Hence we have $\s_1 \rightarrow \t$ and $\s_2 \rightarrow \t$, where $\rho(\s_1,\t) = (4,6)$ and $\rho(\s_2,\t) = (7,8)$. Notice that $(\alpha_6,\alpha_9)$ is not a removable pair in $\t$, because of the existence of $(\alpha_7, \alpha_8)$. But after we remove $(\alpha_7, \alpha_8)$, $(\alpha_6,\alpha_9)$ becomes removable. As an example, set
\begin{align*}
\s_3 & = \left( (1,1), -(1,1), (1,1), -(1,1), (1,1), -(1,1), (1,1), (2,1), -(2,1) \right)\\
& = \left( \emptyset, \ydiag(1), \emptyset, \ydiag(1), \emptyset, \ydiag(1), \emptyset, \ydiag(1), \ydiag(1,1), \ydiag(1) \right) \in \Tud_n(\lambda).
\end{align*}

Then we have $\s_3 \rightarrow \s_2$.
\end{Example}

The next Lemma is obvious by the definition of $\s \rightarrow \t$.

\begin{Lemma} \label{udtab:h2}
Suppose $(\lambda,f) \in \widehat B_n$ and $\t \in \Tud_n(\lambda)$. If there exists an up-down tableau $\s$ such that $\s \rightarrow \t$, then $\s \in \Tud_n(\lambda)$ and $head(\s) = head(t) + 1$.
\end{Lemma}

Suppose $(\lambda,f) \in \widehat B_n$ and $\t \in \Tud_n(\lambda)$ with head $h < f$. Because for any $\s \in \Tud_n(\lambda)$, we have $head(\s) \leq f$. Hence by~\autoref{udtab:h2}, there exists a finite sequence
$$
\t^{(m)} \rightarrow \t^{(m-1)} \rightarrow \ldots \rightarrow \t^{(1)} \rightarrow \t^{(0)} = \t,
$$
where $\t^{(1)}, \ldots, \t^{(m)} \in \Tud_n(\lambda)$ and there is no $\s \in \Tud_n(\lambda)$ such that $\s \rightarrow \t^{(m)}$. We define such sequence to be the \textit{reduction sequence} of $\t$.

\begin{Example} \label{udtab:ex3}
Suppose $\t$ is defined as in~\autoref{udtab:ex1}. Define
\begin{align*}
\u & = \left( (1,1), -(1,1), (1,1), -(1,1), (1,1), -(1,1), (1,1), -(1,1), (1,1)\right)\\
& = \left( \emptyset, \ydiag(1), \emptyset, \ydiag(1), \emptyset, \ydiag(1), \emptyset, \ydiag(1), \emptyset, \ydiag(1) \right) \in \Tud_n(\lambda).
\end{align*}

Hence we have the sequence $\t^{(3)} \rightarrow \t^{(2)} \rightarrow \t^{(1)} \rightarrow \t^{(0)} = \t$, where $\t^{(3)} = \u$,
\begin{align*}
\t^{(2)} & = \left( (1,1), -(1,1), (1,1), -(1,1), (1,1), -(1,1), (1,1), (1,2), -(1,2) \right)\\
& = \left( \emptyset, \ydiag(1), \emptyset, \ydiag(1), \emptyset, \ydiag(1), \emptyset, \ydiag(1), \ydiag(2), \ydiag(1) \right) \in \Tud_n(\lambda),\\
\t^{(1)} & = \left( (1,1), -(1,1), (1,1), -(1,1), (1,1), (1,2), (1,3), -(1,3), -(1,2) \right)\\
& = \left( \emptyset, \ydiag(1), \emptyset, \ydiag(1), \emptyset, \ydiag(1), \ydiag(2), \ydiag(3), \ydiag(2), \ydiag(1) \right) \in \Tud_n(\lambda),\\
\t = \t^{(0)} & = \left( (1,1), -(1,1), (1,1), (2,1), (1,2), -(2,1), (1,3), -(1,3), -(1,2) \right)\\
& = \left( \emptyset, \ydiag(1), \emptyset, \ydiag(1), \ydiag(1,1), \ydiag(2,1), \ydiag(2), \ydiag(3), \ydiag(2), \ydiag(1) \right) \in \Tud_n(\lambda).
\end{align*}
\end{Example}

In general, for $\t \in \Tud_n(\lambda)$ with head $h < f$, there exist more than one reduction sequence of $\t$. For example, as in~\autoref{udtab:ex2} and~\autoref{udtab:ex3}, there exists another sequence $\t^{(3)} \rightarrow \t^{(2)} \rightarrow \t^{(1)} \rightarrow \t$ where $\t^{(3)} = \u$, $\t^{(2)} = \s_3$ and $\t^{(1)} = \s_2$.

\begin{Lemma} \label{udtab:h3}
Suppose $(\lambda,f) \in \widehat B_n$ and $\t = (\alpha_1, \ldots, \alpha_n) \in \Tud_n(\lambda)$ with head $h<f$ and remove pairs $(\alpha_{j_1},\alpha_{k_1}), \ldots, (\alpha_{j_f}, \alpha_{k_f})$. For any reduction sequence $\t^{(m)} \rightarrow \t^{(m-1)} \rightarrow \ldots \rightarrow \t^{(1)} \rightarrow \t^{(0)} = \t$, we have $m = f-h$.
\end{Lemma}

\begin{proof}
By~\autoref{udtab:h2}, $head(\t^{(m)}) = h + m$. If $m > f-h$, we have $head(\t^{(m)}) = h + m > f$. Because $\lambda \vdash n-2f$, $\t^{(m)} \not\in \Tud_n(\lambda)$. Hence we have $m \leq f-h$. If $m < f-h$, then $head(\t^{(m)}) = h + m < f$. By~\autoref{udtab:h1} there exists $\s \in \Tud_n(\lambda)$ such that $\s \rightarrow \t^{(m)}$. Therefore we have $m = f-h$.
\end{proof}

\autoref{udtab:h3} shows that the length of the reduction sequence is determined by $\t$. Moreover, by the construction of the reduction sequence, we have $\t^{(f-h)} = (\beta_1, \ldots, \beta_n) \in \Tud_n(\lambda)$ where $\beta_1 = \beta_3 = \ldots = \beta_{2f-1} = -\beta_2 = -\beta_4 = \ldots = - \beta_{2f} = \alpha_0$, and $(\beta_{2f+1}, \beta_{2f+2}, \ldots, \beta_n)$ is obtained by removing $\alpha_{j_1}, \alpha_{k_1}, \alpha_{j_2}, \alpha_{k_2}, \ldots, \alpha_{j_n}, \alpha_{k_n}$ from $\t = (\alpha_1, \ldots, \alpha_n)$.

\begin{Example}
Suppose $\t$ is defined as in~\autoref{udtab:ex1} and $\u$ is in~\autoref{udtab:ex3}. We can see that by removing $(\alpha_1, \alpha_2), (\alpha_4, \alpha_6), (\alpha_7, \alpha_8), (\alpha_5, \alpha_9)$ from $\t$ and add $f$ number of $(\alpha_0, -\alpha_0)$ at the front of the resulting sequence, we have
$$
\left( (1,1), -(1,1), (1,1), -(1,1), (1,1), -(1,1), (1,1), -(1,1), (1,1)\right) = \u.
$$
\end{Example}

Therefore $\t^{(f-h)}$ is uniquely determined by $\t$, and we denote $h(\t) = \t^{(f-h)}$. Moreover, by the construction of $\t^{(f-h)}$, one can see that if $\t^{(f-h)} = (\beta_1, \ldots, \beta_n)$, by defining $\s = (\beta_{2f+1}, \beta_{2f+2}, \ldots, \beta_n)$, $\s$ is a tableau of shape $\lambda$.

Recall the reduction sequence of $\t$ is not unique in general. A \textit{standard reduction sequence} is a reduction sequence
$$
\t^{(f-h)} \rightarrow \t^{(f-h-1)} \rightarrow \ldots \rightarrow \t^{(1)} \rightarrow \t^{(0)} = \t
$$
such that for any $0 \leq i \leq f-h-1$, $\t^{(i+1)}$ is obtained by removing the first non-trivial removable pairs of $\t^{(i)}$. In more details, suppose $\t^{(i)}$ has head $h + i$ and remove pairs $(\alpha_{j_1},\alpha_{k_1}),\ldots,(\alpha_{j_f},\alpha_{k_f})$, we obtain $\t^{(i+1)}$ by removing $(\alpha_{j_{h+i+1}}, \alpha_{k_{h+i+1}})$, which could be done because of~\autoref{udtab:h1}. As an example, the reduction sequence in~\autoref{udtab:ex3} is a standard reduction sequence. In the rest of this paper, for any $0 \leq i \leq f-h$, $\t^{(i)}$ means the corresponding up-down tableau in the standard reduction sequence of $\t$.

Suppose $\s \rightarrow \t$ where $\s$ has head $h+1$ and $\rho(\s,\t) = (a,b)$. Define
$$
\epsilon_{\s \rightarrow \t} := e(\bi_\s) \epsilon_{2h+2} \epsilon_{2h+3} \ldots \epsilon_a \psi_{a+1} \psi_{a+2} \ldots \psi_{b-1} e(\bi_\t) \in \G{n}.
$$

Fix $(\lambda,f) \in \widehat B_n$. For any $\t \in \Tud_n(\lambda)$ with head $h < f$ and standard reduction sequence
$$
h(\t) = \t^{(f-h)} \rightarrow \t^{(f-h-1)} \rightarrow \ldots \rightarrow \t^{(1)} \rightarrow \t^{(0)} = \t,
$$
define $\epsilon_\t := \epsilon_{\t^{(f-h)} \rightarrow \t^{(f-h-1)}} \epsilon_{\t^{(f-h-1)} \rightarrow \t^{(f-h-2)}} \ldots \epsilon_{\t^{(1)} \rightarrow \t^{(0)}}$; and for $\t \in \Tud_n(\lambda)$ with head $f$, define $\epsilon_\t := e(\bi_\t)$.

Define $\Sym_{2f,n}$ to be the subalgebra of $\Sym_n$ generated by $s_{2f+1}, s_{2f+2}, \ldots, s_n$. Denote $\t^{(\lambda,f)}$ to be the unique up-down tableau with shape $(\lambda,f)$ which is maximal in dominance ordering. Suppose $\t$ is an up-down tableau. Because $h(\t)$ is uniquely determined by $\t$, we abuse the symbol and define $d(\t) = d(h(\t)) \in \Sym_{2f,n}$ to be the reduced word such that $\t^\lambda d(\t) = h(\t)$. Write $d(\t) = s_{k_1} s_{k_2} \ldots s_{k_l}$. We define
$$
\psi_\t := \psi_{d(\t)} = \psi_{k_1} \psi_{k_2} \ldots \psi_{k_l}.
$$

Let $\bi_{(\lambda,f)}$ be the residue sequence of $\t^{(\lambda,f)}$ and define $e_{(\lambda,f)} := e(\bi_\lambda) \epsilon_1 \epsilon_3 \ldots \epsilon_{2f-1} e(\bi_\lambda)$.

\begin{Definition}
Suppose $(\lambda,f) \in \widehat B_n$ and $\s,\t \in \Tud_n(\lambda)$. We define $\psi_{\s\t} := \epsilon_\s^* \psi_\s^* e_{(\lambda,f)} \psi_\t \epsilon_\t$.
\end{Definition}

\begin{Lemma} \label{idem:psi:1}
Suppose $(\lambda,f) \in \widehat B_n$ and $\s,\t \in \Tud_n(\lambda)$. For any $\bi, \bj \in I^n$, we have
$$
e(\bi) \psi_{\s\t} e(\bj) =
\begin{cases}
\psi_{\s\t}, & \text{if $\bi = \bi_\s$ and $\bj = \bi_\t$,}\\
0, & \text{otherwise.}
\end{cases}
$$
\end{Lemma}

\begin{proof}
It is obvious by the definition of $\psi_{\s\t}$ and (\ref{rela:1}).
\end{proof}

For an up-down tableau $\t$, $\psi_\t = \psi_{d(\t)}$ is determined by the choice of reduced expression of $d(\t)$. Here we prove that $\psi_\t$ is actually independent to the choice of reduced expression of $d(\t)$.

Suppose $(\lambda,f) \in \widehat B_n$ and $\t \in \Tud_n(\lambda)$. For $w = s_{r_1} s_{r_2} \ldots s_{r_m} \in \Sym_{2f,n}$, if for any $1 \leq k \leq m$, $\t{\cdot}s_{r_1} s_{r_2} \ldots s_{r_k}$ is an up-down tableau, then we define $s_{r_1} s_{r_2} \ldots s_{r_m}$ to be \textit{semi-reduced correspond to $\t$}. Notice that if we have $\s = \t{\cdot}s_{r_1} s_{r_2} \ldots s_{r_m}$ and $s_{r_1} s_{r_2} \ldots s_{r_m}$ is semi-reduced correspond to $\t$, then $s_{r_m} s_{r_{m-1}} \ldots s_{r_1}$ is semi-reduced correspond to $\s$.

The next Lemma can be easily verified by the definitions of semi-reduced.

\begin{Lemma} \label{semi-red2}
Suppose $(\lambda,f) \in \widehat B_n$ and $\t \in \Tud_n(\lambda)$ with head $f$. Let $s_{r_1} s_{r_2} \ldots s_{r_m} \in \Sym_{2f,n}$ be semi-reduced correspond to $\t$, and $\s = \t{\cdot}s_{r_1} s_{r_2} \ldots s_{r_m} \in \Tud_n(\lambda)$ with head $f$ and $s_{k_1} s_{k_2} \ldots s_{k_l} \in \Sym_{2f,n}$ be semi-reduced correspond to $\s$. Then $s_{r_1} s_{r_2} \ldots s_{r_m} s_{k_1} s_{k_2} \ldots s_{k_l} \in \Sym_{2f,n}$ is semi-reduced correspond to $\t$.
\end{Lemma}

The motivation of semi-reduced is to calculate $e(\bi_\t) \psi_w$ when $w$ is semi-reduced.

\begin{Lemma} \label{semi-red1}
Suppose $(\lambda,f) \in \widehat B_n$ and $\t \in \Tud_n(\lambda)$. If $w = s_{r_1} s_{r_2} \ldots s_{r_m} \in \Sym_{2f,n}$ is semi-reduced correspond to $\t$, then $e(\bi_\t) \psi_{r_1} \psi_{r_2} \ldots \psi_{r_m} = e(\bi_\t) \psi_w$. Equivalently, we have $\psi_{r_m} \psi_{r_{m-1}} \ldots \psi_{r_1} e(\bi_\t) = \psi_{w^{-1}} e(\bi_\t)$.
\end{Lemma}

\proof It is sufficient if we can prove that
$$
e(\bi_\t) \psi_k \psi_{k+1} \psi_k = e(\bi_\t) \psi_{k+1} \psi_k \psi_{k+1},  \qquad e(\bi_\t) \psi_k^2 = e(\bi_\t), \qquad e(\bi_\t) \psi_k \psi_r = e(\bi_\t) \psi_r \psi_k,
$$
under the assumptions of the Lemma, which can be verified directly by checking the relations of $\G{n}$ and the construction of up-down tableaux. \endproof

\begin{Example}
Suppose $\t = (\emptyset, \ydiag(1), \ydiag(2))$. Then $\t{\cdot}s_1^2 = \t$ but $\t{\cdot}s_1$ is not an up-down tableau. So $s_1^2$ is not semi-reduced correspond to $\t$.

Suoose $\t = (\emptyset, \ydiag(1), \ydiag(2), \ydiag(2,1))$. Then $\t{\cdot}s_2^2 = \t$, and $\t{\cdot}s_2 = (\emptyset, \ydiag(1), \ydiag(1,1), \ydiag(2,1))$ is an up-down tableau. Then $s_2^2$ is semi-reduced correspond to $\t$. Moreover, choose $\delta = 1$. Then the residue sequence $\bi_\t = (0,1,2,-1)$. Hence $e(\bi_\t)\psi_2^2 = e(\bi_\t)$.
\end{Example}

The next Corollary is implied by~\autoref{semi-red1}.

\begin{Corollary} \label{A:1}
Suppose $(\lambda,f) \in \widehat B_n$ and $\t \in \Tud_n(\lambda)$. If $w \in \Sym_n$ is semi-reduced correspond to $\t$, then we have $e(\bi_\t) \psi_w \psi_{w^{-1}} = e(\bi_\t)$.
\end{Corollary}

\proof Let $\s = \t{\cdot}w \in \Tud_n(\lambda)$. By the definition of semi-reduced, $w^{-1}$ is semi-reduced correspond to $\s$. Consider $ww^{-1}$ as a word in $\Sym_n$. By~\autoref{semi-red2}, $ww^{-1}$ is semi-reduced correspond to $\t$. Therefore by~\autoref{semi-red1}, we have $e(\bi_\t) \psi_w \psi_{w^{-1}} = e(\bi_\t)$. \endproof

Generally, reduced does not imply semi-reduced. The next Lemma gives a special case when reduced implies semi-reduced. Suppose $(\lambda,f) \in \widehat B_n$ and $\t \in \Tud_n(\lambda)$ with head $f$. Recall $d(\t) \in \Sym_n$ such that $\t^{(\lambda,f)} {\cdot}d(\t) = \t$.

\begin{Lemma} \label{semi-red3}
Suppose $(\lambda,f) \in \widehat B_n$ and $\t \in \Tud_n(\lambda)$ with head $f$. Then $d(\t) = s_{r_1} s_{r_2} \ldots s_{r_m} \in \Sym_{2f,n}$ is reduced implies $d(\t)$ is semi-reduced correspond to $\t^{(\lambda,f)}$.
\end{Lemma}

\proof It is straightforward by the definition of $d(\t)$. \endproof

Suppose $\t$ is an up-down tableau.~\autoref{semi-red1} and ~\autoref{semi-red3} imply $e(\bi_{(\lambda,f)}) \psi_\t$ is independent to the choice of $d(\t)$.


By the definition of $\psi_{\s\t}$, one can see that $\psi_{\s\t}$'s are homogeneous. Next we calculate the degree of $\psi_{\s\t}$'s. Let $(\lambda,f) \in \widehat B_n$ and $\s,\t \in \Tud_n(\lambda)$. First we consider the simplest case, which is $head(\s) = head(\t) = f$.


\begin{Lemma} \label{deg:tab:1}
Suppose $(\lambda,f) \in \widehat B_n$ and $\t \in \Tud_n(\lambda)$ with head $f$. Then $\deg \t = \frac{1}{2} \deg e_{(\lambda,f)} + \deg \psi_{d(\t)} e(\bi_\t) = \frac{1}{2} \deg e_{(\lambda,f)}$.
\end{Lemma}

\proof Suppose $\u \in \Tud_n(\lambda)$ with head $f$ and residue sequence $\bi_\u = (i_1, \ldots, i_n)$. By the construction of $\u$, we have $\u(k) > 0$ for any $2f+1 \leq k \leq n$. If $\s = \u{\cdot}s_k \in \Tud_n(\lambda)$ for some $s_k \in \Sym_{2f,n}$, then $\s(k) = \u(k+1) > 0$ and $\s(k+1) = \u(k) > 0$. Therefore, we have
$$
\deg \s - \deg \u = \deg(\s|_{k-1} \Rightarrow \s|_k) + \deg(\s|_k \Rightarrow \s|_{k+1}) - \deg(\u|_{k-1} \Rightarrow \u|_k) - \deg(\u|_k \Rightarrow \u|_{k+1}) = 0,
$$
by~\autoref{deg:tab:4}. Moreover, as the nodes $\u(k)$ and $\u(k+1)$ are not adjacent by~\autoref{y:h2:8}, we have $|i_k - i_{k+1}| > 1$. Therefore, we have $\deg \s - \deg \u = 0 = \deg \psi_k e(\bi)$. Hence, as $\t^{(\lambda,f)} \in \Tud_n(\lambda)$ with head $f$ and $d(\t) \in \Sym_{2f,n}$, we have $\deg \t - \deg \t^{(\lambda,f)} = 0 = \deg \psi_{d(\t)} e(\bi_\t)$.

It suffices to prove that when $\t = \t^{(\lambda,f)}$ we have $\deg \t = \frac{1}{2} \deg e_{(\lambda,f)}$. By direct calculation, we have $\deg \t^{(\lambda,f)} = 0$ if $\delta \neq 0$ and $\deg \t^{(\lambda,f)} = f$ if $\delta = 0$. Also as
$$
\deg e_{(\lambda,f)} =
\begin{cases}
2f, & \text{if $\delta = 0$,}\\
0, & \text{if $\delta \neq 0$,}
\end{cases}
$$
we have $\deg \t^{(\lambda,f)} = \frac{1}{2} \deg e_{(\lambda,f)}$, which completes the proof. \endproof

\begin{Lemma} \label{deg:tab:2}
Suppose $(\lambda,f) \in \widehat B_n$ and $\s,\t \in \Tud_n(\lambda)$ with $head(\s) = head(\t) = f$. Then we have $\deg \psi_{\s\t} = \deg \s + \deg \t$.
\end{Lemma}

\proof By the definition of $\psi_{\s\t}$'s, we have $\psi_{\s\t} = \psi_\s^* e_{(\lambda,f)} \psi_\t$. Hence, by~\autoref{deg:tab:1}, we have
\begin{align*}
\deg \psi_{\s\t} & = \deg e(\bi_\s)\psi_\s^* + \deg e_{(\lambda,f)} + \deg \psi_\t e(\bi_\t)\\
& = \frac{1}{2}\deg e_{(\lambda,f)} + \deg \psi_\s e(\bi_\s) + \frac{1}{2}\deg e_{(\lambda,f)} + \deg \psi_\s e(\bi_\t)\\
& = \deg \s + \deg \t,
\end{align*}
which proves the Lemma. \endproof

Next we extend~\autoref{deg:tab:1} to arbitrary $\t \in \Tud_n(\lambda)$ to show that
\begin{equation} \label{deg:psi}
\deg \psi_{\s\t} = \deg \s + \deg \t.
\end{equation}

Here we give some examples with the up-down tableau $\t$ in~\autoref{ex:deg:1} and~\ref{ex:deg:2}, and compare the values of $\frac{1}{2}\deg e_{(\lambda,f)} + \deg e(\bi_{(\lambda,f)}) \psi_\t \epsilon_\t$ with $\deg \t$ calculated in~\autoref{ex:deg:1} and~\ref{ex:deg:2}.

\begin{Example} \label{ex:deg:1c}
Let $n = 6$, $\lambda = (1,1)$, $\delta = 1$ and $\t = \left( \emptyset, \ydiag(1), \ydiag(2), \ydiag(2,1), \ydiag(2,2), \ydiag(2,1), \ydiag(1,1) \right) \in \Tud_n(\lambda)$. We have the standard reduction sequence $h(\t) = \t^{(2)} \rightarrow \t^{(1)} \rightarrow \t^{(0)} = \t$ where
\begin{align*}
\t^{(1)} & = \left( \emptyset, \ydiag(1), \emptyset, \ydiag(1), \ydiag(2), \ydiag(2,1), \ydiag(1,1) \right), \\
\t^{(2)} & = \left( \emptyset, \ydiag(1), \emptyset, \ydiag(1), \emptyset, \ydiag(1), \ydiag(1,1) \right),
\end{align*}
and $\rho(\t^{(2)},\t^{(1)}) = (4,6)$, $\rho(\t^{(1)},\t^{(0)}) = (4,5)$. Henceforth we have
$$
e_{(\lambda,f)} \psi_\t \epsilon_\t = e(0,0,0,0,0,-1) \epsilon_1 \epsilon_3 e(0,0,0,0,0,-1) \epsilon_4 \psi_5 e(0,0,0,1,-1,-1) \epsilon_2 \epsilon_3 \epsilon_3 e(0,1,-1,0,0,-1).
$$

By the direct calculations, the degree of elements $\frac{1}{2}\deg e_{(\lambda,f)} + \deg e(\bi_{(\lambda,f)}) \psi_\t \epsilon_\t$ is
$$
\frac{1}{2}\deg e_{(\lambda,f)} + \deg e(\bi_{(\lambda,f)}) \psi_\t \epsilon_\t = \frac{1}{2}\times 0 + 1 -2 + 1 -2 + 1= -1,
$$
which is the same as $\deg \t$.
\end{Example}

\begin{Example} \label{ex:deg:2c}
Let $n = 6$, $\lambda = (1,1)$, $\delta = 0$ and $\t = \left( \emptyset, \ydiag(1), \ydiag(2), \ydiag(2,1), \ydiag(2,2), \ydiag(2,1), \ydiag(1,1) \right) \in \Tud_n(\lambda)$. We have the standard reduction sequence $h(\t) = \t^{(2)} \rightarrow \t^{(1)} \rightarrow \t^{(0)} = \t$ where
\begin{align*}
\t^{(1)} & = \left( \emptyset, \ydiag(1), \emptyset, \ydiag(1), \ydiag(2), \ydiag(2,1), \ydiag(1,1) \right), \\
\t^{(2)} & = \left( \emptyset, \ydiag(1), \emptyset, \ydiag(1), \emptyset, \ydiag(1), \ydiag(1,1) \right),
\end{align*}
and $\rho(\t^{(2)},\t^{(1)}) = (4,6)$, $\rho(\t^{(1)},\t^{(0)}) = (4,5)$. Henceforth we have
$$
e_{(\lambda,f)} \psi_\t \epsilon_\t = e(-\frac{1}{2},\frac{1}{2},-\frac{1}{2},\frac{1}{2},-\frac{1}{2},-\frac{3}{2}) \epsilon_1 \epsilon_3 e(-\frac{1}{2},\frac{1}{2},-\frac{1}{2},\frac{1}{2},-\frac{1}{2},-\frac{3}{2}) \epsilon_4 \psi_5 e(-\frac{1}{2},\frac{1}{2},-\frac{1}{2},\frac{1}{2},-\frac{3}{2},-\frac{1}{2}) \epsilon_2 \epsilon_3 \epsilon_3 e(-\frac{1}{2},\frac{1}{2},-\frac{3}{2},-\frac{1}{2},\frac{1}{2},-\frac{1}{2}).
$$

By the direct calculations, the degree of elements $\frac{1}{2}\deg e_{(\lambda,f)} + \deg e(\bi_{(\lambda,f)}) \psi_\t \epsilon_\t$ is
$$
\frac{1}{2}\deg e_{(\lambda,f)} + \deg e(\bi_{(\lambda,f)}) \psi_\t \epsilon_\t = \frac{1}{2}\times 4 -2 + 1 - 2 + 2 - 1 = 0,
$$
which is the same as $\deg \t$.
\end{Example}

Suppose $(\lambda,f) \in \widehat B_n$ and fix $\t \in \Tud_n(\lambda)$ with head $h \leq f$. Write $\epsilon_\t = e(\bi_{h(\t)}) g_1 g_2 \ldots g_m e(\bi_\t)$ where $g_i \in \set{\psi_k, \epsilon_k | 1 \leq k \leq n-1}$ for $1 \leq i \leq m$. We will extend~\autoref{deg:tab:1} by showing $\deg \t - \deg h(\t) = \deg \epsilon_\t$ using induction on $m$. The base case, $m = 0$, is proved by~\autoref{deg:tab:1}. In order to complete the induction process, we need to show that there exists $\s \in \Tud_n(\lambda)$ such that $\epsilon_\s e(\bi_\s) g_m e(\bi_\t) = \epsilon_\t$.

Let $\t^{(f-h)} \rightarrow \t^{(f-h-1)} \rightarrow \ldots \rightarrow \t^{(0)} = \t$ be the standard reduction sequence of $\t$. Fix $i$ with $0 \leq i \leq f - h -1$ and write $\u = \t^{(i+1)}$ and $\v = \t^{(i)}$. By the definition we have $\u \rightarrow \v$.

\begin{Lemma} \label{deg:help1}
Suppose $\u, \v \in \Tud_n(\lambda)$ are defined as above. If $\rho(\u,\v) = (a,b)$, we have $\v(\l) > 0$ for any $\l$ with $2(i+h) < \l < b$. Moreover, the node $\v(\l)$ is not adjacent to $\v(a)$.
\end{Lemma}

\proof Suppose $\v = (\alpha_1, \ldots, \alpha_n)$ has remove pairs $(\alpha_{j_1}, \alpha_{k_1}), \ldots, (\alpha_{j_f}, \alpha_{k_f})$. By the definition of standard reduction sequence, we have $k_{i+h} = 2(i+h)$ and $k_{i+h+1} = b$. Therefore the Lemma follows because $\alpha_{k_1}, \ldots, \alpha_{k_f}$ are all the negative nodes of $\alpha_1, \ldots, \alpha_n$ and $k_1 < k_2 < \ldots < k_f$.

Now assume $\v(\l)$ is adjacent to $\v(a)$. Then $\v(\l)$ is either below or on the right of $\v(a)$. Hence $\v(a) \not\in \mathscr R(\v_{b-1})$, which contradicts to $\rho(\u,\v) = (a,b)$. \endproof

%

\begin{Lemma} \label{deg:help2}
Suppose $\u, \v \in \Tud_n(\lambda)$ are defined as above and write $\epsilon_{\u\rightarrow\v} = e(\bi_\u) g_1 \ldots g_m e(\bi_\v)$ where $g_i \in \set{\psi_k, \epsilon_k | 1 \leq k \leq n-1}$ for $1 \leq i \leq m$. Then there exists $\s \in \Tud_n(\lambda)$ with $\u \rightarrow \s$ such that $\epsilon_{\u\rightarrow \s} e(\bi_\s) g_m e(\bi_\v) = \epsilon_{\u \rightarrow \v}$. Moreover, the following results hold:
\begin{enumerate}
\item If $g_m = \psi_k$ for some $1 \leq k \leq n-1$, then $\s(k) < 0$, $\s(k+1) > 0$ and $\s = \v{\cdot}s_k$.

\item If $g_m = \epsilon_k$ for some $1 \leq k \leq n-1$, then $\s(k) = -\s(k+1) < 0$ and $\v(k) = -\v(k+1) > 0$.
\end{enumerate}
\end{Lemma}

\proof Suppose $\v = (\alpha_1, \ldots, \alpha_n)$ and $\rho(\u,\v) = (a,b)$. We have $\alpha_a = -\alpha_b > 0$.

(1). When $g_m = \psi_k$, we have $a < b-1 = k$ by the definition of $\epsilon_{\u\rightarrow \v}$. Hence $\alpha_k \neq \alpha_a$ and by~\autoref{deg:help1}, $\alpha_k > 0$. Therefore, we have $\v(k) = \alpha_k > 0$, $\v(k+1) = \alpha_b = -\alpha_a < 0$ and $\v(k) + \v(k+1) \neq 0$. By~\autoref{y:h2:8}, $\v{\cdot}s_k \in \Tud_n(\lambda)$. Let $\s = \v{\cdot}s_k$ and we have $\epsilon_{\u\rightarrow\s} = e(\bi_\u) g_1 \ldots g_{m-1} e(\bi_\s)$ by the definition of $\epsilon_{\u\rightarrow \s}$ and $\s(k) = \alpha_{k+1} < 0$, $\s(k +1) = \alpha_k > 0$.

(2). When $g_m = \epsilon_k$, let $head(\u) = \l$. By the definition of $\epsilon_{\u\rightarrow\v}$ we have $a = b-1 = k$ and $\epsilon_{\u\rightarrow\v} = e(\bi_\u) \epsilon_{2\l} \epsilon_{2\l+1} \ldots \epsilon_k e(\bi_\v)$. Hence we have $\v(k) = -\v(k+1) > 0$. Let $\s = (\alpha_1, \ldots, \alpha_{k-1}, -\alpha_{k-1}, \alpha_{k-1}, \alpha_{k+2}, \ldots, \alpha_n)$, by the construction we have $\epsilon_{\u\rightarrow\s} = e(\bi_\u) \epsilon_{2\l} \epsilon_{2\l+1} \ldots \epsilon_{k-1} e(\bi_\s) = e(\bi_\u) g_1 \ldots g_{m-1} e(\bi_\s)$. By~\autoref{deg:help1} we have $\alpha_{k-1} > 0$. Hence we have $\s(k) = -\s(k+1) < 0$. \endproof

The next Corollary is a direct result of~\autoref{deg:help2}.

\begin{Corollary} \label{deg:help3}
Suppose $(\lambda,f) \in \widehat B_n$ and $\t \in \Tud_n(\lambda)$. Write $\epsilon_\t = e(\bi_{h(\t)}) g_1 g_2 \ldots g_m e(\bi_\t)$ where $g_i \in \set{\psi_k, \epsilon_k | 1 \leq r \leq n-1}$ for $1 \leq i \leq m$. Then there exists $\s \in \Tud_n(\lambda)$ such that $\epsilon_\s e(\bi_\s) g_m e(\bi_\t) = \epsilon_\t$. Moreover, the following results hold:
\begin{enumerate}
\item If $g_m = \psi_k$ for some $1 \leq k \leq n-1$, then $\s(k-1) < 0$, $\s(k) > 0$ and $\s = \t{\cdot}s_k$.

\item If $g_m = \epsilon_k$ for some $1 \leq k \leq n-1$, then $\s(k) = -\s(k+1) < 0$ and $\t(k) = -\t(k+1) > 0$.
\end{enumerate}
\end{Corollary}

\autoref{deg:help3} shows that there exists $\s \in \Tud_n(\lambda)$ such that $\epsilon_\s e(\bi_\s) g_m e(\bi_\t) = \epsilon_\t$ for $g_m \in \set{\psi_k, \epsilon_k | 1 \leq k \leq n-1}$. In order to complete the induction process, we want to prove that
\begin{equation} \label{deg:eq1}
\deg \t - \deg \s = \deg \epsilon_\t - \deg \epsilon_\s = \deg e(\bi_\s) g_m e(\bi_\t).
\end{equation}

Suppose $g_m \in \{\epsilon_k, \psi_k\}$ for some $1 \leq k \leq n-1$. By~\autoref{deg:help3}, we have $\t(r) = \s(r)$ for any $r \neq k,k+1$. Therefore, we can re-write (\ref{deg:eq1}) as
\begin{equation} \label{deg:eq2}
\deg (\t|_{k-1} \Rightarrow \t|_k) + \deg (\t|_k \Rightarrow \t|_{k+1}) - \deg (\s|_{k-1} \Rightarrow \s|_k) - \deg (\s|_k \Rightarrow \s|_{k+1}) = \deg e(\bi_\s) g_m e(\bi_\t).
\end{equation}

First we prove (\ref{deg:eq2}) when $g_m = \epsilon_k$ for some $1 \leq k \leq n-1$. Notice that $\deg e(\bi_\s) \epsilon_k e(\bi_\t) = \deg_k(\bi_\s) + \deg_k(\bi_\t)$. The following results connect $\deg(\t|_{k-1} \Rightarrow \t|_k)$ and $\deg_k(\bi_\t)$.

\begin{Lemma} \label{deg:tab:3}
Suppose $\t$ is an up-down tableau of size $n$ and $1 \leq k \leq n$. If $\t_k \subset \t_{k-1}$, we have
$$
\deg (\t|_{k-1} \Rightarrow \t|_k) = -\deg_k(\bi_\t).
$$
\end{Lemma}

\proof Denote $\lambda = \t_{k-1}$ and $\mu = \t_k$. As $\mu \subset \lambda$, there exists a positive node $\alpha$ such that $\mu = \lambda \backslash \{\alpha\}$. Moreover, if we write $\bi_\t = (i_1, \ldots, i_n)$, we have $i_k = -\res(\alpha) = \res_\lambda(\alpha)$. Because $\bi_\t \in I^n = I_{k,+}^n \sqcup I_{k,-}^n \sqcup I_{k,0}^n$, we prove the Lemma by considering the following cases:

\textbf{Case 1:} $\bi_\t \in I_{k,0}^n$ and $i_k = -\res(\alpha) = \res_\lambda(\alpha) \neq -\frac{1}{2}$.

When $i_k \neq 0, \frac{1}{2}$, we have $h_k(\bi_\t) = -1$. By (\ref{deg:h4:eq2}) we have $\mathscr{AR}_\lambda(i_k) = \{\alpha\} \subset \mathscr R(\lambda)$ and $\mathscr{AR}_\lambda(-i_k) = \emptyset$. As $\mu = \lambda \backslash \{\alpha\}$ and $\res(\alpha) \neq \pm \frac{1}{2}$, by the construction we have $\mathscr {\widehat A}_\t(k) = \mathscr {\widehat R}_\t(k) = \emptyset$. Therefore
$$
\deg (\t|_{k-1} \Rightarrow \t|_k) = |\mathscr {\widehat A}_\t(k)| - |\mathscr {\widehat R}_\t(k)| + \delta_{\res(\alpha),-\frac{1}{2}} = 0 = -\deg_k(\bi_\t).
$$

When $i_k = 0$, we have $h_k(\bi_\t) = 0$. By~\autoref{deg:h4} we have $\mathscr{AR}_\lambda(i_k) = \{\alpha\} \subset \mathscr R(\lambda)$. As $\mu = \lambda \backslash \{\alpha\}$, by the construction we have $\mathscr {\widehat A}_\t(k) = \mathscr {\widehat R}_\t(k) = \emptyset$. Therefore
$$
\deg (\t|_{k-1} \Rightarrow \t|_k) = |\mathscr {\widehat A}_\t(k)| - |\mathscr {\widehat R}_\t(k)| + \delta_{\res(\alpha),-\frac{1}{2}} = 0 = -\deg_k(\bi_\t).
$$

When $i_k = \frac{1}{2}$, we have $\res(\alpha) = -\frac{1}{2}$ and $h_k(\bi_\t) = -1$. Write $\alpha = (i,j) \in [\lambda]$. By (\ref{deg:h4:eq2}) we have $\mathscr{AR}_\lambda(i_k) = \{\alpha\} \subset \mathscr R(\lambda)$. Hence $(i,j+1) \not\in \mathscr A(\lambda)$, which implies $(i-1, j+1) \not \in [\lambda]$. Hence we have $(i,j), (i-1,j+1) \not\in [\mu]$. Therefore $\gamma = (i-1,j) \in \mathscr R(\mu)$, and $\res(\gamma) = \frac{1}{2} = i_k$. Hence $\mathscr {\widehat A}_\t(k) = \emptyset$ and $\mathscr {\widehat R}_\t(k) = \{\gamma\}$, which yields
$$
\deg (\t|_{k-1} \Rightarrow \t|_k) = |\mathscr {\widehat A}_\t(k)| - |\mathscr {\widehat R}_\t(k)| + \delta_{\res(\alpha),-\frac{1}{2}} = 0 = -\deg_k(\bi_\t).
$$

\textbf{Case 2:} $\bi_\t \in I_{k,0}^n$ and $i_k = -\res(\alpha) = \res_\lambda(\alpha) = -\frac{1}{2}$.

As $\bi_\t \in I_{k,0}^n$ and $i_k = -\frac{1}{2}$, we have $h_k(\bi_\t) = -2$. Write $\alpha = (i,j)$ and $\gamma = (i+1,j)$. By~\autoref{deg:h4} we have $\mathscr{AR}_\lambda(i_k) = \{\alpha,\gamma\}$ where $\alpha \in \mathscr R(\lambda)$ and $\gamma \in \mathscr A(\lambda)$. As $\alpha \not\in [\mu]$, one can see that $\gamma \not\in \mathscr{AR}(\mu)$. As $\gamma \in \mathscr A(\lambda)$, we have $(i+1,j-1) \in [\lambda] \cap [\mu]$. Because $\alpha \not\in [\mu]$, we have $(i,j-1) \not\in \mathscr{AR}(\mu)$. Therefore we have $\mathscr {\widehat A}_\t(k) = \mathscr {\widehat R}_\t(k) = \emptyset$, which yields
$$
\deg (\t|_{k-1} \Rightarrow \t|_k) = |\mathscr {\widehat A}_\t(k)| - |\mathscr {\widehat R}_\t(k)| + \delta_{\res(\alpha),-\frac{1}{2}} = 0 = -\deg_k(\bi_\t).
$$

\textbf{Case 3:} $\bi_\t \in I_{k,-}^n$ and $i_k = -\res(\alpha) = \res_\lambda(\alpha) \neq \frac{1}{2}$.

As $\bi_\t \in I_{k,-}^n$, we have $h_k(\bi_\t) = -2$. By~\autoref{deg:h4} we have $\mathscr{AR}_\lambda(i_k) = \{\alpha,\gamma\}$ where $\gamma \in \mathscr A(\lambda)$ with $\res(\gamma) = i_k$, and $\mathscr{AR}_\lambda(-i_k) = \emptyset$. Because $i_k \neq \pm \frac{1}{2}$, we have $\gamma \in \mathscr A(\mu)$, which implies $\mathscr {\widehat A}_\t(k) = \{\gamma\}$ and $\mathscr {\widehat R}_\t(k) = \emptyset$. Therefore
$$
\deg (\t|_{k-1} \Rightarrow \t|_k) = |\mathscr {\widehat A}_\t(k)| - |\mathscr {\widehat R}_\t(k)| + \delta_{\res(\alpha),-\frac{1}{2}} = 1 = -\deg_k(\bi_\t).
$$

\textbf{Case 4:} $\bi_\t \in I_{k,-}^n$ and $i_k = -\res(\alpha) = \res_\lambda(\alpha) = \frac{1}{2}$.

As $\bi_\t \in I_{k,-}^n$, we have $h_k(\bi_\t) = -2$. Write $\alpha = (i,j)$ and $\gamma = (i,j+1)$. By~\autoref{deg:h4} we have $\mathscr{AR}_\lambda(i_k) = \{\alpha,\gamma\}$ where $\alpha \in \mathscr R(\lambda)$ and $\gamma \in \mathscr A(\lambda)$. As $\alpha \not\in [\mu]$, one can see that $\gamma \not\in \mathscr{AR}(\mu)$. As $\gamma \in \mathscr A(\lambda)$, we have $(i-1,j+1) \in [\lambda]\cap[\mu]$. Because $\alpha \not\in [\mu]$, we have $(i-1,j) \not\in \mathscr{AR}(\mu)$. Therefore we have $\mathscr {\widehat A}_\t(k) = \mathscr {\widehat R}_\t(k) = \emptyset$, which implies
$$
\deg (\t|_{k-1} \Rightarrow \t|_k) = |\mathscr {\widehat A}_\t(k)| - |\mathscr {\widehat R}_\t(k)| + \delta_{\res(\alpha),-\frac{1}{2}} = 1 = -\deg_k(\bi_\t).
$$

\textbf{Case 5:} $\bi_\t \in I_{k,+}^n$ and $i_k = -\res(\alpha) = \res_\lambda(\alpha) \neq -\frac{1}{2}$.

As $\bi_\t \in I_{k,+}^n$ and $i_k \neq -\frac{1}{2}$, we have $h_k(\bi_\t) = 0$. Hence by~\autoref{deg:h4} we have $\mathscr{AR}_\lambda(-i_k) = \{\gamma\}$ and $\mathscr{AR}_\lambda(i_k) = \{\alpha\}$, where $\gamma \in \mathscr R(\lambda)$ and $\res(\gamma) = i_k = -\res(\alpha)$. As $i_k \neq \pm \frac{1}{2}$, we have $\gamma \in \mathscr R(\mu)$, which implies $\mathscr {\widehat A}_\t(k) = \emptyset$ and $\mathscr {\widehat R}_\t(k) = \{\gamma\}$. Therefore we have
$$
\deg (\t|_{k-1} \Rightarrow \t|_k) = |\mathscr {\widehat A}_\t(k)| - |\mathscr {\widehat R}_\t(k)| + \delta_{\res(\alpha),-\frac{1}{2}} = -1 = -\deg_k(\bi_\t).
$$

\textbf{Case 6:} $\bi_\t \in I_{k,+}^n$ and $i_k = -\res(\alpha) = \res_\lambda(\alpha) = -\frac{1}{2}$.

As $\bi_\t \in I_{k,+}^n$ and $i_k = \frac{1}{2}$, we have $h_k(\bi_\t) = -1$. Write $\alpha = (i,j) \in [\lambda]$. By (\ref{deg:h4:eq2}) we have $\mathscr{AR}_\lambda(i_k) = \{\alpha\} \subset \mathscr R(\lambda)$ and $\mathscr{AR}_\lambda(-i_k) = \emptyset$. Hence $(i+1,j) \not\in \mathscr A(\lambda)$, which implies $(i+1,j-1) \not\in [\lambda]$. Hence we have $(i,j), (i+1,j-1) \not\in [\mu]$. Therefore $\gamma = (i,j-1) \in \mathscr R(\mu)$, and $\res(\gamma) = -\frac{1}{2} = i_k$. Hence $\mathscr {\widehat A}_\t(k) = \emptyset$ and $\mathscr {\widehat R}_\t(k) = \{\gamma\}$, which implies
$$
\hspace*{3cm} \deg (\t|_{k-1} \Rightarrow \t|_k) = |\mathscr {\widehat A}_\t(k)| - |\mathscr {\widehat R}_\t(k)| + \delta_{\res(\alpha),-\frac{1}{2}} = -1 = -\deg_k(\bi_\t). \hspace*{3cm} \qedhere
$$

\begin{Lemma} \label{deg:tab:4}
Suppose $\t$ is an up-down tableau of size $n$ and $0 \leq k < n$. If $\t_k \subset \t_{k+1}$, we have
$$
\deg (\t|_k \Rightarrow \t|_{k+1}) = 0.
$$
\end{Lemma}

\proof It is easy to verify because $|\mathscr A_\t(k)| = |\mathscr R_\t(k)| = 0$. \endproof

\begin{Lemma} \label{deg:tab:5}
Suppose $\bi \in I^n$ and $1 \leq k \leq n-2$. If $i_k = - i_{k+1} = i_{k+2}$, we have $\deg_k(\bi) = - \deg_{k+1}(\bi)$.
\end{Lemma}

\proof Because $i_k = -i_{k+1}$, by applying (\ref{remark:h:eq2}) to the definitions of $I_{k,0}^n$, $I_{k,-}^n$ and $I_{k,+}^n$, we can see that $\bi \in I_{k,0}^n$ implies $\bi \in I_{k+1,0}^n$, $\bi \in I_{k,-}^n$ implies $\bi \in I_{k+1,+}^n$ and $\bi \in I_{k,+}^n$ implies $\bi \in I_{k+1,-}^n$. Hence we have $\deg_k(\bi) = - \deg_{k+1}(\bi)$ by the definition of $\deg_k(\bi)$. \endproof

\begin{Corollary} \label{deg:tab:6}
Suppose $\t$ is an up-down tableau of size $n$ and $1 \leq k < n$. If $\t(k) + \t(k+1) = 0$, we have
$$
\deg(\t|_{k-1} \Rightarrow \t|_k) + \deg(\t|_k \Rightarrow \t|_{k+1}) =
\begin{cases}
\deg_k(\bi_\t), & \text{if $\t_{k-1} \subset \t_k$,}\\
-\deg_k(\bi_\t), & \text{if $\t_k \subset \t_{k-1}$.}
\end{cases}
$$
\end{Corollary}

\proof When $\t_k \subset \t_{k-1}$, by~\autoref{deg:tab:3} and~\autoref{deg:tab:4}, we have
$$
\deg(\t|_{k-1} \Rightarrow \t|_k) + \deg(\t|_k \Rightarrow \t|_{k+1}) =-\deg_k(\bi_\t).
$$

When $\t_{k-1} \subset \t_k$, because $\t(k) + \t(k+1) = 0$, we have $\t_{k+1} = \t_{k-1} \subset \t_k$. By~\autoref{deg:tab:3} and~\autoref{deg:tab:4}, we have
$$
\deg(\t|_{k-1} \Rightarrow \t|_k) + \deg(\t|_k \Rightarrow \t|_{k+1}) =-\deg_{k+1}(\bi_\t),
$$
and by~\autoref{deg:tab:5} we have $\deg_k(\bi_\t) = - \deg_{k+1}(\bi_\t)$, which completes the proof. \endproof

Now we are ready to prove (\ref{deg:eq2}) when $g_m = \epsilon_k$.

\begin{Lemma} \label{deg:tab:7}
Suppose $(\lambda,f) \in \widehat B_n$ and $\s, \t \in \Tud_n(\lambda)$ such that $\epsilon_\t = \epsilon_\s e(\bi_\s) \epsilon_k e(\bi_\t)$ for $1 \leq k \leq n-1$. Then the equality (\ref{deg:eq2}) holds.
\end{Lemma}

\proof By~\autoref{deg:help3}, we have $\s(k) = -\s(k+1) < 0$ and $\t(k) = -\t(k+1) > 0$. Hence by~\autoref{deg:tab:6}, we have
\begin{align*}
& \deg(\s|_{k-1} \Rightarrow \s|_k) + \deg(\s|_k \Rightarrow \s|_{k+1}) = -\deg_k(\bi_\s),\\
& \deg(\t|_{k-1} \Rightarrow \t|_k) + \deg(\t|_k \Rightarrow \t|_{k+1}) = \deg_k(\bi_\t).
\end{align*}

As $\s(r) = \t(r)$ for any $r \neq k, k+1$, we have
\begin{eqnarray*}
&& \deg (\t|_{k-1} \Rightarrow \t|_k) + \deg (\t|_k \Rightarrow \t|_{k+1}) - \deg (\s|_{k-1} \Rightarrow \s|_k) - \deg (\s|_k \Rightarrow \s|_{k+1})\\
& = & \deg_k(\bi_\s) + \deg_k(\bi_\t) = \deg e(\bi_\s) \epsilon_k e(\bi_\t),
\end{eqnarray*}
which proves the Lemma. \endproof

Then we prove (\ref{deg:eq2}) when $g_m = \psi_k$.

\begin{Lemma} \label{deg:tab:8}
Suppose $(\lambda,f) \in \widehat B_n$ and $\s, \t \in \Tud_n(\lambda)$ such that $\epsilon_\t = \epsilon_\s e(\bi_\s) \psi_k e(\bi_\t)$ for $1 \leq k \leq n-1$. Then the equality (\ref{deg:eq2}) holds.
\end{Lemma}

\proof By~\autoref{deg:help3}, we have $s(k-1) < 0$, $\s(k) > 0$ and $\t = \s{\cdot}s_k$. By~\autoref{deg:tab:3} and~\autoref{deg:tab:4}, we have
\begin{eqnarray} \label{deg:tab:4:eq1}
&& \deg (\t|_{k-1} \Rightarrow \t|_k) + \deg (\t|_k \Rightarrow \t|_{k+1}) - \deg (\s|_{k-1} \Rightarrow \s|_k) - \deg (\s|_k \Rightarrow \s|_{k+1}) \notag \\
& = & \deg(\t|_k \Rightarrow \t|_{k+1}) - \deg(\s|_{k-1} \Rightarrow \s|_k) = \deg_k(\bi_\s) - \deg_{k+1}(\bi_\t).
\end{eqnarray}

Let $\bi_\s = (i_1, \ldots, i_n)$. Notice $\bi_\t = \bi_\s{\cdot}s_k$. If $|i_k - i_{k+1}| > 1$, we have $h_k(\bi_\s) = h_{k+1}(\bi_\t)$, which implies that $\bi_\s \in I_{k,a}^n$ if and only if $\bi_\t \in I_{k+1,a}^n$ for $a \in \{+,-,0\}$. Hence we have $\deg_k(\bi_\s) - \deg_{k+1}(\bi_\t) = 0$. By (\ref{deg:tab:4:eq1}), the Lemma holds.

If $i_k = i_{k+1} \pm 1$, first we exclude some of the cases: when $i_k = 0$ and $\s(k) < 0$, we always $\s(k+1) < 0$ when $i_k = i_{k+1} \pm 1$, which does not satisfy the condition of the Lemma; and when $i_k = \pm \frac{1}{2}$ and $i_{k+1} = - i_k$, we always have $\s(k) + \s(k+1) = 0$ or $\s(k+1) < 0$, which does not satisfy the condition of the Lemma; and when $h_k(\bi_\s) = -2$, by the construction of $\s$ there exists no up-down tableau $\s$ satisfying the conditions of the Lemma. By excluding these cases, we have $h_{k+1}(\bi_\t) = h_k(\bi_\s) - 1$. By direct calculations, we have $\bi_\s \in I_{k,0}^n$ if and only if $\bi_\t \in I_{k+1,-}^n$, and $\bi_\s \in I_{k,+}^n$ if and only if $\bi_\t \in I_{k+1,0}^n$. Hence we have $\deg_k(\bi_\s) - \deg_{k+1}(\bi_\t) = 1$. By (\ref{deg:tab:4:eq1}), the Lemma holds.

If $i_k = i_{k+1}$, by the construction of $\s$ we have $h_k(\bi_\s) = -2$. First we exclude some of the cases: when $i_k = i_{k+1} = \pm \frac{1}{2}$, we have $h_{k+1}(\bi_\s) = h_k(\bi_\s) + 3 = 1$, which implies $\bi_\s \not\in I^n$ by~\autoref{deg:h2}; and when $i_k = 0$, we have $h_k(\bi_\s) = 0$, which contradicts that $h_k(\bi_\s) = -2$. Therefore we have $i_k = i_{k+1} \neq 0, \pm \frac{1}{2}$, which yields $h_{k+1}(\bi_\t) = h_k(\bi_\s) + 2 = 0$ by (\ref{remark:h:eq1}). Hence, we have $\bi_\s \in I_{k,-}^n$, $\bi_\t \in I_{k+1,+}^n$, and $\deg_k(\bi_\s) - \deg_{k+1}(\bi_\t) = -2$. By (\ref{deg:tab:4:eq1}), the Lemma holds.

Therefore, we have considered all the cases, which completes the proof. \endproof

Now we are able to give a proper proof for the equality~\eqref{deg:psi}.


\begin{Proposition} \label{deg:tab}
Suppose $(\lambda,f) \in \widehat B_n$ and $\s,\t \in \Tud_n(\lambda)$. We have $\deg \psi_{\s\t} = \deg \s + \deg \t$.
\end{Proposition}

\proof For any $\t \in \Tud_n(\lambda)$, let $head(\t) = h$. By the definition, we can write $\epsilon_\t = e(\bi_{h(\t)}) g_1 \ldots g_m e(\bi_\t)$ where $g_i \in \set{\psi_k, \epsilon_k | 1 \leq k \leq n-1}$ for $1 \leq i \leq m$. First we prove that
\begin{equation} \label{deg:eq3}
\deg \t - \deg h(\t) = \deg \epsilon_\t.
\end{equation}

We apply induction here. For base case, it is obvious that when $m = 0$ the equality (\ref{deg:eq1}) holds, because when $m = 0$, $\t = h(\t)$ and $\epsilon_\t = e(\bi_\t)$. For induction process, assume the equality holds when $m < m'$. Let $m = m'$ and $\s \in \Tud_n(\lambda)$ such that $\epsilon_\s e(\bi_\s) g_m e(\bi_\t) = \epsilon_\t$ for $g_m \in \set{\psi_k, \epsilon_k | 1 \leq k \leq n-1}$. By induction, we have $\deg \s - \deg h(\s) = \deg \epsilon_\s$. Notice that we have $h(\s) = h(\t)$. Hence, by~\autoref{deg:tab:7} and~\autoref{deg:tab:8}, (\ref{deg:eq2}) holds, which implies (\ref{deg:eq1}) holds. Therefore, we have
$$
\deg \t - \deg h(\t) = \deg \t - \deg \s + \deg \s - \deg h(\s) = \deg e(\bi_\s) g_m e(\bi_\t) + \deg \epsilon_\s = \deg \epsilon_\t,
$$
which completes the induction process. Therefore, the equality (\ref{deg:eq3}) holds.

Now assume $\s,\t \in \Tud_n(\lambda)$. By~\autoref{deg:tab:2}, we have
\begin{align*}
\deg \psi_{\s\t} & = \deg \epsilon_\s^* \psi_{h(\s) h(\t)} \epsilon_\t = \deg \epsilon_\s + \deg \psi_{h(\s) h(\t)} + \deg \epsilon_\t\\
& = \deg \s - \deg h(\s) + \deg h(\s) + \deg h(\t) + \deg \t - \deg h(\t) = \deg \s + \deg \t,
\end{align*}
which completes the proof. \endproof

\subsection{The induction property} \label{sec:ind}

In Section~\ref{sec:basis} we constructed a set of homogeneous elements $\{\psi_{\s\t}\}$ in $\G{n}$. Define
$$
R_n(\delta) := \set{a \in \G{n} | a = \sum_{\s,\t \in \Tud_n(\lambda)} c_{\s\t} \psi_{\s\t} \text{ where $c_{\s\t} \in R$ and $(\lambda,f) \in \widehat B_n$}}.
$$

It is easy to see that $R_n(\delta)$ is a subspaces of $\G{n}$. In this subsection, we prove the induction property of $R_n(\delta)$. The result of this subsection can be directly implied to $\G{n}$ after we prove $R_n(\delta) = \G{n}$ at the end of Section~\ref{sec:span}.

By the definition of $\G{n}$, we can consider $\G{n}$ as a subalgebra of $\G{n+1}$ by identifying $e(\bi) = \sum_{k \in P} e(\bi \vee k)$ for $\bi \in P^n$. Hence we have a sequence
$$
\G{1} \subset \G{2} \subset \G{3} \subset \ldots.
$$

For each $i \in P$, define
$$
e_{n,i} := \sum_{\bj \in P^n} e(\bj \vee i) \in \G{n+1}.
$$

Similar to the cyclotomic Khovnov-Lauda-Rouquier algebras, there is a (non-unital) embedding of $\theta_i^{(n)}:\G{n} \hookrightarrow \G{n+1}$ given by
$$
e(\bj) \mapsto e(\bj \vee i), \quad y_r \mapsto e_{n,i} y_r, \quad \psi_s \mapsto e_{n,i} \psi_s \quad \text{and} \quad \epsilon_s \mapsto e_{n,i} \epsilon_s,
$$
for $\bj \in I^n$, $2f+1 \leq r \leq n$ and $2f+1 \leq s \leq n-1$. As $R_n(\delta)$ is a $R$-subspace of $\G{n}$, we can restrict $\theta_i^{(n)}$ to $R_n(\delta)$.

First we set up the notations and definitions we are going to use. Suppose $0 \leq f \leq \lfloor \frac{n}{2} \rfloor$. Define
\begin{align*}
R_n^{\geq f}(\delta) & := \set{a \in \G{n} | \text{$a = \sum_{\s,\t \in \Tud_n(\mu)} c_{\s\t} \psi_{\s\t}$ where $c_{\s\t} \in R$, $(\mu,m) \in \widehat B_n$ and $m \geq f$}},\\
R_n^{> f}(\delta) & := \set{a \in \G{n} | \text{$a = \sum_{\s,\t \in \Tud_n(\mu)} c_{\s\t} \psi_{\s\t}$ where $c_{\s\t} \in R$, $(\mu,m) \in \widehat B_n$ and $m > f$}}.
\end{align*}

For $a, b \in \G{n}$, we write $a \equiv b \pmod{R_n^{> f}(\delta)}$ if $a - b \in R_n^{> f}(\delta)$, and $a \equiv b \pmod{R_n^{\geq f}(\delta)}$ in a similar way. Next we define a subset of $\bigcup_{i \geq 1} \widehat B_n$.

\begin{Definition}
Define $\widehat {\mathscr B}$ to be the subset of $\bigcup_{i \geq 1} \widehat B_n$ such that $(\lambda,f) \in \widehat {\mathscr B}$ with $\lambda \vdash n - 2f$ if and only if for any $\s,\t \in \Tud_n(\lambda)$ and $a \in \G{n}$, we have
$$
\psi_{\s\t}{\cdot}a \equiv \sum_{\v \in \Tud_n(\lambda)} c_\v \psi_{\s\v} \pmod{R_n^{>f}(\delta)}.
$$
\end{Definition}

\begin{Remark}
By applying the involution $*$ it is easy to see that the above definition is equivalent to say that $(\lambda,f) \in \widehat {\mathscr B}$ with $\lambda \vdash n - 2f$ if and only if for any $\s,\t \in \Tud_n(\lambda)$ and $b \in \G{n}$, we have
$$
b{\cdot}\psi_{\s\t} \equiv \sum_{\u \in \Tud_n(\lambda)} c_\u \psi_{\u\t} \pmod{R_n^{>f}(\delta)}.
$$
\end{Remark}

We define a total ordering on $\bigcup_{i \geq 1} \widehat B_i$ extended by the lexicographic ordering on $\widehat B_n$. Suppose $(\lambda,f), (\mu,m) \in \bigcup_{i \geq 1} \widehat B_i$. We denote $(\mu,m) \geq (\lambda,f)$ if $|\mu| + 2m < |\lambda| + 2f$; or $|\mu| + 2m = |\lambda| + 2f$ and $m > f$; or $|\mu| + 2m = |\lambda| = 2f$ and $m = f$, and $\mu \geq \lambda$. Define $(\mu,m) > (\lambda,f)$ if $(\mu,m) \geq (\lambda,f)$ and $(\mu,m) \neq (\lambda,f)$.

\begin{Definition} \label{def:s}
Define $\mathscr S_n = \set{(\lambda,f) \in \widehat B_n | \text{$(\mu, m) \in \widehat {\mathscr B}$ whenever $(\mu,m) \in \bigcup_{i \geq 1} \widehat B_i$ and $(\mu,m) > (\lambda,f)$}}$.
\end{Definition}

The next Lemma is the key point of $\mathscr S_n$.

\begin{Lemma} \label{two:sided:ideal}
Suppose $0 \leq f \leq \lfloor \frac{n}{2} \rfloor$. If there exists $\sigma \vdash n - 2f$ such that $(\sigma,f) \in \mathscr S_n$, then $R_n^{>f}(\delta)$ is a two-sided $\G{n}$-ideal.
\end{Lemma}

\proof It is directly implied by the definitions of $\mathscr S_n$ and $R_n^{>f}(\delta)$. \endproof

Now we start to prove the induction property of $R_n(\delta)$. Suppose $(\lambda,f) \in \widehat B_{n-1}$. Let $\alpha \in \mathscr A(\lambda)$ and $\mu = \lambda \cup \{\alpha\}$. Define $\t \in \Tud_{n-1}(\lambda)$ and $\s \in \Tud_n(\mu)$ where $\s|_{n-1} = \t$. Write $\mu = (\mu_1, \mu_2, \ldots, \mu_m)$ and $\alpha = (r, \mu_r)$. The next Lemma explores the connection between $\psi_\t\epsilon_\t$ and $\psi_\s \epsilon_\s$.

\begin{Lemma} \label{I:h1:1}
Suppose $\s, \t$ are defined as above. Let $a = 2f + \mu_1 + \mu_2 + \ldots + \mu_r$ and $\res(\alpha) = i \in P$. We have $\psi_\s = \psi_a \psi_{a+1} \ldots \psi_{n-1} \theta_i^{(n-1)}(\psi_\t)$ and $\epsilon_\s = \theta_i^{(n-1)}(\epsilon_\t)$.
\end{Lemma}

\proof By the construction of $\epsilon_\t$, one can see that $\epsilon_\t$ depends on the remove pairs of $\t$ only. Because $\s(n) = \alpha > 0$, the remove pairs of $\t$ and $\s$ are the same. Hence we have $\epsilon_\s = \theta_i^{(n-1)}(\epsilon_\t)$.

We define $\u = h(\t) \in \Tud_{n-1}(\lambda)$ and $\v = h(\s) \in \Tud_n(\mu)$. Let $\w = \t^{(\mu,f)} s_a s_{a+1} \ldots s_{n-1}$. By the construction one can see that $\w(n) = \v(n) = \alpha$ and $\w|_{n-1} = \t^{(\lambda,f)}$. Because $\t = \s|_{n-1}$ and $\s(n) = \alpha > 0$, we have $\u = \v|_{n-1}$ and $\u(n) = \alpha = \w(n)$, which implies
$$
\v = \w{\cdot} d(\u) = \t^{(\lambda,f)} s_a s_{a+1} \ldots s_{n-1} d(\u).
$$

Hence we have $\psi_\s = \psi_a \psi_{a+1} \ldots \psi_{n-1} \theta_i^{(n-1)}(\psi_\t)$, which completes the proof. \endproof


\begin{Example} \label{I:h1:ex1}
Suppose $n = 8$ and $\delta = 1$. Let $(\lambda,f) = (\ \ydiag(2,1), 2) \in \widehat B_{n-1}$, the node $\alpha = (3,1)$ and $i = \res(\alpha) = 2 \in P$. Define $\mu = \lambda \cup \{\alpha\} = \ \ydiag(3,1)$ and the following up-down tableaux
\begin{align*}
\u & = \left( \emptyset, \ydiag(1), \ydiag(2), \ydiag(2,1), \ydiag(2), \ydiag(3), \ydiag(3,1), \ydiag(2,1) \right) \in \Tud_n(\lambda), \\
\v & = \left( \emptyset, \ydiag(1), \ydiag(1,1), \ydiag(2,1), \ydiag(3,1), \ydiag(4,1), \ydiag(3,1), \ydiag(2,1) \right) \in \Tud_n(\lambda), \\
\s & = \left( \emptyset, \ydiag(1), \ydiag(2), \ydiag(2,1), \ydiag(2), \ydiag(3), \ydiag(3,1), \ydiag(2,1), \ydiag(3,1) \right) \in \Tud_{n+1}(\mu), \\
\t & = \left( \emptyset, \ydiag(1), \ydiag(1,1), \ydiag(2,1), \ydiag(3,1), \ydiag(4,1), \ydiag(3,1), \ydiag(2,1), \ydiag(3,1) \right) \in \Tud_{n+1}(\mu).
\end{align*}

One can see that $\s|_{n-1} = \u$, $\t|_{n-1} = \v$ and $\s(n) = \t(n) = \alpha$. By direct calculations, we have $\psi_\u = 1$, $\psi_\v = \psi_6$, $\psi_\s = \psi_7$, $\psi_\t = \psi_7 \psi_6$ and
\begin{align*}
\epsilon_\u & = e(0,0,0,0,0,1,-1) \epsilon_4 \epsilon_5 e(0,0,0,1,-1,1,-1) \epsilon_2 \epsilon_3 \epsilon_4 \epsilon_5 \psi_6 e(0,1,-1,1,2,-1,-2), \\
\epsilon_\v & = e(0,0,0,0,0,-1,1) \epsilon_4 \epsilon_5 \epsilon_6 e(0,0,0,-1,1,2,-2) \epsilon_2 \epsilon_3 \epsilon_4 \epsilon_5 e(0,-1,1,2,3,-3,-2), \\
\epsilon_\s & = e(0,0,0,0,0,1,-1,2) \epsilon_4 \epsilon_5 e(0,0,0,1,-1,1,-1,2) \epsilon_2 \epsilon_3 \epsilon_4 \epsilon_5 \psi_6 e(0,1,-1,1,2,-1,-2,2), \\
\epsilon_\t & = e(0,0,0,0,0,-1,1,2) \epsilon_4 \epsilon_5 \epsilon_6 e(0,0,0,-1,1,2,-2,2) \epsilon_2 \epsilon_3 \epsilon_4 \epsilon_5 e(0,-1,1,2,3,-3,-2,2).
\end{align*}

Because $\mu = (3,1)$, let $a = 2f + \mu_1 = 7$. Hence we have $\psi_\s = \psi_a \psi_{a+1} \ldots \psi_{n-1} \theta^{(n-1)}_i(\psi_\u) = \psi_7 \theta^{(n-1)}_i(\psi_\u)$, $\psi_\t = \psi_a \psi_{a+1} \ldots \psi_{n-1} \theta^{(n-1)}_i(\psi_\v) = \psi_7 \theta^{(n-1)}_i(\psi_\v)$, $\epsilon_\s^* = \theta^{(n-1)}_i(\epsilon_\u^*)$ and $\epsilon_\t = \theta^{(n-1)}_i(\epsilon_\v)$.
\end{Example}

The following two results give the induction property of $R_n(\delta)$.


\begin{Lemma} \label{I:h1}
Suppose $(\lambda,f) \in \widehat B_{n-1}$. If $i \in P$ such that $i = \res(\alpha)$ for some $\alpha \in \mathscr A(\lambda)$, then we have $\theta_i^{(n-1)}(\psi_{\u\v}) = \psi_{\s\t}$ for any $\u,\v \in \Tud_{n-1}(\lambda)$. Moreover, we have $\s|_{n-1} = \u$, $\t|_{n-1} = \v$ and $\s(n) = \t(n) = \alpha$.
\end{Lemma}

\proof First we show that when $\u = \v = \t^{(\lambda,f)}$, we have $\theta_i^{(n-1)}(\psi_{\u\v}) = \psi_{\x\x}$ where $\x|_{n-1} = \t^{(\lambda,f)}$ and $\x(n) = \alpha$.

Set $\mu = \lambda \cup \{\alpha\} = (\mu_1, \mu_2, \ldots, \mu_m)$. One can see that $\alpha = (r, \mu_r)$ for some $1 \leq r \leq m$. Hence let $a = 2f + \mu_1 + \mu_2 + \ldots + \mu_r$. We have a reduced expression $d(\x) = s_a s_{a+1} \ldots s_{n-1} \in \Sym_{2f,n+1}$.

Let $\bi_{(\lambda,f)} = (i_1, i_2, \ldots, i_{n-1})$. By the construction of $\lambda$ and $\mu$, as $i = \res(\alpha)$, we have $|i - i_s| \geq 2$ for any $a \leq s \leq n-1$ and $\bi_{(\mu,f)} = (i_1, \ldots, i_{a-1}, i, i_a, i_{a+1}, \ldots, i_{n-1})$. Hence by (\ref{rela:4}) we have
\begin{align*}
\theta_i^{(n-1)}(\psi_{\t^{(\lambda,f)} \t^{(\lambda,f)}}) & = e(\bi_{(\lambda,f)} \vee i) \epsilon_1 \epsilon_3 \ldots \epsilon_{2f-1} e(\bi_{(\lambda,f)} \vee i)\\
& = e(\bi_{(\lambda,f)} \vee i) \epsilon_1 \epsilon_3 \ldots \epsilon_{2f-1} \psi_n \psi_{n-1} \ldots \psi_a e(\bi_{(\mu,f)}) \psi_a \ldots \psi_{n-1} \psi_n\\
& = \psi_n \psi_{n-1} \ldots \psi_a e(\bi_{(\mu,f)}) \epsilon_1 \epsilon_3 \ldots \epsilon_{2f-1} e(\bi_{(\mu,f)}) \psi_a \ldots \psi_{n-1} \psi_n\\
& = \psi_{d(\x)}^* e_{(\mu,f)} \psi_{d(\x)} = \psi_{\x\x}.
\end{align*}

Now for any $\u,\v \in \Tud_{n-1}(\lambda)$, set $\s, \t$ to be up-down tableaux such that $\s|_{n-1} = \u$, $\t|_{n-1} = \v$ and $\s(n) = \t(n) = \alpha$. By~\autoref{I:h1:1}, we have $\psi_\s = \psi_{d(\x)} \theta_i^{(n-1)}(\psi_\u)$, $\psi_\t = \psi_{d(\x)} \theta_i^{(n-1)}(\psi_\v)$, $\epsilon_\s^* = \theta^{(n-1)}_i(\epsilon_\u^*)$ and $\epsilon_\t = \theta^{(n-1)}_i(\epsilon_\v)$. Hence
\begin{align*}
\theta^{(n-1)}_i(\psi_{\u\v}) & = \theta^{(n-1)}_i(\epsilon_\u^* \psi_\u^* e_{(\lambda,f)} \psi_\v \epsilon_\v) = \theta^{(n-1)}_i(\epsilon_\u^* \psi_\u^*) \theta^{(n-1)}_i(e_{(\lambda,f)}) \theta_i^{(n-1)}(\psi_\v \epsilon_\v) = \theta^{(n-1)}_i(\epsilon_\u^* \psi_\u^*) \psi_{\x\x} \theta_i^{(n-1)}(\psi_\v \epsilon_\v)\\
& = \theta^{(n-1)}_i(\epsilon_\u^*) \theta_i^{(n-1)} (\psi_\u^*) \psi_\x^* e_{(\mu,f)} \psi_\x \theta_i^{(n-1)}(\psi_\v) \theta^{(n-1)}_i(\epsilon_\v) = \epsilon_\s^* \psi_\s^* e_{(\mu,f)} \psi_\t \epsilon_\t = \psi_{\s\t},
\end{align*}
which completes the proof. \endproof

\begin{Lemma} \label{I:h2}
Suppose $(\lambda,f) \in \widehat B_{n-1}$. If there exists $\sigma \vdash n - 2f$ such that $(\sigma,f) \in \mathscr S_n$, for any $\u,\v \in \Tud_{n-1}(\lambda)$ and $i\in P$, we have the following results:
\begin{enumerate}
\item If $\res(\alpha) = i$ for some $\alpha \in \mathscr A(\lambda)$, we have $\theta^{(n-1)}_i(\psi_{\u\v}) y_n \in R_n^{>f}(\delta)$.

\item If $\res(\alpha) \neq i$ for all $\alpha \in \mathscr A(\lambda)$, we have $\theta^{(n-1)}_i(\psi_{\u\v}) \in R_n^{>f}(\delta)$.
\end{enumerate}
\end{Lemma}

\proof Because $\theta^{(n-1)}_i(\psi_{\u\v}) = \theta^{(n-1)}_i(\epsilon_\u^* \psi_\u^*) \theta^{(n-1)}_i(\psi_{\t^{(\lambda,f)} \t^{(\lambda,f)}}) \theta^{(n-1)}_i(\psi_\v \epsilon_\v)$ and $\theta^{(n-1)}_i(\psi_\v \epsilon_\v)$ commutes with $y_n$, by~\autoref{two:sided:ideal} it suffices to prove the Lemma when $\u = \v = \t^{(\lambda,f)}$.


We prove the Lemma by applying induction twice. First we apply induction on $n$. The base step is $n = 1$, which is trivial by (\ref{rela:1}). Assume that there exists $n'$ such that the Lemma holds when $n < n'$ and we set $n = n'$. Then we apply the induction on $f$. The base step is $f = \lfloor \frac{n}{2} \rfloor$, which can be directly verified by~\eqref{rela:4},~\eqref{rela:5:2},~\eqref{rela:5:4} in Case (1) and~\eqref{rela:4} in Case (2). We omit the detailed proof here.

Assume that there exists $f'$ such that the Lemma holds when $f < f'$ and we set $f = f'$. Write $\bi_{(\lambda,f)} = (i_1,i_2,\ldots,i_{n-1})$, $\t^{(\lambda,f)} = (\alpha_1, \ldots, \alpha_{n-1})$ and $\lambda = (\lambda_1, \lambda_2, \ldots, \lambda_m)$.

Here we give two results implied by the induction. These two results will be used in the rest of the proof.

As $(\sigma,f) \in \mathscr S_n$, by~\autoref{two:sided:ideal}, $R_n^{>f}(\delta)$ is a two-sided $\G{n}$-ideal. Hence by induction on $f$, we have $\theta^{(n-1)}_i(R_{n-1}^{>f}(\delta)) \subset R_n^{>f}(\delta)$ for any $i \in P$. Moreover, by induction on $n$, we have $\psi_{\t^{(\lambda,f)} \t^{(\lambda,f)}} y_{n-1} \in R_{n-1}^{>f}(\delta)$, which yields
\begin{equation} \label{I:h2:eq1}
\theta^{(n-1)}_i(\psi_{\t^{(\lambda,f)} \t^{(\lambda,f)}} y_{n-1}) \in R_n^{>f}(\delta).
\end{equation}

In $\G{2f+2}$, by~\eqref{rela:1} we have $\epsilon_1 \epsilon_3 \ldots \epsilon_{2f+1} = \psi_{\t^{(\emptyset,f+1)} \t^{(\emptyset,f+1)}} \in R_{2f+2}^{>f}(\delta)$. Therefore, by induction on $n$ and $f$, we have
$$
\epsilon_1 \epsilon_3 \ldots \epsilon_{2f-1} \epsilon_{2f+1} = \sum_{(i_{2f+3}, \ldots, i_n) \in P^{n-2f-2}} \theta_{i_n}^{(n-1)} \circ \theta_{i_{n-1}}^{(n-2)} \circ \ldots \circ \theta_{i_{2f+3}}^{(2f+2)}(\psi_{\t^{(\emptyset, f+1)} \t^{(\emptyset, f+1)}}) \in R_n^{>f}(\delta),
$$
which can be generalized by~\eqref{rela:7:2} and~\autoref{two:sided:ideal}:
\begin{equation} \label{I:h2:eq2}
\epsilon_1 \ldots \epsilon_{2f-1} \epsilon_\l = \epsilon_\l \epsilon_{\l-1} \ldots \epsilon_{2f+2} {\cdot} \left(\epsilon_1 \epsilon_3 \ldots \epsilon_{2f-1} \epsilon_{2f+1}\right) {\cdot} \epsilon_{2f+2} \ldots \epsilon_{\l-1} \epsilon_\l \in R_n^{>f}(\delta)
\end{equation}
for any $2f + 1 \leq \l \leq n$.

Now we complete the induction process of $f$ by considering different values of $i$. For convenience we write $\bi = (i_1,\ldots,i_n) = \bi_{(\lambda,f)} \vee i \in P^n$.\\

(1). Suppose $\res(\alpha) = i$ for some $\alpha \in \mathscr A(\lambda)$. We consider the following 3 cases. Note that $i = i_{n-1}$ is excluded because it is impossible to find $\alpha \in \mathscr A(\lambda)$ such that $\res(\alpha) = i = i_{n-1}$.

\textbf{Case 1.1:} $i = -i_{n-1}$.


When $i = -\frac{1}{2}$, because $\res(\alpha) = -\frac{1}{2}$ for some $\alpha \in \mathscr A(\lambda)$ and $\res(\alpha_{n-1}) = \frac{1}{2}$ where $\alpha_{n-1} \in \mathscr R(\lambda)$, we have $h_n(\bi) = -2$, which implies $h_{n-1}(\bi) = -1$ by (\ref{remark:h:eq2}). Hence we have $\bi \in I_{n-1,0}^n$. Therefore, by (\ref{rela:5:2}) we have
\begin{align*}
\theta^{(n-1)}_i(\psi_{\t^{(\lambda,f)} \t^{(\lambda,f)}}) y_n & = \theta^{(n-1)}_i(\psi_{\t^{(\lambda,f)} \t^{(\lambda,f)}}) y_n e(\bi) = \theta^{(n-1)}_i(\psi_{\t^{(\lambda,f)} \t^{(\lambda,f)}}) y_{n-1} e(\bi) - 2 \theta^{(n-1)}_i(\psi_{\t^{(\lambda,f)} \t^{(\lambda,f)}}) y_{n-1} e(\bi) \epsilon_{n-1} e(\bi)\\
& = \theta^{(n-1)}_i(\psi_{\t^{(\lambda,f)} \t^{(\lambda,f)}} y_{n-1}) - 2 \theta^{(n-1)}_i(\psi_{\t^{(\lambda,f)} \t^{(\lambda,f)}}y_{n-1}) \epsilon_{n-1} e(\bi),
\end{align*}
which yields $\theta^{(n-1)}_i(\psi_{\t^{(\lambda,f)} \t^{(\lambda,f)}}) y_n \in R_n^{>f}(\delta)$ by (\ref{I:h2:eq1}) and~\autoref{two:sided:ideal}.

When $i = \frac{1}{2}$, following the same argument as $i = -\frac{1}{2}$, we have $h_{n-1}(\bi) = -1$ and $\bi \in I_{n-1,+}^n$. Hence by (\ref{rela:5:1}) we have
\begin{align*}
\theta^{(n-1)}_i(\psi_{\t^{(\lambda,f)} \t^{(\lambda,f)}}) y_n & = \theta^{(n-1)}_i(\psi_{\t^{(\lambda,f)} \t^{(\lambda,f)}}) y_n e(\bi)\\
& = (-1)^{a_{n-1}(\bi)+1} e(\bi) \epsilon_1 \ldots \epsilon_{2f-1} \epsilon_{n-1} e(\bi) + \theta^{(n-1)}_i(\psi_{\t^{(\lambda,f)} \t^{(\lambda,f)}} y_{n-1}),
\end{align*}
which yields $\theta^{(n-1)}_i(\psi_{\t^{(\lambda,f)} \t^{(\lambda,f)}}) y_n \in R_n^{>f}(\delta)$ by (\ref{I:h2:eq1}) and (\ref{I:h2:eq2}).

When $i \neq \pm \frac{1}{2}$, following the same argument as when $i = -\frac{1}{2}$, we have $h_{n-1}(\bi) = 0$ and $\bi \in I_{n-1,+}^n$. By the similar argument as when $i = \frac{1}{2}$, we have $\theta^{(n-1)}_i(\psi_{\t^{(\lambda,f)} \t^{(\lambda,f)}}) y_n \in R_n^{>f}(\delta)$.

\textbf{Case 1.2:} $i = i_{n-1} \pm 1$.

First we consider the case when $i = i_{n-1} - 1$. Notice that when $i = -\frac{1}{2}$, we have $i = -i_{n-1} = -\frac{1}{2}$, which has already been proved in Case 1.1. Hence we set $i \neq -\frac{1}{2}$.

Because $\res(\alpha) = i$ for some $\alpha \in \mathscr A(\lambda)$ and $\res(\alpha_{n-1}) = i+1$ where $\alpha_{n-1} \in [\lambda]$, we have $\lambda_m = 1$ by the construction of $\lambda$. Let $\mu = \lambda|_{n-2}$ and $\bj = \bi_{(\mu,f)} = \bi_{(\lambda,f)}|_{n-2}$. By (\ref{rela:4}) we have
\begin{align*}
\theta^{(n-1)}_i(\psi_{\t^{(\lambda,f)} \t^{(\lambda,f)}}) y_n & = \theta^{(n-1)}_i(\psi_{\t^{(\lambda,f)} \t^{(\lambda,f)}} y_{n-1}) - \psi_{n-1} e(\bj, i, i+1) \epsilon_1 \ldots \epsilon_{2f-1} e(\bj, i, i+1) \psi_{n-1}.
\end{align*}

Because $\lambda_m = 1$, we have $\res(\alpha) \neq i$ for any $\alpha \in \mathscr A(\mu)$. Hence by induction on $n$, we have
$$
e(\bj, i, i+1) \epsilon_1 \ldots \epsilon_{2f-1} e(\bj, i, i+1) = \theta^{(n-1)}_{i+1}( \theta^{(n-2)}_i ( \psi_{\t^{(\mu,f)} \t^{(\mu,f)}}) )\in \theta^{(n-1)}_{i+1}(R_{n-1}^{>f}(\delta)) \subset R_n^{>f}(\delta),
$$
and the Lemma holds by~\eqref{I:h2:eq1}. Following the similar argument, the Lemma holds when $i = i_{n-1} + 1$.

\textbf{Case 1.3:} $|i - i_{n-1}| > 1$.

This case can be verified using (\ref{rela:3:2}) and (\ref{rela:4}) directly and we omit the detailed proof here.\\

(2). Suppose $\res(\alpha) \neq i$ for all $\alpha \in \mathscr A(\lambda)$. We consider the following 5 cases:

\textbf{Case 2.1:} $i = -i_{n-1}$.


When $i = -\frac{1}{2}$, as $\res(\alpha) \neq i$ for any $\alpha \in \mathscr A(\lambda)$, we have $h_n(\bi) = -1$, which implies $h_{n-1}(\bi) = -2$ by~\eqref{remark:h:eq2} and hence $\bi \in I_{n-1,-}^n$. Therefore, by~\eqref{rela:5:4} we have
$$
\theta^{(n-1)}_i(\psi_{\t^{(\lambda,f)} \t^{(\lambda,f)}}) = (-1)^{a_{n-1}(\bi)} \psi_{\t\t}y_{n-1} + (-1)^{a_{n-1}(\bi)} y_{n-1} \psi_{\t\t} \in R_n^{>f}(\delta).
$$

When $i = \frac{1}{2}$, as $\res(\alpha) \neq i$ for any $\alpha \in \mathscr A(\lambda)$, we have $\lambda_{m-1} = \lambda_m$. See the next diagram for $\lambda$:
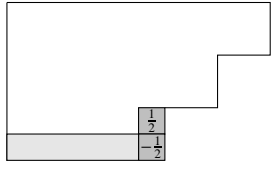
$$
\begin{tikzpicture} [scale = 0.7]
\draw (0,0) -- (5,0) -- (5,-1) -- (4,-1) -- (4,-2) -- (3,-2);
\draw (2.5,-2) [fill = gray!50] rectangle (3,-2.5);
\node [scale = 0.7] at (2.75,-2.25) {$\frac{1}{2}$};
\draw (2.5,-2.5) [fill = gray!50] rectangle (3,-3);
\node [scale = 0.7] at (2.75,-2.75) {$-\frac{1}{2}$};
\draw (2.5,-3) -- (0,-3) -- (0,0);
\draw [fill = gray!20] (0,-2.5) rectangle (2.5,-3);
\end{tikzpicture}
$$

In the above diagram, the entries in the shadowed nodes are their residues and the residues of the light shadowed nodes are less than $-\frac{1}{2}$. Formally, we set $a = n - \lambda$. So the node $\alpha_a$ is the shadowed node with residue $\frac{1}{2}$ and the node $\alpha_{n-1}$ is the shadowed node with residue $-\frac{1}{2}$. Moreover, for any $a < s < n-1$, we have $i_a + i_s \neq 0$ and $|i_a - i_s| > 1$. Therefore, by~\eqref{rela:4} and~\eqref{rela:2:1} we have
$$
\theta^{(n-1)}_i(\psi_{\t^{(\lambda,f)} \t^{(\lambda,f)}}) = \psi_a \ldots \psi_{n-4} \psi_{n-3} e(\bj) \epsilon_1 \ldots \epsilon_{2f-1} e(\bj) \psi_{n-3} \psi_{n-4} \ldots \psi_a,
$$
where $\bj = (i_1,\ldots,i_{a-1}, i_{a+1}, \ldots, i_{n-2}, i_a, i_{n-1}, i)$. We note that $i_a = i = \frac{1}{2}$ and $i_{n-1} = -\frac{1}{2}$. Hence apply~\eqref{rela:5:3} to the above equation and we have
\begin{eqnarray*}
\theta^{(n-1)}_i(\psi_{\t^{(\lambda,f)} \t^{(\lambda,f)}})
& = & (-1)^{a_{n-1}(\bi)} \psi_a \ldots \psi_{n-4} \psi_{n-3} e(\bj) \epsilon_1 \ldots \epsilon_{2f-1} \epsilon_{n-1} e(\bj) \psi_{n-3} \psi_{n-4} \ldots \psi_a\\
&& \ \ - 2(-1)^{a_{n-2}(\bi)} \psi_a \ldots \psi_{n-4} \psi_{n-3} e(\bj) \epsilon_1 \ldots \epsilon_{2f-1} \epsilon_{n-2} e(\bj) \psi_{n-3} \psi_{n-4} \ldots \psi_a\\
&& \ \ + \psi_a \ldots \psi_{n-4} \psi_{n-3} e(\bj) \epsilon_1 \ldots \epsilon_{2f-1} \epsilon_{n-1} \epsilon_{n-2} e(\bj) \psi_{n-3} \psi_{n-4} \ldots \psi_a\\
&& \ \ + \psi_a \ldots \psi_{n-4} \psi_{n-3} e(\bj) \epsilon_1 \ldots \epsilon_{2f-1} \epsilon_{n-2} \epsilon_{n-1} e(\bj) \psi_{n-3} \psi_{n-4} \ldots \psi_a,
\end{eqnarray*}
which yields $\theta^{(n-1)}_i(\psi_{\t^{(\lambda,f)} \t^{(\lambda,f)}}) \in R_n^{>f}(\delta)$ by~\eqref{I:h2:eq2} and~\autoref{two:sided:ideal}.

When $i \neq \pm \frac{1}{2}$, as $\res(\alpha) \neq i$ for any $\alpha \in \mathscr A(\lambda)$, we have $h_{n-1}(\bi) = -1$ and $\bi \in I_{n-1,0}^n$. By (\ref{rela:5:1}) and (\ref{I:h2:eq2}), we have
$$
\theta^{(n-1)}_i(\psi_{\t^{(\lambda,f)} \t^{(\lambda,f)}}) = (-1)^{a_{n-1}(\bi)} e(\bi) \epsilon_1 \ldots \epsilon_{2f-1} \epsilon_{n-1} e(\bi) \in R_n^{>f}(\delta).
$$

\textbf{Case 2.2:} $i = i_{n-1}$.

When $i = i_{n-1} = 0$, we have $\bi \in I_{n-1,0}^n$. Hence following the same argument as in Case 2.1 when $i \neq \pm \frac{1}{2}$, we have $\theta^{(n-1)}_i(\psi_{\t^{(\lambda,f)} \t^{(\lambda,f)}}) \in R_n^{>f}(\delta)$.

When $i = i_{n-1} \neq 0$, it is easy to verify the following equality holds by (\ref{rela:3:1}), (\ref{rela:3:2}) and (\ref{rela:4}):
$$
e(\bi) = - \psi_{n-1} e(\bi) y_{n-1}^2 \psi_{n-1} - e(\bi) y_{n-1} \psi_{n-1} - \psi_{n-1} e(\bi) y_{n-1}.
$$

Hence by (\ref{I:h2:eq1}), we have
$$
\theta^{(n-1)}_i(\psi_{\t^{(\lambda,f)} \t^{(\lambda,f)}}) = - \psi_{n-1} \theta^{(n-1)}_i(\psi_{\t^{(\lambda,f)} \t^{(\lambda,f)}}y_{n-1}^2) \psi_{n-1} - \theta^{(n-1)}_i(\psi_{\t^{(\lambda,f)} \t^{(\lambda,f)}} y_{n-1}) \psi_{n-1} - \psi_{n-1} \theta^{(n-1)}_i(\psi_{\t^{(\lambda,f)} \t^{(\lambda,f)}} y_{n-1}) \in R_n^{>f}(\delta).
$$

\textbf{Case 2.3:} $i = i_{n-1} - 1$.

The case $i = -\frac{1}{2}$ has been proved in Case 2.1. Hence we assume $i \neq -\frac{1}{2}$. As $\res(\alpha) \neq i$ for any $\alpha \in \mathscr A(\lambda)$, we have $\lambda_m > 1$. Hence we have $i_{n-1} = i+1$ and $i_{n-2} = i_{n-1} - 1 = i$. Define $\bj = (i_1, i_2, \ldots, i_{n-3})$. By (\ref{rela:8:7}) we have
\begin{align}
\theta^{(n-1)}_i(\psi_{\t^{(\lambda,f)} \t^{(\lambda,f)}}) & = \psi_{n-2} \psi_{n-1} e(\bj \vee i, i, i+1) \epsilon_1 \ldots \epsilon_{2f-1} e(\bj \vee i, i, i+1) \psi_{n-2} \notag \\
& \ \ - \psi_{n-1} \psi_{n-2} e(\bj \vee i+1, i, i) \epsilon_1 \ldots \epsilon_{2f-1} e(\bj \vee i+1, i, i) \psi_{n-1}. \label{I:h2:eq5}
\end{align}

Let $\mu = \lambda|_{n-2}$ and $\gamma = \lambda|_{n-3}$. One can see that $\res(\alpha) \neq i$ for any $\alpha \in \mathscr A(\mu)$ and $\res(\alpha) \neq i+1$ for any $\alpha \in \mathscr A(\gamma)$ by the construction of $\lambda$. Hence by induction on $n$, we have
\begin{align*}
e(\bj \vee i, i, i+1) \epsilon_1 \ldots \epsilon_{2f-1} e(\bj \vee i, i, i+1) & = \theta^{(n-1)}_{i+1}(\theta^{(n-2)}_i(\psi_{\t^{(\mu,f)} \t^{(\mu,f)}})) \in \theta^{(n-1)}_{i+1}(R_{n-1}^{>f}(\delta)) \subset R_n^{>f}(\delta), \\
e(\bj \vee i+1, i, i) \epsilon_1 \ldots \epsilon_{2f-1} e(\bj \vee i+1, i, i) & = \theta^{(n-1)}_i( \theta^{(n-2)}_i (\theta^{(n-3)}_{i+1}(\psi_{\t^{(\gamma,f)} \t^{(\gamma,f)}}))) \in \theta^{(n-1)}_i(\theta^{(n-2)}_i(R_{n-2}^{>f}(\delta))) \subset R_n^{>f}(\delta).
\end{align*}

Substitute the above equalities to~\eqref{I:h2:eq5}. By~\autoref{two:sided:ideal}, we have $\theta^{(n-1)}_i(\psi_{\t^{(\lambda,f)} \t^{(\lambda,f)}}) \in R_n^{>f}(\delta)$.

\textbf{Case 2.4:} $i = i_{n-1} + 1$.

The case $i = \frac{1}{2}$ has been proved in Case 2.1. Hence we assume $i \neq \frac{1}{2}$. As $\res(\alpha) \neq i$ for any $\alpha \in \mathscr A(\lambda)$, we have $\lambda_{m-1} = \lambda_m$. See the next diagram for $\lambda$:
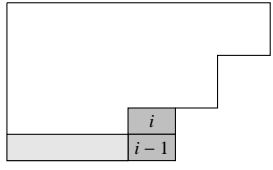
$$
\begin{tikzpicture} [scale = 0.7]
\draw (0,0) -- (5,0) -- (5,-1) -- (4,-1) -- (4,-2) -- (3.2,-2);
\draw (2.3,-2) [fill = gray!50] rectangle (3.2,-2.5);
\node [scale = 0.7] at (2.75,-2.25) {$i$};
\draw (2.3,-2.5) [fill = gray!50] rectangle (3.2,-3);
\node [scale = 0.7] at (2.75,-2.75) {$i-1$};
\draw (2.3,-3) -- (0,-3) -- (0,0);
\draw [fill = gray!20] (0,-2.5) rectangle (2.3,-3);
\end{tikzpicture}
$$

In the above diagram, the entries in the shadowed nodes are their residues and the residues of the light shadowed nodes are less than $i-1$. Formally, we set $a = n - \lambda$. So the node $\alpha_a$ is the shadowed node with residue $i$ and the node $\alpha_{n-1}$ is the shadowed node with residue $i-1$. Moreover, for any $a < s < n-1$, we have $|i_a - i_s| > 1$.

Notice that as $\lambda_{m-1} = \lambda_m$, we have $\t^{(\lambda,f)} s_a s_{a+1} \ldots s_{n-3} \in \Tud_{n-1}(\lambda)$, which is of the form
$$
\left( \ldots\ldots \longrightarrow
\begin{tikzpicture} [baseline=-11mm, scale = 0.7]
\draw (0,0) -- (5,0) -- (5,-1) -- (4,-1) -- (4,-2) -- (2.3,-2) -- (2.3,-3) -- (0,-3) -- (0,0);
\end{tikzpicture}
\longrightarrow
\begin{tikzpicture} [baseline=-11mm, scale = 0.7]
\draw (0,0) -- (5,0) -- (5,-1) -- (4,-1) -- (4,-2) -- (2.3,-2) -- (2.3,-3) -- (0,-3) -- (0,0);
\draw (3.2,-2) -- (3.2,-2.5) -- (2.3,-2.5);
\node [scale = 0.7] at (2.75,-2.25) {$i$};
\end{tikzpicture}
\longrightarrow
\begin{tikzpicture} [baseline=-11mm, scale = 0.7]
\draw (0,0) -- (5,0) -- (5,-1) -- (4,-1) -- (4,-2) -- (2.3,-2) -- (2.3,-3) -- (0,-3) -- (0,0);
\draw (3.2,-2) -- (3.2,-2.5) -- (2.3,-2.5);
\node [scale = 0.7] at (2.75,-2.25) {$i$};
\draw (3.2,-2.5) -- (3.2,-3) -- (2.3,-3);
\node [scale = 0.7] at (2.75,-2.75) {$i-1$};
\end{tikzpicture} \right).
$$

By~\autoref{A:1} and~\eqref{rela:2:1}, we have
$$
\theta^{(n-1)}_i(\psi_{\t^{(\lambda,f)} \t^{(\lambda,f)}}) = \psi_a \psi_{a+1} \ldots \psi_{n-3} e(\bj) \epsilon_1 \ldots \epsilon_{2f-1} e(\bj) \psi_{n-3} \ldots \psi_{a+1} \psi_a,
$$
where $\bj = (i_1,\ldots,i_{a-1}, i_{a+1}, \ldots, i_{n-2}, i_a, i_{n-1}, i)$. We note that $(i_a, i_{n-1}, i) = (i, i-1, i)$. Hence, following the similar argument as in Case 2.3, we have $e(\bj) \epsilon_1 \ldots \epsilon_{2f-1} e(\bj) \in R_n^{>f}(\delta)$. Therefore we have $\theta^{(n-1)}_i(\psi_{\t^{(\lambda,f)} \t^{(\lambda,f)}}) \in R_n^{>f}(\delta)$ by~\autoref{two:sided:ideal}.

\textbf{Case 2.5:} $|i - i_n| > 1$.

This case can be verified by (\ref{rela:4}) directly and we omit the detailed proof here. \endproof

The next Lemma combines the results of~\autoref{I:h1} and~\autoref{I:h2}.

\begin{Lemma} \label{I:h3}
Suppose $0 \leq f \leq \lfloor \frac{n}{2} \rfloor$. If there exists $\sigma \vdash n - 2f$ such that $(\sigma,f) \in \mathscr S_n$, for any $i \in P$ we have $\theta_i^{(n-1)}(R_{n-1}^{\geq f}(\delta)) \subset R_n^{\geq f}(\delta)$ and $\theta_i^{(n-1)}(R_{n-1}^{>f}(\delta)) \subset R_n^{>f}(\delta)$. Moreover, if $\mathscr S_n = \widehat B_n$, then for any $i \in P$ we have $\theta_i^{(n-1)}(R_{n-1}(\delta)) \subseteq R_n(\delta)$.
\end{Lemma}

\proof Suppose $(\mu,m) \in \widehat B_{n-1}$ and $\s,\t \in \Tud_{n-1}(\mu)$. When $m = f$, if $\res(\alpha) = i$ for some $\alpha \in \mathscr A(\mu)$, by~\autoref{I:h1} we have $\theta_i^{(n-1)}(\psi_{\s\t}) \in R_n^{\geq f}(\delta)$; and if $\res(\alpha) \neq i$ for all $\alpha \in \mathscr A(\mu)$, by~\autoref{I:h2} we have $\theta^{(n-1)}_i(\psi_{\s\t}) \in R_n^{>f}(\delta)$. Hence we have
\begin{equation} \label{I:h3:eq1}
\theta^{(n-1)}_i(\psi_{\s\t}) \in R_n^{\geq f}(\delta).
\end{equation}

When $m > f$, by~\autoref{I:h1} and~\autoref{I:h2}, we have
\begin{equation} \label{I:h3:eq2}
\theta^{(n-1)}_i(\psi_{\s\t}) \in R_n^{> f}(\delta).
\end{equation}

By (\ref{I:h3:eq1}) and (\ref{I:h3:eq2}), we have $\theta^{(n-1)}_i(R_{n-1}^{\geq f}(\delta)) \subset R_n^{\geq f}(\delta)$ and $\theta^{(n-1)}_i(R_{n-1}^{>f}(\delta)) \subset R_n^{>f}(\delta)$.

Suppose $\mathscr S_n = \widehat B_n$. Choose $\lambda$ such that $(\lambda,0) \in \widehat B_n$. Hence that $(\lambda,0) \in \mathscr S_n$ and we have $\theta_i^{(n-1)}(R_{n-1}^{\geq 0}(\delta)) \subseteq R_n^{\geq 0}(\delta)$. Notice that $R_n(\delta) = R_n^{\geq 0}(\delta)$. Therefore we have $\theta_i^{(n-1)}(R_{n-1}(\delta)) \subseteq R_n(\delta)$. \endproof

Recall by identifying $e(\bi) = \sum_{i \in P} e(\bi \vee i)$ for $\bi \in P^{n-1}$, we consider $\G{n-1}$ as a subalgebra of $\G{n}$ and obtain a sequence
$$
\G{1} \subset \G{2} \subset \G{3} \subset \ldots.
$$

The key point of~\autoref{I:h3} is that by assuming $\mathscr S_n = \widehat B_n$, we can construct such sequence for $R_n(\delta)$ as well. In~\autoref{I:h3}, as $i$ is chosen arbitrary in $P$, we can consider $R_{n-1}(\delta)$ as a subspace of $R_n(\delta)$ by identifying $e(\bi) = \sum_{i \in P} e(\bi \vee i)$ for $\bi \in P^{n-1}$. Hence we obtain a sequence
$$
R_1(\delta) \subset R_2(\delta) \subset R_3(\delta) \subset \ldots \subset R_n(\delta).
$$

The next Proposition is the most important application of~\autoref{I:h3} in this paper.

\begin{Proposition} \label{I:end}
Suppose $\bigcup_{i=1}^n \widehat B_i \subseteq \widehat{\mathscr B}$. Then $e(\bi) \in R_n(\delta)$ for any $\bi \in P^n$.
\end{Proposition}

\proof Because $\mathscr S_n = \widehat B_n$, by the definition of $\mathscr S_n$ we have $\mathscr S_k = \widehat B_k$ for any $1 \leq k \leq n$. Therefore, if we write $\bi = (i_1, i_2, \ldots, i_n)$, as $i_1 = \frac{\delta-1}{2}$ and $e(i_1) = 1 \in R_1(\delta) \subset \G{1}$, we have
$$
e(\bi) = \theta_{i_n}^{(n-1)} \circ \theta_{i_{n-1}}^{(n-2)} \circ \ldots \circ \theta_{i_2}^{(1)} (1) \in R_n(\delta)
$$
by~\autoref{I:h3}. \endproof

Suppose $\bigcup_{i=1}^n \widehat B_i \subseteq \widehat{\mathscr B}$.~\autoref{I:end} shows that $1 \in R_n(\delta)$. By the definition of $\widehat{\mathscr B}$, $R_n(\delta)$ is a right $\G{n}$-module. Hence, $\bigcup_{i=1}^n \widehat B_i \subseteq \widehat{\mathscr B}$ implies $\set{\psi_{\s\t} | \s,\t \in \Tud_n(\lambda), (\lambda,f) \in \widehat B_n}$ is a $R$-spanning set of $\G{n}$, and it has cellular-like property by the definition of $\widehat{\mathscr B}$.

Finally we introduce some applications of induction property of $R_n(\delta)$.

\begin{Lemma} \label{I:coro:1}
Suppose $0 \leq f \leq \lfloor \frac{n}{2} \rfloor$. If there exists $\sigma \vdash n - 2f$ such that $(\sigma,f) \in \mathscr S_n$, we have $\epsilon_1 \epsilon_3 \ldots \epsilon_{2f-1} \epsilon_k \in R_n^{>f}(\delta)$ for $2f + 1 \leq k \leq n-1$.
\end{Lemma}

\proof It suffices to prove that when $k = 2f+1$ the Lemma holds, because when $k > 2f+1$, by (\ref{rela:7:2}) we have
$$
\epsilon_1 \epsilon_3 \ldots \epsilon_{2f-1} \epsilon_k = \epsilon_k \epsilon_{k-1} \ldots \epsilon_{2f+2} {\cdot} \epsilon_1 \epsilon_3 \ldots \epsilon_{2f-1} \epsilon_{2f+1} {\cdot} \epsilon_{2f+2} \ldots \epsilon_{k-1} \epsilon_k,
$$
and $R_n^{>f}(\delta)$ is a two-sided $\G{n}$-ideal by~\autoref{two:sided:ideal}.

Consider $k = 2f+1$. By (\ref{rela:1}) we have
$$
\epsilon_1 \epsilon_3 \ldots \epsilon_{2f-1} \epsilon_{2f+1} = \sum_{(i_{2f+3}, \ldots, i_n) \in P^{n-2f-2}} \theta_{i_n}^{(n-1)} \circ \theta_{i_{n-1}}^{(n-2)} \circ \ldots \circ \theta_{i_{2f+3}}^{(2f+2)}(\psi_{\t^{(\emptyset, f+1)} \t^{(\emptyset, f+1)}}).
$$

Because $\psi_{\t^{(\emptyset, f+1)} \t^{(\emptyset, f+1)}} \in R_{2f+2}^{>f}(\delta)$, the Lemma follows by~\autoref{I:h3}. \endproof

\begin{Lemma} \label{I:coro:3}
Suppose $0 \leq f \leq \lfloor \frac{n}{2} \rfloor$. If there exists $\sigma \vdash n - 2f$ such that $(\sigma,f) \in \mathscr S_n$, then we have $e(\bi) = 0$ if $h_k(\bi) > 0$ for some $1 \leq k \leq n-1$.
\end{Lemma}

\proof Let $\bj = \bi|_{n-1} \in P^{n-1}$. As $(\sigma,f) \in \mathscr S_n$, we have $\bigcup_{i=1}^{n-1} \widehat B_i \subseteq \widehat{\mathscr B}$, which implies $e(\bj) \in R_{n-1}(\delta)$ by~\autoref{I:end}. If $h_k(\bi) > 0$ for some $1 \leq k \leq n-1$, by~\autoref{deg:h2} we have $\bj \not\in I^{n-1}$. Hence by~\autoref{idem:psi:1} we have $e(\bj) = 0$. Therefore $e(\bi) = \theta_{i_n}(e(\bj)) = 0$. \endproof

\begin{Lemma} \label{remove:tail:1}
Suppose $(\lambda,f) \in \mathscr S_n$ and $\t \in \Tud_n(\lambda)$. Let $1 \leq k \leq n$. If $\t(r) > 0$ for all $k \leq r \leq n$, then for any $a \in \G{k-1}$, we have
$$
\psi_{\t^{(\lambda,f)} \t}{\cdot} a \equiv \sum_{\v \in \Tud_n(\lambda)} c_\v \psi_{\t^{(\lambda,f)} \v} \pmod{R_n^{>f}(\delta)},
$$
where $c_\v \neq 0$ only if $\v(r) = \t(r)$ for any $k \leq r \leq n$.
\end{Lemma}

\proof We only prove the case when $k = n-1$. For smaller $k$ the proof is essentially the same.

Suppose the head of $\t$ is $h$ and $\lambda = (\lambda_1, \ldots, \lambda_m)$. Define $\dot \t = \t|_{n-1}$ and $\mu = \t_{n-1}$. As $\t(n) > 0$, we have $\t(n) = \alpha$ for some $\alpha \in \mathscr A(\mu)$ and $\lambda = \mu \cup \{\alpha\}$. Hence $\alpha = (\l, \lambda_\l)$ for some $1 \leq \l \leq m$. Let $a = \lambda_1 + \lambda_2 + \ldots + \lambda_\l$ and $\res(\alpha) = i_n$. By~\autoref{I:h1:1}, we have $\psi_\t = \psi_a \psi_{a+1} \ldots \psi_{n-1} \theta_{i_n}^{(n-1)}(\psi_{\dot\t})$ and $\epsilon_\t = \theta_{i_n}^{(n-1)}(\epsilon_{\dot\t})$. Hence by~\autoref{I:h1}, we have
$$
\psi_{\t^{(\lambda,f)} \t} a = \psi_a \psi_{a+1} \ldots \psi_{n-1} \theta^{(n-1)}_{i_n}(\psi_{\t^{(\mu,f)} \dot\t} a).
$$

As $(\lambda,f) \in \mathscr S_n$, by~\autoref{I:h3} we have
\begin{equation} \label{remove:tail:1:eq1}
\psi_{\t^{(\lambda,f)} \t} a = \psi_a \psi_{a+1} \ldots \psi_{n-1} \theta^{(n-1)}_{i_n}(\psi_{\t^{(\mu,f)} \dot\t} a) \equiv \sum_{\dot\v \in \Tud_{n-1}(\mu)} c_{\dot\v} \psi_a \psi_{a+1} \ldots \psi_{n-1} \theta_{i_n}^{(n-1)}(\psi_{\t^{(\mu,f)} \dot\v}) \pmod{R_n^{>f}(\delta)}.
\end{equation}

For $\dot\v \in \Tud_{n-1}(\mu)$, define $\v \in \Tud_n(\lambda)$ with $\v|_{n-1} = \dot\v$ and $\v(n) = \alpha$. By~\autoref{I:h1:1} we have $\psi_{\v} = \psi_a \psi_{a+1} \ldots \psi_{n-1} \theta_{i_n}^{(n-1)}(\psi_{\dot\v})$ and $\epsilon_\v = \theta_{i_n}^{(n-1)}(\epsilon_{\dot\v})$. Hence by~\autoref{I:h1}, we have
\begin{equation} \label{remove:tail:1:eq2}
\psi_a \psi_{a+1} \ldots \psi_{n-1} \theta^{(n-1)}_{i_n} (\psi_{\t^{(\mu,f)} \dot\v}) = e_{(\lambda,f)} \psi_\v \epsilon_\v = \psi_{\t^{(\lambda,f)} \v},
\end{equation}
where $\v(n) = \alpha = \t(n)$. The Lemma holds by substituting (\ref{remove:tail:1:eq2}) into (\ref{remove:tail:1:eq1}). \endproof

\subsection{The restriction property} \label{sec:res}


In this subsection we introduce the restriction property of $\G{n}$. Suppose $0 \leq f \leq \lfloor \frac{n}{2} \rfloor$. Define $P_{f,n} = \set{\bi \in P^n | i_1 = i_3 = \ldots = i_{2f+1} = -i_2 = -i_4 = \ldots = -i_{2f} = \frac{\delta-1}{2}}$ and $\G{2f,n}$ as the subalgebra of $\G{n}$ generated by
\begin{eqnarray*}
G_{f,n}(\delta) & = & \set{e(\bi)|\bi \in P_{f,n}} \cup \set{y_k| 2f + 1 \leq k \leq n} \cup \set{\psi_k| 2f + 1\leq k \leq n-1} \cup \set{\epsilon_k| 2f + 1\leq k \leq n-1}.
\end{eqnarray*}

Denote $\epsilon_{1,0} = 1$ and $\epsilon_{1,f} = \epsilon_1 \epsilon_3 \ldots \epsilon_{2f-1}$ for $f > 0$. We define a map $\phi_{f,n}^{(g)}\map{G_{f,n}(\delta)}{G_{n-2f}(\delta)}$ by
$$
e(\bi) \mapsto e(i_{2f+1}, i_{2f+2}, \ldots, i_n),\qquad y_r \mapsto y_{r-2f}, \qquad \psi_s \mapsto \psi_{s-2f}, \qquad \text{and} \qquad \epsilon_s \mapsto \epsilon_{s-2f},
$$
where $\bi = \in P_{f,n}$, $2f+1 \leq r \leq n$ and $2f+1 \leq s \leq n-1$. Extend $\phi^{(g)}_{f,n}$ to a linear map $\phi_{f,n}\map{\epsilon_{1,f}\G{2f,n}}{\G{n-2f}}$ such that for each $a \in \G{2f,n}$, if we can write $a = g_1 g_2 \ldots g_k$ where $g_i \in G_{f,n}(\delta)$ for $1 \leq i \leq k$, then
$$
\phi_{f,n}(\epsilon_{1,f} {\cdot} a) = \phi_{f,n}^{(g)}(g_1)\phi_{f,n}^{(g)}(g_2) \ldots \phi_{f,n}^{(g)}(g_k) \in \G{n-2f}.
$$

\begin{Lemma} \label{A:2}
Suppose $0 \leq f \leq \lfloor \frac{n}{2} \rfloor$. The linear map $\phi_{f,n}\map{\epsilon_{1,f} \G{2f,n}}{\G{n-2f}}$ is well-defined.
\end{Lemma}

\proof When $f = 0$, $\phi_{0,n}$ is the identity map. Hence we assume $f > 0$. It suffices to check the relations of $\epsilon_{1,f} \G{2f,n}$ consist in $\G{n-2f}$ by applying $\phi_{f,n}$. Let $\bi = (i_1, \ldots, i_n) \in P_{f,n}$ and $\bj = (i_{2f+1}, i_{2f+2}, \ldots, i_n) \in P^{n-2f}$. By direct calculation, we have $h_k(\bi) = h_{k-2f}(\bj)$ and $(-1)^{a_k(\bi)} = (-1)^{a_{k-2f}(\bj)}$ for $2f+1 \leq k \leq n$. Hence we only need to check (\ref{rela:5:5}) when $\bi \in P_{k,+}^n$, as all the other relations depend on $h_k(\bi)$, $(-1)^{a_k(\bi)}$, $i_{k-1}$, $i_k$ and $i_{k+1}$.

First we prove the following equality holds: for any $2f+1 \leq k \leq n-1$, we have
\begin{equation} \label{A:2:eq1}
\epsilon_1 \epsilon_3 \ldots \epsilon_{2f-1} (\sum_{\substack{1 \leq r \leq 2f \\ r \in A_{k,1}^\bi}} y_r - 2 \sum_{\substack{1 \leq r \leq 2f \\ r \in A_{k,2}^\bi}} y_r + \sum_{\substack{1 \leq r \leq 2f \\ r \in A_{k,3}^\bi}} y_r - 2\sum_{\substack{1 \leq r \leq 2f \\ r \in A_{k,4}^\bi}} y_r) = 0.
\end{equation}

Because we have $i_1 = i_3 = \ldots = i_{2f-1} = -i_2 = -i_4 = \ldots = -i_{2f} = \frac{\delta-1}{2}$, for any $1 \leq \l \leq 2f$, $2\l - 1 \in A_{k,1}^\bi$ if and only if $2\l \in A_{k,3}^\bi$, and $2\l \in A_{k,1}^\bi$ if and only if $2\l - 1 \in A_{k,3}^\bi$. Similarly, $2\l - 1 \in A_{k,2}^\bi$ if and only if $2\l \in A_{k,4}^\bi$, and $2\l \in A_{k,2}^\bi$ if and only if $2\l - 1 \in A_{k,4}^\bi$. Hence by (\ref{rela:7:2}), (\ref{A:2:eq1}) holds.

Suppose $\bi \in P_{k,+}^n$ and $2f + 1 \leq k \leq n-1$. Recall $(-1)^{a_k(\bi)} = (-1)^{a_{k-2f}(\bj)}$. By (\ref{rela:5:5}) and (\ref{A:2:eq1}) we have
\begin{align*}
\phi_{f,n}(\epsilon_{1,f} \epsilon_k e(\bi) \epsilon_k) & = (-1)^{a_k(\bi)} (1 + \delta_{i_k, -\frac{1}{2}})  \phi_{f,n}(\epsilon_{1,f} (\sum_{r \in A_{k,1}^\bi} y_r - 2 \sum_{r \in A_{k,2}^\bi} y_r + \sum_{r \in A_{k,3}^\bi} y_r - 2\sum_{r \in A_{k,4}^\bi} y_r) \epsilon_k)\\
& = (-1)^{a_k(\bi)} (1 + \delta_{i_k, -\frac{1}{2}})  \phi_{f,n}(\epsilon_{1,f} (\sum_{\substack{2f + 1 \leq r \leq k-1 \\ r \in A_{k,1}^\bi}} y_r - 2 \sum_{\substack{2f + 1 \leq r \leq k-1 \\ r \in A_{k,2}^\bi}} y_r + \sum_{\substack{2f + 1 \leq r \leq k-1 \\ r \in A_{k,3}^\bi}} y_r - 2\sum_{\substack{2f + 1 \leq r \leq k-1 \\ r \in A_{k,4}^\bi}} y_r) \epsilon_k)\\
& = (-1)^{a_{k-2f}(\bj)} (1 + \delta_{j_{k-2f}, -\frac{1}{2}}) (\sum_{r \in A_{k-2f,1}^\bj} y_r - 2 \sum_{r \in A_{k-2f,2}^\bj} y_r + \sum_{r \in A_{k-2f,3}^\bj} y_r - 2\sum_{r \in A_{k-2f,4}^\bj} y_r) \epsilon_{k-2f}\\
& = \epsilon_{k-2f} e(\bj) \epsilon_{k-2f},
\end{align*}
which proves the Lemma. \endproof

The next result is used to prove~\autoref{remove-front}.

\begin{Lemma} \label{idem:psi:2}
Suppose $\bi = (i_1, \ldots, i_n) \in P_{f,n}$. For any $\t \in \Tud_n(\bi)$, we have $\t_{2f+1} = (1)$.
\end{Lemma}

\proof Suppose $\frac{\delta-1}{2} \neq \pm \frac{1}{2}$, the Lemma is obvious by the construction of $\t$. If $\frac{\delta - 1}{2} = \frac{1}{2}$, we prove by induction.

When $f = 0$, the Lemma follows obviously. Suppose when $f < f'$ the Lemma follows. When $f = f'$, by induction we have $\t_{2f-1} = (1)$. Hence by the construction, $\t(2f) = (2,1)$ or $\t(2f) = -(1,1)$. If $\t(2f) = -(1,1)$, then $\t(2f+1) = (1,1)$ and the Lemma holds. If $\t(2f) = (2,1)$, then $\t(2f+1) = -(2,1)$ and the Lemma holds. Hence when $\frac{\delta-1}{2} = \frac{1}{2}$, the Lemma holds. Following the same argument, the Lemma holds when $\frac{\delta-1}{2} = -\frac{1}{2}$. \endproof

\begin{Lemma} \label{remove-front}
The map $\phi_{f,n}$ is a bijection.
\end{Lemma}

\proof By the definition, $\phi_{f,n}$ is surjective. In order to prove $\phi_{f,n}$ is injective, it suffices to show that $\ker \phi_{f,n} = \{0\}$ by checking the relations of $\G{n-2f}$. Following the same argument as in the proof of~\autoref{A:2}, we only need to check the first two relations of (\ref{rela:1}), and (\ref{rela:5:5}) when $\bi \in P_{k,+}^n$.

Suppose $\bi = (i_1, \ldots, i_n) \in P^n$. For the first relation of (\ref{rela:1}), we prove that $e(\bi) \epsilon_{1,f} e(\bi) = 0$ when $i_{2f+1} \neq \frac{\delta-1}{2}$ and $e(\bi) \epsilon_{1,f} y_{2f+1} e(\bi) = 0$. Note that we have $c = \pm 1$ such that $e(\bi) \epsilon_{1,f} e(\bi) = c e(\bi) \epsilon_1 \ldots \epsilon_{2f-1} e(\bi)$. Hence it suffices to prove that $e(\bi) \epsilon_1 \ldots \epsilon_{2f-1} e(\bi) = 0$ and $e(\bi) \epsilon_1 \ldots \epsilon_{2f-1} y_{2f+1} e(\bi) = 0$ under certain conditions.

First we prove $e(\bi) \epsilon_1 \ldots \epsilon_{2f-1} e(\bi) = 0$ when $i_{2f+1} \neq \frac{\delta-1}{2}$. We apply the induction on $f$. When $f = 0$ it is obvious by (\ref{rela:1}). Suppose for $f - 1$ the result holds. By (\ref{rela:7:2}) we have
\begin{align*}
e(\bi) \epsilon_1 \ldots \epsilon_{2f-1} e(\bi) & = e(\bi) \epsilon_1 \ldots \epsilon_{2f-3} \epsilon_{2f-1} \epsilon_{2f} e(\bl) \epsilon_{2f-1} e(\bi)\\
& = e(\bi) \epsilon_{2f-1} \epsilon_{2f} e(\bl) \epsilon_1 \ldots \epsilon_{2f-3} e(\bl) \epsilon_{2f-1} e(\bi)
\end{align*}
where $\bl \in P^n$. If we write $\bl = (\l_1, \l_2, \ldots, \l_n)$, one can see that $\l_{2f-1} = i_{2f+1} \neq \frac{\delta-1}{2}$ by (\ref{rela:1}), which implies $e(\bl) \epsilon_1 \ldots \epsilon_{2f-3} e(\bl) = 0$ by induction. Hence $e(\bi) \epsilon_1 \ldots \epsilon_{2f-1} e(\bi) = 0$ when $i_{2f+1} \neq \frac{\delta-1}{2}$.

Then we prove $e(\bi) \epsilon_1 \ldots \epsilon_{2f-1} y_{2f+1} e(\bi) = 0$. We apply induction on $f$ as well. When $f = 0$ it is obvious by (\ref{rela:1}). Suppose for $f-1$ the result holds. By (\ref{rela:7:2}) we have
\begin{align*}
e(\bi) \epsilon_1 \ldots \epsilon_{2f-1} y_{2f+1} e(\bi) & = e(\bi) \epsilon_1 \ldots \epsilon_{2f-3} \epsilon_{2f-1} \epsilon_{2f} \epsilon_{2f-1} y_{2f+1} e(\bi)\\
& = e(\bi) \epsilon_{2f-1} \epsilon_{2f} \epsilon_1 \ldots \epsilon_{2f-3} y_{2f-1} \epsilon_{2f-1} e(\bi) = 0,
\end{align*}
where the last equality holds because $\epsilon_1 \ldots \epsilon_{2f-3} y_{2f-1} = 0$ by induction. Hence the first relation of (\ref{rela:1}) consists.

For the second relation of (\ref{rela:1}), notice that $\epsilon_{1,f} e(\bi) = 0$ if $e(\bi) \not\in P_{f,n}$ by (\ref{rela:1}). Hence, we have
$$
\phi_{f,n}(\epsilon_{1,f}{\cdot}1) = \phi_{f,n}(\epsilon_{1,f} \sum_{e(\bi) \in P_{f,n}} e(\bi)) = \sum_{\bj \in P^{n-2f}} e(\bj) = 1.
$$

For (\ref{rela:5:5}) when $\bi \in P_{k,+}^n$, the relation consists by following the same argument as in the proof of~\autoref{A:2}. \endproof

\autoref{remove-front} shows that $\epsilon_{1,f} \G{2f,n} \cong \G{n-2f}$ as $R$-space. For $0 \leq f \leq \lfloor \frac{n}{2} \rfloor$, define $R_{2f,n}(\delta)$ to be the subspace of $R_n(\delta)$ spanned by $\psi_{\s\t}$'s, such that $head(\s), head(\t) \geq f$. By the definition of $\psi_{\s\t}$'s, we have $R_{f,n}(\delta) \subseteq \epsilon_{1,f} \G{2f,n}$. Therefore, we restrict $\phi_{f,n}$ to $R_{f,n}(\delta)$.

\begin{Lemma} \label{res:3}
Suppose $\s,\t \in \Tud_n(\lambda)$ with $head(\s), head(\t) \geq f$. Then we have $\phi_{f,n}(\psi_{\s\t}) = \psi_{\u\v} \in R_{n-2f}(\delta)$. Moreover, if we write $\s = (\alpha_1, \ldots, \alpha_n)$ and $\t = (\beta_1, \ldots, \beta_n)$, then we have $\u = (\alpha_{2f+1}, \ldots, \alpha_n)$ and $\v = (\beta_{2f+1}, \ldots, \beta_n)$.
\end{Lemma}

\proof Because $head(\s), head(\t) \geq f$, we have
\begin{align*}
& \alpha_1 = \alpha_3 = \ldots = \alpha_{2f-1} = -\alpha_2 = -\alpha_4 = \ldots = -\alpha_{2f} = \alpha_0,\\
& \beta_1 = \beta_3 = \ldots = \beta_{2f-1} = -\beta_2 = -\beta_4 = \ldots = -\beta_{2f} = \alpha_0.
\end{align*}

Hence the Lemma follows by direct calculations. \endproof

By~\autoref{res:3}, one can see that $\phi_{f,n}(R_{f,n}(\delta)) = R_{n-2f}(\delta)$. Hence we have a sequence
$$
R_n(\delta) = R_{0,n}(\delta) \supseteq R_{1,n}(\delta) \supseteq R_{2,n}(\delta) \supseteq \ldots.
$$

Finally we introduce some application of restriction property of $R_n(\delta)$.

\begin{Lemma} \label{I:coro:2}
Suppose $(\sigma,f) \in \mathscr S_n$. Then for any $(\lambda,f) \in \widehat B_{n-1}$ and $k \in P$, if $\res(\alpha) \neq k$ for all $\alpha \in \mathscr A(\lambda)$, we have $\theta^{(n-1)}_k(\psi_{\t^{(\lambda,f)} \t^{(\lambda,f)}}) = \sum c_{\u\v} \psi_{\u\v} \in R_n^{>f}(\delta)$, where $c_{\u\v} \neq 0$ only if $head(\u) \geq f$ and $head(\v) \geq f$.
\end{Lemma}

\proof If $f = 0$, the Lemma follows by~\autoref{I:h2}. If $f > 0$, recall $\phi_{f,n}\map{\epsilon_{1,f} \G{2f,n}}{\G{n-2f}}$ defined in~\autoref{A:2}. By the definition of $\theta^{(n-1)}_k$ and $\phi_{f,n}$ one can see that $\phi_{f,n}\circ\theta_k^{(n-1)} = \theta_k^{(n-2f-1)}\circ\phi_{f,n-1}$. By~\autoref{remove-front}, $\phi_{f,n}$ is a bijection. Hence by~\autoref{I:h2} we have
$$
\theta^{(n-1)}_k(\psi_{\t^{(\lambda,f)} \t^{(\lambda,f)}}) = \phi^{-1}_{f,n} (\theta^{(n-2f-1)}_k(\phi_{f,n-1}(\psi_{\t^{(\lambda,f)} \t^{(\lambda,f)}}))) = \phi^{-1}_{f,n}(\theta^{(n-2f-1)}_k(\psi_{\t^{(\lambda,0)} \t^{(\lambda,0)}})) = \sum_{\s,\t}c_{\s\t} \phi^{-1}_{f,n}(\psi_{\s\t}),
$$
where $\s,\t \in \Tud_{n-2f}(\mu)$ with $\mu \vdash n-2f-2m$ and $m > 0$. Hence the Lemma follows by~\autoref{res:3}. \endproof

\begin{Lemma} \label{remove:head}
Suppose $(\lambda,f) \in \mathscr S_n$ and $\t \in \Tud_n(\lambda)$ with head $h > 0$. Then for any $a \in \G{2h,n}$, we have
$$
\psi_{\t^{(\lambda,f)} \t} a \equiv \sum_{\v \in \Tud_n(\lambda)} c_\v \psi_{\t^{(\lambda,f)} \v} \pmod{R_n^{>f}(\delta)},
$$
where $c_\v \neq 0$ only if $head(\v) \geq h$.
\end{Lemma}

\proof Suppose $\t = (\alpha_1, \ldots, \alpha_n)$. Because $head(\t) = h \leq f$, we have $\alpha_1 = \alpha_3 = \ldots = \alpha_{2h-1} = -\alpha_2 = -\alpha_4 = \ldots = -\alpha_{2h} = \alpha_0$. Define $\dot\t = (\alpha_{2h+1}, \alpha_{2h+2}, \ldots, \alpha_n)$ and we have $\dot\t \in \Tud_{n-2h}(\lambda)$.

Recall $\phi_{h,n}\map{\epsilon_{1,h} \G{2h,n}}{\G{n-2h}}$ defined in~\autoref{A:2}. Because $\psi_{\t^{(\lambda,f)} \t}a \in \epsilon_{1,h} \G{2h,n}$, by~\autoref{res:3}, we have $\phi_{h,n}(\psi_{\t^{(\lambda,f)} \t}a) = \psi_{\t^{(\lambda,f-h)} \dot\t}a'$, for some $a' \in \G{n-2h}$.

Because $(\lambda,f) \in \mathscr S_n$ and $h > 0$, we have $(\lambda,f-h) \in \widehat B_{n-2h} \subset \widehat{\mathscr B}$. Hence by the definition of $\widehat {\mathscr B}$, we have
$$
\phi_{h,n}(\psi_{\t^{(\lambda,f)} \t}a) = \psi_{\t^{(\lambda,f-h)} \dot\t}a' \equiv \sum_{\dot\v \in \Tud_{n-2h}(\lambda)} c_{\dot\v} \psi_{\t^{(\lambda,f-h)} \dot\v} \pmod{R_{n-2h}^{>f-h}(\delta)},
$$
which implies $\psi_{\t^{(\lambda,f)} \t} a \equiv \sum_{\dot\v \in \Tud_{n-2h}(\lambda)} c_{\dot\v} \phi_{h,n}^{-1}(\psi_{\t^{(\lambda,f-h)} \dot\v}) \equiv \sum_{\v \in \Tud_n(\lambda)} c_\v \psi_{\t^{(\lambda,f)} \v} \pmod{R_n^{>f}(\delta)}$, where $c_\v \neq 0$ only if $head(\v) \geq h$. \endproof

\section{A spanning set of $\G{n}$} \label{sec:span}

In this section we prove $\set{\psi_{\s\t} | \s,\t \in \Tud_n(\lambda), (\lambda,f) \in \widehat B_n}$ is a $R$-spanning set of $\G{n}$. The main idea is to prove $\bigcup_{i=1}^n \widehat B_i \subseteq \widehat {\mathscr B}$ by induction on $\mathscr S_n$. We show that if $(\lambda,f) \in \mathscr S_n$, we have
\begin{equation} \label{Induction:1}
\psi_{\t^{(\lambda,f)} \t} a \equiv \sum_{\v \in \Tud_n(\lambda)} c_{\v} \psi_{\t^{(\lambda,f)}\v} \pmod{R_n^{>f}(\delta)},
\end{equation}
for any $\t \in \Tud_n(\lambda)$ and $a \in \G{n}$. As a byproduct, (\ref{Induction:1}) shows the cellular-like property of $\psi_{\s\t}$'s, which will directly apply that $\set{\psi_{\s\t} | \s,\t \in \Tud_n(\lambda), (\lambda,f) \in \widehat B_n}$ is a graded cellular basis of $\G{n}$ after we prove $\G{n} \cong \B$.

\subsection{The base case}

Fix $(\lambda,f) \in \mathscr S_n$. We start by proving (\ref{Induction:1}) in the most simple case: when $\t \in \Tud_n(\lambda)$ with head $f$, which will be used in the following subsections for computational purposes when we prove more complicated cases. It suffices to prove (\ref{Induction:1}) when $a$ is one of the generators of $\G{n}$. For $e(\bi)$ with $\bi \in P^n$ we have $\psi_{\t^{(\lambda,t)} \t} e(\bi) = \delta_{\bi,\bi_\t} \psi_{\t^{(\lambda,t)} \t}$. Therefore it left us to consider $y_k, \psi_s$ and $\epsilon_s$ with $1 \leq k \leq n$ and $1 \leq s \leq n-1$.

\begin{Lemma} \label{commute:1:h1}
Suppose $(\lambda,f) \in \mathscr S_n$ and $\t \in \Tud_n(\lambda)$ with head $f$. For $2f+1 \leq k \leq n-1$, we have $\psi_{\t^{(\lambda,f)} \t} \epsilon_k \in R_n^{>f}(\delta)$.
\end{Lemma}

\proof We have $\psi_{\t^{(\lambda,f)} \t} = e_{(\lambda,f)} \psi_{d(\t)}$ with $d(\t) \in \Sym_{2f,n}$. Hence for any $\bk \in P^n$, by~\autoref{I:coro:1} and~\autoref{two:sided:ideal} we have $\psi_{\t^{(\lambda,f)} \t} \epsilon_k e(\bk) = e(\bi_{(\lambda,f)}) \psi_{d(\t)} \epsilon_1 \epsilon_3 \ldots \epsilon_{2f-1} \epsilon_k e(\bk) \in R_n^{>f}(\delta)$, which proves the Lemma. \endproof

\begin{Lemma} \label{base:y:1}
Suppose $(\lambda,f) \in \mathscr S_n$ and $\t \in \Tud_n(\lambda)$ with head $f$. Then for any $1 \leq k \leq n$, we have $\psi_{\t^{(\lambda,f)} \t} y_k \in R_n^{>f}(\delta)$.
\end{Lemma}

\proof First we prove that for
\begin{equation} \label{base:y:1:eq1}
\psi_{\t^{(\lambda,f)} \t^{(\lambda,f)}} y_k = e_{(\lambda,f)} y_k \in R_n^{>f}(\delta),
\end{equation}
the Lemma holds. When $f = 0$, the Lemma follows by~\autoref{I:h2} and~\autoref{I:h3}. When $f \geq 1$, if $k = 1$, the Lemma follows by (\ref{rela:1}); and if $k = 2$, by (\ref{rela:7:2}) we have $e_{(\lambda,f)}y_2 = - e{(\lambda,f)} y_1$, and the Lemma follows; and if $k \geq 3$ the Lemma follows by~\autoref{remove:head}.

For arbitrary $\t \in \Tud_n(\lambda)$ with head $f$, suppose $d(\t) = s_{r_1} s_{r_2} \ldots s_{r_\l} \in \Sym_{2f,n}$ is a reduced expression of $d(\t)$. We prove the Lemma by induction. As the base step, when $\l = 0$ the Lemma follows by (\ref{base:y:1:eq1}). For the induction step, we assume that when $\l < \l'$ the Lemma holds. When $\l = \l'$, set $\s = \t^{(\lambda,f)} s_{r_1} \ldots s_{r_{\l-1}} \in \Tud_n(\lambda)$ and $r_\l = r$. One can see that $\s = \t{\cdot}s_r$. If we write $\bi_\t = (i_1, \ldots, i_n)$, by~\autoref{deg:help1} we have $|i_r - i_{r+1}| > 1$. Hence by (\ref{rela:3:1}) and (\ref{rela:3:2}), we have
$$
\psi_{\t^{(\lambda,f)} \t} y_k = \psi_{\t^{(\lambda,f)} \s} y_{s_r(k)} \psi_r \pm\psi_{\t^{(\lambda,f)} \s} \epsilon_r e(\bi_\t) \in R_n^{>f}(\delta),
$$
where the first term is in $R_n^{>f}(\delta)$ by induction, and the second term is in $R_n^{>f}(\delta)$ by~\autoref{commute:1:h1}. Hence we proves the Lemma. \endproof

\begin{Lemma} \label{base:psi:1}
Suppose $(\lambda,f) \in \mathscr S_n$ and $\t \in \Tud_n(\lambda)$ with head $f$. Then for any $1 \leq k \leq n - 1$, we have
$$
\begin{cases}
\psi_{\t^{(\lambda,f)} \t} \psi_k = \psi_{\t^{(\lambda,f)} \v} & \text{if $\v = \t {\cdot}s_k \in \Tud_n(\lambda)$,}\\
\psi_{\t^{(\lambda,f)} \t} \psi_k \in R_n^{>f}(\delta) & \text{if $\t{\cdot}s_k$ is not an up-down tableau.}
\end{cases}
$$
\end{Lemma}

\proof We prove the Lemma by considering consider the following different cases.

\textbf{Case 1:} $1 \leq k \leq 2f$.

In this case, as $\t$ has head $f$, $\t{\cdot}s_k$ is not an up-down tableau. Write $\bi_\t = (i_1, i_2, \ldots, i_n)$. Then we have $i_1 = i_3 = \ldots = i_{2f-1} = \frac{\delta-1}{2}$ and $i_2 = i_4 = \ldots = i_{2f} = -\frac{\delta-1}{2}$.

When $\frac{\delta-1}{2} = 0$, we have $\bi_\t \in P_{k,0}^n$. By (\ref{rela:6:1}) we have $e(\bi_\t) \psi_k \epsilon_k e(\bi_\t) = e(\bi_\t) \psi_k e(\bi_\t) \epsilon_k e(\bi_\t) = 0$ and by (\ref{rela:5:1}) we have $e(\bi_\t) \epsilon_k e(\bi_\t) = (-1)^{a_k(\bi_\t)} e(\bi_\t)$. Hence, we have $e(\bi_\t) \psi_k e(\bi_\t) = e(\bi_\t) \psi_k = 0$, which implies $\psi_{\t^{(\lambda,f)} \t} \psi_k = 0$ by~\autoref{idem:psi:1}.

When $\frac{\delta-1}{2} \neq 0$, we have $h_k(\bi_\t) < 0$ by~\autoref{deg:h1}, which implies $h_k(\bi_\t{\cdot}s_k) > 0$. Hence $e(\bi{\cdot}s_k) = 0$ by~\autoref{I:coro:3}. Therefore by (\ref{rela:2:1}), we have $e(\bi_\t) \psi_k = \psi_k e(\bi_\t{\cdot}s_k) = 0$, which implies $\psi_{\t^{(\lambda,f)} \t} \psi_k = 0$ by~\autoref{idem:psi:1}. So the Lemma follows when $1 \leq k \leq 2f$.

\textbf{Case 2:} $2f+1 \leq k \leq n-1$ and $\v = \t{\cdot}s_k \in \Tud_n(\lambda)$.

In this case, as $\t \in \Tud_n(\lambda)$ with head $f$, we have $\psi_{\t^{(\lambda,f)} \t} = e_{(\lambda,f)} \psi_{d(\t)}$. By~\autoref{semi-red2} and~\autoref{semi-red3}, we have $d(\t){\cdot}s_k$ is semi-reduced correspond to $\t^{(\lambda,f)}$. Hence by~\autoref{semi-red1}, we have $e(\bi_{(\lambda,f)})\psi_{d(\t)} \psi_k = e(\bi_{(\lambda,f)}) \psi_{d(\t){\cdot}s_k}$. Furthermore, as $\v = \t{\cdot}s_k$, we have $d(\v) = d(\t){\cdot}s_k$, which implies $\psi_{d(\t){\cdot}s_k} = \psi_{d(\v)}$. Therefore, $\psi_{\t^{(\lambda,f)} \t}\psi_k = \psi_{\t^{(\lambda,f)} \v}$, where the Lemma holds.

\textbf{Case 3:} $2f+1 \leq k \leq n-1$ and $\t{\cdot}s_k$ is not an up-down tableau.

Write $\mu = \t_{k-1}$. As $\t{\cdot}s_k$ is not an up-down tableau, we have $\res(\alpha) \neq i_{k+1}$ for all $\alpha \in \mathscr A(\mu)$. Therefore, $e(\bi_\t{\cdot}s_k) \epsilon_1 \ldots \epsilon_{2f-1} e(\bi_\t{\cdot}s_k) \in R_k^{>f}(\delta)$ by~\autoref{I:h2} and~\autoref{I:h3}. Because we have
\begin{align*}
\psi_{\t^{(\lambda,f)} \t} \psi_k & = e(\bi_{(\lambda,f)}) \epsilon_1 \ldots \epsilon_{2f-1} e(\bi_{(\lambda,f)}) \psi_{d(\t)} \psi_k \\
& = e(\bi_{(\lambda,f)}) \psi_{d(\t)} \psi_k e(\bi_\t{\cdot}s_k) \epsilon_1 \ldots \epsilon_{2f-1} e(\bi_\t{\cdot}s_k),
\end{align*}
the Lemma follows by~\autoref{two:sided:ideal}. \endproof

\begin{Lemma} \label{base:ep:1}
Suppose $(\lambda,f) \in \mathscr S_n$ and $\t \in \Tud_n(\lambda)$ with head $f$. Then for any $1 \leq k \leq n-1$, we have
$$
\begin{cases}
\psi_{\t^{(\lambda,f)} \t} \epsilon_k \equiv \sum_{\v \in \Tud_n(\lambda)} c_\v \psi_{\t^{(\lambda,f)} \v} \pmod{R_n^{>f}(\delta)}, & \text{if $1 \leq k \leq 2f$,}\\
\psi_{\t^{(\lambda,f)} \t} \epsilon_k \in R_n^{>f}, & \text{if $2f+1 \leq k \leq n-1$.}
\end{cases}
$$
\end{Lemma}

\proof Write $\bi_{(\lambda,f)} = (i_1, \ldots, i_n)$. We consider different values of $k$.

\textbf{Case 1:} $k = 2\l - 1$ for $1 \leq \l \leq f$.

We have $h_k(\bi_{(\lambda,f)}) = -1$. Because $\psi_{d(\t)} \in \Sym_{2f,n}$, by (\ref{rela:5:5}), we have
$$
\psi_{\t^{(\lambda,f)} \t} \epsilon_{2\l-1} = e_{(\lambda,f)} \epsilon_{2\l-1} \psi_{d(\t)} =
\begin{cases}
(-1)^{a_{2\l-1}(\bi_{(\lambda,f)})} \psi_{\t^{(\lambda,f)} \t}, & \text{if $i_{2\l - 1} \neq -\frac{1}{2}$,}\\
2(-1)^{a_{2\l-1}(\bi_{(\lambda,f)})} \psi_{\t^{(\lambda,f)} \t} f(y_1, \ldots, y_{2\l-2}), & \text{if $i_{2\l - 1} = -\frac{1}{2}$,}
\end{cases}
$$
where $f(y_1, \ldots, y_{2\l-2})$ is a polynomial of $y_1, \ldots, y_{2\l-2}$. Hence, by~\autoref{base:y:1}, the Lemma holds when $k = 2\l - 1$ with $1 \leq \l \leq f$.

\textbf{Case 2:} $k = 2\l$ with $1 \leq \l \leq f$.

Write $\t = (\alpha_1, \ldots, \alpha_n)$. Let $\beta = (2,1)$ and $\gamma = (1,2)$, and define
\begin{align*}
\u & = (\alpha_1, \ldots, \alpha_{2\l-1}, \beta, -\beta, \alpha_{2\l+2}, \ldots, \alpha_n),\\
\v & = (\alpha_1, \ldots, \alpha_{2\l-1}, \gamma, -\gamma, \alpha_{2\l+2}, \ldots, \alpha_n).
\end{align*}

If $\frac{\delta-1}{2} = \pm \frac{1}{2}$, we have $\bi_\u = \bi_\t$ or $\bi_\v = \bi_\t$. Hence, by (\ref{rela:1}), we have $\psi_{\t^{(\lambda,f)} \t} \epsilon_k = \psi_{\t^{(\lambda,f)} \t} \epsilon_k (e(\bi_\u) + e(\bi_\v))$. By directly comparing both sides of the equation, we have $\psi_{\t^{(\lambda,f)} \t} \epsilon_k = \psi_{\t^{(\lambda,f)} \u} + \psi_{\t^{(\lambda,f)} \v}$, because $\psi_{\t^{(\lambda,f)} \t} \epsilon_k e(\bi_\u) = \psi_{\t^{(\lambda,f)} \u}$ and $\psi_{\t^{(\lambda,f)} \t} \epsilon_k e(\bi_\v) = \psi_{\t^{(\lambda,f)} \v}$.

If $\frac{\delta-1}{2} \neq \pm \frac{1}{2}$, we have $\bi_\t \neq \bi_\u, \bi_\v$. Hence, by (\ref{rela:1}), we have $\psi_{\t^{(\lambda,f)} \t} \epsilon_k = \psi_{\t^{(\lambda,f)} \t} \epsilon_k (e(\bi_\t) + e(\bi_\u) + e(\bi_\v))$. Because $h_k(\bi_\t) = -1$, by (\ref{rela:5:1}) we have $\psi_{\t^{(\lambda,f)} \t} \epsilon_k e(\bi_\t) = \pm \psi_{\t^{(\lambda,f)} \t}$. Following the similar argument as above, we have $\psi_{\t^{(\lambda,f)} \t} \epsilon_k = \pm \psi_{\t^{(\lambda,f)} \t} + \psi_{\t^{(\lambda,f)} \u} + \psi_{\t^{(\lambda,f)} \v}$, which proves the Lemma when $1 \leq k \leq 2f$.

\textbf{Case 3:} $2f + 1 \leq k \leq n-1$.

In this case the Lemma follows by ~\autoref{commute:1:h1}. \endproof

\subsection{A weak version} \label{sec:induction:1}

In this subsection we prove (\ref{Induction:1}) for $a \in \G{n-1}$ by considering $\G{n-1}$ as a subalgebra of $\G{n}$. We separate the question into two cases by considering $\t(n) > 0$ and $\t(n) < 0$. For $\t \in \Tud_n(\lambda)$ with $\t(n) > 0$, (\ref{Induction:1}) for $a \in \G{n-1}$ is directly implied by~\autoref{remove:tail:1}. Hence we have the following Lemma.

\begin{Lemma} \label{y:first:1}
Suppose $(\lambda,f) \in \mathscr S_n$ and $\t \in \Tud_n(\lambda)$ with $\t(n) > 0$. Then the equality (\ref{Induction:1}) holds when $a \in \G{n-1}$ and $\t \in \Tud_n(\lambda)$ with $\t(n) > 0$.
\end{Lemma}

To prove (\ref{Induction:1}) for $a \in \G{n-1}$ and $\t \in \Tud_n(\lambda)$ with $\t(n) < 0$, first we introduce a commutation rule of $\G{n}$ as a technical result.

Suppose $(\lambda,f) \in \widehat B_n$ and $\t \in \Tud_n(\lambda)$ with head $f-1$ and $\t(n) < 0$. Let $\s = h(\t) \rightarrow \t$. It is easy to see that $\rho(\s,\t) = (a,n)$ because $\t(n) < 0$ and $head(\t) = f-1$. The following Lemma gives the commutation rule of $\psi_k$ and $\epsilon_{\s \rightarrow \t}$ when $2f+2 \leq k \leq n-1$.

\begin{Lemma} \label{y:h2:2}
Suppose $(\lambda,f) \in \widehat B_n$ and $\t \in \Tud_n(\lambda)$ with head $f-1$ and $\t(n) < 0$. Let $\s = h(\t) \rightarrow \t$ and $\rho(\s,\t) = (a,n)$. If $\s{\cdot}s_k \in \Tud_n(\lambda)$ for some $2f + 2 \leq k \leq n-1$, there exists $w \in \Sym_{2f-2, n-1}$ such that $w$ is semi-reduced correspond to $\t$, and $\psi_k \epsilon_{\s \rightarrow \t} = \epsilon_{\s{\cdot}s_k \rightarrow \t{\cdot}w} \psi_{w^{-1}}$. In more details, we have
$$
w =
\begin{cases}
s_{k-2}, & \text{if $2f+2 \leq k \leq a$,}\\
s_{k-1}s_{k-2}, & \text{if $k = a+1$,}\\
s_{k-1}, & \text{if $a + 1 < k \leq n-1$.}
\end{cases}
$$
\end{Lemma}

\proof As $\rho(\s,\t) = (a,n)$ and $\s \in \Tud_n(\lambda)$ with head $f$, we can write $\epsilon_{\s\rightarrow \t} = e(\bi_\s) \epsilon_{2f}\epsilon_{2f+1} \ldots \epsilon_a \psi_{a+1} \ldots \psi_{n-1} e(\bi_\t)$. Write $\t = (\alpha_1, \ldots, \alpha_n)$ and $\s = (\alpha_0, -\alpha_0, \alpha_1, \ldots, \alpha_{a-1}, \alpha_{a+1}, \ldots, \alpha_{n-1})$. We consider different values of $k$.

\textbf{Case 1:} $2f+2 \leq k \leq a$.

By (\ref{rela:12}) we have
$$
\psi_k \epsilon_{\s \rightarrow \t} = e(\bi_\s{\cdot}s_k) \epsilon_{2f} \epsilon_{2f+1} \ldots \epsilon_{k-2} \epsilon_{k-1} \epsilon_k \ldots \epsilon_a \psi_{a+1} \ldots \psi_{n-1} e(\bi_\s{\cdot}s_{k-2}) \psi_{k-2}.
$$

We write
\begin{align*}
\s{\cdot}s_k & = (\alpha_0, -\alpha_0, \alpha_1, \ldots, \alpha_{k-3}, \alpha_{k-1}, \alpha_{k-2}, \alpha_k, \ldots, \alpha_{a-1}, \alpha_{a+1}, \ldots, \alpha_{n-1}),\\
\t{\cdot}s_{k-2} & = (\alpha_1, \ldots, \alpha_{k-3}, \alpha_{k-1}, \alpha_{k-2}, \alpha_k, \ldots, \alpha_n).
\end{align*}

As $\s{\cdot}s_k \in \Tud_n(\lambda)$, by~\autoref{y:h2:8}, $\alpha_{k-2}$ and $\alpha_{k-1}$ are not adjacent, which implies $\t{\cdot}s_{k-2} \in \Tud_n(\lambda)$ by~\autoref{y:h2:8}. By the above expression, we have $\s{\cdot}s_k \rightarrow \t{\cdot}s_{k-2}$ and
$$
\epsilon_{\s{\cdot}s_k \rightarrow \t{\cdot}s_{k-2}} = e(\bi_\s{\cdot}s_k) \epsilon_{2f} \epsilon_{2f+1} \ldots \epsilon_a \psi_{a+1} \ldots \psi_{n-1} e(\bi_\t{\cdot}s_{k-2}).
$$

Hence by setting $w = s_{k-2}$, the Lemma follows.

\textbf{Case 2:} $k = a+1$.

By (\ref{rela:14}) we have
$$
\psi_k \epsilon_{\s \rightarrow \t} = \psi_{a+1} \epsilon_{\s \rightarrow \t} = e(\bi_\s{\cdot}s_{a+1}) \epsilon_{2f} \epsilon_{2f+1} \ldots \epsilon_{a-1} \epsilon_a \epsilon_{a+1} \psi_{a+2} \ldots \psi_{n-1} e(\bi_\t{\cdot} s_a s_{a-1}) \psi_{a-1} \psi_a.
$$

We write
\begin{align*}
\s{\cdot}s_k = \s{\cdot}s_{a+1} & = (\alpha_0, -\alpha_0, \alpha_1, \ldots, \alpha_{a-2}, \alpha_{a+1}, \alpha_{a-1}, \alpha_{a+2}, \ldots, \alpha_{n-1}),\\
\t{\cdot}s_{k-1} = \t{\cdot}s_a & = (\alpha_1, \ldots, \alpha_{a-2}, \alpha_{a-1}, \alpha_{a+1}, \alpha_a, \alpha_{a+2},\ldots, \alpha_n),\\
\t{\cdot}s_{k-1}s_{k-2} = \t{\cdot}s_a s_{a-1} & = (\alpha_1, \ldots, \alpha_{a-2}, \alpha_{a+1}, \alpha_{a-1}, \alpha_a, \alpha_{a+2},\ldots, \alpha_n).
\end{align*}

By~\autoref{deg:help1}, $\alpha_{a+1}$ is not adjacent to $\alpha_a$, which implies $\t{\cdot}s_{k-1} \in \Tud_n(\lambda)$ by~\autoref{y:h2:8}. As $\s{\cdot}s_k \in \Tud_n(\lambda)$, by~\autoref{y:h2:8} we have $\alpha_{a-1}$ and $\alpha_{a+1}$ are not adjacent, which implies $\t{\cdot}s_{k-1}s_{k-2} \in \Tud_n(\lambda)$ by~\autoref{y:h2:8}. Hence $s_{k-1}s_{k-2}$ is semi-reduced correspond to $\t$. By the above expression, we have $\s{\cdot}s_k \rightarrow \t{\cdot}s_{k-1}s_{k-2}$ and
$$
\epsilon_{\s{\cdot}s_k \rightarrow \t{\cdot}s_{k-1}s_{k-2}} = e(\bi_\s{\cdot}s_{a+1}) \epsilon_{2f} \epsilon_{2f+1} \ldots \epsilon_{a-1} \epsilon_a \epsilon_{a+1} \psi_{a+2} \ldots \psi_{n-1} e(\bi_\t{\cdot} s_a s_{a-1}).
$$

By setting $w = s_{k-1} s_{k-2}$, the Lemma follows.

\textbf{Case 3:} $a + 1 < k \leq n-1$.

Because $\psi_k$ commutes with $\epsilon_1 \epsilon_3 \ldots \epsilon_{2f-1}$, we have
$$
\psi_k \epsilon_{\s \rightarrow \t} = e(\bi_\s{\cdot}s_k) \epsilon_{2f} \epsilon_{2f+1} \ldots \epsilon_a \psi_{a+1} \ldots \psi_{k-2} \psi_k \psi_{k-1} \psi_k e(\bi_\t{\cdot} s_{n-1} s_{n-2} \ldots s_{k+1}) \psi_{k+1} \ldots \psi_{n-1}.
$$

If we write $\bi_\t = (i_1, \ldots, i_n)$, the relation of $\psi_k \psi_{k-1} \psi_k  e(\bi_\t{\cdot} s_{n-1} s_{n-2} \ldots s_{k+1})$ is determined by $i_{k-1}, i_k$ and $i_n$. Notice that $\alpha_r > 0$ and $i_r = \res(\alpha_r)$ for $a \leq r \leq n-1$, and $i_n = -\res(\alpha_a)$.

Because $\alpha_n = -\alpha_a$, we have $\alpha_a \in \mathscr R(\t_{n-1})$. By~\autoref{deg:help1}, we have $\alpha_r > 0$ for all $a < r < n$. Hence, by the construction of up-down tableau we have $\res(\alpha_r) \neq \res(\alpha_a) = i_a = -i_n$ for any $a < r < n$. Therefore we have $i_{k-1} + i_n \neq 0$ and $i_k + i_n \neq 0$. Hence (\ref{rela:8:1}), (\ref{rela:8:2}), (\ref{rela:8:5}) and (\ref{rela:8:6}) will not apply here; if $i_{k-1} + i_k = 0$, by~\autoref{deg:help1}, $\alpha_{k-1}$ and $\alpha_k$ are not adjacent to $\alpha_a$. $\alpha_k$ not adjacent to $\alpha_a$ implies $i_k \neq -i_n \pm 1$ and $\alpha_{k-1}$ not adjacent to $\alpha_a$ implies $i_k = -i_{k-1} \neq i_n \pm 1$. Hence (\ref{rela:8:3}) and (\ref{rela:8:4}) will not apply here; as $\s{\cdot}s_k \in \Tud_n(\lambda)$, by~\autoref{y:h2:8}, $\alpha_{k-1}$ and $\alpha_k$ are not adjacent, which implies $|i_{k-1} - i_k| > 1$. Hence (\ref{rela:8:7}) and (\ref{rela:8:8}) will not apply here. Therefore we have $\psi_k \psi_{k-1} \psi_k  e(\bi_\t{\cdot} s_{n-1} s_{n-2} \ldots s_{k+1}) = \psi_{k-1} \psi_k \psi_{k-1} e(\bi_\t{\cdot} s_{n-1} s_{n-2} \ldots s_{k+1})$, which implies
$$
\psi_k \epsilon_{\s \rightarrow \t} = e(\bi_\s{\cdot}s_k) \epsilon_{2f} \epsilon_{2f+1} \ldots \epsilon_a \psi_{a+1} \ldots \psi_{k-2} \psi_{k-1} \psi_k \psi_{k+1} \ldots \psi_{n-1} e(\bi_\t{\cdot} s_{k-1}) \psi_{k-1}.
$$

We write
\begin{align*}
\s{\cdot}s_k & = (\alpha_0, -\alpha_0, \alpha_1, \ldots, \alpha_{a-1}, \alpha_{a+1}, \ldots, \alpha_{k-2}, \alpha_k, \alpha_{k-1},\alpha_{k+1}, \ldots, \alpha_{n-1}),\\
\t{\cdot}s_{k-1} & = (\alpha_1, \ldots, \alpha_{k-2}, \alpha_k, \alpha_{k-1}, \alpha_{k+1}, \ldots, \alpha_n).
\end{align*}

As $\s{\cdot} s_k \in \Tud_n(\lambda)$, by~\autoref{y:h2:8}, $\alpha_{k-1}$ and $\alpha_k$ are not adjacent, which implies $\t{\cdot}s_{k-1} \in \Tud_n(\lambda)$ by~\autoref{y:h2:8}. By the above expression, we have $\s{\cdot}s_k \rightarrow \t{\cdot}s_{k-1}$ and
$$
\epsilon_{\s{\cdot}s_k \rightarrow \t{\cdot}s_{k-1}} = e(\bi_\s{\cdot}s_k) \epsilon_{2f} \epsilon_{2f+1} \ldots \epsilon_a \psi_{a+1} \ldots \psi_{k-2} \psi_{k-1} \psi_k \psi_{k+1} \ldots \psi_{n-1} e(\bi_\t{\cdot} s_{k-1})
$$

By setting $w = s_{k-1}$, the Lemma follows. \endproof

Suppose $\s$ and $\t$ are defined as in~\autoref{y:h2:2}. Let $\mu = \t_{n-1}$ and $\u \in \Tud_n(\lambda)$ with $\u|_{n-1} = \t^{(\mu,f-1)}$. It is easy to see that $\u(n) = \t(n)$. We abuse the symbol and define $d(\t) \in \Sym_{2f-2,n-1}$ such that $\t = \u{\cdot}d(\t)$. Moreover, by the definition of $\u$, we have $\t^{(\lambda,f)} \rightarrow \u$.

\begin{Example}
Suppose $(\lambda,f) = ((2),2)$ and $n = 6$. Hence we have $(\lambda,f) \in \widehat B_n$. Let
$$
\t = \left(\emptyset, \ydiag(1), \emptyset, \ydiag(1), \ydiag(1,1), \ydiag(2,1), \ydiag(2) \right).
$$

Then we have $\t \in \Tud_n(\lambda)$ with head $f-1$ and $\t(n) < 0$. We can find a unique $\s \in \Tud_n(\lambda)$ such that $\s = h(\t) \rightarrow \t$, which is
$$
\s = \left(\emptyset, \ydiag(1), \emptyset, \ydiag(1), \emptyset, \ydiag(1), \ydiag(2) \right).
$$

We have $\mu = \t_{n-1} = (2,1)$, and define
$$
\u = \left(\emptyset, \ydiag(1), \emptyset, \ydiag(1), \ydiag(2), \ydiag(2,1), \ydiag(2) \right).
$$

Then we have $\u \in \Tud_n(\lambda)$ where $\u|_{n-1} = \t^{(\mu,f-1)}$, and $\u(n) = \t(n)$. We have $d(\t) = s_4 \in \Sym_{2f-1,n-1}$ and $\t = \u{\cdot}d(\t)$. Also, we have $\t^{(\lambda,f)} \rightarrow \u$.
\end{Example}

The next Lemma is an extension of~\autoref{y:h2:2}.

\begin{Lemma} \label{y:h2}
Suppose $\s,\t$ and $\u$ are defined as above. Then we have $\psi_{d(\s)} \epsilon_{\s \rightarrow \t} = \epsilon_{\t^{(\lambda,f)} \rightarrow \u} \psi_{d(\t)}$.
\end{Lemma}

\proof By~\autoref{y:h2:2}, there exists $w \in \Sym_{2f-2,n-1}$ such that $\w$ is semi-reduced correspond to $\t$ and $\psi_{d(\s)} \epsilon_{\s \rightarrow \t} = \epsilon_{\t^{(\lambda,f)} \rightarrow \t{\cdot}w} \psi_{w^{-1}}$.

Write $\u = (\alpha_1, \ldots, \alpha_n)$ and $\rho(\t^{(\lambda,f)},\u) = (a,n)$. Define $\v = \t{\cdot}w$. As $w \in \Sym_{n-1}$, we have $\v(n) = \t(n) = \u(n)$. Hence $\t^{(\lambda,f)} \rightarrow \v$ and we write $\rho(\t^{(\lambda,f)},\v) = (b,n)$. By the construction of $\u$ and $\v$, we have $b \geq a$.

Define $d(\v)$ such that $\v = \u{\cdot}d(\v)$. Because $b \geq a$, we have $d(\v) = s_a s_{a+1} \ldots s_{b-1}$. By~\autoref{deg:help1}, for any $a < r < n$, $\alpha_r$ is not adjacent to $\alpha_a$ and $\alpha_r > 0$. Hence by~\autoref{y:h2:8}, for any $a \leq k \leq b-1$, $\u{\cdot}s_a s_{a+1} \ldots s_k \in \Tud_n(\lambda)$. Therefore $d(\v)$ is semi-reduced correspond to $\u$. As $w$ is semi-reduced correspond to $\t$, $w^{-1}$ is semi-reduced correspond to $\v = \t{\cdot}w$. By~\autoref{semi-red2}, $d(\v) w^{-1}$ is semi-reduced correspond to $\u$. Hence we have $e(\bi_\u) \psi_{d(\v)}\psi_{w^{-1}} = e(\bi_\u) \psi_{d(\t)}$ by~\autoref{semi-red1} as $\u{\cdot}d(\v)w^{-1} = \u{\cdot}d(\t) = \t$.

Because $d(\v)$ is semi-reduced correspond to $\u$, by~\autoref{A:1} we have $e(\bi_\v) = \psi_{d(\v)^{-1}} e(\bi_\u) \psi_{d(\v)}$. Therefore by (\ref{rela:9}), we have
\begin{align*}
\epsilon_{\t^{(\lambda,f)} \rightarrow \v} \psi_{w^{-1}} & = e(\bi_{(\lambda,f)}) \epsilon_{2f} \epsilon_{2f+1} \ldots, \epsilon_b \psi_{b+1} \ldots \psi_{n-1} e(\bi_\v) \psi_{w^{-1}}\\
& = e(\bi_{(\lambda,f)}) \epsilon_{2f} \epsilon_{2f+1} \ldots, \epsilon_b \psi_{b+1} \ldots \psi_{n-1} \psi_{b-1} \psi_{b-2} \ldots \psi_a e(\bi_\u) \psi_{d(\v)} \psi_{w^{-1}}\\
& = e(\bi_{(\lambda,f)}) \epsilon_{2f} \epsilon_{2f+1} \ldots, \epsilon_a \psi_{a+1} \ldots \psi_{n-1} e(\bi_\u) \psi_{d(\v)} \psi_{w^{-1}} = e_{\t^{(\lambda,f)} \rightarrow \u} \psi_{d(\t)},
\end{align*}
which completes the proof. \endproof

The key point of~\autoref{y:h2} is it implies the next Lemma, which allow us to apply induction to prove (\ref{Induction:1}) for $a \in \G{n-1}$ and $\t \in \Tud_n(\lambda)$ with $\t(n) < 0$.

In the rest of this subsection we fix $\t \in \Tud_n(\lambda)$ with $\t(n) = -\alpha < 0$, $\mu = \lambda \cup \{\alpha\}$ and $\u \in \Tud_n(\lambda)$ such that $\u|_{n-1} = \t^{(\mu,f-1)}$ and $\u(n) = \t(n)$.

\begin{Lemma} \label{y:h4}
Suppose $(\lambda,f) \in \mathscr S_n$. Then
\begin{enumerate}
\item we have
\begin{equation} \label{y:h4:1}
\psi_{\t^{(\lambda,f)} \t} = \psi_{\t^{(\lambda,f)} \u} \epsilon_2 \epsilon_4 \ldots \epsilon_{2f-2} \theta_k^{(n-1)}(\psi_{\t^{(\mu,f-1)} \dot\t}),
\end{equation}
where $k = -\res(\alpha) \in P$ and $\dot\t = \t|_{n-1} \in \Tud_{n-1}(\mu)$.

\item for any $a \in \G{n-1}$ we have
\begin{equation} \label{y:h4:2}
\psi_{\t^{(\lambda,f)} \t} a \equiv \sum_{\v \in \Tud_n(\lambda)} c_{\v} \psi_{\t^{(\lambda,f)} \v} + \sum_{\substack{\dot\x, \dot\y \in \Tud_{n-1}(\gamma) \\ (\gamma,f) \in \widehat B_{n-1}}} c_{\dot\x \dot\y} \psi_{\t^{(\lambda,f)} \u} \epsilon_2 \epsilon_4 \ldots \epsilon_{2f-2} \theta_k^{(n-1)}(\psi_{\dot\x \dot\y}) \pmod{R_n^{>f}(\delta)},
\end{equation}
where $c_{\dot\x \dot\y} \neq 0$ only if $\bi_{\dot\x} = \bi_{(\mu,f-1)}$.

\item for any $a \in \G{n-1}$ we have
\begin{equation} \label{y:h4:3}
\psi_{\t^{(\lambda,f)} \t} a \equiv \sum_{\v \in \Tud_n(\lambda)} c_\v \psi_{\t^{(\lambda,f)} \v} + \sum_{\substack{\x, \y \in \Tud_n(\sigma) \\ (\sigma,f) \in \widehat B_n}} c_{\x \y} \psi_{\t^{(\lambda,f)} \u} \epsilon_2 \epsilon_4 \ldots \epsilon_{2f-2} \psi_{\x \y} \pmod{R_n^{>f}(\delta)},
\end{equation}
where $c_{\x\y} \neq 0$ only if $\x(n) = \y(n) > 0$ and $\bi_\x = \bi_\u$.
\end{enumerate}
\end{Lemma}

\proof (1). As $\t(n) < 0$, in the standard reduction sequence of $\t$
$$
h(\t) = \t^{(m)} \rightarrow \t^{(m-1)} \rightarrow \ldots \rightarrow \t^{(1)} \rightarrow \t^{(0)} = \t,
$$
we have $\t^{(m-1)}(n) = \t(n)$. Because $\dot\t = \t|_{n-1}$, by the definition, the standard reduction sequence of $\dot\t$ is
$$
h(\dot\t) = \dot\t^{(m-1)} \rightarrow \dot\t^{(m-2)} \rightarrow \ldots \rightarrow \dot\t^{(1)} \rightarrow \dot\t^{(0)} = \dot\t,
$$
where $\dot\t^{(i)} = \t^{(i)}|_{n-1}$ and $\rho(\t^{(i)}, \t^{(i-1)}) = \rho(\dot\t^{(i)}, \dot\t^{(i-1)})$ for $1 \leq i \leq m-1$. Hence we have
\begin{equation} \label{y:h4:eq1}
\theta_k^{(n-1)}(\epsilon_{\dot\t}) = \theta_k^{(n-1)}(\epsilon_{\dot\t^{(m-1)} \rightarrow \dot\t^{(m-2)}} \ldots \epsilon_{\dot\t^{(1)} \rightarrow \dot\t^{(0)}}) = \epsilon_{\t^{(m-1)} \rightarrow \t^{(m-2)}} \ldots \epsilon_{\t^{(1)} \rightarrow \t^{(0)}}.
\end{equation}

By (\ref{rela:7:2}) and (\ref{rela:1}), we have
\begin{equation} \label{y:h4:eq2}
e_{(\lambda,f)} = e(\bi_{(\lambda,f)}) \epsilon_1 \epsilon_3 \ldots \epsilon_{2f-1} e(\bi_{(\lambda,f)}) = e_{(\lambda,f)} \epsilon_2 \epsilon_4 \ldots \epsilon_{2f-2} {\cdot} \epsilon_1 \epsilon_3 \ldots \epsilon_{2f-3} e(\bi_{(\lambda,f)}).
\end{equation}

Because $h(\dot\t) = \t^{(m-1)}|_{n-1}$ and $\t^{(m-1)}(n) = \t(n) = \u(n)$, we have $\u{\cdot}d(h(\dot\t)) = \t^{(m-1)}$. Hence, by~\autoref{y:h2}, we have
$$
\psi_{\t^{(\lambda,f)} \t} = e_{(\lambda,f)} \psi_{d(h(\t))} \epsilon_{h(\t) \rightarrow \t^{(m-1)}} \epsilon_{\t^{(m-1)} \rightarrow \t^{(m-2)}} \ldots \epsilon_{\t^{(1)} \rightarrow \t^{(0)}} = e_{(\lambda,f)} \epsilon_{\t^{(\lambda,f)} \rightarrow \u} \psi_{d(h(\dot\t))} \epsilon_{\t^{(m-1)} \rightarrow \t^{(m-2)}} \ldots \epsilon_{\t^{(1)} \rightarrow \t^{(0)}}.
$$

Finally, because $\epsilon_{\t^{(\lambda,f)} \rightarrow \u} \in \G{2f,n}$, which commutes with $\epsilon_2 \epsilon_4 \ldots \epsilon_{2f-2} \epsilon_1 \epsilon_3 \ldots \epsilon_{2f-3}$, by (\ref{y:h4:eq1}) and (\ref{y:h4:eq2}) we have
\begin{align*}
\psi_{\t^{(\lambda,f)} \t} & = e_{(\lambda,f)} \epsilon_{\t^{(\lambda,f)} \rightarrow \u} \psi_{d(h(\dot\t))} \epsilon_{\t^{(m-1)} \rightarrow \t^{(m-2)}} \ldots \epsilon_{\t^{(1)} \rightarrow \t^{(0)}}\\
& = e_{(\lambda,f)} \epsilon_{\t^{(\lambda,f)} \rightarrow \u} \epsilon_2 \epsilon_4 \ldots \epsilon_{2f-2} e(\bi_\u) \epsilon_1 \epsilon_3 \ldots \epsilon_{2f-3} e(\bi_\u) \psi_{d(h(\dot\t))} \theta_k^{(n-1)}(\epsilon_{\dot\t})\\
& = \psi_{\t^{(\lambda,f)} \u} \epsilon_2 \epsilon_4 \ldots \epsilon_{2f-2} \theta_k^{(n-1)}(e_{(\mu,f-1)}) \theta_k^{(n-1)}(\psi_{d(h(\dot\t))}) \theta_k^{(n-1)}(\epsilon_{\dot\t})\\
& = \psi_{\t^{(\lambda,f)} \u} \epsilon_2 \epsilon_4 \ldots \epsilon_{2f-2} \theta_k^{(n-1)}(\psi_{\t^{(\mu,f-1)} \dot\t}),
\end{align*}
where $k = -\res(\alpha) \in P$, which proves part (1).

(2). Because $(\lambda,f) \in \mathscr S_n$, by the definition of $\mathscr S_n$ we have
$$
\psi_{\t^{(\mu,f-1)} \dot\t} a \equiv \sum_{\dot\v \in \Tud_{n-1}(\mu)} c_{\dot\v} \psi_{\t^{(\mu,f-1)} \dot\v} + \sum_{\substack{\dot\x,\dot\y \in \Tud_{n-1}(\gamma) \\ (\gamma,f) \in \widehat B_{n-1}}} c_{\dot\x \dot\y} \psi_{\dot\x \dot\y} \pmod{R_{n-1}^{>f}(\delta)},
$$
where $c_{\dot\x \dot\y} \neq 0$ only if $\bi_{\dot\x} = \bi_{(\mu,f-1)}$. By~\autoref{I:h3} we have $\theta_k^{(n-1)}(R_{n-1}^{>f}(\delta)) \subseteq R_n^{>f}(\delta)$. Therefore by substituting the above equation into (\ref{y:h4:1}) we have
$$
\psi_{\t^{(\lambda,f)} \t} a \equiv \sum_{\dot\v \in \Tud_{n-1}(\mu)} c_{\dot\v} \psi_{\t^{(\lambda,f)} \u} \epsilon_2 \epsilon_4 \ldots \epsilon_{2f-2} \theta_k^{(n-1)}(\psi_{\t^{(\mu,f-1)} \dot\v}) + \sum_{\substack{\dot\x, \dot\y \in \Tud_{n-1}(\gamma) \\ (\gamma,f) \in \widehat B_{n-1}}} c_{\dot\x \dot\y} \psi_{\t^{(\lambda,f)} \u} \epsilon_2 \epsilon_4 \ldots \epsilon_{2f-2} \theta_k^{(n-1)}(\psi_{\dot\x \dot\y}) \pmod{R_n^{>f}(\delta)},
$$
where $c_{\dot\x \dot\y} \neq 0$ only if $\bi_{\dot\x} = \bi_{(\mu,f-1)}$. By (\ref{y:h4:1}), the first term of the equality equals $\sum_{\v \in \Tud_n(\lambda)} c_{\v} \psi_{\t^{(\lambda,f)} \v}$, where $\v \in \Tud_n(\lambda)$ with $\v|_{n-1} = \dot\v$ and $c_\v = c_{\dot\v}$, which proves part (2).

(3). For $(\gamma,f) \in \widehat B_{n-1}$ and $\dot\x,\dot\y \in \Tud_{n-1}(\gamma)$, if $\res(\beta) \neq k$ for any $\beta \in \mathscr A(\gamma)$, we have $\theta_k^{(n-1)}(\psi_{\dot\x \dot\y}) \in R_n^{>f}(\delta)$ by~\autoref{I:h2}. Then by~\autoref{two:sided:ideal}, we have
\begin{equation} \label{y:h4:eq3}
\psi_{\t^{(\lambda,f)} \u} \epsilon_2 \epsilon_4 \ldots \epsilon_{2f-2} \theta_k^{(n-1)}(\psi_{\dot\x \dot\y}) \in R_n^{>f}(\delta).
\end{equation}

If there exists $\beta \in \mathscr A(\gamma$ with $\res(\beta) = k$, set $\sigma = \gamma \cup\{\beta\}$. By~\autoref{I:h1}, we have
\begin{equation} \label{y:h4:eq4}
\theta_k^{(n-1)}(\psi_{\dot\x \dot\y}) = \psi_{\x\y},
\end{equation}
where $\x,\y \in \Tud_n(\sigma)$ such that $\x|_{n-1} = \dot\x$ and $\y|_{n-1} = \dot\x$. This implies that $\x(n) = \y(n) = \beta > 0$ and $\bi_\x = \bi_{\dot\x}\vee k = \bi_{(\mu,f-1)}\vee k = \bi_\u$. Hence we prove part (3) by substituting (\ref{y:h4:eq3}) and (\ref{y:h4:eq4}) into (\ref{y:h4:2}). \endproof

The following Lemmas simplify the expression of $\psi_{\t^{(\lambda,f)} \u} \epsilon_2 \epsilon_4 \ldots \epsilon_{2f-2} \psi_{\x\y}$ on the RHS of (\ref{y:h4:3}).

\begin{Lemma} \label{y:first:h1}
Suppose $(\sigma,f) \in \widehat B_n$ and $\x,\y \in \Tud_n(\sigma)$ with $\bi_\x = \bi_\u$ and $\x(n) = \y(n) > 0$. Then we have
$$
\psi_{\t^{(\lambda,f)} \u} \epsilon_2 \epsilon_4 \ldots \epsilon_{2f-2} \psi_{\x\y} \equiv c_\sigma \psi_{\t^{(\lambda,f)} \u} \left(\psi_{b-1} \ldots \psi_{a+1} \epsilon_a \ldots \epsilon_{2f} \epsilon_{2f-1} \right) \psi_\x^* e(\bi_{(\sigma,f)}) \psi_\y \epsilon_\y \pmod{R_n^{>f}(\delta)},
$$
where $c_\sigma \in R$.
\end{Lemma}

\proof Recall $\u \in \Tud_n(\lambda)$ with head $f-1$. Hence if we write $\bi_\x = (i_1,\ldots,i_n)$, $\bi_\u = \bi_\x$ implies $i_1 = i_3 = \ldots = i_{2f-3} = -i_2 = -i_4 = \ldots = -i_{2f-2} = \frac{\delta-1}{2}$. By~\autoref{idem:psi:2}, we have $\x_{2f-1} = (1)$. When $\frac{\delta-1}{2} \neq \pm\frac{1}{2}$, one can see that for $1 \leq \l \leq f-1$, we have $\x_{2\l-1} = (1)$ and $\x_{2\l} = \emptyset$; and when $\frac{\delta-1}{2} = \frac{1}{2}$, for $1 \leq \l \leq f-1$ we have $\x_{2\l-1} = (1)$ and $\x_{2\l} = \emptyset$ or $(1,1)$; and when $\frac{\delta-1}{2} = -\frac{1}{2}$, for $1 \leq \l \leq f-1$ we have $\x_{2\l-1} = (1)$ and $\x_{2\l} = \emptyset$ or $(2)$.

Because $\Shape(\x|_{2f-1}) = ((1),f-1)$ and $\x \in \Tud_n(\sigma)$ where $(\sigma,f) \in \widehat B_n$, by the construction of $\epsilon_\x$ we have
\begin{equation} \label{y:first:h1:eq1}
\epsilon_\x = e(\bi_{h(\x)}) g_1 g_2 \ldots g_{f-1} \epsilon_{2f}\epsilon_{2f+1} \ldots \epsilon_a \psi_{a+1} \ldots \psi_{b-1} e(\bi_\x),
\end{equation}
where $2f \leq a < b \leq n$, and $g_\l = 1$ if $\x_{2\l} = \emptyset$ and $g_\l = \epsilon_{2\l}$ if $\x_{2\l} \neq \emptyset$ for $1 \leq \l \leq f-1$. By (\ref{rela:5:5}) and (\ref{rela:7:2}), we have
\begin{equation} \label{y:first:h1:eq2}
\epsilon_{2\l-1} \epsilon_{2\l} g_\l \epsilon_{2\l-1} \ldots \epsilon_3 \epsilon_1 =
f_\l \epsilon_1 \epsilon_3 \ldots \epsilon_{2\l-1},
\end{equation}
where $f_\l \in R$ if $g_\l = 1$ and $f_\l$ is a polynomial of $y_1,y_2,\ldots,y_{2\l-1}, y_{2\l+2}$'s if $g_\l = \epsilon_{2\l}$. As $\psi_\u, \psi_\x^* \in \G{2f,n}$, they commute with $\epsilon_1 \epsilon_3 \ldots \epsilon_{2f-1}$; and as $\epsilon_\u \in \G{2f-1,n}$ which is implied by $head(\u) = f-1$, it commutes with $\epsilon_1 \epsilon_3 \ldots \epsilon_{2f-3}$. Therefore, by (\ref{y:first:h1:eq1}) and (\ref{y:first:h1:eq2}) we have
\begin{eqnarray*}
&& \psi_{\t^{(\lambda,f)} \u} \epsilon_2 \epsilon_4 \ldots \epsilon_{2f-2} \psi_{\x\y}\\
& = & e(\bi_{(\lambda,f)}) \epsilon_{2f-1} \psi_\u \epsilon_\u \epsilon_1 \epsilon_3 \ldots \epsilon_{2f-3} \epsilon_\x^* \epsilon_{2f-1} \epsilon_{2f-3} \ldots \epsilon_3 \epsilon_1 {\cdot} \psi_\x^*  e(\bi_{(\sigma,f)}) \psi_\y \epsilon_\y\\
& = & e(\bi_{(\lambda,f)}) \epsilon_{2f-1} \psi_\u \epsilon_\u \epsilon_1 \epsilon_3 \ldots \epsilon_{2f-3} \psi_{b-1} \ldots \psi_{a+1} \epsilon_a \ldots \epsilon_{2f} \epsilon_{2f-1} g_{f-1} \ldots g_2 g_1 \epsilon_{2f-3}  \ldots \epsilon_3 \epsilon_1 {\cdot} \psi_\x^* e(\bi_{(\sigma,f)}) \psi_\y \epsilon_\y\\
& = & e(\bi_{(\lambda,f)}) \epsilon_{2f-1} \psi_\u \epsilon_\u \psi_{b-1} \ldots \psi_{a+1} \epsilon_a \ldots \epsilon_{2f} \epsilon_{2f-1} {\cdot} \left( \epsilon_1 \epsilon_3 \ldots \epsilon_{2f-3} g_{f-1} \ldots g_2 g_1 \epsilon_{2f-3} \ldots \epsilon_3 \epsilon_1 \right) {\cdot} \psi_\x^*  e(\bi_{(\sigma,f)}) \psi_\y \epsilon_\y.
\end{eqnarray*}

By applying (\ref{y:first:h1:eq2}) on the term in the bracket recursively, from $\l = 1$ to $f-1$, we have
\begin{align*}
\psi_{\t^{(\lambda,f)} \u} \epsilon_2 \epsilon_4 \ldots \epsilon_{2f-2} \psi_{\x\y} & = \left(\prod_{\l=1}^{f-1} f_\l \right) e(\bi_{(\lambda,f)}) \epsilon_{2f-1} \psi_\u \epsilon_\u \psi_{b-1} \ldots \psi_{a+1} \epsilon_a \ldots \epsilon_{2f} \epsilon_{2f-1} {\cdot} \epsilon_1 \epsilon_3 \ldots \epsilon_{2f-3} {\cdot} \psi_\x^* e(\bi_{(\sigma,f)}) \psi_\y \epsilon_\y\\
& = \left(\prod_{\l=1}^{f-1} f_\l \right) \left( e(\bi_{(\lambda,f)}) \epsilon_1 \epsilon_3 \ldots \epsilon_{2f-1} \psi_\u \epsilon_\u \right) \left(\psi_{b-1} \ldots \psi_{a+1} \epsilon_a \ldots \epsilon_{2f} \epsilon_{2f-1} \right) \psi_\x^* e(\bi_{(\sigma,f)}) \psi_\y \epsilon_\y\\
& = \left(\prod_{\l=1}^{f-1} f_\l \right) \psi_{\t^{(\lambda,f)} \u} \left(\psi_{b-1} \ldots \psi_{a+1} \epsilon_a \ldots \epsilon_{2f} \epsilon_{2f-1} \right) \psi_\x^* e(\bi_{(\sigma,f)}) \psi_\y \epsilon_\y.
\end{align*}

Finally, by~\autoref{base:y:1}, we have $y_s \psi_{\t^{(\lambda,f)} \u} \in R_n^{>f}(\delta)$ for any $1 \leq s \leq n$. As $f_\l \in R$ if $g_\l = 1$ and $f_\l$ is a polynomial of $y_1,y_2,\ldots,y_{2\l-1}, y_{2\l+2}$'s if $g_\l = \epsilon_{2\l}$, we have
$$
\left(\prod_{\l=1}^{f-1} f_\l \right) \psi_{\t^{(\lambda,f)} \u} = c_\sigma \psi_{\t^{(\lambda,f)} \u}
$$
for some $c_\sigma \in R$, which completes the proof of the Lemma. \endproof


The following two Lemmas are technical results which will be used for computational purposes.

\begin{Lemma} \label{commute:1}
Suppose $(\lambda,f) \in \mathscr S_n$ and $\t \in \Tud_n(\lambda)$ with head $f-1$. Let $\s = h(\t) \rightarrow \t$ and $\rho(\s,\t) = (a,b)$ for some $2f \leq a < b \leq n$. Then for any $k$ with $2f-1 \leq k \leq n-1$, we have
$$
\psi_{\t^{(\lambda,f)} \t} \epsilon_k \epsilon_{k-1} \ldots \epsilon_{2f-1} \equiv c_w{\cdot} \psi_{\t^{(\lambda,f)} \s} \psi_w f_w \pmod{R_n^{>f}(\delta)},
$$
where $c_w \in R$, $w \in \Sym_{2f,n}$ and $f_w$ is a polynomial of $y_1, \ldots, y_{k+1}$.
\end{Lemma}

\proof It suffices to prove that for any $\bk \in P^n$, we have
$$
\psi_{\t^{(\lambda,f)} \t} \epsilon_k \epsilon_{k-1} \ldots \epsilon_{2f-1} e(\bk) \equiv c_\w{\cdot}\psi_{\t^{(\lambda,f)} \s} \psi_w f_w e(\bk) \pmod{R_n^{>f}(\delta)},
$$
where $c_w \in R$, $w \in \Sym_{a,b+1}$ and $f_w$ is a polynomial of $y_1, \ldots, y_{k+1}$. In the proof we omit the $e(\bk)$, but readers have to remember that there is always some $e(\bk)$ on the right of each element.

Because $\rho(\s,\t) = (a,b)$, we can write $\psi_{\t^{(\lambda,f)} \t} = \psi_{\t^{(\lambda,f)} \s} \epsilon_{2f} \epsilon_{2f+1} \ldots \epsilon_a \psi_{a+1} \ldots \psi_{b-1} e(\bi_\t)$. Because $head(\s) = f$, we have $d(\s) \in \Sym_{2f+1,n}$, which implies $\psi_\s$ commutes with $\epsilon_{2f}$. Therefore, by (\ref{rela:7:2}), we have
\begin{equation} \label{commute:1:eq1}
\psi_{\t^{(\lambda,f)} \s} \epsilon_{2f} \epsilon_{2f-1} = e(\bi_{(\lambda,f)}) \epsilon_1 \ldots \epsilon_{2f-1} e(\bi_{(\lambda,f)}) \psi_\s \epsilon_{2f} \epsilon_{2f-1} = e(\bi_{(\lambda,f)}) \epsilon_1 \ldots \epsilon_{2f-1} \epsilon_{2f} \epsilon_{2f-1} e(\bi_{(\lambda,f)}) \psi_\s = \psi_{\t^{(\lambda,f)} \s}.
\end{equation}

We consider the following cases for different values of $k$, $a$ and $b$.

\textbf{Case 1:} $k < a-1$.

By (\ref{rela:13}), we have
$$
\psi_{\t^{(\lambda,f)} \t} \epsilon_k = \psi_{\t^{(\lambda,f)} \s} \epsilon_{2f} \epsilon_{2f+1} \ldots \epsilon_a \psi_{a+1} \ldots \psi_{b-1} e(\bi_\t) \epsilon_k = \psi_{\t^{(\lambda,f)} \s}\epsilon_{k+2} {\cdot}\epsilon_{2f} \epsilon_{2f+1} \ldots \epsilon_a \psi_{a+1} \ldots \psi_{b-1}.
$$

By~\autoref{commute:1:h1} we have $\psi_{\t^{(\lambda,f)} \s}\epsilon_{k+2} \in R_n^{>f}(\delta)$. Hence by~\autoref{two:sided:ideal}, we have $\psi_{\t^{(\lambda,f)} \t} \epsilon_k \epsilon_{k-1} \ldots \epsilon_{2f-1} \in R_n^{>f}(\delta)$.

\textbf{Case 2:} $k = a-1$.

By (\ref{rela:7:2}) we have
\begin{align*}
\psi_{\t^{(\lambda,f)} \t} \epsilon_{a-1} \ldots \epsilon_{2f-1} & = \psi_{\t^{(\lambda,f)} \s} \epsilon_{2f} \epsilon_{2f+1} \ldots \epsilon_a \psi_{a+1} \ldots \psi_{b-1} e(\bi_\t) \epsilon_{a-1} \ldots \epsilon_{2f-1}\\
&= \psi_{\t^{(\lambda,f)} \s} \epsilon_{2f} \epsilon_{2f-1} \psi_{a+1} \ldots \psi_{b-1}.
\end{align*}

Hence, by (\ref{commute:1:eq1}) we have $\psi_{\t^{(\lambda,f)} \t} \epsilon_{a-1} \ldots \epsilon_{2f-1} = \psi_{\t^{(\lambda,f)} \s} \psi_{a+1} \ldots \psi_{b-1}$, where the Lemma holds by setting $c_w = 1$, $w = s_{a+1} s_{a+2} \ldots s_{b-1}$ and $f_w = 1$.

\textbf{Case 3:} $k = a$.

If $k = a < b-1$, by (\ref{rela:15}) and (\ref{rela:7:2}) we have
\begin{align*}
\psi_{\t^{(\lambda,f)} \t} \epsilon_a \epsilon_{a-1} \ldots \epsilon_{2f-1} & = \psi_{\t^{(\lambda,f)} \s} \epsilon_{2f} \epsilon_{2f+1} \ldots \epsilon_a \psi_{a+1} \ldots \psi_{b-1} e(\bi_\t) \epsilon_a \epsilon_{a-1} \ldots \epsilon_{2f-1}\\
& = \pm \psi_{\t^{(\lambda,f)} \s} \epsilon_{2f} \epsilon_{2f+1} \ldots \epsilon_{a-1} \epsilon_a \psi_{a+2} \ldots \psi_{b-1} e(\bi_\t) \epsilon_{a-1} \ldots \epsilon_{2f-1}\\
& = \pm \psi_{\t^{(\lambda,f)} \s} \epsilon_{2f}\epsilon_{2f-1} \psi_{a+2} \ldots \psi_{b-1}.
\end{align*}

Hence by (\ref{commute:1:eq1}) we have $\psi_{\t^{(\lambda,f)} \t} \epsilon_a \epsilon_{a-1} \ldots \epsilon_{2f-1} = \pm \psi_{\t^{(\lambda,f)} \s} \psi_{a+2} \ldots \psi_{b-1}$, where the Lemma holds by setting $c_w = \pm 1$, $w = s_{a+2} \ldots s_{b-1}$ and $f_w = 1$.

If $k = a = b-1$, by (\ref{rela:5:5}), (\ref{rela:16}) and (\ref{rela:7:2}), we have
\begin{align*}
\psi_{\t^{(\lambda,f)} \t} \epsilon_k \epsilon_{k-1} \ldots \epsilon_{2f-1} & = \psi_{\t^{(\lambda,f)} \s} \epsilon_{2f} \epsilon_{2f+1} \ldots \epsilon_k e(\bi_\t) \epsilon_k \epsilon_{k-1} \ldots \epsilon_{2f-1}\\
& = \psi_{\t^{(\lambda,f)} \s} \epsilon_{2f} \epsilon_{2f+1} \ldots \epsilon_{k-1} f_w(y_1, \ldots, y_{k-1}) \epsilon_k \epsilon_{k-1} \ldots \epsilon_{2f-1}\\
& = \psi_{\t^{(\lambda,f)} \s} \epsilon_{2f} \epsilon_{2f+1} \ldots \epsilon_{k-1} \epsilon_k \epsilon_{k-1} \ldots \epsilon_{2f-1} f_w(y_1, \ldots, y_{2f-2}, y_{2f+1}, \ldots, y_{k+1})\\
& = \psi_{\t^{(\lambda,f)} \s} \epsilon_{2f} \epsilon_{2f-1} f_w(y_1, \ldots, y_{2f-2}, y_{2f+1}, \ldots, y_{k+1}),
\end{align*}
where $f_w$ is a polynomial of $y_1, \ldots, y_{k+1}$. Hence by (\ref{commute:1:eq1}) we have
$$
\psi_{\t^{(\lambda,f)} \t} \epsilon_k \epsilon_{k-1} \ldots \epsilon_{2f-1} = \psi_{\t^{(\lambda,f)} \s} f_w(y_1, \ldots, y_{2f-2}, y_{2f+1}, \ldots, y_{k+1}),
$$
where the Lemma holds by setting $c_w = 1$, $w = 1$.

\textbf{Case 4:} $k = a+1$.

When $k = a+1 = b$, by (\ref{rela:7:2}) we have
\begin{align*}
\psi_{\t^{(\lambda,f)} \s} \epsilon_{2f} \epsilon_{2f+1} \ldots \epsilon_a e(\bi_\t) \epsilon_k \epsilon_{k-1} \ldots \epsilon_{2f-1}
& = \psi_{\t^{(\lambda,f)} \s} \epsilon_{2f} \epsilon_{2f+1} \ldots \epsilon_a \epsilon_{a+1} \epsilon_a \ldots \epsilon_{2f-1} = \psi_{\t^{(\lambda,f)} \s} \epsilon_{2f} \epsilon_{2f-1},
\end{align*}
and the Lemma holds by (\ref{commute:1:eq1}).

When $k = a + 1 = b-1$, by (\ref{rela:6:1}) we have
$$
\psi_{\t^{(\lambda,f)} \t} \epsilon_k \epsilon_{k-1} \ldots \epsilon_{2f-1} = \psi_{\t^{(\lambda,f)} \s} \epsilon_{2f} \epsilon_{2f+1} \ldots \epsilon_a \psi_{a+1} e(\bi_\t) \epsilon_k \epsilon_{k-1} \ldots \epsilon_{2f-1} = \pm \psi_{\t^{(\lambda,f)} \s} \epsilon_{2f} \epsilon_{2f+1} \ldots \epsilon_a \epsilon_{a+1} \epsilon_a \ldots \epsilon_{2f-1},
$$
and following the same argument as when $k = b$, the Lemma follows.

When $k = a + 1 < b-1$, we have $a < b-2$. Hence by (\ref{rela:7:1}), (\ref{rela:4}), (\ref{rela:7:2}) and (\ref{rela:16}), we have
\begin{eqnarray}
\psi_{\t^{(\lambda,f)} \t} \epsilon_k \epsilon_{k-1} \ldots \epsilon_{2f-1}  & = & \psi_{\t^{(\lambda,f)} \s} \epsilon_{2f} \epsilon_{2f+1} \ldots \epsilon_a \psi_{a+1} \ldots \psi_{b-1} e(\bi_\t) \epsilon_k \epsilon_{k-1} \ldots \epsilon_{2f-1} \notag \\
& = & \psi_{\t^{(\lambda,f)} \s} \epsilon_{2f} \epsilon_{2f+1} \ldots \epsilon_a \psi_{a+1} \psi_{a+2} \epsilon_{a+1} \epsilon_a \ldots \epsilon_{2f-1} \psi_{a+3} \ldots \psi_{b-1} \notag \\
& = & \psi_{\t^{(\lambda,f)} \s} \epsilon_{2f} \epsilon_{2f+1} \ldots \epsilon_a \psi_{a+1}^2 \epsilon_{a+2} \epsilon_{a+1} \epsilon_a \ldots \epsilon_{2f-1} \psi_{a+3} \ldots \psi_{b-1} \notag \\
& = & \psi_{\t^{(\lambda,f)} \s} \epsilon_{2f} \epsilon_{2f+1} \ldots \epsilon_a f(y_{a+1}, y_{a+2}) \epsilon_{a+2} \epsilon_{a+1} \epsilon_a \ldots \epsilon_{2f-1} \psi_{a+3} \ldots \psi_{b-1}, \label{commute:1:eq3}
\end{eqnarray}
where $f(y_{a+1}, y_{a+2}) = 0$, $1$ or $\pm (y_{a+2} - y_{a+1})$. By~\autoref{commute:1:h1} and~\autoref{two:sided:ideal}, for any $f(y_{a+1}, y_{a+2})$ we have that (\ref{commute:1:eq3}) is an element of $R_n^{>f}(\delta)$.

\textbf{Case 5:} $k > a+1$.

When $k > b$, $\epsilon_k$ commutes with $\epsilon_{\s\rightarrow\t}$. Hence by~\autoref{commute:1:h1} we have $\psi_{\t^{(\lambda,f)} \t} \epsilon_k \epsilon_{k-1} \ldots \epsilon_{2f-1} \in R_n^{>f}(\delta)$ and the Lemma follows.

When $k = b$, by (\ref{rela:9}) and (\ref{rela:7:2}) we have
\begin{align*}
\psi_{\t^{(\lambda,f)} \t} \epsilon_k \epsilon_{k-1} \ldots \epsilon_{2f-1} & =\psi_{\t^{(\lambda,f)} \s} \epsilon_{2f} \epsilon_{2f+1} \ldots \epsilon_a \psi_{a+1} \ldots \psi_{b-1} e(\bi_\t) \epsilon_k \epsilon_{k-1} \ldots \epsilon_{2f-1}\\
& = \psi_{\t^{(\lambda,f)} \s} \epsilon_{2f} \epsilon_{2f+1} \ldots \epsilon_a \psi_{a+1} \ldots \psi_{b-1} \epsilon_b \epsilon_{b-1} \ldots \epsilon_{2f-1}\\
& = \psi_{\t^{(\lambda,f)} \s} \epsilon_{2f} \epsilon_{2f+1} \ldots \epsilon_a \psi_b \psi_{b-1} \ldots \psi_{a+2} \epsilon_{a+1} \epsilon_a \ldots \epsilon_{2f-1}\\
& = \psi_{\t^{(\lambda,f)} \s} \epsilon_{2f} \epsilon_{2f-1} \psi_b \psi_{b-1} \ldots \psi_{a+2}.
\end{align*}

Hence by (\ref{commute:1:eq1}) we have $\psi_{\t^{(\lambda,f)} \t} \epsilon_k \epsilon_{k-1} \ldots \epsilon_{2f-1} = \psi_{\t^{(\lambda,f)} \s} \psi_b \psi_{b-1} \ldots \psi_{a+2}$, and the Lemma follows by setting $c_w = 1$, $w = s_b s_{b-1} \ldots s_{a+2}$ and $f_w = 1$.

When $k = b-1$, by (\ref{rela:6:1}) we have
\begin{align*}
\psi_{\t^{(\lambda,f)} \s} \epsilon_{2f} \epsilon_{2f+1} \ldots \epsilon_a \psi_{a+1} \ldots \psi_{b-1} e(\bi_\t) \epsilon_k \epsilon_{k-1} \ldots \epsilon_{2f-1}
& = \psi_{\t^{(\lambda,f)} \s} \epsilon_{2f} \epsilon_{2f+1} \ldots \epsilon_a \psi_{a+1} \ldots \psi_{b-1} \epsilon_{b-1} \epsilon_{b-2} \ldots \epsilon_{2f-1}\\
& = \pm \psi_{\t^{(\lambda,f)} \s} \epsilon_{2f} \epsilon_{2f+1} \ldots \epsilon_a \psi_{a+1} \ldots \psi_{b-2} \epsilon_{b-1} \epsilon_{b-2} \ldots \epsilon_{2f-1},
\end{align*}
and following the same argument as when $k = b$, the Lemma holds.

When $k < b-1$, following the similar argument as in Case 4 when $k < b-1$, the Lemma follows. \endproof

\begin{Lemma} \label{commute:3}
Suppose $(\lambda,f) \in \widehat B_n$ with $f \geq 1$ and $\alpha \in \mathscr A(\lambda)$. Define $\mu = \lambda \cup \{\alpha\}$ and $\u \in \Tud_n(\lambda)$ such that $\u|_{n-1} = \t^{(\mu,f-1)}$ and $\u(n) = -\alpha$. For any $2f \leq k < \l \leq n-1$, we have
$$
\psi_{\t^{(\lambda,f)} \u} \psi_{\l-1} \psi_{\l-2} \ldots \psi_{k+1} \epsilon_k \ldots \epsilon_{2f} \epsilon_{2f-1} \equiv \sum_{w \in \Sym_{2f,n}} c_w{\cdot} \psi_{\t^{(\lambda,f)} \t^{(\lambda,f)}} \psi_w f_w \pmod{R_n^{>f}(\delta)},
$$
where $c_w \in R$, $w \in \Sym_{2f,n}$ and $f_w$'s are polynomials of $y_1, \ldots, y_n$.
\end{Lemma}

\proof Because $\u \in \Tud_n(\lambda)$ and $\u|_{n-1} = \t^{(\mu,f-1)}$, we have $\t^{(\lambda,f)} \rightarrow \u$ and $\rho(\t^{(\lambda,f)},\u) = (a,n)$. Hence we can write
$$
\psi_{\t^{(\lambda,f)} \u} = e(\bi_{(\lambda,f)})\epsilon_1 \epsilon_3 \ldots \epsilon_{2f-1} e(\bi_{(\lambda,f)}) \epsilon_{2f} \epsilon_{2f+1} \ldots \epsilon_a \psi_{a+1} \ldots \psi_{n-1} e(\bi_\u).
$$

We consider different values of $\l$ and $a$.

\textbf{Case 1:} $\l < a$.

By (\ref{rela:12}) and (\ref{rela:13}), we have
\begin{eqnarray*}
&& \psi_{\t^{(\lambda,f)} \u} \psi_{\l-1} \psi_{\l-2} \ldots \psi_{k+1} \epsilon_k \ldots \epsilon_{2f} \epsilon_{2f-1}\\
& = & e(\bi_{(\lambda,f)})\epsilon_1 \epsilon_3 \ldots \epsilon_{2f-1} \epsilon_{2f} \epsilon_{2f+1} \ldots \epsilon_a \psi_{a+1} \ldots \psi_{n-1} \psi_{\l-1} \psi_{\l-2} \ldots \psi_{k+1} \epsilon_k \ldots \epsilon_{2f} \epsilon_{2f-1}\\
& = & e(\bi_{(\lambda,f)})\psi_{\l+1} \psi_{\l} \ldots \psi_{k+3} \epsilon_1 \epsilon_3 \ldots \epsilon_{2f-1} \epsilon_{2f} \epsilon_{2f+1} \ldots \epsilon_a \psi_{a+1} \ldots \psi_{n-1} \epsilon_k \ldots \epsilon_{2f} \epsilon_{2f-1}\\
& = & e(\bi_{(\lambda,f)})\psi_{\l+1} \psi_{\l} \ldots \psi_{k+3} \epsilon_1 \epsilon_3 \ldots \epsilon_{2f-1} \epsilon_{k+2} \epsilon_{2f} \epsilon_{2f+1} \ldots \epsilon_a \psi_{a+1} \ldots \psi_{n-1} \epsilon_{k-1} \ldots \epsilon_{2f} \epsilon_{2f-1}.
\end{eqnarray*}

By~\autoref{I:coro:1}, we have $\epsilon_1 \epsilon_3 \ldots \epsilon_{2f-1} \epsilon_{k+2} \in R_n^{>f}(\delta)$. Hence $\psi_{\t^{(\lambda,f)} \u} \psi_{\l-1} \psi_{\l-2} \ldots \psi_{k+1} \epsilon_k \ldots \epsilon_{2f} \epsilon_{2f-1} \in R_n^{>f}(\delta)$ by~\autoref{two:sided:ideal}.

\textbf{Case 2:} $\l = a$.

By (\ref{rela:9}) and (\ref{rela:7:2}) we have
\begin{eqnarray*}
&& \psi_{\t^{(\lambda,f)} \u} \psi_{\l-1} \psi_{\l-2} \ldots \psi_{k+1} \epsilon_k \ldots \epsilon_{2f} \epsilon_{2f-1}\\
& = & e(\bi_{(\lambda,f)})\epsilon_1 \epsilon_3 \ldots \epsilon_{2f-1} \epsilon_{2f} \epsilon_{2f+1} \ldots \epsilon_a \psi_{a+1} \ldots \psi_{n-1} \psi_{\l-1} \psi_{\l-2} \ldots \psi_{k+1} \epsilon_k \ldots \epsilon_{2f} \epsilon_{2f-1}\\
& = & e(\bi_{(\lambda,f)}) \epsilon_1 \epsilon_3 \ldots \epsilon_{2f-1} \epsilon_{2f} \epsilon_{2f+1} \ldots \epsilon_{k+1} \psi_{k+2} \ldots \psi_{n-1} \epsilon_k \ldots \epsilon_{2f} \epsilon_{2f-1}\\
& = & e(\bi_{(\lambda,f)}) \epsilon_1 \epsilon_3 \ldots \epsilon_{2f-1} \psi_{k+2} \ldots \psi_{n-1} = \psi_{\t^{(\lambda,f)} \t^{(\lambda,f)}} \psi_w,
\end{eqnarray*}
where $w = s_{k+2} s_{k+3} \ldots s_{n-1}$.

\textbf{Case 3:} $a < \l < n-1$.

By the assumption, we have $\l \leq n-2$ and $a \leq n-3$. Hence we have $\u(n-1) > 0$, $\u(n) < 0$ and $\u(n-1) + \u(n) \neq 0$. By~\autoref{deg:tab:4}, we have $\u{\cdot}s_{n-1} \in \Tud_n(\lambda)$. Set $\v = \u{\cdot}s_{n-1}$ and we have $\psi_{\t^{(\lambda,f)} \u} = \psi_{\t^{(\lambda,f)} \v} \psi_{n-1}$. Therefore, we have
\begin{equation} \label{commute:3:eq1}
\psi_{\t^{(\lambda,f)} \u} \psi_{\l-1} \psi_{\l-2} \ldots \psi_{k+1} \epsilon_k \ldots \epsilon_{2f} \epsilon_{2f-1} = \psi_{\t^{(\lambda,f)} \v} \psi_{\l-1} \psi_{\l-2} \ldots \psi_{k+1} \epsilon_k \ldots \epsilon_{2f} \epsilon_{2f-1} \psi_{n-1}.
\end{equation}

By the construction, $\v$ has head $f-1$ and $\v(n) = \u(n-1) > 0$. By~\autoref{remove:tail:1} and~\autoref{remove:head}, we have
\begin{equation} \label{commute:3:eq2}
\psi_{\t^{(\lambda,f)} \v} \psi_{\l-1} \psi_{\l-2} \ldots \psi_{k+1} \equiv \sum_{\y \in \Tud_n(\lambda)} c_\y \psi_{\t^{(\lambda,f)} \y} \pmod{R_n^{>f}(\delta)},
\end{equation}
where $c_\y \neq 0$ only if $head(\y) \geq f-1$ and $\y(n) = \v(n) > 0$. Therefore, we have $h(\y) \rightarrow \y$ and $\rho(h(\y),\y) = (s,m)$ where $s < m \leq n-1$. Because $k < \l \leq n-2$, we have $k+1 \leq n-2$. Hence by~\autoref{commute:1}, we have
$$
\psi_{\t^{(\lambda,f)} \y} \epsilon_k \epsilon_{k-1} \ldots \epsilon_{2f-1} \psi_{n-1} \equiv c_w \psi_{\t^{(\lambda,f)} \t^{(\lambda,f)}} \psi_w f_w \psi_{n-1} \pmod{R_n^{>f}(\delta)},
$$
where $w \in \Sym_{2f,n}$ and $f_w$ is a polynomial of $y_1, \ldots, y_{n-2}$. Hence $\psi_{n-1}$ commutes with $f_w$. By~\autoref{base:psi:1} and~\autoref{two:sided:ideal}, we have
$$
\psi_{\t^{(\lambda,f)} \t^{(\lambda,f)}} \psi_w \psi_{n-1}
\begin{cases}
= \psi_{\t^{(\lambda,f)} \t^{(\lambda,f)}} \psi_{w'}, & \text{if $w' = w{\cdot}s_{n-1}$ is semi-reduced corresponding to $\t^{(\lambda,f)}$,}\\
\in R_n^{>f}(\delta), & \text{if $w' = w{\cdot}s_{n-1}$ is not semi-reduced corresponding to $\t^{(\lambda,f)}$.}
\end{cases}
$$

Therefore, by~\autoref{two:sided:ideal}, we have
$$
\psi_{\t^{(\lambda,f)} \y} \epsilon_k \epsilon_{k-1} \ldots \epsilon_{2f-1} \psi_{n-1} \equiv c_w \psi_{\t^{(\lambda,f)} \t^{(\lambda,f)}} \psi_w \psi_{n-1} f_w \equiv c_{w'} \psi_{\t^{(\lambda,f)} \t^{(\lambda,f)}} \psi_{w'} f_w \pmod{R_n^{>f}(\delta)},
$$
where $c_{w'} \in R$, $w' \in \Sym_{2f,n}$ and $f_{w'}$ is a polynomial of $y_1, \ldots, y_{n-2}$. The Lemma holds by substituting the above equality and (\ref{commute:3:eq2}) into (\ref{commute:3:eq1}).

\textbf{Case 4:} $\l = n-1$.

If $a = n-1$, then we have $\l = a$, which is proved in Case 2. Assume $a < n-1$. Following the same argument as in Case 3, we have $\u(n-1) > 0$ and $\u{\cdot}s_{n-1} \in \Tud_n(\lambda)$. Denote $\v = \u{\cdot}s_{n-1}$. By (\ref{rela:10}) we have
\begin{align*}
\psi_{\t^{(\lambda,f)} \u} \psi_{\l-1} \psi_{\l-2} \ldots \psi_{k+1} \epsilon_k \ldots \epsilon_{2f} \epsilon_{2f-1}
& = \psi_{\t^{(\lambda,f)} \v} \psi_{n-1} \psi_{n-2} \ldots \psi_{k+1} \epsilon_k \ldots \epsilon_{2f} \epsilon_{2f-1}\\
& = \psi_{\t^{(\lambda,f)} \v} \psi_k \psi_{k+1} \ldots \psi_{n-2} \epsilon_{n-1} \ldots \epsilon_{2f} \epsilon_{2f-1}.
\end{align*}

Because $head(\v) = f-1$, by~\autoref{remove:tail:1} and~\autoref{remove:head}, we have
$$
\psi_{\t^{(\lambda,f)} \v} \psi_k \psi_{k+1} \ldots \psi_{n-2} \epsilon_{n-1} \ldots \epsilon_{2f} \epsilon_{2f-1} \equiv \sum_{\y \in \Tud_n(\lambda)} c_\y \psi_{\t^{(\lambda,f)} \y} \epsilon_{n-1} \ldots \epsilon_{2f} \epsilon_{2f-1} \pmod{R_n^{>f}(\delta)},
$$
where $c_\y \neq 0$ only if $head(\y) \geq f-1$. Therefore, by~\autoref{commute:1}, we have
$$
\psi_{\t^{(\lambda,f)} \y} \epsilon_{n-1} \ldots \epsilon_{2f} \epsilon_{2f-1} \equiv c_w \psi_{\t^{(\lambda,f)} \t^{(\lambda,f)}} \psi_w f_w \pmod{R_n^{>f}(\delta)},
$$
where $c_w \in R$, $w \in \Sym_{2f,n}$ and $f_w$ is a polynomial of $y_1, \ldots, y_n$. By combining the above two equalities, the Lemma holds. \endproof

Now we are ready to prove (\ref{Induction:1}) when $a \in \G{n-1}$.

\begin{Lemma} \label{y:first}
Suppose $(\lambda,f) \in \mathscr S_n$ and $\t \in \Tud_n(\lambda)$. For any $a \in \G{n-1}$, the equality (\ref{Induction:1}) holds.
\end{Lemma}

\proof When $\t(n) > 0$, the Lemma follows by~\autoref{y:first:1}. Hence we only have to consider the case when $\t(n) < 0$.

Suppose $\t(n) < 0$. By~\autoref{y:h4}, we have
\begin{equation} \label{y:first:eq1}
\psi_{\t^{(\lambda,f)} \t} a \equiv \sum_{\v \in \Tud_n(\lambda)} c_\v \psi_{\t^{(\lambda,f)} \v} + \sum_{\substack{\x, \y \in \Tud_n(\sigma) \\ (\sigma,f) \in \widehat B_n}} c_{\x \y} \psi_{\t^{(\lambda,f)} \u} \epsilon_2 \epsilon_4 \ldots \epsilon_{2f-2} \psi_{\x \y} \pmod{R_n^{>f}(\delta)},
\end{equation}
where $c_{\x\y} \neq 0$ only if $\x(n) = \y(n) > 0$ and $\bi_\x = \bi_\u$. For the second term of (\ref{y:first:eq1}), by~\autoref{y:first:h1} and~\autoref{commute:3}, we have
\begin{eqnarray}
\psi_{\t^{(\lambda,f)} \u} \epsilon_2 \epsilon_4 \ldots \epsilon_{2f-2} \psi_{\x\y}
& \equiv & c_\sigma \psi_{\t^{(\lambda,f)} \u} \left(\psi_{b-1} \ldots \psi_{a+1} \epsilon_a \ldots \epsilon_{2f} \epsilon_{2f-1} \right) \psi_\x^* e(\bi_{(\sigma,f)}) \psi_\y \epsilon_\y \notag \\
& \equiv & \sum_{w \in \Sym_{2f,n}} c_\w c_\sigma{\cdot} \psi_{\t^{(\lambda,f)} \t^{(\lambda,f)}} \psi_w f_w \psi_\x^* e(\bi_{(\sigma,f)}) \psi_\y \epsilon_\y \pmod{R_n^{>f}(\delta)}, \label{y:first:eq2}
\end{eqnarray}
where $c_\sigma, c_w \in R$, $w \in \Sym_{2f,n}$ and $f_w$'s are polynomials of $y_1,\ldots,y_n$. Hence, by~\autoref{base:y:1} and~\autoref{base:psi:1}, we have
$$
\psi_{\t^{(\lambda,f)} \t^{(\lambda,f)}} \psi_w f_w \psi_\x^* e(\bi_{(\sigma,f)}) \psi_\y \epsilon_\y \equiv c{\cdot} \psi_{\t^{(\lambda,f)} \w} \epsilon_\y \pmod{R_n^{>f}(\delta)},
$$
where $c \in R$ and $\w \in \Tud_n(\lambda)$ with head $f$. We note that $c \neq 0$ only if $\w = \t^{(\lambda,f)} w d(h(\x))^* d(h(\y))$ and $f_w = 1$. Notice that $\t(n) < 0$ forces $f > 0$. Hence we have $\epsilon_y \in \G{2,n}$ unless $\epsilon_y = e(\bi_\y)$. Because $head(w) = head(\t^{(\lambda,f)}) = f \geq 1$, by~\autoref{remove-front} we have
\begin{equation} \label{y:first:eq3}
\psi_{\t^{(\lambda,f)} \t^{(\lambda,f)}} \psi_w f_w \psi_\x^* e(\bi_{(\sigma,f)}) \psi_\y \epsilon_\y \equiv c{\cdot}\psi_{\t^{(\lambda,f)} \w} \epsilon_\y \equiv \sum_{\v' \in \Tud_n(\lambda)} c_{\v'} \psi_{\t^{(\lambda,f)} \v'} \pmod{R_n^{>f}(\delta)},
\end{equation}
where $c_{\v'} \in R$. Substitute (\ref{y:first:eq3}) into (\ref{y:first:eq2}) implies
$$
\psi_{\t^{(\lambda,f)} \u} \epsilon_2 \epsilon_4 \ldots \epsilon_{2f-2} \psi_{\x\y} \equiv \sum_{\v' \in \Tud_n(\lambda)} c_{\v'} \psi_{\t^{(\lambda,f)} \v'} \pmod{R_n^{>f}(\delta)}.
$$

Finally, the Lemma follows by substituting the above equality into (\ref{y:first:eq1}). \endproof

\subsection{The spanning set of $\G{n}$} \label{sec:induction:2}

In the previous subsection we have proved that (\ref{Induction:1}) holds for any $a \in \G{n-1}$. In this subsection we extend this result to arbitrary $a \in \G{n}$ by showing that (\ref{Induction:1}) holds when $a \in \{y_n, \psi_{n-1}, \epsilon_{n-1}\}$. Then using (\ref{Induction:1}) we prove that $\psi_{\s\t}$'s form a spanning set of $\G{n}$.

We start by proving the equality (\ref{Induction:1}) holds when $a = y_n$.

\begin{Lemma} \label{base:y:2}
Suppose $(\lambda,f) \in \mathscr S_n$ and $\t \in \Tud_n(\lambda)$ with $\t(n) < 0$ and head $f-1$. Then the equality (\ref{Induction:1}) holds when $a = y_n$.
\end{Lemma}

\proof As $\t \in \Tud_n(\lambda)$ with $\t(n) < 0$ and head $f-1$, write $\s = h(\t) \rightarrow \t$  and we have $\rho(\s,\t) = (a,n)$ for $2f \leq a < n$. When $a = n-1$, we have
$$
\psi_{\t^{(\lambda,f)} \t} y_n = e_{(\lambda,f)} \psi_{d(\s)} \epsilon_{2f} \epsilon_{2f+1} \ldots \epsilon_{n-1} e(\bi_\t) y_n = - e_{(\lambda,f)} \psi_{d(\s)} \epsilon_{2f} \epsilon_{2f+1} \ldots \epsilon_{n-1} e(\bi_\t) y_{n-1} = -\psi_{\t^{(\lambda,f)} \t} y_{n-1},
$$
by (\ref{rela:7:2}), and the Lemma holds by~\autoref{y:first}.

When $a < n-1$, we have $\t(n-1) > 0$, $\t(n) < 0$ and $\t(n-1) + \t(n) \neq 0$ by~\autoref{deg:help1}. By~\autoref{y:h2:8} we have $\u = \t{\cdot}s_{n-1} \in \Tud_n(\lambda)$. Hence the constructions of $\t$ and $\u$, we have $\s \rightarrow \u$ and $\epsilon_{\s \rightarrow \u} = e(\bi_\s) \epsilon_{2f} \epsilon_{2f+1} \ldots \epsilon_a \psi_{a+1} \ldots \psi_{n-2} e(\bi_\u)$.

Write $\bi_\t = (i_1, \ldots, i_n)$. Because $\t(n-1) > 0$, $\t(n) < 0$ and $\t(n-1) + \t(n) \neq 0$, we have $i_{n-1} + i_n \neq 0$. If $i_{n-1} \neq i_n$, by (\ref{rela:3:1}) and~\autoref{y:first} we have
\begin{equation} \label{base:y:2:eq1}
\psi_{\t^{(\lambda,f)} \t} y_n = \psi_{\t^{(\lambda,f)} \u} y_{n-1} \psi_{n-1} \equiv \sum_{\v \in \Tud_n(\lambda)} c_{\w} \psi_{\t^{(\lambda,f)} \w}\psi_{n-1} \pmod{R_n^{>f}},
\end{equation}
with $c_\w \in R$. Moreover, because $y_{n-1} \in \G{n-1}$ and $\u(n) = \t(n-1) > 0$, $c_\w \neq 0$ only if $\w(n) > 0$; and because $\t(n) > 0$ forces $f \geq 1$ and $y_{n-1} \in \G{2(f-1),n}$, by~\autoref{remove:head} we have $c_\w \neq 0$ only if $head(\w) \geq f-1$.

For any $\w \in \Tud_n(\lambda)$ with $c_\w \neq 0$ in (\ref{base:y:2:eq1}), by~\autoref{idem:psi:1} we have $\bi_\w = \bi_\u$. If $\w(n-1) < 0$, $head(\w) \geq f-1$ forces $\w = \u$, which implies $\psi_{\t^{(\lambda,f)} \w} \psi_{n-1} = \psi_{\t^{(\lambda,f)} \t}$; and if $\w(n-1) > 0$, we have $\epsilon_\w \in \G{n-2}$ because $\w(n) > 0$, which implies
$$
\psi_{\t^{(\lambda,f)} \w} \psi_{n-1} = \psi_{\t^{(\lambda,f)} \s} \epsilon_\w \psi_{n-1} = \psi_{\t^{(\lambda,f)} \s} \psi_{n-1} \epsilon_\w \equiv \sum_{\v \in \Tud_n(\lambda)} c_{\v} \psi_{\t^{(\lambda,f)}\v} \pmod{R_n^{>f}(\delta)},
$$
by (\ref{rela:2:3}),~\autoref{base:psi:1} and~\autoref{y:first}. Hence the Lemma holds when $i_{n-1} \neq i_n$.

If $i_{n-1} = i_n$, by (\ref{rela:3:1}) we have $\psi_{\t^{(\lambda,f)} \t} y_n = \psi_{\t^{(\lambda,f)} \u} \psi_{n-1} y_n = \psi_{\t^{(\lambda,f)} \u} y_{n-1} \psi_{n-1} + \psi_{\t^{(\lambda,f)} \u}$. Hence the Lemma holds when $i_{n-1} = i_n$ by following the same argument as when $i_{n-1} \neq i_n$. \endproof

\begin{Lemma} \label{y:end}
Suppose $(\lambda,f) \in \mathscr S_n$ and $\t \in \Tud_n(\lambda)$. Then the equality (\ref{Induction:1}) holds when $a = y_n$.
\end{Lemma}

\proof Suppose $head(\t) = h$. We have the standard reduction sequence of $\t$:
$$
h(\t) = \t^{(f-h)} \rightarrow \t^{(f-h-1)} \rightarrow \ldots \rightarrow \t^{(1)} \rightarrow \t^{(0)} = \t.
$$

If $\t(n) > 0$, we have $\epsilon_\t \in \G{n-1}$. Hence by (\ref{rela:2:2}) we have $\psi_{\t^{(\lambda,f)} \t} y_n = \psi_{\t^{(\lambda,f)} h(\t)} y_n \epsilon_\t$. By~\autoref{base:y:1} and~\autoref{y:first}, the Lemma holds.

If $\t(n) < 0$, denote $\s = \t^{(f-h-1)}$. By the construction, we have $\s(n) = \t(n) < 0$ and $\s \in \Tud_n(\lambda)$ with head $f-1$. Because $\epsilon_{\t^{(f-h-1)} \rightarrow \t^{(f-h-2)}} \ldots \epsilon_{\t^{(1)} \rightarrow \t^{(0)}} \in \G{n-1}$, by (\ref{rela:2:2}),~\autoref{base:y:2} and~\autoref{y:first} we have
$$
\psi_{\t^{(\lambda,f)} \t} y_n = \psi_{\t^{(\lambda,f)} \s} y_n \epsilon_{\t^{(f-h-1)} \rightarrow \t^{(f-h-2)}} \ldots \epsilon_{\t^{(1)} \rightarrow \t^{(0)}} \equiv \sum_{\v \in \Tud_n(\lambda)} c_{\v} \psi_{\t^{(\lambda,f)} \v} \pmod{R_n^{>f}(\delta)},
$$
which completes the proof. \endproof

Next we prove that when $a = \epsilon_{n-1}$, the equality (\ref{Induction:1}) holds. We separate the question by considering different $\t(n-1)$ and $\t(n)$. The next Lemma shows (\ref{Induction:1}) holds when $\t(n-1) > 0$ and $\t(n) > 0$.


\begin{Lemma} \label{ep:h1}
Suppose $(\lambda,f) \in \mathscr S_n$ and $\t \in \Tud_n(\lambda)$ with $\t(n-1) > 0$ and $\t(n) > 0$. Then the equality (\ref{Induction:1}) holds when $a = \epsilon_{n-1}$.
\end{Lemma}

\proof In this case we have $\epsilon_\t, \psi_\t \in \G{n-2}$, which commute with $\epsilon_{n-1}$. Hence we have
$$
\psi_{\t^{(\lambda,f)} \t} \epsilon_{n-1} = e_{(\lambda,f)} \psi_\t \epsilon_\t \epsilon_{n-1} = e(\bi_{(\lambda,f)}) \epsilon_1 \epsilon_3 \ldots \epsilon_{2f-1} \epsilon_{n-1} \psi_\t \epsilon_\t \in R_n^{>f}(\delta),
$$
by~\autoref{I:coro:1} and~\autoref{two:sided:ideal}. \endproof

Before proceeding further, we introduce a technical result~\autoref{ep:h3}, which will be used to prove (\ref{Induction:1}) when at least one of $\t(n-1)$ and $\t(n)$ is negative for $a = \epsilon_{n-1}$, and $a = \psi_{n-1}$. The following two Lemmas will be used to prove~\autoref{ep:h3}.

\begin{Lemma} \label{ep:h2}
Suppose $(\lambda,f-1) \in \widehat B_{n-2}$ and $(\mu,f) \in \widehat B_n$. We have either (1) $\mu = \lambda$, or (2) $\u(n-1) > 0$ and $\u(n) > 0$ if there exist $\s \in \Tud_{n-2}(\lambda)$ and $\u \in \Tud_n(\mu)$ such that the following conditions are satisfied for some $i \in P$:
\begin{enumerate}
\item $head(\s) = head(\u) = f-1$.

\item $\bi_\u = (\bi_\s \vee i, -i)$.

\item $\res(\alpha) \neq i$ for all $\alpha \in \mathscr A(\lambda)$.
\end{enumerate}
\end{Lemma}

\proof Because $\u \in \Tud_n(\lambda)$ with $(\lambda,f) \in \widehat B_n$ and $head(\u) = f-1$, there exists a unique $k$ with $2f+1 \leq k \leq n$ such that $\u(k) < 0$. If $k \leq n-2$, we have $\u(n-1) > 0$ and $\u(n) > 0$.

Whenever $n-1 \leq k \leq n$, let $\x = \u|_{n-2}$. Then we have $\Shape(\x) = (\gamma,f-1) \in \widehat B_{n-2}$ for some partition $\gamma$ and $head(\x) = head(\u) = f-1$. Moreover, because $\bi_\u = (\bi_\s\vee i,-i)$, we have $\bi_\x = \bi_\s$. By~\autoref{I:h1} it forces $\x = \s$ and $\u_{n-2} = \gamma = \lambda$.

Because $\res(\alpha) \neq i$ for all $\alpha \in \mathscr A(\lambda)$, it forces $|\mathscr {AR}_\lambda(i)| = 1$. Therefore, by~\eqref{deg:h4:eq1} -~\eqref{deg:h4:eq3} we have $h_k(\bi_\u) = 0$ or $-1$, which implies $\u(n-1) + \u(n) = 0$ by~\autoref{deg:h4:2}. Hence, we have $\u_{n-2} = \mu$, which yields $\mu = \lambda$ as we have shown in the last paragraph that $\u_{n-2} = \lambda$. \endproof

\begin{Lemma} \label{ep:h4}
Suppose $(\lambda,f) \in \widehat B_n$ and $\t, \x, \y \in \Tud_n(\lambda)$. If $head(\x) \geq f-1$, we have
$$
e(\bi_{(\lambda,f)}) \psi_{d(\t)} \epsilon_{2f-1} \epsilon_{2f} \ldots \epsilon_{n-1} \psi_{\x\y} \equiv c{\cdot} \psi_{\t^{(\lambda,f)} \y} \pmod{R_n^{>f}(\delta)},
$$
for some $c \in R$.
\end{Lemma}

\proof Because $\x \in \Tud_n(\lambda)$ with $(\lambda,f) \in \widehat B_n$ and $head(\x) \geq f-1$, we have $head(\x) = f-1$ or $f$. When $head(\x) = f$, by~\autoref{base:ep:1} we have $\epsilon_{n-1} \psi_{\x\y} \in R_n^{>f}(\delta)$ and the Lemma follows by~\autoref{two:sided:ideal}.

When $head(\x) = f-1$, let $\w = h(\x) \rightarrow \x$. By~\autoref{commute:1}, we have
$$
e(\bi_{(\lambda,f)}) \psi_{d(\t)} \epsilon_{2f-1} \epsilon_{2f} \epsilon_{2f+1} \ldots \epsilon_{n-1} \psi_{\x\y} \equiv c_w e(\bi_{(\lambda,f)}) \psi_{d(\t)} f_w \psi_w \psi_{\w\y} \pmod{R_n^{>f}(\delta)},
$$
where $c_w \in R$, $w \in \Sym_n$ and $f_w$ is a polynomial of $y_1, \ldots, y_n$. As $head(\w) = f$, by~\autoref{base:psi:1} and~\autoref{base:y:1}, we have
\begin{align*}
e(\bi_{(\lambda,f)}) \psi_{d(\t)} \epsilon_{2f-1} \epsilon_{2f} \epsilon_{2f+1} \ldots \epsilon_{n-1} \psi_{\x\y} & \equiv c_w e(\bi_{(\lambda,f)}) \psi_{d(\t)} f_w \psi_w \psi_{\w\y} \equiv c_w{\cdot} \psi_{\u \y} \pmod{R_n^{>f}(\delta)},
\end{align*}
for some $\u \in \Tud_n(\lambda)$ with head $f$. By~\autoref{idem:psi:1}, we have $\bi_\u = \bi_{(\lambda,f)}$. As $\u$ has head $f$, by~\autoref{I:h1}, it forces $\u = \t^{(\lambda,f)}$. Hence the Lemma holds as $c_w \in R$. \endproof

\begin{Lemma} \label{ep:h3}
Suppose $(\lambda,f) \in \mathscr S_n$ and $\t \in \Tud_n(\lambda)$ with head $f$. Then we have
$$
\psi_{\t^{(\lambda,f)} \t} \epsilon_{2f} \epsilon_{2f+1} \ldots \epsilon_{n-1} \equiv \sum_{\v \in \Tud_n(\lambda)} c_{\v} \psi_{\t^{(\lambda,f)}\v} \pmod{R_n^{>f}(\delta)},
$$
where $c_\v \in R$.
\end{Lemma}

\proof It suffices to prove the Lemma by showing
\begin{equation} \label{ep:h3:eq}
\psi_{\t^{(\lambda,f)} \t} \epsilon_{2f} \epsilon_{2f+1} \ldots \epsilon_{n-1} e(\bj) \equiv \sum_{\v \in \Tud_n(\lambda)} c_{\v} \psi_{\t^{(\lambda,f)}\v} \pmod{R_n^{>f}(\delta)},
\end{equation}
for any $\bj \in P^n$. Notice that if we write $\bi_\t = (i_1, \ldots, i_n)$, $e(\bi_\t) \epsilon_{2f}\epsilon_{2f+1} \ldots \epsilon_{n-1} e(\bj) \neq 0$ only if $j_r = i_r$ for $1 \leq r \leq 2f-1$, $j_r = i_{r+2}$ for $2f \leq r \leq n-2$ and $j_{n-1} + j_n = 0$ by (\ref{rela:1}). In the rest of the proof, we assume $\bj$ has such property.

Suppose there exists $\alpha \in \mathscr A(\lambda)$ such that $\res(\alpha) = j_{n-1}$. Write $\t = (\alpha_1, \ldots, \alpha_n)$ and define
$$
\u = (\alpha_1, \alpha_2, \ldots, \alpha_{2f-1}, \alpha_{2f+2}, \alpha_{2f+3}, \ldots, \alpha_n, \alpha, -\alpha).
$$

Hence we have $\bi_\u = \bj$. Moreover, we have $\t \rightarrow \u$ and $\epsilon_{\t\rightarrow \u} = e(\bi_\t) \epsilon_{2f} \epsilon_{2f+1} \ldots \epsilon_{n-1} e(\bj)$. Therefore, by the definition of $\psi_{\s\t}$'s we have $\psi_{\t^{(\lambda,f)} \t} \epsilon_{2f} \epsilon_{2f+1} \ldots \epsilon_{n-1} e(\bj) = \psi_{\t^{(\lambda,f)} \u}$, which proves that~\eqref{ep:h3:eq} holds.

It left us to show that~\eqref{ep:h3:eq} holds when $\res(\alpha) \neq j_{n-1}$ for all $\alpha \in \mathscr A(\lambda)$. In this case, we have
\begin{equation} \label{ep:h3:eq1}
\psi_{\t^{(\lambda,f)} \t} \epsilon_{2f} \epsilon_{2f+1} \ldots \epsilon_{n-1} e(\bj) = e(\bi_{(\lambda,f)}) \psi_{d(\t)} \epsilon_{2f-1} \epsilon_{2f} \epsilon_{2f+1} \ldots \epsilon_{n-1} e(\bj) \epsilon_1 \epsilon_3 \ldots \epsilon_{2f-3} e(\bj).
\end{equation}

Because $j_r = i_r$ for $1 \leq r \leq 2f-1$, we have
$$
e(\bj) \epsilon_1 \epsilon_3 \ldots \epsilon_{2f-3} e(\bj) = \theta^{(n-1)}_{j_n} \circ \theta^{(n-2)}_{j_{n-1}} \circ \ldots \circ\theta_{j_{2f-1}}^{(2f-2)}(\psi_{\t^{(\emptyset,f-1)} \t^{(\emptyset,f-1)}}).
$$

By applying~\autoref{I:h1} and~\autoref{I:coro:2} recursively to the above equality, we have
\begin{equation} \label{ep:h3:eq3}
e(\bj) \epsilon_1 \epsilon_3 \ldots \epsilon_{2f-3} e(\bj) \equiv \sum_{\substack{(\mu,f) \in \widehat B_n \\ \x,\y \in \Tud_n(\mu)}} c_{\x\y} \psi_{\x\y} \pmod{R_n^{>f}(\delta)},
\end{equation}
where $c_{\x\y} \neq 0$ only if $head(\x) \geq f-1$ and $head(\y) \geq f-1$. Moreover, by~\autoref{idem:psi:1}, we have $\bi_\x = \bi_\y = \bj$ if $c_{\x\y} \neq 0$. Hence, by~\autoref{ep:h2}, we have either $\x,\y \in \Tud_n(\lambda)$ or $\x(n-1) > 0$ and $\x(n) > 0$ if $c_{\x\y} \neq 0$.

Substituting (\ref{ep:h3:eq3}) into (\ref{ep:h3:eq1}) yields
\begin{equation} \label{ep:h3:eq4}
\psi_{\t^{(\lambda,f)} \t} \epsilon_{2f} \epsilon_{2f+1} \ldots \epsilon_{n-1} e(\bj) \equiv \sum_{\substack{(\mu,f) \in \widehat B_n \\ \x,\y \in \Tud_n(\mu)}} c_{\x\y} e(\bi_{(\lambda,f)}) \psi_{d(\t)} \epsilon_{2f-1} \epsilon_{2f} \epsilon_{2f+1} \ldots \epsilon_{n-1} \psi_{\x\y} \pmod{R_n^{>f}(\delta)}.
\end{equation}

If $\x(n-1) > 0$ and $\x(n) > 0$, by~\autoref{ep:h1} and~\autoref{two:sided:ideal} we have
\begin{equation} \label{ep:h3:eq5}
e(\bi_{(\lambda,f)}) \psi_{d(\t)} \epsilon_{2f-1} \epsilon_{2f} \epsilon_{2f+1} \ldots \epsilon_{n-1} \psi_{\x\y} \in R_n^{>f}(\delta).
\end{equation}

If $\x,\y \in \Tud_n(\lambda)$, recall $head(\x) \geq f-1$. Hence by~\autoref{ep:h4}, we have
\begin{equation} \label{ep:h3:eq6}
e(\bi_{(\lambda,f)}) \psi_{d(\t)} \epsilon_{2f-1} \epsilon_{2f} \epsilon_{2f+1} \ldots \epsilon_{n-1} \psi_{\x\y}
\equiv c{\cdot} \psi_{\t^{(\lambda,f)} \y} \pmod{R_n^{>f}(\delta)},
\end{equation}
where $c \in R$. Hence~\eqref{ep:h3:eq} follows by substituting (\ref{ep:h3:eq5}) and (\ref{ep:h3:eq6}) into (\ref{ep:h3:eq4}). \endproof

Then we consider the cases when at least one of $\t(n-1)$ and $\t(n)$ is negative.

\begin{Lemma} \label{ep:first:1}
Suppose $(\lambda,f) \in \mathscr S_n$ and $\t \in \Tud_n(\lambda)$ with $\t(n-1) < 0$ and $\t(n) < 0$. Then the equality (\ref{Induction:1}) holds when $a = \epsilon_{n-1}$.
\end{Lemma}

\proof Because $\t(n-1) < 0$ and $\t(n) < 0$, the standard reduction sequence $\t$ is
$$
\s = \t^{(m)} \rightarrow \t^{(m-1)} \rightarrow \ldots \t^{(1)} \rightarrow \t^{(0)} = \t,
$$
where $\rho(\t^{(m)}, \t^{(m-1)}) = (a,n)$ and $\rho(\t^{(m-1)},\t^{(m-2)}) = (b,n-1)$ for $a \leq n-1$ and $b \leq n-2$. Hence we can write
\begin{align*}
\epsilon_{\t^{(m)} \rightarrow \t^{(m-1)}} & = e(\bi_{\t^{(m)}}) \epsilon_{2f} \epsilon_{2f+1} \ldots \epsilon_a \psi_{a+1} \ldots \psi_{n-1} e(\bi_{\t^{(m-1)}}),\\
\epsilon_{\t^{(m-1)} \rightarrow \t^{(m-2)}} & = e(\bi_{\t^{(m-1)}}) \epsilon_{2f-2} \epsilon_{2f-1} \ldots \epsilon_b \psi_{b+1} \ldots \psi_{n-2} e(\bi_{\t^{(m-2)}}).
\end{align*}

By (\ref{rela:9}), we have
\begin{equation} \label{ep:first:1:eq1}
\epsilon_{\t^{(m-1)} \rightarrow \t^{(m-2)}} \epsilon_{n-1} = e(\bi_{\t^{(m-1)}}) \epsilon_{2f-2} \epsilon_{2f-1} \ldots \epsilon_{n-1} (\psi_{n-3} \psi_{n-4} \ldots \psi_b);
\end{equation}
and by (\ref{rela:12}) and (\ref{rela:13}) we have
\begin{equation} \label{ep:first:1:eq2}
\epsilon_{\t^{(m)} \rightarrow \t^{(m-1)}} \epsilon_{2f-2} \epsilon_{2f-1} \ldots \epsilon_{n-1} = e(\bi_{\t^{(m)}}) \epsilon_{2f} \epsilon_{2f+1} \ldots \epsilon_{n-1} \left(\epsilon_{2f-2} \epsilon_{2f-1} \ldots \epsilon_{a-2} \psi_{a-1} \ldots \psi_{n-3}\right).
\end{equation}

Set $x = \epsilon_{2f-2} \epsilon_{2f-1} \ldots \epsilon_{a-2} \psi_{a-1} \ldots \psi_{n-3} {\cdot}\psi_{n-3} \psi_{n-4} \ldots \psi_b$ and $y = \epsilon_{\t^{(m-2)} \rightarrow \t^{(m-3)}} \ldots \epsilon_{\t^{(1)} \rightarrow \t} \in \G{n-2}$. Because $a \leq n-1$ and $b \leq n-2$, we have $x \in \G{n-1}$. Then by (\ref{rela:2:3}), (\ref{ep:first:1:eq1}) and (\ref{ep:first:1:eq2}), we have
\begin{align*}
\psi_{\t^{(\lambda,f)} \t} \epsilon_{n-1} & = \psi_{\t^{(\lambda,f)} \s} \epsilon_{\t^{(m)} \rightarrow \t^{(m-1)}} \epsilon_{\t^{(m-1)} \rightarrow \t^{(m-2)}} y \epsilon_{n-1} = \psi_{\t^{(\lambda,f)} \s} \epsilon_{\t^{(m)} \rightarrow \t^{(m-1)}} \epsilon_{\t^{(m-1)} \rightarrow \t^{(m-2)}} \epsilon_{n-1} y\\
& = \psi_{\t^{(\lambda,f)} \s} \epsilon_{2f} \epsilon_{2f+1} \ldots \epsilon_{n-1} {\cdot}xy.
\end{align*}

Because $\s = h(\t) \in \Tud_n(\lambda)$, the head of $\s$ is $f$. Then because $xy \in \G{n-1}$, the Lemma follows by~\autoref{ep:h3} and~\autoref{y:first}. \endproof

\begin{Lemma} \label{ep:first:2}
Suppose $(\lambda,f) \in \mathscr S_n$ and $\t \in \Tud_n(\lambda)$ with $\t(n-1) < 0$, $\t(n) > 0$ or $\t(n-1) > 0$, $\t(n) < 0$. Then the equality (\ref{Induction:1}) holds when $a = \epsilon_{n-1}$.\end{Lemma}

\proof Suppose $\bi_\t = (i_1, \ldots, i_n)$. We assume $i_{n-1} + i_n = 0$, otherwise by (\ref{rela:1}) we have $\psi_{\t^{(\lambda,f)} \t} \epsilon_{n-1} = 0$. Hence, as $\t(n-1) < 0$, $\t(n) > 0$ or $\t(n-1) > 0$, $\t(n) < 0$, by the construction of up-down tableaux, we have $\t(n-1) + \t(n) = 0$.

Let the standard reduction sequence of $\t$ be
$$
\s = \t^{(m)} \rightarrow \t^{(m-1)} \rightarrow \ldots \t^{(1)} \rightarrow \t^{(0)} = \t,
$$
and denote $x = \epsilon_{\t^{(m-1)} \rightarrow \t^{(m-2)}} \epsilon_{\t^{(m-2)} \rightarrow \t^{(m-3)}} \ldots \epsilon_{\t^{(1)} \rightarrow \t}$. Notice that because $\t(n-1) > 0$, $\t(n) < 0$ or $\t(n-1) < 0$, $\t(n) > 0$, we have $x \in \G{n-2}$ in either cases. Hence by (\ref{rela:2:3}), $x$ commutes with $\epsilon_{n-1}$.

When $\t(n-1) > 0$ and $\t(n) < 0$, we have $\rho(\t^{(m)}, \t^{(m-1)}) = (n-1,n)$. Then by (\ref{rela:5:5}) and (\ref{rela:16}), we have
\begin{align*}
\psi_{\t^{(\lambda,f)} \t} \epsilon_{n-1}  & = \psi_{\t^{(\lambda,f)} \s} \epsilon_{\t^{(m)} \rightarrow \t^{(m-1)}} x \epsilon_{n-1}  = \psi_{\t^{(\lambda,f)} \s} \epsilon_{\t^{(m)} \rightarrow \t^{(m-1)}} \epsilon_{n-1} x \\
& = \psi_{\t^{(\lambda,f)} \s} \epsilon_{2f} \epsilon_{2f+1} \ldots \epsilon_{n-1} e(\bi_{\t^{(m-1)}}) \epsilon_{n-1} x \\
& = \psi_{\t^{(\lambda,f)} \s} \epsilon_{2f} \epsilon_{2f+1} \ldots \epsilon_{n-1} f(y_1, \ldots, y_{n-2}) x ,
\end{align*}
where $f(y_1, \ldots, y_{n-2})$ is a polynomial of $y_1, \ldots, y_{n-2}$.

Because $\s = h(\t) \in \Tud_n(\lambda)$ with head $f$ and $f(y_1, \ldots, y_{n-2}) x \in \G{n-2}$, by~\autoref{ep:h3} and~\autoref{y:first}, the Lemma holds when $\t(n-1) > 0$ and $\t(n) < 0$.

When $\t(n-1) < 0$ and $\t(n) > 0$, we have $\rho(\t^{(m)}, \t^{(m-1)}) = (a,n-1)$ with $a \leq n-2$. Then by (\ref{rela:9}) we have
\begin{align*}
\psi_{\t^{(\lambda,f)} \t} \epsilon_{n-1}  & = \psi_{\t^{(\lambda,f)} \s} \epsilon_{\t^{(m)} \rightarrow \t^{(m-1)}} x \epsilon_{n-1}  = \psi_{\t^{(\lambda,f)} \s} \epsilon_{\t^{(m)} \rightarrow \t^{(m-1)}} \epsilon_{n-1} x \\
& = \psi_{\t^{(\lambda,f)} \s} \epsilon_{2f} \epsilon_{2f+1} \ldots \epsilon_a \psi_{a+1} \ldots \psi_{n-2} \epsilon_{n-1} x  \\
& = \psi_{\t^{(\lambda,f)} \s} \epsilon_{2f} \epsilon_{2f+1} \ldots \epsilon_{n-1} \psi_{n-3} \psi_{n-4} \ldots \psi_a x .
\end{align*}

Because $\s = h(\t) \in \Tud_n(\lambda)$ with head $f$ and $\psi_{n-3} \psi_{n-4} \ldots \psi_a x \in \G{n-2}$, by~\autoref{ep:h3} and~\autoref{y:first}, the Lemma holds when $\t(n-1) < 0$ and $\t(n) > 0$. \endproof

Combining~\autoref{ep:h1},~\autoref{ep:first:1} and~\autoref{ep:first:2}, the next Lemma follows.

\begin{Lemma} \label{ep:end}
Suppose $(\lambda,f) \in \mathscr S_n$ and $\t \in \Tud_n(\lambda)$. Then the equality (\ref{Induction:1}) holds when $a = \epsilon_{n-1}$.
\end{Lemma}

Finally, we prove that when $a = \psi_{n-1}$, the equality (\ref{Induction:1}) holds. Similar as before, we separate the question by considering different $\t(n-1)$ and $\t(n)$.

\begin{Lemma} \label{psi:first:1}
Suppose $(\lambda,f) \in \mathscr S_n$ and $\t \in \Tud_n(\lambda)$ with $\t(n) > 0$. Then the equality (\ref{Induction:1}) holds when $a = \psi_{n-1}$.
\end{Lemma}

\proof When $\t(n-1) > 0$, we have $\epsilon_\t \in \G{n-2}$, which commutes with $\psi_{n-1}$. Denote $\s = h(\t) \in \Tud_n(\lambda)$, and we have $head(\s) = f$. By~\autoref{base:psi:1} and~\autoref{y:first}, we have
$$
\psi_{\t^{(\lambda,f)} \t} \psi_{n-1}  = \psi_{\t^{(\lambda,f)} \s} \psi_{n-1} \epsilon_\t  \equiv \sum_{\v \in \Tud_n(\lambda)} c_\v \psi_{\t^{(\lambda,f)} \v} \pmod{B_n^{>f}(\delta)}.
$$

When $\t(n-1) < 0$, write $\t(n-1) = -\alpha$ and $\t(n) = \beta$. If $\alpha \neq \beta$, by~\autoref{y:h2:8} we have $\u = \t{\cdot}s_{n-1} \in \Tud_n(\lambda)$ and by the definition of $\psi_{\s\t}$'s, we have $\psi_{\t^{(\lambda,f)} \t} \psi_{n-1} = \psi_{\t^{(\lambda,f)} \u}$. Hence the Lemma holds.

If $\alpha = \beta$, let the standard reduction sequence of $\t$ be
$$
\s = \t^{(m)} \rightarrow \t^{(m-1)} \rightarrow \ldots \t^{(1)} \rightarrow \t^{(0)} = \t,
$$
and denote $x = \epsilon_{\t^{(m-1)} \rightarrow \t^{(m-2)}} \epsilon_{\t^{(m-2)} \rightarrow \t^{(m-3)}} \ldots \epsilon_{\t^{(1)} \rightarrow \t}$. Notice that as $\t(n-1) < 0$ and $\t(n) > 0$, we have $x \in \G{n-2}$, which commutes with $\psi_{n-1}$. Let $\rho(\t^{(m)},\t^{(m-1)}) = (a,n-1)$ with $a \leq n-2$. By (\ref{rela:9}) we have
\begin{align*}
\psi_{\t^{(\lambda,f)} \t} \psi_{n-1}  & = \psi_{\t^{(\lambda,f)} \s} \epsilon_{\t^{(m)} \rightarrow \t^{(m-1)}} x \psi_{n-1}  = \psi_{\t^{(\lambda,f)} \s} \epsilon_{\t^{(m)} \rightarrow \t^{(m-1)}} \psi_{n-1} x \\
& = \psi_{\t^{(\lambda,f)} \s} \epsilon_{2f} \epsilon_{2f+1} \ldots \epsilon_a \psi_{a+1} \ldots \psi_{n-2} \psi_{n-1} x \\
& = \psi_{\t^{(\lambda,f)} \s} \epsilon_{2f} \epsilon_{2f+1} \ldots \epsilon_{n-1} \left(\psi_{n-2} \psi_{n-3} \ldots \psi_a x \right).
\end{align*}

Because $\s = h(\t) \in \Tud_n(\lambda)$ with head $f$ and $\psi_{n-2} \psi_{n-3} \ldots \psi_a x \in \G{n-1}$, the Lemma follows by~\autoref{ep:h3} and~\autoref{y:first}. \endproof

\begin{Lemma} \label{psi:first:2}
Suppose $(\lambda,f) \in \mathscr S_n$ and $\t \in \Tud_n(\lambda)$ with $\t(n-1) > 0$ and $\t(n) < 0$. Then the equality (\ref{Induction:1}) holds when $a = \psi_{n-1}$.
\end{Lemma}

\proof Let the standard reduction sequence of $\t$ be
$$
\s = \t^{(m)} \rightarrow \t^{(m-1)} \rightarrow \ldots \t^{(1)} \rightarrow \t^{(0)} = \t,
$$
and denote $x = \epsilon_{\t^{(m-1)} \rightarrow \t^{(m-2)}} \epsilon_{\t^{(m-2)} \rightarrow \t^{(m-3)}} \ldots \epsilon_{\t^{(1)} \rightarrow \t}$. Notice that as $\t(n-1) > 0$ and $\t(n) < 0$, we have $x \in \G{n-2}$, which commutes with $\psi_{n-1}$.

Suppose $\t(n-1) = \beta$ and $\t(n) = -\alpha$. If $\alpha \neq \beta$, by~\autoref{y:h2:8} we have $\u = \t{\cdot}s_{n-1} \in \Tud_n(\lambda)$ and $\psi_{\t^{(\lambda,f)} \t} = \psi_{\t^{(\lambda,f)} \u} \psi_{n-1}$ by the construction of $\psi_{\t^{(\lambda,f)} \t}$ and $\psi_{\t^{(\lambda,f)} \u}$. Hence by (\ref{rela:4}) we have
$$
\psi_{\t^{(\lambda,f)} \t} \psi_{n-1}  = \psi_{\t^{(\lambda,f)} \u} \psi_{n-1}^2  = \psi_{\t^{(\lambda,f)} \u} f(y_{n-1}, y_n) ,
$$
where $f(y_{n-1}, y_n)$ is a polynomial of $y_{n-1}$ and $y_n$ determined by $\bi_\u$. Hence the Lemma holds by~\autoref{y:end}.

If $\alpha = \beta$, we have $\rho(\t^{(m)}, \t^{(m-1)}) = (n-1,n)$. By (\ref{rela:6:2}), we have
\begin{align*}
\psi_{\t^{(\lambda,f)} \t} \psi_{n-1}  & = \psi_{\t^{(\lambda,f)} \s} \epsilon_{\t^{(m)}, \t^{(m-1)}} x \psi_{n-1}  = \psi_{\t^{(\lambda,f)} \s} \epsilon_{\t^{(m)}, \t^{(m-1)}} \psi_{n-1} x \\
& = \psi_{\t^{(\lambda,f)} \s} \epsilon_{2f} \epsilon_{2f+1} \ldots \epsilon_{n-1} e(\bi_\t) \psi_{n-1} x  \\
& = c{\cdot} \psi_{\t^{(\lambda,f)} \s} \epsilon_{2f} \epsilon_{2f+1} \ldots \epsilon_{n-1} x ,
\end{align*}
where $c \in R$. Hence, as $\s = h(\t) \in \Tud_n(\lambda)$ with head $f$ and $x \in \G{n-2}$, by~\autoref{ep:h3} and~\autoref{y:first}, the Lemma holds. \endproof

\begin{Lemma} \label{psi:first:3}
Suppose $(\lambda,f) \in \mathscr S_n$ and $\t \in \Tud_n(\lambda)$ with $\t(n-1) < 0$ and $\t(n) < 0$. Then the equality (\ref{Induction:1}) holds when $a = \psi_{n-1}$.
\end{Lemma}

\proof Let the standard reduction sequence of $\t$ be
$$
\s = \t^{(m)} \rightarrow \t^{(m-1)} \rightarrow \ldots \t^{(1)} \rightarrow \t^{(0)} = \t,
$$
and denote $x = \epsilon_{\t^{(m-2)} \rightarrow \t^{(m-3)}} \epsilon_{\t^{(m-3)} \rightarrow \t^{(m-4)}} \ldots \epsilon_{\t^{(1)} \rightarrow \t}$. Notice that as $\t(n-1) < 0$ and $\t(n) < 0$, we have $x \in \G{n-2}$, which commutes with $\psi_{n-1}$.

Suppose $\rho(\t^{(m)}, \t^{(m-1)}) = (a,n)$ and $\rho(\t^{(m-1)}, \t^{(m-2)}) = (b,n-1)$ where $a \leq n-1$ and $b \leq n-2$. Hence we can write
\begin{align*}
\epsilon_{\t^{(m)} \rightarrow \t^{(m-1)}} & = e(\bi_{\t^{(m)}}) \epsilon_{2f} \epsilon_{2f+1} \ldots \epsilon_a \psi_{a+1} \ldots \psi_{n-1} e(\bi_{\t^{(m-1)}}),\\
\epsilon_{\t^{(m-1)} \rightarrow \t^{(m-2)}} & = e(\bi_{\t^{(m-1)}}) \epsilon_{2f-2} \epsilon_{2f-1} \ldots \epsilon_b \psi_{b+1} \ldots \psi_{n-2} e(\bi_{\t^{(m-2)}}).
\end{align*}

By (\ref{rela:9}), (\ref{rela:12}) and (\ref{rela:13}), we have
\begin{align*}
\psi_{\t^{(\lambda,f)} \t} \psi_{n-1}  & = \psi_{\t^{(\lambda,f)} \s} \epsilon_{\t^{(m)} \rightarrow \t^{(m-1)}} \epsilon_{\t^{(m-1)} \rightarrow \t^{(m-2)}} x \psi_{n-1}  = \psi_{\t^{(\lambda,f)} \s} \epsilon_{\t^{(m)} \rightarrow \t^{(m-1)}} \epsilon_{\t^{(m-1)} \rightarrow \t^{(m-2)}} \psi_{n-1} x \\
& = \psi_{\t^{(\lambda,f)} \s} \left(\epsilon_{2f} \epsilon_{2f+1} \ldots \epsilon_a \psi_{a+1} \ldots \psi_{n-1}\right){\cdot} \left(\epsilon_{2f-2} \epsilon_{2f-1} \ldots \epsilon_b \psi_{b+1} \ldots \psi_{n-1}\right) x  \\
& = \psi_{\t^{(\lambda,f)} \s} \left(\epsilon_{2f} \epsilon_{2f+1} \ldots \epsilon_a \psi_{a+1} \ldots \psi_{n-1}\right) {\cdot} \left(\epsilon_{2f-2} \epsilon_{2f-1} \ldots \epsilon_{n-1}\right){\cdot} \left(\psi_{n-2} \psi_{n-3} \ldots \psi_b \right) x  \\
& = \psi_{\t^{(\lambda,f)} \s} \epsilon_{2f} \epsilon_{2f+1} \ldots \epsilon_{n-1} {\cdot} \left(\epsilon_{2f-2} \epsilon_{2f-1} \ldots \epsilon_{a-2} \psi_{a-1} \ldots \psi_{n-3}\right){\cdot} \left(\psi_{n-2} \psi_{n-3} \ldots \psi_b\right) x .
\end{align*}

As $\s = h(\t) \in \Tud_n(\lambda)$ with head $f$ and $\epsilon_{2f-2} \epsilon_{2f-1} \ldots \epsilon_{a-2} \psi_{a-1} \ldots \psi_{n-3} \psi_{n-2} \psi_{n-3} \ldots \psi_b x \in \G{n-1}$, by~\autoref{ep:h3} and~\autoref{y:first}, the Lemma holds. \endproof

Combining~\autoref{psi:first:1},~\autoref{psi:first:2} and~\autoref{psi:first:3}, the next Lemma follows.

\begin{Lemma} \label{psi:end}
Suppose $(\lambda,f) \in \mathscr S_n$ and $\t \in \Tud_n(\lambda)$. Then the equality (\ref{Induction:1}) holds when $a = \psi_{n-1}$.
\end{Lemma}

Therefore, we have proved that (\ref{Induction:1}) holds. Combining all the results of this section, we are ready to give the first main result of this paper.

\begin{Proposition} \label{basis:main}
We have $\bigcup_{n \geq 1} \widehat B_n = \widehat{\mathscr B}$.
\end{Proposition}

\proof By the definition we have $\widehat {\mathscr B} \subseteq \bigcup_{n\geq 1}\widehat B_n$. It suffices to show that $\bigcup_{n\geq 1}\widehat B_n \subseteq \widehat {\mathscr B}$.

We have defined a total ordering $<$ on $\bigcup_{n\geq 1}\widehat B_n$. Hence we can list all the elements of $\bigcup_{n\geq 1}\widehat B_n$ in decreasing order as:
$$
(\lambda_1, f_1) > (\lambda_2,f_2) > (\lambda_3,f_3) > \ldots.
$$

We prove the Proposition by induction. For the base step it is easy to see that the maximal element in $\bigcup_{n \geq 1} \widehat B_n$ is $(\lambda_1, f_1) = ((1), 0) \in \widehat B_1$. Because $\G{1} \cong R$, it is obvious that for any $\s,\t \in \Tud_1(\lambda_1)$ and $a \in \G{1}$ we have
$$
\psi_{\s \t} a \equiv \sum_{\v \in \Tud_1(\lambda_1)} c_\v \psi_{\s \v} \pmod{R_1^{>f_1}(\delta)},
$$
which implies $(\lambda_1, f_1) \in \widehat {\mathscr B}$.

For induction step, assume $k > 1$ and $(\lambda_i,f_i) \in \widehat {\mathscr B}$ for any $1 \leq i \leq k-1$. Let $(\lambda,f) = (\lambda_k,f_k) \in \widehat B_n$. By the definition we have $(\lambda,f) \in \mathscr S_n$. For any $\s,\t \in \Tud_n(\lambda)$ and $a \in \G{n}$, by~\autoref{y:first},~\autoref{y:end},~\autoref{ep:end} and~\autoref{psi:end}, we have
$$
\psi_{\t^{(\lambda,f)} \t} a \equiv \sum_{\v \in \Tud_n(\lambda)} c_\v \psi_{\t^{(\lambda,f)} \v} \pmod{R_n^{>f}(\delta)}.
$$

Multiply $\epsilon_\s^* \psi_\s^*$ from left to the above equation. By~\autoref{two:sided:ideal} we have
$$
\psi_{\s \t} a \equiv \sum_{\v \in \Tud_n(\lambda)} c_\v \psi_{\s \v} \pmod{R_n^{>f}(\delta)},
$$
which implies $(\lambda,f) = (\lambda_k,f_k) \in \widehat {\mathscr B}$. Therefore we completes the induction process. Hence, for any $(\lambda,f) \in \bigcup_{n\geq 1}\widehat B_n$, we have $(\lambda,f) \in \widehat {\mathscr B}$, which proves $\bigcup_{n\geq 1}\widehat B_n \subseteq \widehat {\mathscr B}$. \endproof

\begin{Theorem} \label{main:1}
The algebra $\G{n}$ is spanned by $\set{\psi_{\s\t} | (\lambda,f) \in \widehat B_n, \s,\t \in \Tud_n(\lambda)}$.
\end{Theorem}

\proof \autoref{basis:main} implies that for any $(\lambda,f) \in \widehat B_n$ and $\s,\t \in \Tud_n(\lambda)$, we have $\psi_{\s\t} \G{n} \subseteq R_n(\delta)$. Hence we have $R_n(\delta) \G{n} \subseteq R_n(\delta)$. By~\autoref{basis:main} we have $\mathscr S_n = \widehat B_n$, which implies $e(\bi) \in R_n(\delta)$ for any $\bi \in P^n$ by~\autoref{I:end}. As $1 = \sum_{\bi\in P^n} e(\bi)$ by (\ref{rela:1}), we have $1 \in R_n(\delta)$, which yields $\G{n} \subseteq R_n \G{n} \subseteq R_n(\delta)$.

Because $R_n(\delta) \subseteq \G{n}$ by the definition, we have $R_n(\delta) = \G{n}$, which proves the Theorem. \endproof

\autoref{main:1} shows that $\G{n}$ is a finite-dimensional $R$-space with $\dim \G{n} \leq (2n-1)!!$. The following results are directly implied.

Recall $I^n$ is the subset of $P^n$ containing all the residue sequence of up-down tableaux.

\begin{Corollary} \label{idem:nonzero}
Suppose $\bi \in P^n$. We have $e(\bi) = 0$ if $\bi \not\in I^n$.
\end{Corollary}

\proof If $\bi \not\in I^n$, we have $\Tud_n(\bi) = \emptyset$. By~\autoref{idem:psi:1}, we have $\psi_{\s\t} e(\bi) = 0$ for any up-down tableaux $\s$ and $\t$. Hence by~\autoref{main:1}, we have $e(\bi) = e(\bi)^2 = \sum_{\s,\t} c_{\s\t} \psi_{\s\t} e(\bi) = 0$.  \endproof

\begin{Corollary} \label{y:nil}
The elements $y_k \in \G{n}$ are nilpotent for all $1 \leq k \leq n$.
\end{Corollary}

\proof By~\autoref{main:1}, we have $\G{n} = R_n(\delta)$. Because there are finite elements in $\set{\psi_{\s\t} | (\lambda,f) \in \widehat B_n, \s,\t \in \Tud_n(\lambda)}$, there exists $m$ such that $\deg \psi_{\s\t} \leq m$. Hence for any homogeneous element $a \in \G{n}$, we have $\deg a \leq m$. Choose $N = \lfloor\frac{m}{2}\rfloor + 1$. For any $1 \leq k \leq n$ we have $\deg y_k^N > m$, which forces $y_k^N = 0$. \endproof

\begin{Remark}
\autoref{main:1} shows that $\set{\psi_{\s\t} | (\lambda,f) \in \widehat B_n, \s,\t \in \Tud_n(\lambda)}$ spans $\G{n}$, and~\autoref{basis:main} shows that $\psi_{\s\t}$'s have cellular-like property. Therefore, $\set{\psi_{\s\t} | (\lambda,f) \in \widehat B_n, \s,\t \in \Tud_n(\lambda)}$ is a potential cellular basis of $\G{n}$. We will prove this result in Section~\ref{sec:rela:main}.
\end{Remark}

\section{A generating set of $\B$} \label{sec:gene}

In the rest of this paper, we are going to prove that $\G{n} \cong \B$. The first step of the proof is to define
$$
G_n(\delta) = \set{e(\bi)|\bi \in P^n} \cup \set{y_k| 1 \leq k \leq n} \cup \set{\psi_k| 1\leq k \leq n-1} \cup \set{\epsilon_k| 1\leq k \leq n-1}
$$
in $\B$ and show that $G_n(\delta)$ generates $\B$. Here we abuse the symbols and use $e(\bi)$, $y_k$, $\psi_r$ and $\epsilon_r$ as elements in both $\G{n}$ and $\mathscr B_n(\delta)$. Then we construct a mapping from $\G{n}$ to $\B$ by sending generators to generators and show this mapping is a surjective homomorphism. In this section, we construct the elements $e(\bi)$, $y_k$, $\psi_r$ and $\epsilon_r$ in $\B$ and show these elements generate $\B$.

First we recall the definitions and notations we need for the rest of the paper, which have been introduced in Section~\ref{sec:semi}. Recall $R$ is a field with characteristic $0$ and fix $\delta \in R$. Define $\F = R(x)$ to be the rational field with indeterminate $x$ and $\O = R[x]_{(x - \delta)} = R[[x - \delta]]$. Let $\m = (x - \delta)\O \subset \O$. Then $\m$ is a maximal ideal of $\O$ and $R \cong \O/\m$.

Let $\Bx$ and $\BOx$ be the Brauer algebras over $\F$ and $\O$, respectively. Then $\Bx = \BOx \otimes_\O \F$ and $\B \cong \BOx \otimes_\O R \cong \BOx/(x - \delta)\BOx$. In order to avoid confusion we will write the generators of $\BOx$ and $\Bx$ as $s_k^\O$ and $e_k^\O$ and generators of $\B$ as $s_k$ and $e_k$. Hence for any element $w \in \B$, we write $w^\O = w \otimes_R 1_\O \in \BOx$, so that $w = w^\O \otimes_\O 1_R$.

Because $\B \cong \BOx\otimes_\O R \cong \BOx/(x-\delta)\BOx$, if $x,y \in \BOx$ and we have $x \equiv y \pmod{(x-\delta)\BOx}$, then $x\otimes_\O 1_R = y\otimes_\O 1_R$ as elements of $\B$. This observation will give us a way to extend the results of $\BOx$ to $\B$.

\subsection{Gelfand-Zetlin algebra of $\B$} \label{sec:L}

Following Okounkov-Vershik~\cite{OV:GelfandZetlin}, define \textit{Gelfand-Zetlin subalgebra} of $\B$ to be the algebra $\mathscr L_n$ generated by $L_1, L_2, \ldots, L_n$. By the definition one can see that $\mathscr L_n$ is a commutative subalgebra of $\B$. Similarly we define $\mathscr L_n(\O)$ in $\BOx$. In this subsection we define the idempotents $e(\bi)$ with $\bi \in I^n$ and nilpotency elements $y_k$ with $1 \leq k \leq n$ in $\B$.

Let $M$ be a finite dimensional $\mathscr B_n(\delta)$-module. Similarly as in Brundan and Kleshchev~\cite[Section 3.1]{BK:GradedKL}, the eigenvalues of each $L_k$ on $M$ belongs to $P$. So $M$ decomposes as the direct sum $M = \bigoplus_{\bi \in P^n} M_\bi$ of weight spaces
$$
M_\bi = \set{v \in M | \text{$(L_k - i_k)^N v = 0$ for all $k = 1, 2, \ldots, n$ and $N \gg 0$}}.
$$

We deduce that there is a system $\set{e(\bi) | \bi \in P^n}$ of mutually orthogonal idempotents in $\mathscr B_n(\delta)$ such that $M e(\bi) = M_\bi$ for each finite dimensional module $M$. In fact, $e(\bi)$ lies in $\mathscr L_n$.

Hu-Mathas~\cite[Proposition 4.8]{HuMathas:GradedCellular} proved the following result in the cyclotomic Hecke algebras. Their result can be directly extended to $\B$ following the same proof.

\begin{Lemma}[\protect{Hu-Mathas~\cite[Proposition 4.8]{HuMathas:GradedCellular}}] \label{idem-semi}
Suppose that $e(\bi) \neq 0$ for some $\bi \in P^n$ and let
$$
e(\bi)^{\mathscr O} = \sum_{\t \in \mathscr T^{ud}_n(\bi)} \f{\t} \in \mathscr B_n^\F(x).
$$

Then $e(\bi)^{\mathscr O} \in \mathscr B_n^{\mathscr O}(x)$ and $e(\bi) = e(\bi)^{\mathscr O} \otimes_{\mathscr O} 1_R$.
\end{Lemma}

By~\autoref{idem-semi}, it is straightforward that $e(\bi) \neq 0$ only if $\bi \in I^n$, i.e. $\bi$ is the residue sequence of an up-down tableau. Therefore, by defining $e(\bi) = 0$ for $\bi \not\in I^n$, we construct a set of orthogonal elements $\set{e(\bi)^\O | \bi \in P^n}$ in $\BOx$ and $\set{e(\bi) | \bi \in P^n}$ in $\B$, such that $\sum_{\bi \in P^n} e(\bi)^\O = \sum_{\bi \in I^n} e(\bi)^\O = 1_\O$ and $\sum_{\bi \in P^n} e(\bi) = \sum_{\bi \in I^n} e(\bi) = 1_R$ by~\autoref{idem:semi:1} and~\autoref{idem-semi}.


For an integer $k$ with $1 \leq k \leq n$, define $y_k^\O := \sum_{\bi \in I^n} (L_k^\O - i_k) e(\bi)^\O \in \mathscr L_n(\O)$ and
\begin{equation} \label{def:y}
y_k := y_k^\O \otimes_\O 1_R = \sum_{\bi \in I^n} (L_k - i_k) e(\bi) \in \mathscr L_n.
\end{equation}

By~\autoref{B:basis}, one can see that $y_k$ is nilpotent in $\B$. Hence for any polynomial $\phi(L_1, L_2, \ldots, L_n) e(\bi) \in \mathscr L_n$ with $\phi(i_1, i_2, \ldots, i_n) \neq 0$, we can define a element $\phi(L_1, \ldots, L_n)^{-1} \in \mathscr L_n$ such that
$$
\phi(L_1, L_2, \ldots, L_n) \phi(L_1, L_2, \ldots, L_n)^{-1}e(\bi) = \phi(L_1, L_2, \ldots, L_n)^{-1} \phi(L_1, L_2, \ldots, L_n)e(\bi) = e(\bi).
$$

Suppose $\phi(x_1, \ldots, x_n) \in R[x_1, \ldots, x_n]$ is a polynomial. Define $R^\times = \set{r \in R | r \neq 0}$. The next result is the extended version of Hu-Mathas~\cite[Proposition 4.6]{HuMathas:SemiQuiver}, followed by the same proof.

\begin{Lemma} \label{L:inv:1}
Suppose $\bi \in I^n$ and $\phi(x_1, \ldots, x_n) \in R[x_1, \ldots, x_n]$ is a polynomial. If $\phi(i_1, \ldots, i_n) \in R^\times$, we have
$$
\sum_{\t \in \mathscr T^{ud}_n(\bi)} \frac{1}{\phi(c_\t(1), \ldots, c_\t(n))} \f{\t} \in \mathscr L_n(\O).
$$
\end{Lemma}

\begin{Lemma} \label{L:inv:2}
Suppose $\bi \in I^n$ and $\phi(x_1, \ldots, x_n) \in R(x_1, \ldots, x_n)$ is a rational function. If $\phi(i_1, \ldots, i_n) \in R$, we have
$$
\sum_{\t \in \mathscr T^{ud}_n(\bi)} \phi(c_\t(1), \ldots, c_\t(n)) \f{\t} \in \mathscr L_n(\O).
$$
\end{Lemma}

\proof Because $\phi$ is a rational function, there are two polynomials $\phi_1$ and $\phi_2$ such that $\phi = \phi_1 / \phi_2$. It is obvious that $\phi(i_1, \ldots, i_n) \in R$ if and only if $\phi_2(i_1, \ldots, i_n) \in R^\times$. Hence by~\autoref{L:inv:1}, we have
$$
\sum_{\t \in \mathscr T^{ud}_n(\bi)} \frac{1}{\phi_2(c_\t(1), \ldots, c_\t(n))} \f{\t} \in \L_n(\O).
$$

Hence, because $\phi_1$ is a polynomial, we have $\phi_1(L_1^\O,\ldots,L_n^\O) \in \L_n(\O)$. So
$$
\phi_1(L_1^\O,\ldots,L_n^\O) \sum_{\t \in \mathscr T^{ud}_n(\bi)} \frac{1}{\phi_2(c_\t(1), \ldots, c_\t(n))} \f{\t} = \sum_{\t \in \mathscr T^{ud}_n(\bi)} \phi(c_\t(1), \ldots, c_\t(n)) \f{\t} \in \L_n(\O),
$$
which completes the proof.
\endproof

Let $\phi$ be a polynomial in $R[x_1, \ldots, x_n]$ satisfying the assumptions of~\autoref{L:inv:1}. Then
$$
\phi(L_1^\O, \ldots, L_n^\O) \sum_{\t \in \mathscr T^{ud}_n(\bi)} \frac{1}{\phi(c_\t(1), \ldots, c_\t(n))} \f{\t} = e(\bi)^\O = \sum_{\t \in \mathscr T^{ud}_n(\bi)} \frac{1}{\phi(c_\t(1), \ldots, c_\t(n))} \f{\t} \phi(L_1^\O, \ldots, L_n^\O).
$$

Abusing notations, in this situation we write
$$
\frac{1}{\phi(L_1^\O,\ldots,L_n^\O)} e(\bi)^\O = e(\bi)^\O \frac{1}{\phi(L_1^\O,\ldots,L_n^\O)} = \sum_{\t \in \mathscr T^{ud}_n(\bi)} \frac{1}{\phi(c_\t(1), \ldots, c_\t(n))} \f{\t} \in \L_n(\O).
$$

Similarly, let $\phi$ be a rational function in $R(x_1, \ldots, x_n)$ satisfying the assumptions of~\autoref{L:inv:2}. Then we write
$$
\phi(L_1^\O, \ldots, L_n^\O) e(\bi)^\O = e(\bi)^\O \phi(L_1^\O, \ldots, L_n^\O) = \sum_{\t \in \mathscr T^{ud}_n(\bi)} \phi(c_\t(1), \ldots, c_\t(n)) \f{\t}.
$$

\subsection{Modification terms $P_k(\bi)$, $Q_k(\bi)$ and $V_k(\bi)$} \label{sec:PQV}

In this subsection we define three elements $P_k(\bi)$, $Q_k(\bi)$ and $V_k(\bi)$ in $\B$, which are essential when we define $\psi_k$ and $\epsilon_k$ in $\B$. These terms are defined so that the actions of $\psi_k$ and $\epsilon_k$'s on seminormal forms of $\B$ are well-behaved (cf.~\autoref{semi:1} and~\autoref{semi:2}).

Let $x_1, \ldots, x_n$ be invariants. For each $\bi \in I^n$ and rational function $\phi$, we say $\phi(x_1,\ldots,x_{k-1},\frac{x_k - x_{k+1}}{2}) \in R(x_1,\ldots,x_{k+1})$ is \textit{invertible over $\bi$} if $\phi(i_1,i_2,\ldots,i_{k-1},(i_k - i_{k+1})/2) \in R^\times$. It is obvious by the definition that $\phi_1, \phi_2 \in R(x_1, \ldots, x_{k+1})$ invertible over $\bi$ implies $\phi_1{\cdot}\phi_2\in R(x_1, \ldots, x_{k+1})$ invertible over $\bi$.

First we define the elements $P_k(\bi)$ and $Q_k(\bi)$. For $1 \leq k \leq n-1$ and $1 \leq r \leq k-1$, define
\begin{align*}
L_{k,r} & = \left\{\frac{x_k - x_{k+1}}{2} + x_r + 1, \frac{1}{\frac{x_k - x_{k+1}}{2} - x_r + 1}, -(\frac{x_k - x_{k+1}}{2} - x_r), -(\frac{1}{\frac{x_k - x_{k+1}}{2} + x_r})\right\} \subset R(x_1,\ldots,x_{k+1}),\\
R_{k,r} & = \left\{- (\frac{x_k - x_{k+1}}{2} + x_r - 1), -\frac{1}{\frac{x_k - x_{k+1}}{2} - x_r - 1}, \frac{x_k - x_{k+1}}{2} - x_r, \frac{1}{\frac{x_k - x_{k+1}}{2} + x_r}\right\} \subset R(x_1,\ldots,x_{k+1}),\\
S_k & = \left\{- \frac{1}{\frac{x_k - x_{k+1}}{2} - x_1}, -\frac{1}{x_k - x_{k+1}}, x_k - x_{k+1} + 1\right\} \subset R(x_1,x_k, x_{k+1}), \\
T_k & = \left\{\frac{x_k - x_{k+1}}{2} + x_1 \right\} \subset R(x_1,x_k, x_{k+1}).
\end{align*}

For any $\bi \in P^n$, we define
\begin{align*}
L_{k,r}(\bi) & = \set{w\in L_{k,r}|\text{$w$ is invertible over $\bi$}}, &
S_k(\bi) & = \set{w\in S_k|\text{$w$ is invertible over $\bi$}},\\
R_{k,r}(\bi) & = \set{w\in R_{k,r}|\text{$w$ is invertible over $\bi$}}, &
T_k(\bi) & = \set{w\in T_k|\text{$w$ is invertible over $\bi$}}.
\end{align*}

Let
\begin{align*}
P_k^\bi (x_1,\ldots,x_{k-1},\frac{x_k - x_{k+1}}{2}) & := \prod_{w \in S_k(\bi)}w \prod_{r = 1}^{k-1} \left( \prod_{w \in L_{k,r}(\bi)}w \right) \in R(x_1,\ldots,x_{k+1}), \\
Q_k^\bi (x_1,\ldots,x_{k-1},\frac{x_k - x_{k+1}}{2}) & := \prod_{w \in T_k(\bi)}w \prod_{r = 1}^{k-1} \left( \prod_{w \in R_{k,r}(\bi)}w \right) \in R(x_1,\ldots,x_{k+1}).
\end{align*}

By the definitions, one can see that $P_k^\bi, Q_k^\bi \in R(x_1, \ldots, x_{k+1})$ are invertible over $\bi$. Hence
\begin{equation} \label{PQ:eq1}
\begin{array} {ll}
P_k^\bi(i_1, \ldots, i_{k-1}, \frac{i_k - i_{k+1}}{2}) \in R^\times, & Q_k^\bi(i_1, \ldots, i_{k-1}, \frac{i_k - i_{k+1}}{2}) \in R^\times, \\
P_k^\bi(i_1, \ldots, i_{k-1}, \frac{i_k - i_{k+1}}{2})^{-1} \in R^\times, & Q_k^\bi(i_1, \ldots, i_{k-1}, \frac{i_k - i_{k+1}}{2})^{-1} \in R^\times.
\end{array}
\end{equation}

For any $\t \in \Tud_n(\bi)$, we define $P_k(\t), Q_k(\t) \in \F$ by
\begin{align*}
P_k(\t) & = P_k^{\bi}(c_\t(1),\ldots,c_\t(k-1),(c_\t(k) - c_\t(k+1))/2),\\
Q_k(\t) & = Q_k^{\bi}(c_\t(1),\ldots,c_\t(k-1),(c_\t(k) - c_\t(k+1))/2).
\end{align*}

Define $P_k^\O(\bi), Q_k^\O(\bi), P_k^\O(\bi)^{-1}, Q_k^\O(\bi)^{-1} \in \Bx$ by
\begin{equation} \label{PQ:eq2}
\begin{array} {ll}
P_k^\O (\bi) = \sum_{\t \in \mathscr T^{ud}_n(\bi)} P_k(\t) \f{\t}, & Q_k^\O (\bi) = \sum_{\t \in \mathscr T^{ud}_n(\bi)} Q_k(\t) \f{\t},\\
P_k^\O (\bi)^{-1} = \sum_{\t \in \mathscr T^{ud}_n(\bi)} \frac{1}{P_k(\t)} \f{\t}, & Q_k^\O (\bi)^{-1} = \sum_{\t \in \mathscr T^{ud}_n(\bi)} \frac{1}{Q_k(\t)} \f{\t},
\end{array}
\end{equation}
and by (\ref{PQ:eq1}) and~\autoref{L:inv:2}, we have $P_k^\O(\bi), Q_k^\O(\bi), P_k^\O(\bi)^{-1}, Q_k^\O(\bi)^{-1} \in \L_n(\O)$.

\begin{Definition} \label{def:PQ}
Suppose $\bi \in I^n$ and $1 \leq k \leq n-1$. We define elements $P_k(\bi), Q_k(\bi), P_k(\bi)^{-1}, Q_k(\bi)^{-1} \in \mathscr L_n$ by
\begin{align*}
& P_k (\bi) = P_k^\O (\bi) \otimes_\O 1_R, && Q_k (\bi) = Q_k^\O (\bi) \otimes_\O 1_R,\\
& P_k (\bi)^{-1} = P_k^\O (\bi)^{-1} \otimes_\O 1_R, && Q_k (\bi)^{-1} = Q_k^\O (\bi)^{-1} \otimes_\O 1_R.
\end{align*}
\end{Definition}

The next Lemma is directly implied by (\ref{PQ:eq2}).

\begin{Lemma} \label{PQ:9}
Suppose $\bi \in I^n$ and $1 \leq k \leq n-1$. We have
\begin{align*}
& P_k(\bi) P_k(\bi)^{-1} = P_k(\bi)^{-1} P_k(\bi) = e(\bi), \\
& Q_k(\bi) Q_k(\bi)^{-1} = Q_k(\bi)^{-1} Q_k(\bi) = e(\bi).
\end{align*}
\end{Lemma}

Next, for each $1 \leq k \leq n-1$ and $\bi \in I^n$ with $i_k = i_{k+1}$, we define element $V_k(\bi)$. First we define a rational function
$$
V_k^\bi (x_1,\ldots,x_{k-1},\frac{x_k - x_{k+1}}{2}) := \frac{P_k^\bi Q_k^\bi - 1}{x_k - x_{k+1}} \in R(x_1,\ldots,x_{k+1}).
$$

The next Lemma shows that $V_k^\bi$ is invertible over $\bi$.

\begin{Lemma} \label{V:1}
Suppose $1 \leq k \leq n-1$ and $\bi = (i_1, \ldots, i_n) \in I^n$ with $i_k = i_{k+1}$. Then we have $V_k^\bi(i_1, \ldots, i_{k-1}, \frac{i_k - i_{k+1}}{2}) \in R$.
\end{Lemma}

\proof As $i_k = i_{k+1}$, we have $\frac{i_k - i_{k+1}}{2} = 0$. Because all the factors of $P_k^\bi$ and $Q_k^\bi$ are invertible over $\bi$, we can write
$$
P_k^\bi Q_k^\bi = \sum_{i = 0}^\infty c_i(x_1, \ldots, x_{k-1}) \left( \frac{x_k - x_{k+1}}{2} \right)^i,
$$
where $c_i(x_1, \ldots, x_{k-1}) \in R(x_1, \ldots, x_{k-1})$ and $c_i(i_1, \ldots, i_{k-1}) \in R$. Hence
$$
V_k^\bi(i_1,\ldots,i_{k-1},\frac{x_k - x_{k+1}}{2}) = \frac{c_1(i_1, \ldots, i_{k-1}) - 1}{x_k - x_{k+1}} + \sum_{i = 0}^\infty \frac{c_{i+1}(i_1, \ldots, i_{k-1})}{2} \left( \frac{x_k - x_{k+1}}{2} \right)^i,
$$
and $V_k^\bi(i_1, \ldots, i_{k-1}, \frac{i_k - i_{k+1}}{2}) \in R$ if and only if $c_1(i_1, \ldots, i_{k-1}) = 1$. By the definitions of $P_k^\bi$ and $Q_k^\bi$, we have
$$
c_1(i_1, \ldots, i_{k-1}) = P_k^\bi(i_1, \ldots, i_{k-1}, 0) Q_k^\bi(i_1, \ldots, i_{k-1}, 0) = 1,
$$
which completes the proof. \endproof

Suppose $\t \in \mathscr T^{ud}_n(\bi)$. Define $V_k(\t) = V_k^\bi(c_\t(1), \ldots, c_\t(k-1), \frac{c_\t(k) - c_\t(k+1)}{2})$. We define $V_k^\O(\bi) \in \Bx$ by
$$
V_k^\O(\bi) = \sum_{\t \in \mathscr T^{ud}_n(\bi)} V_k(\t) \f{\t},
$$
and by~\autoref{L:inv:2} and~\autoref{V:1}, we have $V_k^\O(\bi) \in \L_n(\O)$.

\begin{Definition} \label{def:V}
Suppose $1 \leq k \leq n-1$ and $\bi = (i_1, \ldots, i_n) \in I^n$ with $i_k = i_{k+1}$. We define the element $V_k(\bi) \in \mathscr L_n$ by $V_k(\bi):= V_k^\O(\bi) \otimes_\O 1_R$.
\end{Definition}

The next Lemma shows the connections of $V_k(\t)$ with $P_k(\t)$ and $Q_k(\t)$.

\begin{Lemma} \label{V:2}
Suppose $1 \leq k \leq n-1$ and $\bi = (i_1, \ldots, i_n) \in I^n$ with $i_k = i_{k+1}$. For any $\t \in \Tud_n(\bi)$, we have $P_k(\t) Q_k(\t) = (c_\t(k) - c_t(k+1)) V_k(\t) + 1$.
\end{Lemma}

\proof By the definition of $V_k^\bi$, we have $P_k^\bi Q_k^\bi = (x_k - x_{k+1})V_k^\bi + 1$. Hence the Lemma follows straightforward. \endproof


In the rest of this subsection we introduce some of the properties of $P_k(\bi)$, $Q_k(\bi)$ and $V_k(\bi)$, or more precisely, properties of $P_k(\t)$, $Q_k(\t)$ and $V_k(\t)$ for up-down tableaux $\t$. These results will be used frequently in the rest of this paper when we derive the relations of the generators of $\B$. These properties make the actions of the $\psi_k$'s and $\epsilon_k$'s on the seminormal basis well-behaved and they are the reasons why we define the modification terms in such a way.

For any rational function $w \in R(x_1, \ldots, x_n)$, we denote $w(\t) = w(c_\t(1), \ldots, c_\t(n))$.

\begin{Lemma} \label{PQ:1}
Suppose $\bi = (i_1, \ldots, i_n) \in I^n$ and $\t \in \Tud_n(\bi)$. For $1 \leq k \leq n-1$, if $\u = \t{\cdot}s_k$ exists, then we have
$$
P_k(\t)^{-1} Q_k(\u)^{-1} =
\begin{cases}
 \frac{1}{1 - c_\u(k) + c_\u(k+1)}, & \text{if $i_{k+1} = i_k$,}\\
 c_\u(k) - c_\u(k+1), & \text{if $i_{k+1} = i_k - 1$,}\\
 \frac{c_\u(k) - c_\u(k+1)}{1 - c_\u(k) + c_\u(k+1)}, & \text{if $i_{k+1} \neq i_k,i_k - 1$.}
\end{cases}
$$
\end{Lemma}

\proof Because $\t = \u{\cdot}s_k$, we have $c_\u(r) = c_\t(r)$ for $1 \leq r \leq k-1$ and $c_\u(k) - c_\u(k+1) = -(c_\t(k) - c_\t(k+1)$. Hence for any $1 \leq \l \leq k-1$, we have
$$
\left( \prod_{w \in L_{k,r}(\bi{\cdot}s_k)}w(\t)\right) \left( \prod_{w \in R_{k,r}(\bi)}w(\u) \right) = 1.
$$

Therefore, by the definition of $P_k(\t)$ and $Q_k(\u)$, we have
\begin{align*}
P_k(\t) Q_k(\u) & = \prod_{w \in S_k(\bj)}w(\t)
    \prod_{w \in T_k(\bi)}w(\u)
=
\begin{cases}
1 - c_\u(k) + c_\u(k+1), & \text{if $i_{k+1} = i_k$,}\\
\frac{1}{c_\u(k) - c_\u(k+1)}, & \text{if $i_{k+1} = i_k - 1$,}\\
\frac{1 - c_\u(k) + c_\u(k+1)}{c_\u(k) - c_\u(k+1)}, & \text{otherwise.}
\end{cases}
\end{align*}

Hence the Lemma follows. \endproof

\begin{Lemma} \label{PQ:2}
Suppose $\bi \in I^n$ and $\t \in \Tud_n(\bi)$. If $\t(k-1) = \t(k+1) = -\t(k)$, we have
$$
Q_k(\t) P_{k-1}(\t) = P_k(\t) Q_{k-1}(\t) = 1,
$$
for $2 \leq k \leq n-1$.
\end{Lemma}

\proof By the construction of $\t$, we have $c_\t(k-1) = c_\t(k+1) = -c_\t(k)$. Hence we have $c_\t(k) = (c_\t(k) - c_\t(k+1))/2$ and $c_\t(k-1) = (c_\t(k-1) - c_\t(k))/2$. Because $c_\t(k-1) = -c_\t(k)$, for $1 \leq r \leq k-2$, we have
$$
\left( \prod_{w \in L_{k,r}(\bi)}w(\t)\right) \left( \prod_{w \in R_{k-1,r}(\bi)}w(\t) \right) = 1.
$$

Hence, by the definition of $P_k(\t)$ and $Q_{k-1}(\t)$, we have
$$
P_k(\t) Q_{k-1}(\t) = \prod_{w \in S_k(\bi)}w(\t) \prod_{w \in T_{k-1} (\bi)}w(\t) \left( \prod_{w \in L_{k,k-1}(\bi)}w(\t)\right).
$$

Then $P_k(\t) Q_{k-1}(\t)$ is the product of the non-invertible elements of
$$
\left\{ -\frac{1}{c_\t(k) - c_\t(1)}, -\frac{1}{2c_\t(k)}, 2c_\t(k) + 1, 1, \frac{1}{2c_\t(k) + 1}, -2c_\t(k), -c_\t(k) + c_\t(1) \right\}.
$$

By direct calculation, we have $P_k(\t) Q_{k-1}(\t) = 1$. $Q_k(\t) P_{k-1}(\t) = 1$ follows by the similar argument. \endproof

Recall that
$$
a_k(\bi) =
\begin{cases}
\#\set{1 \leq r \leq k-1 | i_r \in \{-1,1\} } + 1 + \delta_{\frac{i_k - i_{k+1}}{2}, \frac{\delta-1}{2}}, & \text{if $(i_k - i_{k+1})/2 = 0$,}\\
\#\set{1 \leq r \leq k-1 | i_r \in \{-1,1\} } + \delta_{\frac{i_k - i_{k+1}}{2}, \frac{\delta-1}{2}}, & \text{if $(i_k - i_{k+1})/2 = 1$,}\\
\delta_{\frac{i_k - i_{k+1}}{2}, \frac{\delta-1}{2}}, & \text{if $(i_k - i_{k+1})/2 = 1/2$,}\\
\#\set{1 \leq r \leq k-1 | i_r \in \{\frac{i_k - i_{k+1}}{2}, \frac{i_k - i_{k+1}}{2} - 1, -\frac{i_k - i_{k+1}}{2}, -\frac{i_k - i_{k+1}}{2} + 1\} } + \delta_{\frac{i_k - i_{k+1}}{2}, \frac{\delta-1}{2}}, & \text{otherwise.}
\end{cases}
$$

For $1 \leq k \leq n-1$ and $1 \leq r \leq k-1$, define
\begin{align*}
\widehat L_{k,r} & = \left\{\frac{x_k - x_{k+1}}{2} + x_r + 1, \frac{1}{\frac{x_k - x_{k+1}}{2} - x_r + 1}, \frac{x_k - x_{k+1}}{2} - x_r, \frac{1}{\frac{x_k - x_{k+1}}{2} + x_r}\right\} \subset R(x_1,\ldots,x_{k+1}),\\
\widehat R_{k,r} & = \left\{\frac{x_k - x_{k+1}}{2} + x_r - 1, \frac{1}{\frac{x_k - x_{k+1}}{2} - x_r - 1}, \frac{x_k - x_{k+1}}{2} - x_r, \frac{1}{\frac{x_k - x_{k+1}}{2} + x_r}\right\} \subset R(x_1,\ldots,x_{k+1}),\\
\widehat S_k & = \left\{\frac{1}{\frac{x_k - x_{k+1}}{2} - x_1}, \frac{1}{x_k - x_{k+1}}, x_k - x_{k+1} + 1\right\} \subset R(x_1,x_k, x_{k+1}), \\
\widehat T_k & = \left\{\frac{x_k - x_{k+1}}{2} + x_1 \right\} \subset R(x_1,x_k, x_{k+1}).
\end{align*}

For any residue sequence $\bi \in I^n$, define
\begin{align*}
\widehat L_{k,r}(\bi) & = \set{w\in \widehat L_{k,r}|\text{$w$ is invertible over $\bi$}}, &
\widehat S_k(\bi) & = \set{w\in \widehat S_k|\text{$w$ is invertible over $\bi$}},\\
\widehat R_{k,r}(\bi) & = \set{w\in \widehat R_{k,r}|\text{$w$ is invertible over $\bi$}}, &
\widehat T_k(\bi) & = \set{w\in \widehat T_k|\text{$w$ is invertible over $\bi$}}.
\end{align*}

For any $\t \in \Tud_n(\bi)$, let
$$
\widehat P_k(\t) := \prod_{w \in \widehat S_k(\bi)}w(\t) \prod_{r = 1}^{k-1} \left( \prod_{w \in \widehat L_{k,r}(\bi)}w(\t) \right), \qquad \text{and} \qquad
\widehat Q_k(\t) := \prod_{w \in \widehat T_k(\bi)}w(\t) \prod_{r = 1}^{k-1} \left( \prod_{w \in \widehat R_{k,r}(\bi)}w(\t) \right).
$$

\begin{Lemma} \label{PQ:3:h}
Suppose $\bi \in I^n$ and $1 \leq k \leq n-1$. For any $\t \in \Tud_n(\bi)$, we have $P_k(\t) Q_k(\t) = (-1)^{a_k(\bi)} \widehat P_k(\t) \widehat Q_k(\t)$.
\end{Lemma}

\proof For $1 \leq r \leq k-1$, define $b_{k,r}(\bi) =
\delta_{i_r, \frac{i_k - i_{k+1}}{2}} + \delta_{i_r, \frac{i_k - i_{k+1}}{2}-1} + \delta_{i_r, -\frac{i_k - i_{k+1}}{2}} + \delta_{i_r, -\frac{i_k - i_{k+1}}{2}+1}$. By comparing $a_k(\bi)$ and $\sum_{r=1}^{k-1}b_{k,r}(\bi)$, it is easy to see that
\begin{equation} \label{PQ:3:h:eq1}
(-1)^{a_k(\bi)} = (-1)^{\sum_{r=1}^{k-1}b_{k,r}(\bi) + \delta_{\frac{i_k - i_{k+1}}{2}, 0} + \delta_{\frac{i_k - i_{k+1}}{2}, \frac{\delta-1}{2}}}.
\end{equation}

By the definitions of $L_{k,r}(\bi), \widehat L_{k,r}(\bi), R_{k,r}(\bi)$ and $\widehat R_{k,r}(\bi)$, for any $1 \leq r \leq k-1$, we have
$$
\prod_{w \in L_{k,r}(\bi)}w(\t)\prod_{w \in R_{k,r}(\bi)}w(\t) = (-1)^{b_{k,r}(\bi)}\prod_{w \in \widehat L_{k,r}(\bi)}w(\t)\prod_{w \in \widehat R_{k,r}(\bi)}w(\t),
$$
which implies
\begin{equation} \label{PQ:3:h:eq2}
\prod_{r = 1}^{k-1} \left(\prod_{w \in L_{k,r}(\bi)}w(\t)\prod_{w \in R_{k,r}(\bi)}w(\t) \right) = (-1)^{\sum_{r=1}^{k-1} b_{k,r}(\bi)}\prod_{r = 1}^{k-1} \left(\prod_{w \in \widehat L_{k,r}(\bi)}w(\t)\prod_{w \in \widehat R_{k,r}(\bi)}w(\t) \right).
\end{equation}

By the definitions of $S_k(\bi), \widehat S_k(\bi), T_k(\bi)$ and $\widehat T_k(\bi)$, we have
\begin{equation} \label{PQ:3:h:eq3}
\prod_{w \in S_k(\bi)}w(\t) \prod_{w \in T_k(\bi)}w(\t) = (-1)^{\delta_{\frac{i_k - i_{k+1}}{2}, 0} + \delta_{\frac{i_k - i_{k+1}}{2}, \frac{\delta-1}{2}}} \prod_{w \in \widehat S_k(\bi)}w(\t) \prod_{w \in \widehat T_k(\bi)}w(\t).
\end{equation}

Combining (\ref{PQ:3:h:eq2}) and (\ref{PQ:3:h:eq3}), we have
$$
P_k(\t) Q_k(\t) = (-1)^{\sum_{r=1}^{k-1}b_{k,r}(\bi) + \delta_{\frac{i_k - i_{k+1}}{2}, 0} + \delta_{\frac{i_k - i_{k+1}}{2}, \frac{\delta-1}{2}}} \widehat P_k(\t) \widehat Q_k(\t).
$$

Hence the Lemma follows by (\ref{PQ:3:h:eq1}). \endproof

\begin{Lemma} \label{PQ:3}
Suppose $\bi = (i_1, \ldots, i_n) \in I^n$ with $i_k + i_{k+1} = 0$ for $1 \leq k \leq n-2$ and $\t \in \Tud_n(\bi)$ with $\t(k) + \t(k+1) = 0$. Then we have
$$
P_k(\t) Q_k(\t) =
\begin{cases}
(-1)^{a_k(\bi)} {\cdot} e_k(\t,\t), & \text{if $\bi \in I_{k,0}^n$,}\\
(-1)^{a_k(\bi)} {\cdot} 2(c_\t(k) - i_k) e_k(\t,\t), & \text{if $\bi \in I_{k,-}^n$,}\\
(-1)^{a_k(\bi)} {\cdot} \frac{1}{2(c_\t(k) - i_k)}e_k(\t,\t), & \text{if $\bi \in I_{k,+}^n$.}
\end{cases}
$$
\end{Lemma}

\proof By~\autoref{Naz}, we have
\begin{equation} \label{PQ:3:eq1}
(2u+1) \sum_{\substack{\s \overset{k}\sim \t\\ \s \neq \t}} \frac{u + c_\s(k)}{u - c_\s(k)} = (2u+1) \frac{u - c_\t(k)}{u + c_\t(k)} \frac{u + c_\t(1)}{u - c_\t(1)} \prod_{r = 1}^{k-1} \frac{(u + c_\t(r))^2 - 1}{(u - c_\t(r))^2 - 1} \frac{(u - c_\t(r))^2}{(u + c_\t(r))^2}.
\end{equation}

Define $f(x_1, x_2, \ldots, x_{k-1}, \frac{x_k - x_{k+1}}{2})$ to be a rational function in $R(x_1, \ldots, x_{k+1})$ obtained by removing all factors which are non-invertible over $\bi$ from
\begin{equation} \label{PQ:3:eq4}
(x_k - x_{k+1} +1) \frac{1}{x_k - x_{k+1}} \frac{\frac{x_k - x_{k+1}}{2} + x_1}{\frac{x_k - x_{k+1}}{2} - x_1} \prod_{r = 1}^{k-1} \frac{(\frac{x_k - x_{k+1}}{2} + x_r) + 1}{(\frac{x_k - x_{k+1}}{2} - x_r) + 1} \frac{(\frac{x_k - x_{k+1}}{2} + x_r) - 1}{(\frac{x_k - x_{k+1}}{2} - x_r) - 1} \frac{(\frac{x_k - x_{k+1}}{2} - x_r)^2}{(\frac{x_k - x_{k+1}}{2} + x_r)^2}.
\end{equation}

Because $c_\t(k) + c_\t(k+1) = 0$, we have $\frac{c_\t(k) - c_\t(k+1)}{2} = c_\t(k)$. If $w$ is a factor of (\ref{PQ:3:eq4}) which is non-invertible over $\bi$, one can see that $w(\t) \in \{ 2(c_\t(k) - i_k), \frac{1}{2(c_\t(k) - i_k)}, 0\}$. Hence, by (\ref{PQ:3:eq1}) we have
\begin{equation} \label{PQ:3:eq2}
f(c_\t(1),\ldots, c_\t(k-1), \frac{c_\t(k) - c_\t(k+1)}{2}) = f(c_\t(1),\ldots, c_\t(k)) = (2(c_\t(k) - i_k))^\l e_k(\t,\t)
\end{equation}
for some $\l \in \Z$.

By the definition of $f$, we can see that $f(x_1,\ldots,x_{k-1},\frac{x_k - x_{k+1}}{2}) = \widehat P_k^\bi (x_1,\ldots,x_{k-1},\frac{x_k - x_{k+1}}{2}) \widehat Q_k^\bi (x_1,\ldots,x_{k-1},\frac{x_k - x_{k+1}}{2})$. Hence by~\autoref{PQ:3:h} and (\ref{PQ:3:eq2}), we have
\begin{equation} \label{PQ:3:eq3}
P_k(\t) Q_k(\t) = (-1)^{a_k(\bi)} f(c_\t(1), \ldots, c_\t(k)) = (-1)^{a_k(\bi)} (2(c_\t(k) - i_k))^\l e_k(\t,\t),
\end{equation}
where $\l \in \Z$.

By~\autoref{e_k:h1}, we have
$$
\sum_{\substack{\s \overset{k}\sim \t\\ \s \neq \t}} \frac{c_\t(k) + c_\s(k)}{c_\t(k) - c_\s(k)} = (2(c_\t(k) - i_k))^{|\mathscr{AR}_\lambda(-i_k)| - |\mathscr{AR}_\lambda(-i_k)| + 1} v,
$$
for some $v$ invertible in $\O$. Hence, by the definitions of $I_{k,0}^n, I_{k,-}^n, I_{k,+}^n$ and (\ref{deg:h4:eq1}) - (\ref{deg:h4:eq3}), we have
$$
e_k(\t,\t) = (2c_\t(k) + 1) \sum_{\substack{\s \overset{k}\sim \t\\ \s \neq \t}} \frac{c_\t(k) + c_\s(k)}{c_\t(k) - c_\s(k)} =
\begin{cases}
v, & \text{if $\bi \in I_{k,0}^n$,}\\
\frac{1}{2(c_\t(k) - i_k)} v, & \text{if $\bi \in I_{k,-}^n$,}\\
2(c_\t(k) - i_k) v, & \text{if $\bi \in I_{k,+}^n$,}
\end{cases}
$$
for some $v$ invertible in $\O$. Hence as $P_k(\t), Q_k(\t)$ are invertible in $\O$, by (\ref{PQ:3:eq3}) we complete the proof. \endproof

\begin{Lemma} \label{PQ:4}
Suppose $1 \leq r, k \leq n-1$ with $r < k-1$ and $\t$ is an up-down tableau. If $\s = \t{\cdot}s_r$ exists, we have $P_k(\t) = P_k(\s)$ and $Q_k(\t) = Q_k(\s)$.
\end{Lemma}

\proof Because $\s = \t{\cdot}s_r$, we have $\t(\l) = \s(\l)$ for $1 \leq \l \leq n$ and $\l \neq r,r+1$, and $\t(r) = \s(r+1)$, $\t(r+1) = \s(r)$. Let $\bi_\t$ and $\bi_\s$ be the residue sequences of $\t$ and $\s$, respectively. Hence, we have $\bi_\s = \bi_\t{\cdot}s_r$ because $\s = \t{\cdot}s_r$. By the construction, we have $\prod_{w\in L_{k,\l}(\bi_\t)} w(\t) = \prod_{w\in L_{k,\l}(\bi_\s)} w(\s)$ when $\l \neq r,r+1$, and
$$
\prod_{w\in L_{k,r}(\bi_\t)} w(\t) = \prod_{w\in L_{k,r+1}(\bi_\s)} w(\s), \qquad
\prod_{w\in L_{k,r+1}(\bi_\t)} w(\t) = \prod_{w\in L_{k,r}(\bi_\s)} w(\s).
$$

As $\s = \t{\cdot}s_r$ is an up-down tableau, we have $r > 1$. Hence we have $\prod_{w\in S_k(\bi_\t)} w(\t) = \prod_{w\in S_k(\bi_\s)} w(\s)$. Therefore, by combining the above results, we have $P_k(\t) = P_k(\s)$. We can prove $Q_k(\t) = Q_k(\s)$ following by the same process. \endproof

\begin{Lemma} \label{PQ:5}
Suppose $1 \leq r, k \leq n-1$ with $r < k-1$ and $\t$ is an up-down tableau with $\t(r) + \t(r+1) = 0$. For any $\s \overset{r}\sim \t$, we have $P_k(\t) = P_k(\s)$ and $Q_k(\t) = Q_k(\s)$.
\end{Lemma}

\proof Because $\s \overset{r}\sim \t$, we have $\t(\l) = \s(\l)$ for $1 \leq \l \leq n$ and $\l \neq r, r+1$, and $c_\s(r) + c_\s(r+1) = c_\t(r) + c_\t(r+1) = 0$. Let $\bi_\t$ and $\bi_\s$ be the residue sequences of $\t$ and $\s$, respectively. By the construction, when $\l \neq r, r+1$, we have $\prod_{w\in L_{k,\l}(\bi_\t)} w(\t) = \prod_{w\in L_{k,\l}(\bi_\s)} w(\s)$. As $c_\s(r) + c_\s(r+1) = c_\t(r) + c_\t(r+1) = 0$, we have
$$
\prod_{w\in L_{k,r}(\bi_\t)} w(\t) \prod_{w\in L_{k,r+1}(\bi_\t)} w(\t) = \prod_{w\in L_{k,r}(\bi_\s)} w(\s) \prod_{w\in L_{k,r+1}(\bi_\s)} w(\s) = 1.
$$

As $c_\t(1) = c_\s(1)$, we have $\prod_{w\in S_k(\bi_\t)} w(\t) = \prod_{w\in S_k(\bi_\s)} w(\s)$. Therefore, by combining the above results, we have $P_k(\t) = P_k(\s)$. We can prove $Q_k(\t) = Q_k(\s)$ following by the same process. \endproof

The next Lemma is used to prove~\autoref{V:3}, and also will be used in the next section when we prove the essential commutation relations of $\B$.

\begin{Lemma} \label{esscom:1}
Suppose $\bi \in I^n$ with $i_k + i_{k+1} = 0$ for $1 \leq k \leq n-1$ and $h_k(\bi) = 0$. Then $(-1)^{a_k(\bi)} = 1$ when $i_k = 0$ and $(-1)^{a_k(\bi) + a_k(\bi{\cdot}s_k)} = 1$ when $i_k \neq 0$.
\end{Lemma}

\proof Suppose $i_k = 0$. We have $a_k(\bi) = \#\set{1 \leq r \leq k-1 | i_r \in \{-1,1\}} + 1 + \delta_{\frac{i_k - i_{k+1}}{2}, \frac{\delta - 1}{2}}$ by the definition of $a_k(\bi)$. Suppose $\t \in \Tud_n(\bi)$ and write $\t_{k-1} = \lambda$. By (\ref{deg:h4:eq3}), we have $|\mathscr{AR}_\lambda(0)| = 1$. Hence by~\autoref{deg:h4:4}, we have
$$
\begin{cases}
\text{$\#\set{1 \leq r \leq k-1 | i_r \in \{-1,1\}}$ is odd if $\frac{\delta - 1}{2} \neq 0$,}\\
\text{$\#\set{1 \leq r \leq k-1 | i_r \in \{-1,1\}}$ is even if $\frac{\delta - 1}{2} = 0$,}
\end{cases}
$$
which implies that $a_k(\bi)$ is even. Hence $(-1)^{a_k(\bi)} = 1$ when $i_k = 0$.

Suppose $i_k \neq 0$. Because $i_k - i_{k+1} = 2i_k = -2i_{k+1}$, we have
\begin{align*}
a_k(\bi) & = \#\set{1 \leq r \leq k-1 | i_r \in \{i_k, i_k - 1, i_{k+1}, i_{k+1} + 1\}} + \delta_{i_k, (\delta-1)/2}, \\
a_k(\bi{\cdot}s_k) & = \#\set{1 \leq r \leq k-1 | i_r \in \{i_{k+1}, i_{k+1} - 1, i_k, i_k + 1\}} + \delta_{i_{k+1}, (\delta-1)/2}.
\end{align*}

Suppose $\t \in \Tud_n(\bi)$ and write $\t_{k-1} = \lambda$. Because $h_k(\bi) = 0$, we have $i_k \neq \pm \frac{1}{2}$ by~\autoref{deg:h3}. As $i_k = -i_{k+1}$, we have $\{i_k \pm 1\} \cap \{i_{k+1} \pm 1\} = \emptyset$. Hence by the definition of $a_k(\bi)$ and the construction of $\lambda$, we have
\begin{align}
(-1)^{a_k(\bi) + a_k(\bi{\cdot}s_k)} & = (-1)^{\#\set{1 \leq r \leq k-1 | i_r \in \{i_k \pm 1, i_{k+1} \pm 1\}} + \delta_{i_k, (\delta-1)/2} + \delta_{i_{k+1}, (\delta-1)/2}} \notag \\
& = (-1)^{\#\set{\alpha \in [\lambda] | \res(\alpha) \in \{i_k\pm 1,i_{k+1}\pm 1\}} + \delta_{i_k, (\delta-1)/2} + \delta_{i_{k+1}, (\delta-1)/2}}. \label{esscom:1:eq1}
\end{align}

Because $h_k(\bi) = 0$ and $i_k \neq 0$, by~\autoref{deg:h4}, there exist $\beta$ and $\gamma$ such that $\res(\beta) = -\res(\gamma) = i_k$ and either $\beta,\gamma \in \mathscr A(\lambda)$ or $\beta,\gamma \in \mathscr R(\lambda)$. Therefore, by the construction of $\lambda$, we have
$$
\begin{cases}
\text{$\#\set{\alpha \in [\lambda] | \res(\alpha) \in \{i_m-1,i_m+1\}}$ is odd if $i_m \neq \frac{\delta-1}{2}$,}\\
\text{$\#\set{\alpha \in [\lambda] | \res(\alpha) \in \{i_m-1,i_m+1\}}$ is even if $i_m = \frac{\delta-1}{2}$,}
\end{cases}
$$
where $m \in \{k,k+1\}$. Therefore, by (\ref{esscom:1:eq1}) and $i_k \neq i_{k+1}$, we have $(-1)^{a_k(\bi) + a_k(\bi{\cdot}s_k)} = 1$. \endproof

\begin{Lemma} \label{V:3}
Suppose $\bi = (i_1, \ldots, i_n) \in I^n$ and $\t \in \Tud_n(\bi)$. If $i_k = i_{k+1}$ and $\t(k) + \t(k+1) = 0$ for some $1 \leq k \leq n-1$, we have $V_k(\t) = s_k(\t,\t)$.
\end{Lemma}

\proof When $i_k = i_{k+1}$ and $\t(k) + \t(k+1) = 0$, we have $i_k = i_{k+1} = 0$, which implies $\bi \in I_{k,0}^n$. By~\autoref{PQ:3} and~\autoref{esscom:1}, we have $P_k(\t) Q_k(\t) = e_k(\t,\t)$. Therefore
$$
V_k(\t) = \frac{P_k(\t) Q_k(\t) - 1}{c_\t(k) - c_\t(k+1)} = \frac{e_k(\t,\t) - 1}{c_\t(k) + c_\t(k)} = s_k(\t,\t),
$$
which completes the proof. \endproof

\begin{Lemma} \label{V:4}
Suppose $1 \leq k \leq n-1$ and $\t$ is an up-down tableau. For up-down tableau $\s$, we have $V_k(\t) = V_k(\s)$ if one of the following conditions holds:

(a). If $\s = \t{\cdot}s_r$ for some $1 \leq r < k-1$.

(b). If $\t(r) + \t(r+1) = \s(r) + \s(r+1) = 0$ and $\t \overset{r}\sim \s$ for some $1 \leq r < k-1$.
\end{Lemma}

\proof Because $r < k-1$, in both (a) and (b), we have $c_\t(k) = c_\s(k)$ and $c_\t(k+1) = c_\s(k+1)$. Hence we have $c_\t(k) - c_\t(k+1) = c_\s(k) - c_\s(k+1)$. By~\autoref{PQ:4} and~\autoref{PQ:5}, we have $P_k(\t) = P_k(\s)$ and $Q_k(\t) = Q_k(\s)$. Hence, we have
$$
\hspace*{4cm} V_k(\t) = \frac{P_k(\t) Q_k(\t) - 1}{c_\t(k) - c_\t(k+1)} = \frac{P_k(\s) Q_k(\s) - 1}{c_\s(k) - c_\s(k+1)} = V_k(\s). \hspace*{4cm}\qedhere
$$
\endproof

\subsection{Generators of $\B$} \label{sec:generators}

In this subsection, we define
$$
G_n(\delta) = \set{e(\bi)|\bi \in P^n} \cup \set{y_k| 1 \leq k \leq n} \cup \set{\psi_k| 1\leq k \leq n-1} \cup \set{\epsilon_k| 1\leq k \leq n-1}
$$
in $\B$ and show that $G_n(\delta)$ is a generating set of $\B$.

By~\autoref{idem-semi} we have defined a set $\set{e(\bi) | \bi \in P^n} \subset \L_n$ where $e(\bi) \neq 0$ only if $\bi \in I^n$, and by (\ref{def:y}) we have a set of elements $y_k \in \L_n$ for $1 \leq k \leq n$. Recall that we have showed that $\sum_{\bi \in I^n} e(\bi)^\O = \sum_{\bi \in P^n} e(\bi)^\O = 1_\F$, $\sum_{\bi \in I^n} e(\bi) = \sum_{\bi \in P^n} e(\bi) = 1_R$ and $y_k$ is nilpotent for $1 \leq k \leq n$.

It left us to define $\set{\psi_r, \epsilon_r | 1 \leq r \leq n-1}$. Suppose $\l$ is a positive integer and $0 \leq r \leq \l$. Define
\begin{equation} \label{crl}
c_r^{(\l)} =
\begin{cases}
0, & \text{if $r = \l/2$, $r \neq 0$ and $4\ |\ \l$,}\\
2, & \text{if $k = \l/2$, $r \neq 0$ and $2\ |\ \l$ but $4\ \not|\ \l$,}\\
1, & \text{otherwise;}
\end{cases}
\end{equation}
and $z_\l = \sum_{r = 0}^\l c_r^{(\l)}$. We have $z_\l > 0$ because $c_0^{(\l)} = 1$ and $c_r^{(\l)} \geq 0$ for any $r$ and $\l$.

Suppose $\bi = (i_1, \ldots, i_n) \in I^n$ and $1 \leq k \leq n-1$. If $i_k \neq i_{k+1}$, by~\autoref{L:inv:1} we define
$$
\frac{1}{L_k^\O - L_{k+1}^\O} e(\bi)^\O := \sum_{\t \in \mathscr T^{ud}_n(\bi)} \frac{1}{c_\t(k) - c_\t(k+1)} \f{\t} \in \L_n(\O),
$$
and $\frac{1}{L_k - L_{k+1}} e(\bi) = \frac{1}{L_k^\O - L_{k+1}^\O} e(\bi)^\O \otimes_\O 1_R \in \mathscr L_n$.

\begin{Lemma} \label{idem:2}
Suppose $1 \leq k \leq n-1$ and $\bi = (i_1, \ldots, i_n) \in I^n$ with $i_k + i_{k+1} \neq 0$. Then we have $e(\bi)^\O e_k^\O = e_k^\O e(\bi)^\O = 0$ in $\BOx$ and $e(\bi)e_k = e_k e(\bi) = 0$ in $\B$.
\end{Lemma}

\proof For any $\t \in \Tud_n(\bi)$, we have $\t(k) + \t(k+1) \neq 0$ because $i_k + i_{k+1} \neq 0$. Hence we have $\f{\t} e_k^\O = e_k^\O \f{\t} = 0$ in $\BOx$, which implies $e(\bi)^\O e_k^\O = e_k^\O e(\bi)^\O = 0$. Then we have $e(\bi) e_k = e(\bi)^\O e_k^\O \otimes_\O 1_R = 0$ and $e_k e(\bi) = e_k^\O e(\bi)^\O \otimes_\O 1_R = 0$ in $\B$. \endproof

Suppose $1 \leq k \leq n-1$ and $\bi, \bj \in I^n$. Define
\begin{align*}
e(\bi)^{\mathscr O} \psi_k^{\mathscr O} e(\bj)^{\mathscr O} & :=
\begin{cases}
e(\bi)^\O (s_r^\O + e(\bi)^\O\frac{1}{L_k^\O - L_{k+1}^\O}e(\bj)^\O - \frac{1}{i_k+j_k} e_k^\O  & \\
\hspace*{1cm} - \frac{1}{i_k+j_k} \sum_{\l = 1}^\infty (-\frac{2}{i_k + j_k})^\l \frac{1}{z_\l} \left(\sum_{r=0}^\l c_r^{(\l)} (L_k^\O - i_k)^{\l - r} e_k^\O (L_k^\O - j_k)^r\right))e(\bj)^\O, & \text{if $\bj \neq \bi{\cdot}s_k$,}\\
e(\bi)^{\mathscr O} P_k^{\mathscr O}(\bi)^{-1} (s_k^{\mathscr O} - V_k^\O(\bi)) Q_k^{\mathscr O}(\bj)^{-1} e(\bj)^{\mathscr O}, & \text{if $\bj = \bi{\cdot}s_k$;}
\end{cases}\\
e(\bi)^{\mathscr O} \epsilon_k^{\mathscr O} e(\bj)^{\mathscr O} & := e(\bi)^{\mathscr O} P_k^{\mathscr O}(\bi)^{-1} e_k^{\mathscr O} Q_k^{\mathscr O}(\bj)^{-1} e(\bj)^{\mathscr O}.
\end{align*}

The above definition is well-defined. It is obvious that $e(\bi)^{\mathscr O} \epsilon_k^{\mathscr O} e(\bj)^{\mathscr O}$ is well-defined. For $e(\bi)^{\mathscr O} \psi_k^{\mathscr O} e(\bj)^{\mathscr O}$ with $\bj \neq \bi{\cdot}s_k$, by~\autoref{idem:2}, $e(\bi)^\O e_k^\O e(\bj)^\O \neq 0$ implies $i_k + i_{k+1} = j_k + j_{k+1} = 0$. Hence we have $i_k + j_k \neq 0$ because $\bj \neq \bi{\cdot}s_k$, which implies $e(\bi)^{\mathscr O} \psi_k^{\mathscr O} e(\bj)^{\mathscr O}$ with $\bj \neq \bi{\cdot}s_k$ is well-defined. For $e(\bi)^{\mathscr O} \psi_k^{\mathscr O} e(\bj)^{\mathscr O}$ with $\bj = \bi{\cdot}s_k$, $V_k(\bi)$ is only defined when $i_k = i_{k+1}$. But it only exists if $\bi = \bj$, otherwise $e(\bi)^\O P_k^\O(\bi)^{-1} V_k^\O(\bi) Q_k^\O(\bj)^{-1} e(\bj)^\O = 0$. When $\bi = \bj$, because $\bj = \bi{\cdot}s_k$, we have $i_k = i_{k+1}$. Therefore $e(\bi)^{\mathscr O} \psi_k^{\mathscr O} e(\bj)^{\mathscr O}$ with $\bj = \bi{\cdot}s_k$ is well-defined.

Then, define
$$
\psi_k^\O = \sum_{\bi \in I^n} \sum_{\bj \in I^n} e(\bi)^\O \psi_k^\O e(\bj)^\O \in \BOx, \qquad
\epsilon_k^\O = \sum_{\bi \in I^n} \sum_{\bj \in I^n} e(\bi)^\O \epsilon_k^\O e(\bj)^\O \in \BOx,
$$
and for any $\bi,\bj \in I^n$ and $1 \leq k \leq n-1$, we define $e(\bi) \psi_k e(\bj) = e(\bi)^{\mathscr O} \psi_k^{\mathscr O} e(\bj)^{\mathscr O} \otimes_\O 1_R$ and $e(\bi) \epsilon_k e(\bj) = e(\bi)^{\mathscr O} \epsilon_k^{\mathscr O} e(\bj)^{\mathscr O} \otimes_\O 1_R$, and
$$
\psi_k = \sum_{\bi \in I^n} \sum_{\bj \in I^n} e(\bi) \psi_k e(\bj) \in \B, \qquad
\epsilon_k = \sum_{\bi \in I^n} \sum_{\bj \in I^n} e(\bi) \epsilon_k e(\bj) \in \B.
$$

\begin{Remark} \label{remark:generators}
By the definitions of $e(\bi)^{\mathscr O} \psi_k^{\mathscr O} e(\bj)^{\mathscr O}$ and $e(\bi)^{\mathscr O} \epsilon_k^{\mathscr O} e(\bj)^{\mathscr O}$, it is easy to see that it is equivalently to define $e(\bi) \psi_k e(\bj)$ and $e(\bi) \epsilon_k e(\bj)$ by
\begin{align*}
e(\bi) \psi_k e(\bj) & :=
\begin{cases}
e(\bi) (s_r + e(\bi)\frac{1}{L_k - L_{k+1}}e(\bj) - \frac{1}{i_k+j_k} e_k  & \\
\hspace*{1cm} - \frac{1}{i_k+j_k} \sum_{\l = 1}^\infty (-\frac{2}{i_k + j_k})^\l \frac{1}{z_\l} \left(\sum_{r=0}^\l c_r^{(\l)} (L_k - i_k)^{\l - r} e_k (L_k - j_k)^r\right))e(\bj), & \text{if $\bj \neq \bi{\cdot}s_k$,}\\
e(\bi) P_k(\bi)^{-1} (s_k - V_k(\bi)) Q_k(\bj)^{-1} e(\bj), & \text{if $\bj = \bi{\cdot}s_k$.}
\end{cases}\\
e(\bi) \epsilon_k e(\bj) & := e(\bi) P_k(\bi)^{-1} e_k Q_k(\bj)^{-1} e(\bj).
\end{align*}
\end{Remark}

\begin{Proposition} \label{B:span}
The elements
$$
G_n(\delta) = \set{e(\bi) | \bi \in P^n} \cup \set{y_k | 1\leq k \leq n} \cup \set{\psi_k | 1 \leq k \leq n-1} \cup \set{\epsilon_k | 1 \leq k \leq n-1}
$$
generates $\mathscr B_n(\delta)$.
\end{Proposition}

\proof Suppose $S \subseteq \B$ is generated by $G_n(\delta)$. It is sufficient to prove that $\set{s_k, e_k | 1 \leq k \leq n-1}$ is contained in $S$.

For any $1 \leq k \leq n-1$ and $\bi, \bj \in I^n$, we have $e(\bi) e_k e(\bj) = e(\bi) P_k(\bi) \epsilon_k Q_k(\bi) e(\bj) \in S$. Hence, by $\sum_{\bi \in I^n} e(\bi) = 1$, we have $e_k = \sum_{\bi \in I^n} \sum_{\bj \in I^n} e(\bi) e_k e(\bj) \in S$.

For any $1 \leq k \leq n-1$ and $\bi, \bj \in I^n$, if $\bj \neq \bi{\cdot}s_k$, we have
$$
e(\bi) s_k e(\bj) = e(\bi) \left(\psi_k - \frac{1}{L_k - L_{k+1}} +  \frac{1}{i_k+j_k} e_k + \frac{1}{i_k+j_k} \sum_{\l = 1}^\infty (-\frac{2}{i_k + j_k})^\l \frac{1}{z_\l} \left(\sum_{r=0}^\l c_r^{(\l)} (L_k - i_k)^{\l - r} e_k (L_k - j_k)^r\right)\right)e(\bj) \in S;
$$
and for $\bj = \bi{\cdot}s_k$, we have
$$
e(\bi) s_k e(\bj) = e(\bi)P_k(\bi)(\psi_k + V_k(\bi))Q_k(\bj) e(\bj) \in S,
$$
which implies that $s_k = \sum_{\bi \in I^n} \sum_{\bj \in I^n} e(\bi) s_k e(\bj) \in S$ by $\sum_{\bi \in I^n} e(\bi) = 1$. Hence the Proposition holds. \endproof

Reader may notice that the definitions of $e(\bi)^\O \psi_k^\O e(\bj)^\O$ and $e(\bi) \psi_k e(\bj)$ when $\bj \neq \bi{\cdot}s_k$ are comparatively complicated. In the rest of this subsection, we simplify the definitions of these two elements, and the results we obtained will be used in the next section.

Suppose $\t$ is an up-down tableau with $\t(k) + \t(k+1) = 0$ and write $\lambda = \t_{k-1} = \t_{k+1}$, $\mu = \t_k$ and $\alpha = \lambda\ominus\mu$. We say $\t$ is \textit{$k$-added} if $\mu = \lambda\cup\{\alpha\}$, and $\t$ is \textit{$k$-removed} if $\mu = \lambda\backslash\{\alpha\}$. Define $A_k$ to be the set of all $k$-added up-down tableaux and $R_k$ to be the set of all $k$-removed up-down tableaux.

\begin{Lemma} \label{com:1:9}
Suppose $\l$ is a positive integer, and $c_r^{(\l)}$ and $z_\l$ are defined as in (\ref{crl}). We have $z_\l > 0$ and $\sum_{r = 0}^\l (-1)^r c_r^{(\l)} = \sum_{r = 0}^\l (-1)^{\l-r} c_r^{(\l)} = 0$.
\end{Lemma}

\begin{proof}
It is obvious that $z_\l > 0$ because $c_r^{(\l)} \geq 0$ and $c_1^\l = 1$ for any $\l$. The following diagram gives the values of $(-1)^r c_r$:
$$
\begin{array} {cccccccc}
\l = 1: & 1 & -1 &&&&& \\
\l = 2: & 1 & -2 & 1 &&&& \\
\l = 3: & 1 & -1 & 1 & -1 &&& \\
\l = 4: & 1 & -1 & 0 & 1 & -1 && \\
\l = 5: & 1 & -1 & 1 & -1 & 1 & -1 & \\
\l = 6: & 1 & -1 & 1 & -2 & 1 & -1 & 1 \\
\ldots & \ldots & \ldots & \ldots & \ldots & \ldots & \ldots & \ldots
\end{array}
$$

By direct calculation, we have $\sum_{r = 0}^\l (-1)^r c_r^{(\l)} = 0$. Similarly, we have $\sum_{r = 0}^\l (-1)^{\l-r} c_r^{(\l)} = 0$, and the Lemma follows.
\end{proof}

For convenience, we set $d = \frac{x - \delta}{2}$.

\begin{Lemma} \label{com:1:1}
Suppose $\t$ is an up-down tableau with residue sequence $\bi$ and $1 \leq k \leq n-1$. Then
$$
y_k^\O f_{\t\t} = f_{\t\t} y_k^\O =
\begin{cases}
d f_{\t\t}, & \text{if $\t \in A_k$,}\\
-d f_{\t\t}, & \text{if $\t \in R_k$.}
\end{cases}
$$
\end{Lemma}

\begin{proof}
This Lemma is a direct application of the definition of $y_k^\O$.
\end{proof}

\begin{Lemma} \label{com:1:2}
Suppose $(\lambda,f) \in \widehat B_n$ and $\s, \t \in \Tud_n(\lambda)$ with $\s \overset{k}\sim \t$. For any integer $\l > 0$, we have
$$
\f{\s} ( \sum_{r = 0}^\l c_r^{(\l)} (y_k^{\mathscr O})^{\l - r} e_k^{\mathscr O} (y_k^{\mathscr O})^r) \f{\t} =
\begin{cases}
z_\l d^\l \f{\s} e_k^{\mathscr O} \f{\t}, & \text{if $\s,\t \in A_k$,}\\
z_\l (-d)^\l \f{\s} e_k^{\mathscr O} \f{\t}, & \text{if $\s,\t \in R_k$,}\\
0, & \text{otherwise.}
\end{cases}
$$
\end{Lemma}

\proof Suppose $\s, \t \in A_k$. By~\autoref{com:1:1} we have
$$
\f{\s} ( \sum_{r = 0}^\l c_r^{(\l)} (y_k^{\mathscr O})^{\l - r} e_k^{\mathscr O} (y_k^{\mathscr O})^r) \f{\t} = \f{\s} ( \sum_{r = 0}^\l c_r^{(\l)} d^{\l - r} e_k^{\mathscr O} d^r) \f{\t} = z_\l d^\l \f{\s} e_k^{\mathscr O} \f{\t}.
$$

Suppose $\s, \t \in R_k$. By~\autoref{com:1:1} we have
$$
\f{\s} ( \sum_{r = 0}^\l c_r^{(\l)} (y_k^{\mathscr O})^{\l - r} e_k^{\mathscr O} (y_k^{\mathscr O})^r) \f{\t} = \f{\s} ( \sum_{r = 0}^\l c_r^{(l)} (-d)^{\l - r} e_k^{\mathscr O} (-d)^r) \f{\t} = z_\l (-d)^\l \f{\s} e_k^{\mathscr O} \f{\t}.
$$

Suppose $\s \in R_k$ and $\t \in A_k$. By~\autoref{com:1:1} and~\autoref{com:1:9} we have
$$
\f{\s} ( \sum_{r = 0}^\l c_r^{(\l)} (y_k^{\mathscr O})^{\l - r} e_k^{\mathscr O} (y_k^{\mathscr O})^r) \f{\t} = \f{\s} ( \sum_{r = 0}^\l (-1)^{\l - r} c_r^{(\l)} d^{\l - r} e_k^{\mathscr O} d^r) \f{\t} = 0.
$$

For $\s \in A_k$ and $\t \in R_k$, following the same argument we have $\f{\s} ( \sum_{r = 0}^\l c_r (y_k^{\mathscr O})^{\l - r} e_k^{\mathscr O} (y_k^{\mathscr O})^r) \f{\t} = 0$, which completes the proof. \endproof

Suppose $\s, \t$ are up-down tableaux with residue sequences $\bi = (i_1, \ldots, i_n)$ and $\bj = (j_1, \ldots, j_n)$, respectively, and $\s \overset{k}\sim \t$. Then $s_k(\t,\s) = \frac{1}{c_\s(k) + c_\t(k)} (e_k(\t,\s) - \delta_{\s,\t})$. One can see that $i_k + i_{k+1} = j_k + j_{k+1} = 0$. If $\bj \neq \bi{\cdot}s_k$, i.e. $i_k + j_k \neq 0$, we have $\frac{1}{c_\s(k) + c_\t(k)} \in \mathscr O$. As $\mathscr O = R[[x-\delta]]$,~\autoref{com:1:2} gives us a method to express $e(\bi)^{\mathscr O} s_k^{\mathscr O} e(\bj)^{\mathscr O}$ using $e(\bi)^{\mathscr O} e_k^{\mathscr O} e(\bj)^{\mathscr O}$ and $y_k^{\mathscr O}$'s.

\begin{Lemma} \label{com:1:3}
Suppose $\s, \t$ are up-down tableaux with residue sequences $\bi$ and $\bj$, respectively, and $\s \overset{k}\sim \t$. If $\bj \neq \bi{\cdot}s_k$, we have
$$
\f{\s} s_k^{\mathscr O} \f{\t} = \f{\s} (\frac{1}{i_k + j_k} (e_k^{\mathscr O} + \sum_{\l = 1}^\infty (- \frac{2}{i_k + j_k})^\l \frac{1}{z_\l} (\sum_{r = 0}^\l c_r^{(\l)} (y_k^{\mathscr O})^{\l - r} e_k^{\mathscr O} (y_k^{\mathscr O})^r))) \f{\t} - e(\bi)^\O\frac{1}{L_k^{\mathscr O} - L_{k+1}^{\mathscr O}} e(\bj)^\O \f{\t}.
$$
\end{Lemma}

\proof As $s_k(\t,\s) = \frac{1}{c_\s(k) + c_\t(k)} (e_k(\t,\s) - \delta_{\s,\t})$, we have
\begin{equation} \label{com:1:3:eq1}
\f{\s} s_k^{\mathscr O} \f{\t} = \f{\s} \frac{e_k^{\mathscr O} - \delta_{\s\t}}{c_\s(k) + c_\t(k)} \f{\t} = \f{\s} \frac{e_k^{\mathscr O}}{c_\s(k) + c_\t(k)} \f{\t} - \f{\s} \frac{\delta_{\s\t}}{c_\s(k) + c_\t(k)} \f{\t}.
\end{equation}

If $\s = \t$, we have $\bi = \bj$. Hence, by the definition of $\frac{1}{L_k^\O - L_{k+1}^\O} e(\bj)^\O$, we have
\begin{equation} \label{com:1:3:eq2}
\f{\s}\frac{\delta_{\s,\t}}{c_\s(k) + c_\t(k)} \f{\t} = \frac{\delta_{\s,\t}}{c_\t(k) - c_\t(k+1)} \f{\t} = \f{\s} e(\bi)^\O\frac{1}{L_k^{\mathscr O} - L_{k+1}^{\mathscr O}} e(\bj)^\O \f{\t}.
\end{equation}

By (\ref{com:1:3:eq1}) and (\ref{com:1:3:eq2}), it is sufficient to show that
\begin{equation} \label{com:1:3:eq3}
\f{\s} \frac{e_k^{\mathscr O}}{c_\s(k) + c_\t(k)} \f{\t} = \f{\s} (\frac{1}{i_k + j_k} (e_k^{\mathscr O} + \sum_{\l = 1}^\infty (- \frac{2}{i_k + j_k})^\l \frac{1}{z_\l} (\sum_{r = 0}^\l c_r^{(\l)} (y_k^{\mathscr O})^{\l - r} e_k^{\mathscr O} (y_k^{\mathscr O})^r))) \f{\t}.
\end{equation}

Suppose $\s,\t \in A_k$. We have
$$
\f{\s} \frac{e_k^\O}{c_\s(k) + c_\t(k)} \f{\t} = \f{\s} \frac{e_k^\O}{2d + i_k + j_k} \f{\t} = \frac{1}{i_k + j_k}\sum_{\l = 0}^\infty (-\frac{2}{i_k + j_k})^\l d^\l \f{\s} e_k^\O \f{\t}.
$$

Then by~\autoref{com:1:2}, we have
\begin{eqnarray*}
&& \f{\s} (\frac{1}{i_k + j_k} (e_k^{\mathscr O} + \sum_{\l = 1}^\infty (- \frac{2}{i_k + j_k})^\l \frac{1}{z_\l} (\sum_{r = 0}^\l c_r^{(\l)} (y_k^{\mathscr O})^{\l - r} e_k^{\mathscr O} (y_k^{\mathscr O})^r))) \f{\t}\\
& = & \frac{1}{i_k + j_k} \left(1 + \sum_{\l = 1}^\infty (- \frac{2}{i_k + j_k})^\l d^\l \right) \f{\s} e_k^{\mathscr O} \f{\t} = \f{\s} \frac{e_k^{\mathscr O}}{c_\s(k) + c_\t(k)} \f{\t}.
\end{eqnarray*}

Suppose $\s,\t \in R_k$. We have
$$
\f{\s} \frac{e_k^\O}{c_\s(k) + c_\t(k)} \f{\t} = \f{\s} \frac{e_k^\O}{-2d + i_k + j_k} \f{\t} = \frac{1}{i_k + j_k}\sum_{\l = 0}^\infty (-\frac{2}{i_k + j_k})^\l (-d)^\l \f{\s} e_k^\O \f{\t}.
$$

Then by~\autoref{com:1:2}, we have
\begin{eqnarray*}
&& \f{\s} (\frac{1}{i_k + j_k} (e_k^{\mathscr O} + \sum_{\l = 1}^\infty (- \frac{2}{i_k + j_k})^\l \frac{1}{z_\l} (\sum_{r = 0}^\l c_r^{(\l)} (y_k^{\mathscr O})^{\l - r} e_k^{\mathscr O} (y_k^{\mathscr O})^r))) \f{\t}\\
& = & \frac{1}{i_k + j_k} \left(1 + \sum_{\l = 1}^\infty (- \frac{2}{i_k + j_k})^\l (-d)^\l \right) \f{\s} e_k^{\mathscr O} \f{\t} = \f{\s} \frac{e_k^{\mathscr O}}{c_\s(k) + c_\t(k)} \f{\t}.
\end{eqnarray*}

Suppose $\s \in A_k$ and $\t \in R_k$. We have $c_\s(k) = d + i_k$ and $c_\t(k) = -d + j_k$. Then
$$
\f{\s} \frac{e_k^\O}{c_\s(k) + c_\t(k)} \f{\t} = \frac{1}{i_k + j_k} \f{\s} e_k^\O \f{\t}.
$$

Then by~\autoref{com:1:2}, we have
$$
\f{\s} (\frac{1}{i_k + j_k} (e_k^{\mathscr O} + \sum_{\l = 1}^\infty (- \frac{2}{i_k + j_k})^\l \frac{1}{z_\l} (\sum_{r = 0}^\l c_r^{(\l)} (y_k^{\mathscr O})^{\l - r} e_k^{\mathscr O} (y_k^{\mathscr O})^r))) \f{\t} = \frac{1}{i_k + j_k} \f{\s} e_k^{\mathscr O} \f{\t} = \f{\s} \frac{e_k^\O}{c_\s(k) + c_\t(k)} \f{\t}.
$$

Follow the same argument, we can show that (\ref{com:1:3:eq3}) holds for $\s \in R_k$ and $\t \in A_k$. Therefore (\ref{com:1:3:eq3}) holds, which proves the Lemma. \endproof

\begin{Lemma} \label{com:1:4}
Suppose $1 \leq k \leq n-1$ and $\s$ is an up-down tableau with residue sequence $\bi$ and $\s(k) + \s(k+1) \neq 0$. Then we have
$$
\f{\s} s_k^\O \f{\s} = -\f{\s} e(\bi)^\O\frac{1}{L_k^{\mathscr O} - L_{k+1}^{\mathscr O}} e(\bj)^\O \f{\s}.
$$
\end{Lemma}

\proof By~\autoref{semi:B}, we have $\f{\s} s_k^\O \f{\s} = \frac{1}{c_\s(k+1) - c_\s(k)} \f{\s}$. The Lemma follows because $\f{\s} e(\bi)^\O\frac{1}{L_k^{\mathscr O} - L_{k+1}^{\mathscr O}} e(\bj)^\O \f{\s} = \frac{1}{c_\s(k) - c_\s(k+1)} \f{\s}$. \endproof

By~\autoref{com:1:3} and~\autoref{com:1:4}, we can simplify the definition of $e(\bi)^\O \psi_k^\O e(\bj)^\O$ and $e(\bi) \psi_k e(\bj)$ when $\bj \neq \bi{\cdot}s_k$.

\begin{Corollary} \label{com:1}
Suppose $\bi, \bj \in I^n$. We have $e(\bi)^\O \psi_k^\O e(\bj)^\O = 0$ and $e(\bi) \psi_k e(\bj) = 0$ if $\bi \neq \bj{\cdot}s_k$.
\end{Corollary}

\proof Suppose $\s \in \Tud_n(\bi)$. If $\s(k) + \s(k+1) \neq 0$, by~\autoref{semi:B}, we have $f_{\s\s} s_k^\O e(\bj)^\O \neq 0$ only if $\bj = \bi$ or $\bj = \bi{\cdot}s_k$. Hence we assume $\bj = \bi$. Then by~\autoref{com:1:4}, we have
\begin{equation} \label{com:1:eq1}
f_{\s\s} s_k^\O e(\bj)^\O = f_{\s\s} s_k^\O \f{\s} = -f_{\s\s} e(\bi)^\O\frac{1}{L_k^{\mathscr O} - L_{k+1}^{\mathscr O}} e(\bj)^\O \f{\s}.
\end{equation}

If $\s(k) + \s(k+1) = 0$, by~\autoref{idem-semi} and~\autoref{com:1:3}, we have
\begin{align}
f_{\s\s} s_k^{\mathscr O} e(\bj)^\O & = f_{\s\s} s_k^\O \left(\sum_{\t \in \Tud_n(\bj)} \f{\t} \right) = \sum_{\substack{\t \in \Tud_n(\bj) \\ \s \overset{k} \sim \t}} f_{\s\s} s_k^\O \f{\t} \notag \\
& = f_{\s\s} (\frac{1}{i_k + j_k} (e_k^{\mathscr O} + \sum_{\l = 1}^\infty (- \frac{2}{i_k + j_k})^\l \frac{1}{z_\l} (\sum_{r = 0}^\l c_r^{(\l)} (y_k^{\mathscr O})^{\l - r} e_k^{\mathscr O} (y_k^{\mathscr O})^r))) e(\bj)^\O - f_{\s\s}\frac{1}{L_k^{\mathscr O} - L_{k+1}^{\mathscr O}} e(\bj)^\O. \label{com:1:eq2}
\end{align}

Therefore, as $e(\bi)^\O y_k^\O = e(\bi)^\O (L_k^\O - i_k)$ and $y_k^\O e(\bj)^\O = (L_k^\O - j_k) e(\bj)^\O$, (\ref{com:1:eq1}) - (\ref{com:1:eq2}) implies $e(\bi)^\O \psi_k^\O e(\bj)^\O = 0$. By lifting the element of $\BOx$ into $\B$, we have $e(\bi) \psi_k e(\bj) = e(\bi)^\O \psi_k^\O e(\bj)^\O \otimes_\O 1_R = 0$. \endproof

By~\autoref{com:1}, we re-write the definitions of $e(\bi)^\O \psi_k^\O e(\bj)^\O$ and $e(\bi) \psi_k e(\bj)$ as
\begin{align*}
e(\bi)^{\mathscr O} \psi_k^{\mathscr O} e(\bj)^{\mathscr O} & :=
\begin{cases}
0, & \text{if $\bj \neq \bi{\cdot}s_k$,}\\
e(\bi)^{\mathscr O} P_k^{\mathscr O}(\bi)^{-1} (s_k^{\mathscr O} - V_k^\O(\bi)) Q_k^{\mathscr O}(\bj)^{-1} e(\bj)^{\mathscr O}, & \text{if $\bj = \bi{\cdot}s_k$;}
\end{cases}\\
e(\bi)^{\mathscr O} \epsilon_k^{\mathscr O} e(\bj)^{\mathscr O} & := e(\bi)^{\mathscr O} P_k^{\mathscr O}(\bi)^{-1} e_k^{\mathscr O} Q_k^{\mathscr O}(\bj)^{-1} e(\bj)^{\mathscr O},
\end{align*}
and
\begin{align*}
e(\bi) \psi_k e(\bj) & :=
\begin{cases}
0, & \text{if $\bj \neq \bi{\cdot}s_k$,}\\
e(\bi) P_k(\bi)^{-1} (s_k - V_k(\bi)) Q_k(\bj)^{-1} e(\bj), & \text{if $\bj = \bi{\cdot}s_k$.}
\end{cases}\\
e(\bi) \epsilon_k e(\bj) & := e(\bi) P_k(\bi)^{-1} e_k Q_k(\bj)^{-1} e(\bj).
\end{align*}

Moreover,~\autoref{com:1} forces the elements $\psi_k$ to be the intertwining elements of $\B$, i.e. for any $\bi \in P^n$ we have $e(\bi) \psi_k = \psi_k e(\bi{\cdot}s_k)$.

\section{Grading of Brauer algebras} \label{sec:grade}

In this section we are going to prove that the elements of $G_n(\delta)$ in $\B$ follow the same relations associated in $\G{n}$, which implies that we can define a surjective homomorphism $\G{n} \longrightarrow \mathscr B_n(\delta)$ by sending
$$
e(\bi) \mapsto e(\bi), \qquad y_r \mapsto y_r, \qquad \psi_k \mapsto \psi_k, \qquad \epsilon_k \mapsto \epsilon_k
$$
with $\bi \in P^n$, $1 \leq r \leq n$ and $1 \leq k \leq n-1$. Throughout of this section we are going to work in $\BOx$ and $\B$. So $e(\bi)$, $y_r$, $\psi_k$ and $\epsilon_k$ we used in this section will be elements in $\B$.

The goal of this section is to prove that relations (\ref{rela:1}) - (\ref{rela:8:9}) holds in $\B$. Note that by~\autoref{idem-semi}, we have $e(\bi) = 0$ if $\bi \not\in I^n$. By~\autoref{idem:nonzero}, we have similar property in $\G{n}$. Hence in this section we assume $\bi, \bj, \bk \in I^n$ except when we prove relation~\eqref{rela:5:5}, because apart from~\eqref{rela:5:5}, when any of $\bi,\bj,\bk$ involved in the relations is not a residue sequence, i.e. any of $e(\bi), e(\bj), e(\bk)$ involved in the relations equals $0$, then both sides of the relations will be $0$ and there is nothing to prove. In relation~\eqref{rela:5:5}, whenever $\bi \not\in I^n$, i.e. $e(\bi) = 0$, the left hand side of the relation is $0$. But it is not obvious that the right hand side of the relation equals $0$. Hence when we prove~\eqref{rela:5:5}, we will allow $\bi$ chosen from $P^n$ rather than assuming $\bi \in I^n$.

\subsection{Actions of generators on seminormal forms} \label{sec:semiB}

Most of the calculations in this section are in $\BOx$. Hence our first step is to calculate the actions of generators on seminormal forms of $\B$.

Fix $(\lambda,f) \in \widehat B_n$ and $\s,\t \in \Tud_n(\lambda)$. The actions of $e(\bi)^\O$ and $y_k^\O$ on $f_{\s\t}$ from right and $f_{\t\s}$ from left can be easily calculated by the definitions.

\begin{Lemma} \label{semi:y}
Suppose $\t \in \Tud_n(\lambda)$ with residue sequence $\bj = (j_1, \ldots, j_n)$ and $1 \leq k \leq n$. Then for any $\s \in \Tud_n(\lambda)$ and $\bi \in P^n$, we have
\begin{align*}
& f_{\s\t} e(\bi)^\O = \delta_{\bi,\bj} f_{\s\t}, \qquad \qquad f_{\s\t} y_k^\O = (c_\t(k) - j_k) f_{\s\t},\\
& e(\bi)^\O f_{\t\s} = \delta_{\bi,\bj} f_{\t\s}, \qquad \qquad y_k^\O f_{\t\s} = (c_\t(k) - j_k) f_{\t\s}.
\end{align*}
\end{Lemma}

Next we calculate the actions of $\psi_k^\O$ and $\epsilon_k^\O$ on $f_{\s\t}$ from right and $f_{\t\s}$ from left when $\t(k) + \t(k+1) \neq 0$.

\begin{Lemma} \label{semi:1}
Suppose $1 \leq k \leq n - 1$ and $\t \in \Tud_n(\lambda)$ with residue sequence $\bi = (i_1, \ldots, i_n)$ and $\t(k) + \t(k+1) \neq 0$. For any $\s \in \Tud_n(\lambda)$,
\begin{enumerate}
\item if $\t{\cdot}s_k$ does not exist, we have $f_{\s\t} \psi_k^\O = \psi_k^\O f_{\t\s} = 0$; and if $\u = \t{\cdot}s_k$ exists, we have
\begin{align}
f_{\s\t} \psi_k^\O & = \begin{cases}
\frac{1}{c_\t(k+1) - c_\t(k)} f_{\s\t} + \frac{s_k(\t,\u)}{1 - c_\u(k) + c_\u(k+1)} f_{\s\u}, & \text{if $i_k = i_{k+1}$,}\\
s_k(\t,\u)(c_\u(k) - c_\u(k+1)) f_{\s\u}, & \text{if $i_k = i_{k+1} - 1$,}\\
s_k(\t,\u)\frac{c_\u(k) - c_\u(k+1)}{1 - c_\u(k) + c_\u(k+1)} f_{\s\u}, & \text{otherwise;}
\end{cases} \label{semi:1:eq1} \\
\psi_k^\O f_{\t\s} & = \begin{cases}
\frac{1}{c_\t(k+1) - c_\t(k)} f_{\t\s} + \frac{s_k(\t,\u)}{1 - c_\t(k) + c_\t(k+1)} f_{\u\s}, & \text{if $i_k = i_{k+1}$,}\\
s_k(\t,\u)(c_\t(k) - c_\t(k+1)) f_{\u\s}, & \text{if $i_k = i_{k+1} + 1$,}\\
s_k(\t,\u)\frac{c_\t(k) - c_\t(k+1)}{1 - c_\t(k) + c_\t(k+1)} f_{\s\u}, & \text{otherwise.}
\end{cases} \label{semi:1:eq2}
\end{align}

\item we have $f_{\s\t} \epsilon_k^\O = \epsilon_k^\O f_{\t\s} = 0$.
\end{enumerate}
\end{Lemma}

\proof (1). Suppose $\t{\cdot}s_k \not\in \Tud_n(\lambda)$. By~\autoref{com:1}, we have $f_{\s\t} \psi_k^\O = f_{\s\t} \psi_k^\O e(\bi{\cdot}s_k)^\O$. Let $\bj = \bi{\cdot}s_k$. By~\autoref{semi:B} and the definition of $\psi_k^\O$, we have $f_{\s\t}\psi_k^\O e(\bj)^\O = a f_{\s\t} e(\bj)^\O$ for some $a \in R(x)$. Hence, we have $f_{\s\t}\psi_k^\O e(\bj)^\O \neq 0$ only if $\bi = \bj$ by~\autoref{semi:y}.

Because $\bj = \bi = \bi{\cdot}s_k$, we have $i_k = i_{k+1}$. As $\t(k) + \t(k+1) \neq 0$, by the construction of $\t$, we have $i_k = i_{k+1} \neq 0$. By the definition of $h_k$, we have $h_{k+1}(\bi) \geq h_k(\bi) + 2$. As $\bi \in I^n$, we have $-2 \leq h_k(\bi), h_{k+1}(\bi) \leq 0$ by~\autoref{deg:h2}, which forces $h_k(\bi) = -2$. By~\autoref{deg:h4:h}, we have either $\t(k) > 0$ and $\t(k+1) < 0$, or $\t(k) < 0$ and $\t(k+1) > 0$. Therefore, $\t{\cdot}s_k \in \Tud_n(\lambda)$ by~\autoref{y:h2:8}. This implies $f_{\s\t} \psi_k^\O = 0$ when $\t{\cdot}s_k \not\in \Tud_n(\lambda)$. Following the similar argument, we have $\psi_k^\O f_{\t\s} = 0$ when $\t{\cdot}s_k \not\in \Tud_n(\lambda)$.

Suppose $\t{\cdot}s_k \in \Tud_n(\lambda)$. By~\autoref{com:1}, we have $f_{\s\t} \psi_k^\O = f_{\s\t} \psi_k^\O e(\bi{\cdot}s_k)^\O$. When $i_k = i_{k+1}$, we have $\bj = \bi = \bi{\cdot}s_k$. By the definition of $\psi_k^\O$,~\autoref{PQ:1} and~\autoref{V:2}, we have
\begin{align*}
f_{\s\t} \psi_k^\O & = P_k(\t)^{-1} Q_k(\t)^{-1} (\frac{1}{c_\t(k+1) - c_\t(k)} - V_k(\t))f_{\s\t} + P_k(\t)^{-1} Q_k(\u)^{-1} s_k(\t,\u)f_{\s\u}\\
& = \frac{1}{(c_\t(k) - c_\t(k+1))V_k(\t) + 1} \frac{1 + (c_\t(k) - c_\t(k+1))V_k(\t)}{c_\t(k+1) - c_\t(k)} f_{\s\t} + \frac{s_k(\t,\u)}{1 - c_\u(k) + c_\u(k+1)} f_{\s\u}\\
& = \frac{1}{c_\t(k+1) - c_\t(k)} f_{\s\t} + \frac{s_k(\t,\u)}{1 - c_\u(k) + c_\u(k+1)} f_{\s\u}.
\end{align*}

When $i_k = i_{k+1} - 1$, by the definition of $\psi_k^\O$ and~\autoref{PQ:1}, we have
$$
f_{\s\t} \psi_k^\O = P_k(\t)^{-1} Q_k(\u)^{-1} s_k(\t,\u) f_{\s\u} = s_k(\t,\u) (c_\u(k) - c_\u(k+1)) f_{\s\u}.
$$

For the other cases, by the definition of $\psi_k^\O$ and~\autoref{PQ:1}, we have
$$
f_{\s\t} \psi_k^\O = P_k(\t)^{-1} Q_k(\u)^{-1} s_k(\t,\u) f_{\s\u} = s_k(\t,\u) \frac{c_\u(k) - c_\u(k+1)}{1 - c_\u(k) + c_\u(k+1)} f_{\s\u}.
$$

Therefore, (\ref{semi:1:eq1}) holds. Following similar argument, we can prove (\ref{semi:1:eq2}).

(2). By~\autoref{idem:2} and~\autoref{semi:y}, we have $f_{\s\t} \epsilon_k^\O = f_{\s\t} e(\bi)^\O \epsilon_k^\O = 0$. Similarly, we have $\epsilon_k^\O f_{\t\s} = 0$.\endproof

Notice that the actions of $e(\bi)^\O$, $\psi_k^\O$ and $y_k^\O$ on $f_{\s\t}$ from right and $f_{\t\s}$ from left with $\t(k) + \t(k+1) \neq 0$ are the same as in the KLR algebras. See Hu-Mathas~\cite[Lemma 4.23]{HuMathas:SemiQuiver}. Following the same process as Hu-Mathas~\cite[Proposition 4.28, 4.29]{HuMathas:SemiQuiver}, we have the following Corollary.

\begin{Corollary} \label{KLR:1}
Suppose $(\lambda,f) \in \widehat B_n$, $\t \in \Tud_n(\lambda)$ with residue sequence $\bi = (i_1, \ldots, i_n)$ and $1 \leq k \leq n-1$. If $\t(k) + \t(k+1) \neq 0$, for any $\s \in \Tud_n(\lambda)$, we have
$$
f_{\s\t} (\psi_k^\O)^2 =
\begin{cases}
0, & \text{if $i_k = i_{k+1}$,}\\
f_{\s\t} (y_k^\O - y_{k+1}^\O), & \text{if $i_k = i_{k+1} - 1$,}\\
f_{\s\t} (y_{k+1}^\O - y_k^\O), & \text{if $i_k = i_{k+1} + 1$,}\\
f_{\s\t}, & \text{if $|i_k - i_{k+1}| > 1$.}
\end{cases}
$$

Similarly, for $1 \leq k \leq n-2$, if $\t(k) + \t(k+1) \neq 0$, $\t(k) + \t(k+2) \neq 0$ and $\t(k+1) + \t(k+2) \neq 0$, we have
$$
f_{\s\t} (\psi_k^\O \psi_{k+1}^\O \psi_k^\O - \psi_{k+1}^\O \psi_k^\O \psi_{k+1}^\O) =
\begin{cases}
f_{\s\t}, & \text{if $i_{k-1} = i_{k+1} = i_k - 1$,}\\
-f_{\s\t}, & \text{if $i_{k-1} = i_{k+1} = i_k + 1$,}\\
0, & \text{otherwise.}
\end{cases}
$$
\end{Corollary}

Finally we calculate the actions of $\psi_k^\O$ and $\epsilon_k^\O$ on $f_{\s\t}$ from right and $f_{\t\s}$ from left when $\t(k) + \t(k+1) = 0$.

\begin{Lemma} \label{semi:2}
Suppose $1 \leq k \leq n-1$ and $\t \in \Tud_n(\lambda)$ with residue sequence $\bi = (i_1, \ldots, i_n)$ and $\t(k) + \t(k+1) = 0$. For any $\s \in \Tud_n(\lambda)$, we have
\begin{align}
f_{\s\t} \psi_k^\O & =
\begin{cases}
0, & \text{if $i_k = i_{k+1} = 0$ or $h_k(\bi) \neq 0$,}\\
\frac{1}{c_\t(k) + c_\u(k)} P_k(\t)^{-1} e_k(\t,\u) Q_k(\u)^{-1} f_{\s\u}, & \text{otherwise;}
\end{cases} \label{semi:2:eq1} \\
\psi_k^\O f_{\t\s} & =
\begin{cases}
0, & \text{if $i_k = i_{k+1} = 0$ or $h_k(\bi) \neq 0$,}\\
\frac{1}{c_\t(k) + c_\u(k)} P_k(\u)^{-1} e_k(\t,\u) Q_k(\t)^{-1} f_{\u\s}, & \text{otherwise,}
\end{cases} \label{semi:2:eq3}
\end{align}
where $\u$ is the unique up-down tableau with residue sequence $\bi{\cdot}s_k$ and $\u \overset{k}\sim \t$; and
\begin{align}
f_{\s\t} \epsilon_k^\O & = \sum_{\v \overset{k}\sim \t} P_k(\t)^{-1} e_k(\t,\v) Q_k(\v)^{-1} f_{\s\v}, \label{semi:2:eq2}\\
\epsilon_k^\O f_{\t\s} & = \sum_{\v \overset{k}\sim \t} P_k(\v)^{-1} e_k(\t,\v) Q_k(\t)^{-1} f_{\v\s}. \label{semi:2:eq4}
\end{align}
\end{Lemma}

\proof \eqref{semi:2:eq2} and \eqref{semi:2:eq4} can be easily derived by the definition of $\epsilon_k^\O$ and~\autoref{semi:B}. For \eqref{semi:2:eq1}, by~\autoref{com:1} we have $f_{\s\t} \psi_k^\O = f_{\s\t} \psi_k^\O e(\bi{\cdot}s_k)^\O$. Hence, if $h_k(\bi) \neq 0$, by~\autoref{A:9} we have $\bi{\cdot}s_k \not\in I^n$, which implies $f_{\s\t}\psi_k^\O = f_{\s\t} \psi_k^\O e(\bi{\cdot}s_k)^\O = 0$ by~\autoref{idem-semi}. If $i_k = i_{k+1} = 0$, by~\autoref{V:3} we have $f_{\s\t}(s_k^\O - V_k(\bi)^\O) = 0$, which implies $f_{\s\t}\psi_k^\O = 0$. Therefore we have proved that $f_{\s\t}\psi_k^\O = 0$ if $i_k = i_{k+1} = 0$ or $h_k(\bi) \neq 0$.

For the other cases, by~\autoref{deg:h4:1} there exists a unique $\u \in \Tud_n(\bi{\cdot}s_k)$ such that $\u \overset{k}\sim \t$. Because $i_k \neq 0$, we have $\bi \neq \bi{\cdot}s_k$, and $\t \neq \u$. Therefore, we have
$$
f_{\s\t} \psi_k^\O = f_{\s\t} \psi_k^\O e(\bi{\cdot}s_k)^\O = P_k(\t)^{-1} s_k(\t,\u) Q_k(\u)^{-1} f_{\s\u} = \frac{1}{c_\t(k) + c_\u(k)} P_k(\t)^{-1} e_k(\t,\u) Q_k(\u)^{-1} f_{\s\u}.
$$

Hence, (\ref{semi:2:eq1}) holds. Following the same argument, (\ref{semi:2:eq3}) holds. \endproof

\subsection{Idempotent and (essential) commutation relations} \label{sec:rela:1}

In this subsection we are going to prove the idempotent relations, the commutation relations and the essential commutation relations hold in $\B$. First we prove the idempotent relations.

\begin{Lemma} \label{idem:1}
Suppose $\bi \in I^n$. We have $y_1^{\delta_{i_1,\frac{\delta-1}{2}}} e(\bi) = 0$.
\end{Lemma}

\proof Because $\bi \in I^n$, we have $i_1 = \frac{\delta-1}{2}$. Hence by~\autoref{idem-semi} and~\autoref{semi:y}, we have
$$
y_1^{\mathscr O} e(\bi)^{\mathscr O} = \sum_{\t \in \mathscr T^{ud}_n(\bi)} y_1^{\mathscr O} \f{\t} = \frac{x - \delta}{2} \sum_{\t \in \mathscr T^{ud}_n(\bi)} \frac{1}{\gamma_\t} f_{\t\t} = \frac{x - \delta}{2} e(\bi)^{\mathscr O} \in \mathscr (x - \delta)B_n^{\mathscr O}(\delta).
$$

Hence we have $y_1^{\delta_{i_1,\frac{\delta-1}{2}}} e(\bi) = y_1 e(\bi) = y_1^{\mathscr O} e(\bi)^{\mathscr O} \otimes_\O 1_R = 0$. \endproof

\begin{Proposition} \label{idem}
In $\mathscr B_n(\delta)$, the idempotent relations hold.
\end{Proposition}

\proof By~\autoref{idem-semi} we have $\sum_{\bi \in P^n} e(\bi)^\O = \sum_{\bi \in I^n} e(\bi)^\O = 1$ and $e(\bi)^\O e(\bj)^\O = \delta_{\bi,\bj} e(\bi)^\O \in \BOx$, which implies $\sum_{\bi \in P^n} e(\bi) = 1$ and $e(\bi) e(\bj) = \delta_{\bi,\bj} e(\bi)$ in $\B$. By~\autoref{idem:2}, we have $e(\bi)e_k = e_k e(\bi) = 0$ if $i_k + i_{k+1} \neq 0$. Hence by~\autoref{idem:1}, we complete the proof. \endproof

Next we prove the commutation relations. First we prove that $y_k$'s commute with $e(\bi)$'s.

\begin{Lemma} \label{com:2}
Suppose $\bi \in I^n$ and $1 \leq k \leq n$. We have $y_k e(\bi) = e(\bi) y_k$.
\end{Lemma}

\proof By~\autoref{idem-semi} and~\autoref{semi:y}, we have
$$
y_k^\O e(\bi)^\O = \sum_{\t \in \mathscr T^{ud}_n(\bi)} y_1^\O \f{\t} = \sum_{\t \in \mathscr T^{ud}_n(\bi)} c_\t(k) \f{\t} = \sum_{\t \in \mathscr T^{ud}_n(\bi)} \f{\t} y_k^\O = e(\bi)^\O y_k^\O \in \BOx,
$$
which implies $y_k e(\bi) = e(\bi) y_k$ by lifting the elements into $\B$. \endproof

Hence by~\autoref{com:1} and~\autoref{com:2}, we have shown that (\ref{rela:2:1}) holds in $\B$. Now we prove the rest of commutation relations hold in $\B$ as well.

\begin{Lemma} \label{com:3}
For $1 \leq k,r \leq n$ and $1 \leq m \leq n-1$, we have $y_k y_r = y_r y_k$, and if $|k - m| > 1$, we have $y_k \psi_m = \psi_m y_k$ and $y_k \epsilon_m = \epsilon_m y_k$.
\end{Lemma}

\proof By the definition of $y_k$, we have $y_k, y_r \in \L$. Because $\L$ is a commutative subalgebra of $\B$, we have $y_k y_r = y_r y_k$ for $1 \leq k,r \leq n$.

Suppose $1 \leq k \leq n$ and $1 \leq m \leq n-1$, and $|k - m| > 1$. For any $\bi \in I^n$ and $\t \in \Tud_n(\bi)$, if $\t(m) + \t(m+1) \neq 0$, by~\autoref{semi:1}, without loss of generality, we have $\psi_m^\O f_{\t\t} = a f_{\t\t} + b f_{\u\t}$, where $a,b \in R(x)$ and $\u = \t{\cdot}s_m$; and if $\t(m) + \t(m+1) = 0$, by~\autoref{semi:2} we have $\psi_m^\O f_{\t\t} = a f_{\u\t}$, where $a \in R(x)$ and $\u \overset{m}\sim \t$. In either case, for $|k - m| > 1$, we have $y_k^\O \psi_m^\O f_{\t\t} = \psi_m^\O f_{\t\t} y_k^\O$ because $c_\u(k) = c_\t(k)$, which implies
$$
y_k^\O \psi_m^\O = \sum_{\substack{\bi \in I^n \\ \t \in \Tud_n(\bi)}} y_k^\O \psi_m^\O \f{\t} = \sum_{\substack{\bi \in I^n \\ \t \in \Tud_n(\bi)}} \psi_m^\O \f{\t} y_k^\O = \psi_m^\O y_k^\O \in \BOx.
$$

Hence we have $y_k \psi_m = \psi_m y_k$ by lifting the elements into $\B$.

Suppose $1 \leq k \leq n$ and $1 \leq m \leq n-1$, and $|k - m| > 1$. For any $\bi \in I^n$ and $\t \in \Tud_n(\bi)$ with $\t(m) + \t(m+1) = 0$, by~\autoref{semi:2}, we have $\epsilon_k^\O f_{\t\t} = \sum_{\u \overset{m}\sim \t} a_\u f_{\u\t}$, where $a_\u \in R(x)$. Hence, because $|k-m| > 1$, we have $y_k^\O \epsilon_m^\O f_{\t\t} = \epsilon_m^\O f_{\t\t} y_k^\O$ because $c_\u(k) = c_\t(k)$ for any $\u \overset{m}\sim \t$, which implies
$$
y_k^\O \epsilon_m^\O = \sum_{\substack{\bi \in I^n \\ \t \in \Tud_n(\bi)}} y_k^\O \epsilon_m^\O \f{\t} = \sum_{\substack{\bi \in I^n \\ \t \in \Tud_n(\bi)}} \epsilon_m^\O \f{\t} y_k^\O = \epsilon_m^\O y_k^\O \in \BOx.
$$

Hence we have $y_k \epsilon_m = \epsilon_m y_k$ by lifting the elements into $\B$. \endproof

\begin{Lemma} \label{com:4:1}
Suppose $\bi, \bj \in I^n$ and $1 \leq r, k \leq n-1$. If $|k - r| > 1$, we have
\begin{align*}
& e(\bi) P_k(\bi)^{-1} s_r e(\bj) = e(\bi) s_r P_k(\bj)^{-1} e(\bj), && e(\bi) Q_k(\bi)^{-1} s_r e(\bj) = e(\bi) s_r Q_k(\bj)^{-1} e(\bj), && e(\bi) V_k(\bi) s_r e(\bj) = e(\bi) s_r V_k(\bj) e(\bj)\\
& e(\bi) P_k(\bi)^{-1} e_r e(\bj) = e(\bi) e_r P_k(\bj)^{-1} e(\bj), && e(\bi) Q_k(\bi)^{-1} e_r e(\bj) = e(\bi) e_r Q_k(\bj)^{-1} e(\bj), && e(\bi) V_k(\bi)^{-1} e_r e(\bj) = e(\bi) e_r V_k(\bj)^{-1} e(\bj).
\end{align*}
\end{Lemma}

\proof We only prove $e(\bi) P_k(\bi)^{-1} s_r e(\bj) = e(\bi) s_r P_k(\bj)^{-1} e(\bj)$ here. The rest of the equalities follow by the similar argument, except that we use~\autoref{V:4} instead of~\autoref{PQ:4} and~\autoref{PQ:5} when we prove
$$
e(\bi) V_k(\bi) s_r e(\bj) = e(\bi) s_r V_k(\bj) e(\bj) \qquad \text{and} \qquad  e(\bi) V_k(\bi)^{-1} e_r e(\bj) = e(\bi) e_r V_k(\bj)^{-1} e(\bj).
$$

Suppose $k > r$. Because $|k - r| > 1$, we have $k - 1 > r$. Choose arbitrary $\t \in \Tud_n(\bi)$. If $\t(r) + \t(r+1) \neq 0$, by~\autoref{semi:B}, without loss of generality, we can write
$$
f_{\t\t} s_r^\O e(\bj)^\O = a_\s f_{\t\s} + a_\t f_{\t\t},
$$
where $\s = \t{\cdot}s_r$, and $a_\s, a_\t \in R(x)$. Note we set $a_\s = 0$ if $\s$ is not an up-down tableau. Hence, by~\autoref{PQ:4}, we have
$$
f_{\t\t} P_k^\O(\bi)^{-1} s_r^\O e(\bj)^\O = a_\s P_k(\t)^{-1} f_{\t\s} + a_\t P_k(\t)^{-1} f_{\t\t} = a_\s  f_{\t\s} P_k(\s)^{-1} + a_\t f_{\t\t} P_k(\t)^{-1} = f_{\t\t} s_r^\O P^\O_k(\bj)^{-1} e(\bj)^\O,
$$
which implies
\begin{equation} \label{com:4:1:eq1}
f_{\t\t} P_k^\O(\bi)^{-1} s_r^\O e(\bj)^\O = f_{\t\t} s_r^\O P_k(\bj)^{-1} e(\bj)^\O,
\end{equation}
when $\t(r) + \t(r+1) \neq 0$.

If $\t(r) + \t(r+1) = 0$, by~\autoref{semi:B}, we have $f_{\t\t} s_r^\O e(\bj)^\O = \sum_{\substack{\s \in \Tud_n(\bj) \\ \s \overset{r}\sim \t}} a_\s f_{\t\s}$, where $a_\s \in R(x)$. Hence, by~\autoref{PQ:5}, we have
$$
f_{\t\t} P_k^\O(\bi)^{-1} s_r^\O e(\bj)^\O = \sum_{\substack{\s \in \Tud_n(\bj) \\ \s \overset{r}\sim \t}} a_\s P_k(\t)^{-1} f_{\t\s} = \sum_{\substack{\s \in \Tud_n(\bj) \\ \s \overset{r}\sim \t}} a_\s f_{\t\s} P_k(\s)^{-1} = f_{\t\t} s_r^\O P^\O_k(\bj)^{-1} e(\bj)^\O,
$$
which implies
\begin{equation} \label{com:4:1:eq2}
f_{\t\t} P_k^\O(\bi)^{-1} s_r^\O e(\bj)^\O = f_{\t\t} s_r^\O P_k(\bj)^{-1} e(\bj)^\O,
\end{equation}
when $\t(r) + \t(r+1) = 0$.

By (\ref{com:4:1:eq1}) and (\ref{com:4:1:eq2}), for arbitrary $\t \in \Tud_n(\bi)$, we have $f_{\t\t} P_k^\O(\bi)^{-1} s_r^\O e(\bj)^\O = f_{\t\t} s_r^\O P_k(\bj)^{-1} e(\bj)^\O$. Therefore, by~\autoref{idem-semi}, we have
\begin{align*}
e(\bi)^\O P^\O_k(\bi)^{-1} s_r^\O e(\bj)^\O & = \sum_{\t \in \Tud_n(\bi)} f_{\t\t} P^\O_k(\bi)^{-1} s_r^\O e(\bj)^\O\\
& = \sum_{\t \in \Tud_n(\bi)} f_{\t\t} s_r^\O e(\bj)^\O P^\O_k(\bj)^{-1} = e(\bi)^\O s_r^\O P_k(\bj)^{-1} e(\bj)^\O \in \BOx.
\end{align*}

By lifting the elements into $\B$, we have $e(\bi) P_k(\bi)^{-1} s_r e(\bj) = e(\bi) s_r P_k(\bj)^{-1} e(\bj)$ when $k > r$.

Suppose $k < r$. Because $|k - r| > 1$, we have $k < r-1$. Choose arbitrary $\t \in \Tud_n(\bi)$. By~\autoref{semi:B}, without loss of generality, we have $f_{\t\t} s_r^\O e(\bj)^\O = \sum_{\substack{ \s \in \Tud_n(\bj) \\ \s|_{r-1} = \t|_{r-1}}} a_\s f_{\t\s}$, where $a_\s \in R(x)$. By the definition, $P_k(\s)$ depends on $\s(1), \s(2), \ldots, \s(k+1)$. Because $k < r-1$, for any $\s \in \Tud_n(\bj)$ with $a_\s \neq 0$, we have $P_k(\s)^{-1} = P_k(\t)^{-1}$. Therefore, we have
\begin{align*}
f_{\t\t} P_k^\O(\bi)^{-1}s_r^\O e(\bj)^\O & = \sum_{\substack{ \s \in \Tud_n(\bj) \\ \s|_{r-1} = \t|_{r-1}}} a_\s P_k(\t)f_{\t\s} = \sum_{\substack{ \s \in \Tud_n(\bj) \\ \s|_{r-1} = \t|_{r-1}}} a_\s f_{\t\s} P_k(\s) = f_{\t\t} s_r^\O P_k^\O(\bj)^{-1} e(\bj)^\O.
\end{align*}

Because $\t$ is chosen arbitrary, following the same argument as when $k > r$, we have $e(\bi) P_k(\bi)^{-1} s_r e(\bj) = e(\bi) s_r P_k(\bj)^{-1} e(\bj)$ when $k < r$. \endproof

\begin{Lemma} \label{com:4}
For $1 \leq k, r \leq n-1$ where $|k-r| > 1$, we have $\psi_k \psi_r = \psi_r \psi_k$, $\psi_k \epsilon_r = \epsilon_r \psi_k$ and $\epsilon_k \epsilon_r = \epsilon_r \epsilon_k$.
\end{Lemma}

\proof Without loss of generality, we assume $k > r$. We will only prove $\psi_k \psi_r = \psi_r \psi_k$, and the rest of the equalities follow by the same argument.

Choose arbitrary $\bi \in I^n$. By the definitions, we have $P_k^\O(\bi)^{-1} e(\bi)^\O = P_k^\O(\bi)^{-1}$ and $Q_k^\O(\bi)^{-1} e(\bi)^\O = Q_k^\O(\bi)^{-1}$, which implies $P_k(\bi)^{-1} e(\bi) = P_k(\bi)^{-1}$ and $Q_k(\bi)^{-1} e(\bi) = Q_k(\bi)^{-1}$ by lifting the elements into $\B$. Because $|k - r| > 1$, $s_r$ and $s_k$ commutes. Hence, by~\autoref{com:1} and~\autoref{com:4:1}, we have
\begin{align*}
\psi_k \psi_r e(\bi) & = e(\bi{\cdot}s_k s_r) \psi_k e(\bi{\cdot}s_r) \psi_r e(\bi)\\
& = e(\bi{\cdot}s_k s_r) P_k(\bi{\cdot}s_k s_r)^{-1} (s_k - V_k(\bi{\cdot}s_r)) Q_k(\bi{\cdot}s_r)^{-1} P_r(\bi{\cdot}s_r)^{-1}(s_r - V_r(\bi)) Q_r(\bi)^{-1} e(\bi)\\
& = e(\bi{\cdot}s_k s_r) P_r(\bi{\cdot}s_k s_r)^{-1} (s_r - V_r(\bi{\cdot}s_k)) Q_r(\bi{\cdot}s_k)^{-1} P_k(\bi{\cdot}s_k)^{-1}(s_k - V_k(\bi)) Q_k(\bi)^{-1} e(\bi)\\
& = e(\bi{\cdot}s_ks_r) \psi_r e(\bi{\cdot}s_k) \psi_k e(\bi) = \psi_r \psi_k e(\bi).
\end{align*}

As $\bi$ is chosen arbitrary, we have $\psi_k \psi_r = \psi_r \psi_k$, which completes the proof. \endproof

\begin{Proposition} \label{com}
In $\B$, the commutation relations hold.
\end{Proposition}

\proof By~\autoref{com:1}, we have $e(\bi) \psi_k e(\bj) = 0$ if $\bi \neq \bj{\cdot}s_k$. Hence, we have $e(\bi) \psi_k = e(\bi) \psi_k e(\bi{\cdot}s_k) = \psi_k e(\bi{\cdot}s_k)$. Therefore, (\ref{rela:2:1}) holds by~\autoref{com:2}. The relation (\ref{rela:2:2}) holds by~\autoref{com:3} and (\ref{rela:2:3}) holds by~\autoref{com:4}. \endproof

In the rest of this subsection, we prove that the essential commutation relations hold in $\B$. First we introduce the following results, which will be used.

\begin{Lemma} \label{esscom:2}
Suppose $1 \leq k \leq n-1$ and $\bi \in I_{k,0}^n$ with $i_k = -i_{k+1} \neq \pm \frac{1}{2}$. Then we have $e(\bi) \epsilon_k e(\bi) = (-1)^{a_k(\bi)} e(\bi)$.
\end{Lemma}

\proof As $\bi \in I_{k,0}^n$ with $i_k = -i_{k+1} \neq \frac{1}{2}$, when $i_k = -i_{k+1} \neq 0$, we have $h_k(\bi) = -1$, which implies $h_{k+1}(\bi) = -1$. Choose arbitrary $\t \in \Tud_n(\bi)$ and let $\t_{k-1} = \lambda$ and $\t_k = \mu$. We have a unique $\alpha \in \mathscr{AR}_\lambda(i_k)$ and unique $\beta \in  \mathscr{AR}_\mu(i_{k+1})$ because $h_k(\bi) = h_{k+1}(\bi) = -1$. Therefore, we have $\alpha = \beta$. Similarly, when $i_k = -i_{k+1} = 0$, choose arbitrary $\t \in \Tud_n(\bi)$, we have $\t(k) + \t(k+1) = 0$.

Choose arbitrary $\t \in \Tud_n(\bi)$. By~\autoref{PQ:3}, we have
$$
\f{\t} \epsilon_k^\O \f{\t} = \f{\t} P^\O_k(\bi)^{-1} e_k^\O Q^\O_k(\bi)^{-1} \f{\t} = P_k(\t)^{-1} Q_k(\t)^{-1} e_k(\t,\t) \f{\t} = (-1)^{a_k(\bi)}\f{\t}.
$$

As $\bi \in I_{k,0}^n$ with $i_k = -i_{k+1} \neq \pm\frac{1}{2}$, when $i_k = -i_{k+1} \neq 0$, we have $h_k(\bi) = -1$. Hence by~\autoref{deg:h4:1}, for any $\s \in \mathscr T^{ud}_n(\bi)$ with $\s \overset{k}\sim \t$, we have $\s = \t$. When $i_k = -i_{k+1} = 0$, we have $h_k(\bi) = 0$. Hence by~\autoref{deg:h4:1}, there exists an unique up-down tableau $\s \in \bi{\cdot}s_k$ such that $\s \overset{k}\sim \t$. Because $\bi{\cdot}s_k = \bi$ and $\t \overset{k}\sim \t$, it forces $\s = \t$. Therefore, we conclude that for any $\s \in \mathscr T^{ud}_n(\bi)$ with $\s \overset{k}\sim \t$, we have $\s = \t$. Hence,
$$
e(\bi)^\O \epsilon_k^\O e(\bi)^O = \left(\sum_{\s \in \Tud_n(\bi)} \f{\s}\right) \epsilon_k^\O \left(\sum_{\s \in \Tud_n(\bi)} \f{\s}\right) = \sum_{\s \in \Tud_n(\bi)} \f{\s} \epsilon_k^\O \f{\s} = (-1)^{a_k(\bi)} \sum_{\s \in \Tud_n(\bi)} \f{\s} = (-1)^{a_k(\bi)} e(\bi)^\O \in \BOx,
$$
and we have $e(\bi) \epsilon_k e(\bi) = (-1)^{a_k(\bi)} e(\bi)$ by lifting elements to $\B$. \endproof

\begin{Lemma} \label{esscom:3}
Suppose $1 \leq k \leq n-1$ and $\bi \in I^n$ with $i_k = i_{k+1} = 0$. Then we have $e(\bi)^\O \psi_k^\O e(\bi)^\O = 0$ and $e(\bi) \psi_k e(\bi) = 0$.
\end{Lemma}

\proof Suppose $\t \in \mathscr T^{ud}_n(\bi)$. Because $i_k = i_{k+1} = 0$, we have $c_\t(k) + c_\t(k+1) = 0$. By~\autoref{V:3}, we have $V_k(\t) = s_k(\t,\t)$. Hence we have
$$
\f{\t} \psi_k^\O \f{\t} = P_k(\t) Q_k(\t) (s_k(\t,\t) - V_k(\t)) \f{\t} = 0.
$$

Because $i_k = 0$, we have $h_k(\bi) = 0$. By~\autoref{deg:h4:1}, there exists an unique $\s \in \Tud_n(\bi{\cdot}s_k)$ such that $\s \overset{k}\sim \t$. Because $\bi = \bi{\cdot}s_k$, we have $\s = \t$. Therefore,
$$
e(\bi)^\O \psi_k^\O e(\bi)^\O = \left( \sum_{\t \in \mathscr T^{ud}_n(\bi)} \f{\t} \right) \psi_k^\O \left( \sum_{\t \in \mathscr T^{ud}_n(\bi)} \f{\t} \right) = \sum_{\t \in \mathscr T^{ud}_n(\bi)} \f{\t} \psi_k^\O \f{\t} = 0,
$$
which implies $e(\bi) \psi_k e(\bi) = 0$ by lifting the elements into $\B$. \endproof


Recall that $L_k s_k - s_k L_{k+1} = s_k L_k - L_{k+1} s_k = e_k - 1$. The next Proposition shows that the essential commutation relations hold in $\B$.

\begin{Proposition} \label{esscom}
In $\B$, the essential commutation relations hold.
\end{Proposition}

\proof Suppose $i_k = i_{k+1} = 0$. We have $\bi{\cdot}s_k = \bi$. Then by~\autoref{esscom:1},~\autoref{esscom:2} and~\autoref{esscom:3}, we have
$$
e(\bi) y_k \psi_k = 0 = e(\bi) \psi_k y_{k+1} + e(\bi) \epsilon_k e(\bi) - e(\bi).
$$

Suppose $i_k = i_{k+1} \neq 0$. Then we have $\bi{\cdot}s_k = \bi$. Therefore, by~\autoref{V:2}, we have
\begin{align*}
e(\bi) y_k \psi_k & = e(\bi) P_k(\bi)^{-1} (s_k - V_k(\bi)) Q_k(\bi)^{-1} y_{k+1} + e(\bi) \epsilon_k e(\bi) - e(\bi) P_k(\bi)^{-1} Q_k(\bi)^{-1} \left(V_k(\bi)(L_k - L_{k+1}) + 1\right)\\
& = e(\bi) \psi_k y_{k+1} + e(\bi) \epsilon_k e(\bi) - e(\bi).
\end{align*}

Suppose $i_k \neq i_{k+1}$ and let $\bj = \bi{\cdot}s_k$. Notice that $\bj \neq \bi$. Hence we have $e(\bi) y_k \psi_k = e(\bi) \psi_k y_{k+1} + e(\bi) \epsilon_k e(\bj)$ by direct calculation. Hence the relation (\ref{rela:3:1}) holds. By applying the same method as above, (\ref{rela:3:2}) holds, which completes the proof. \endproof

\subsection{Inverse relations, essential idempotent relations and untwist relations} \label{sec:rela:2}

In this subsection, we are going to prove the inverse relations, essential idempotent relations and untwist relations hold in $\B$. First we prove the inverse relations of $\B$.

\begin{Lemma} \label{inv:1}
Suppose $\bi \in I^n$ with $|i_k - i_{k+1}| > 1$ and $\t \in \Tud_n(\bi)$ with $\t(k) + \t(k+1) = 0$ for $1 \leq k \leq n-1$. Then we have $f_{\t\t} (\psi_k^\O)^2 = f_{\t\t}$ if $h_k(\bi) = 0$.
\end{Lemma}

\proof We have $i_k + i_{k+1} = 0$ because $\t(k) + \t(k+1) = 0$. Then, as $|i_k - i_{k+1}| > 1$ and $h_k(\bi) = 0$, we have $\bi \in I_{k,+}^n$.

By~\autoref{deg:h4:1}, there exists a unique $\u \in \Tud_n(\bi{\cdot}s_k)$ such that $\u \overset{k}\sim \t$ and $c_\t(k) - i_k = c_\u(k) - i_{k+1}$. Therefore,
\begin{align*}
f_{\t\t} (\psi_k^\O)^2 & = \frac{P_k(\t)^{-1} e_k(\t,\u) Q_k(\u)^{-1} P_k(\u)^{-1} e_k(\u,\t) Q_k(\t)^{-1}}{(c_\t(k) + c_\u(k))^2}f_{\t\t}\\
& = \frac{P_k(\t)^{-1}Q_k(\t)^{-1} Q_k(\u)^{-1} P_k(\u)^{-1} e_k(\t,\t) e_k(\u,\u) }{4(c_\t(k) - i_k) (c_\u(k) - i_{k+1})}f_{\t\t} = (-1)^{a_k(\bi) + a_k(\bi{\cdot}s_k)}f_{\t\t},
\end{align*}
by~\autoref{semi:B},~\autoref{semi:2} and~\autoref{PQ:3}. As $h_k(\bi) = 0$, by~\autoref{esscom:1}, we have $(-1)^{a_k(\bi) + a_k(\bi{\cdot}s_k)} = 1$, which proves the Lemma. \endproof

\begin{Lemma} \label{inv:2}
Suppose $\bi \in I^n$ and $\t \in \Tud_n(\bi)$. For any $1 \leq k \leq n-1$, we have
$$
f_{\t\t} (\psi_k^\O)^2 =
\begin{cases}
0, & \text{if $i_k = i_{k+1}$, or $i_k + i_{k+1} = 0$ and $h_k(\bi) \neq 0$,}\\
f_{\t\t} (y_k^\O - y_{k+1}^\O), & \text{if $i_k = i_{k+1} - 1$ and $i_k + i_{k+1} \neq 0$,}\\
f_{\t\t} (y_{k+1}^\O - y_k^\O), & \text{if $i_k = i_{k+1} + 1$ and $i_k + i_{k+1} \neq 0$,}\\
f_{\t\t}, & \text{otherwise.}
\end{cases}
$$
\end{Lemma}

\proof We prove the Lemma by considering each case.

\textbf{Case 1:} $i_k = i_{k+1}$, or $i_k + i_{k+1} = 0$ and $h_k(\bi) \neq 0$.

In this case, when $i_k = i_{k+1} \neq 0$, we have $i_k + i_{k+1} \neq 0$, which implies $\t(k) + \t(k+1) \neq 0$. Hence we have $f_{\t\t} (\psi_k^\O)^2 = 0$ by~\autoref{KLR:1}. When $i_k = i_{k+1} = 0$, by~\autoref{esscom:3}, we have $e(\bi)^\O \psi_k^\O = 0$, which implies $f_{\t\t} (\psi_k^\O)^2 = 0$. When $i_k + i_{k+1} = 0$ and $h_k(\bi) \neq 0$, we have $\bi{\cdot}s_k \not\in I^n$ by~\autoref{A:9}, which implies $f_{\t\t} (\psi_k^\O)^2 = f_{\t\t} \psi_k^\O e(\bi{\cdot}s_k)^\O \psi_k^\O = 0$ because $e(\bi)^\O \psi_k^\O = \psi_k^\O e(\bi{\cdot}s_k)^\O$ by~\autoref{com:1} and $e(\bi{\cdot}s_k)^\O = 0$ by~\autoref{idem-semi}.

\textbf{Case 2:} $i_k = i_{k+1} - 1$ and $i_k + i_{k+1} \neq 0$.

In this case, $i_k + i_{k+1} \neq 0$ forces $\t(k) + \t(k+1) \neq 0$. Hence we have $f_{\t\t}(\psi_k^\O)^2 = f_{\t\t} (y_k^\O - y_{k+1}^\O)$ by~\autoref{KLR:1}.

\textbf{Case 3:} $i_k = i_{k+1} + 1$ and $i_k + i_{k+1} \neq 0$.

Following the same argument as in Case 2, we have $f_{\t\t} (\psi_k^\O)^2 = f_{\t\t} (y_{k+1}^\O - y_k^\O)$.

\textbf{Case 4:} $|i_k - i_{k+1}| > 1$ with either $i_k + i_{k+1} \neq 0$, or $i_k + i_{k+1} = 0$ and $h_k(\bi) = 0$.

First we show this case contains all $\bi$ which does not satisfy Case 1 - 3. Choose arbitrary $\bi \in I^n$ does not satisfy Case 1 - 3. We have $i_k \neq i_{k+1}$ because $\bi$ does not satisfy Case 1. Assume $i_k \neq i_{k+1} \pm 1$. Then $i_k + i_{k+1} = 0$ because $\bi$ does not satisfy Case 2 - 3, and if $i_k + i_{k+1} = 0$, we have $h_k(\bi) \neq 0$ by~\autoref{deg:h3}, which contradicts that $\bi$ does not satisfy Case 1. Hence we always have $|i_k - i_{k+1}| > 1$. It is easy to see we have either $i_k + i_{k+1} \neq 0$, or $i_k + i_{k+1} = 0$ and $h_k(\bi) = 0$ because $\bi$ does not satisfy Case 1. This proves that Case 4 contains all $\bi$ which does not satisfy Case 1 - 3.

We separate this case further by considering $\t(k) + \t(k+1) \neq 0$ and $\t(k) + \t(k+1) = 0$. When $\t(k) + \t(k+1) \neq 0$, as $|i_k - i_{k+1}| > 1$, we have $f_{\t\t} (\psi_k^\O)^2 = f_{\t\t}$ by~\autoref{KLR:1}. When $\t(k) + \t(k+1) = 0$, it forces $i_k + i_{k+1} = 0$. Hence we have $h_k(\bi) = 0$ and we have $f_{\t\t} (\psi_k^\O)^2 = f_{\t\t}$ by~\autoref{inv:1}. \endproof

\begin{Proposition} \label{inv}
In $\B$, the inverse relations hold.
\end{Proposition}

\proof In~\autoref{inv:2}, as $\t$ is chosen arbitrary, we have
$$
e(\bi)^\O (\psi_k^\O)^2 =
\begin{cases}
0, & \text{if $i_k = i_{k+1}$ or $i_k + i_{k+1} = 0$ and $h_k(\bi) \neq 0$,}\\
e(\bi)^\O (y_k^\O - y_{k+1}^\O), & \text{if $i_k = i_{k+1} + 1$ and $i_k + i_{k+1} \neq 0$,}\\
e(\bi)^\O (y_{k+1}^\O - y_k^\O), & \text{if $i_k = i_{k+1} - 1$ and $i_k + i_{k+1} \neq 0$,}\\
e(\bi)^\O, & \text{otherwise,}
\end{cases}
$$
by~\autoref{idem-semi}, and the Proposition follows by lifting the elements into $\B$. \endproof

Then we prove the essential idempotent relations. Recall that~\autoref{esscom:2} proved the first relation of (\ref{rela:5:1}). The next two Lemmas prove the rest of relations of (\ref{rela:5:1}).

\begin{Lemma} \label{essidem:1}
Suppose $1 \leq k \leq n-1$ and $\bi \in I_{k,+}^n$ with $i_k = -i_{k+1} \neq -\frac{1}{2}$. Then we have
$$
e(\bi) \epsilon_k e(\bi) = (-1)^{a_k(\bi) + 1} (y_{k+1} - y_k) e(\bi).
$$
\end{Lemma}

\proof Choose arbitrary $\t \in \Tud_n(\bi)$. As $\bi \in I_{k,+}^n$ and $i_k \neq -\frac{1}{2}$, we have $h_k(\bi) = 0$ and $i_k \neq 0$. By~\eqref{remark:h:eq2} we have $h_{k+1}(\bi) = -2$. Hence, by~\autoref{deg:h4:2}, we have either $\t(k) + \t(k+1) = 0$, or $\t(k) + \t(k+1) \neq 0$ and $c_\t(k) - i_k = c_\t(k+1) - i_{k+1}$.

Suppose $\t(k) + \t(k+1) \neq 0$ and $c_\t(k) - i_k = c_\t(k+1) - i_{k+1}$. We have
$$
f_{\t\t} \epsilon_k^\O e(\bi)^\O = 0 = (-1)^{a_k(\bi) + 1} ((c_\t(k+1) - i_{k+1}) - (c_\t(k) - i_k))f_{\t\t} = (-1)^{a_k(\bi) + 1} (y_{k+1}^\O - y_k^\O) f_{\t\t},
$$
by~\autoref{semi:y} and~\autoref{semi:1}.

Suppose $\t(k) + \t(k+1) = 0$. We have $c_\t(k) - i_k = - (c_\t(k+1) - i_{k+1})$. By~\autoref{deg:h4:1}, for any $\s \in \Tud_n(\bi)$ with $\s \overset{k}\sim \t$, we have $\s = \t$. Therefore, by~\autoref{PQ:3}, we have
$$
f_{\t\t} \epsilon_k^\O e(\bi)^\O = f_{\t\t} P_k(\t)^{-1} e_k(\t,\t) Q_k(\t)^{-1} \f{\t} = (-1)^{a_k(\bi)} 2(c_\t(k) - i_k) f_{\t\t} = (-1)^{a_k(\bi)} (y_{k+1}^\O - y_k^\O) f_{\t\t}.
$$

As $\t$ is chosen arbitrary, we have $e(\bi)^\O \epsilon_k^\O e(\bi)^\O = (-1)^{a_k(\bi)} (y_{k+1}^\O - y_k^\O) e(\bi)^\O \in \BOx$ by~\autoref{idem-semi}, which proves the Lemma by lifting the elements into $\B$ and $(y_k + y_{k+1}) \epsilon_k = 0$. \endproof

\begin{Lemma} \label{essidem:2}
Suppose $1 \leq k \leq n-1$ and $\bi \in I_{k,+}^n$ with $i_k = -i_{k+1} = -\frac{1}{2}$. Then we have
$$
e(\bi) \epsilon_k e(\bi) = (-1)^{a_k(\bi) + 1} (y_{k+1} - y_k) e(\bi).
$$
\end{Lemma}

\proof Choose arbitrary $\t \in \Tud_n(\bi)$. As $\bi = I_{k,+}^n$ and $i_k = -\frac{1}{2}$, we have $h_k(\bi) = -1$. By~\eqref{remark:h:eq2} we have $h_{k+1}(\bi) = -2$. Hence, by~\autoref{deg:h4:2}, we have either $\t(k) + \t(k+1) = 0$, or $\t(k) + \t(k+1) \neq 0$ and $c_\t(k) - i_k = c_\t(k+1) - i_{k+1}$. Hence, following the same argument as in the proof of~\autoref{essidem:2}, the Lemma holds. \endproof

Next we prove relations (\ref{rela:5:2}) - (\ref{rela:5:4}).

\begin{Lemma} \label{essidem:3}
Suppose $1 \leq k \leq n-1$ and $\bi \in I_{k,0}^n$ with $i_k = -i_{k+1} = \frac{1}{2}$. Then we have
$$
y_{k+1} e(\bi) = y_k e(\bi) - 2 y_k e(\bi) \epsilon_k e(\bi) = y_k e(\bi) - 2 e(\bi) \epsilon_k e(\bi) y_k.
$$
\end{Lemma}

\proof Choose arbitrary $\t \in \Tud_n(\bi)$. As $\bi \in I_{k,0}^n$ and $i_k = \frac{1}{2}$, we have $h_k(\bi) = -1$. By~\eqref{remark:h:eq2} we have $h_{k+1}(\bi) = -2$. Hence, by~\autoref{deg:h4:2}, we have either $\t(k) + \t(k+1) = 0$, or $\t(k) + \t(k+1) \neq 0$ and $c_\t(k) - i_k = c_\t(k+1) - i_{k+1}$.

Suppose $\t(k) + \t(k+1) \neq 0$ and $c_\t(k) - i_k = c_\t(k+1) - i_{k+1}$. By~\autoref{semi:y} and~\autoref{semi:1}, we have
\begin{equation} \label{essidem:3:eq1}
f_{\t\t} y_{k+1}^\O = f_{\t\t} y_k^\O = f_{\t\t} y_k^\O - 2y_k^\O f_{\t\t} \epsilon_k^\O e(\bi)^\O.
\end{equation}

Suppose $\t(k) + \t(k+1) = 0$. We have $c_\t(k) - i_k = -(c_\t(k+1) - i_{k+1})$. As $h_k(\bi) = -1$, by~\autoref{deg:h4:1}, if $\s \in \Tud_n(\bi)$ and $\s \overset{k}\sim \t$, then $\s = \t$. Hence, by~\autoref{PQ:3} and~\autoref{semi:2}, we have
$$
f_{\t\t} \epsilon_k^\O e(\bi)^\O = P_k(\t)^{-1} Q_k(\t)^{-1} e_k(\t,\t) f_{\t\t} = (-1)^{a_k(\bi)} f_{\t\t}.
$$

Therefore, by~\autoref{semi:y}, we have
\begin{align}
f_{\t\t} (y_{k+1}^\O - y_k^\O) & = ((c_\t(k+1) - i_{k+1}) - (c_\t(k) - i_k)) f_{\t\t} = -2(c_\t(k) - i_k) f_{\t\t} \notag \\
& = -2 f_{\t\t} y_k^\O = -2(-1)^{a_k(\bi)} f_{\t\t}y_k^\O \epsilon_k^\O e(\bi)^\O. \label{essidem:3:eq2}
\end{align}

Hence, by~\autoref{idem-semi}, (\ref{essidem:3:eq1}) and (\ref{essidem:3:eq2}) implies
$$
e(\bi)^\O y_{k+1}^\O = e(\bi)^\O y_k^\O - 2y_k^\O e(\bi)^\O \epsilon_k^\O e(\bi)^\O \in \BOx,
$$
and we have $e(\bi)^\O y_{k+1}^\O = e(\bi)^\O y_k^\O - 2 e(\bi)^\O \epsilon_k^\O e(\bi)^\O y_k^\O$ following the same argument. The Lemma follows by lifting the elements into $\B$. \endproof

\begin{Lemma} \label{a_ki}
Suppose $2 \leq k \leq n-1$ and $\bi \in I^n$ with $i_{k-1} = -i_k = i_{k+1}$. We have $(-1)^{a_k(\bi)} = (-1)^{a_{k-1}(\bi) + 1}$ if $\bi \in I_{k,-}^n \cup I_{k,+}^n$ and $(-1)^{a_k(\bi)} = (-1)^{a_{k-1}(\bi)}$ if $\bi \in I_{k,0}^n$.
\end{Lemma}

\begin{proof}
By~\autoref{deg:tab:2}, we have $\deg_{k-1}(\bi) = -\deg_k(\bi)$. Hence $\bi \in I_{k-1,0}^n$ if and only if $\bi \in I_{k,0}^n$, $\bi \in I_{k-1,-}^n$ if and only if $\bi \in I_{k,+}^n$ and $\bi \in I_{k-1,+}^n$ if and only if $\bi \in I_{k,-}^n$.

Choose any $\t \in \mathscr T^{ud}_n(\bi)$ with $c_\t(k-1) = -c_\t(k) = c_\t(k+1)$. Because $i_{k-1} = -i_k = i_{k+1}$, such $\t$ exists. By $e_k^\O = e_k^\O e_{k-1}^\O e_k^\O$, we have $e_k(\t,\t) = e_k(\t,\t) e_{k-1}(\t,\t) e_k(\t,\t)$, which implies $e_k(\t,\t) e_{k-1}(\t,\t) = 1$. By~\autoref{PQ:2} and~\autoref{PQ:3}, we have
$$
1 = P_k(\t) Q_{k-1}(\t) Q_k(\t) P_{k-1}(\t) = P_k(\t) Q_k(\t) P_{k-1}(\t) Q_{k-1}(\t) =
\begin{cases}
-(-1)^{a_k(\bi)} (-1)^{a_{k-1}(\bi)}, & \text{if $\bi \in I_{k,-}^n \cup I_{k,+}^n$,}\\
(-1)^{a_k(\bi)} (-1)^{a_{k-1}(\bi)}, & \text{if $\bi \in I_{k,0}^n$,}
\end{cases}
$$
which proves the Lemma.
\end{proof}

\begin{Lemma} \label{essidem:4}
Suppose $2 \leq k \leq n-1$ and $\bi \in I_{k,0}^n$ with $i_{k-1} = -i_k = i_{k+1} = \frac{1}{2}$. Then we have
$$
e(\bi) = (-1)^{a_k(\bi)}e(\bi) \epsilon_k e(\bi) + 2(-1)^{a_{k-1}(\bi)} e(\bi) \epsilon_{k-1} e(\bi) - e(\bi) \epsilon_{k-1} \epsilon_k e(\bi) - e(\bi) \epsilon_k \epsilon_{k-1} e(\bi).
$$
\end{Lemma}

\proof Because $\bi \in I_{k,0}^n$ and $i_k = -\frac{1}{2}$, we have $h_k(\bi) = -2$ and $h_{k+1}(\bi) = h_{k-1}(\bi) = -h_k(\bi) - 3 = -1$. As $h_{k-1}(\bi) = -1$ and $i_{k-1} = \frac{1}{2}$, we have $\bi \in I_{k-1,0}^n$.

For any $\t \in \Tud_n(\bi)$, write $\lambda = \t_{k-1}$ and $\mu = \t_k$. Let $\alpha = \lambda \ominus \mu$. Because $i_k = -i_{k+1}$, it is easy to see that $\alpha \in \mathscr{AR}_\mu(i_{k+1})$. As $h_k(\bi) = -1$, by \eqref{deg:h4:eq2}, we have $\mathscr{AR}_\mu(i_{k+1}) = \{\alpha\}$. Hence, we have $\t(k+1) = \pm \alpha$, which implies $\t(k) + \t(k+1) = 0$.

Because $h_k(\bi) = -2$, by~\autoref{deg:h4:1}, there exists a unique $\s \in \Tud_n(\bi)$ such that $\s \overset{k}\sim \t$ and $\s \neq \t$. Let $\beta$ and $\gamma$ be positive nodes such that $\s(k) = \pm \beta$ and $\t(k-1) = \s(k-1) = \pm \gamma$. By the construction of up-down tableaux, we have $\alpha, \beta, \gamma \in \mathscr{AR}_\lambda(i_k)$.

Because $h_k(\bi) = -2$, we have $\mathscr{AR}_\lambda(i_k) = 2$ by \eqref{deg:h4:eq1}. Hence, as $\s \neq \t$, we have $\alpha \neq \beta$, which forces $\gamma = \alpha$ or $\gamma = \beta$. Therefore, we have either $\s(k-1) + \s(k) = 0$ or $\t(k-1) + \t(k) = 0$.

If $\t(k-1) + \t(k) \neq 0$, then we have $\s(k-1) + \s(k) = 0$. As $\bi \in I_{k-1,0}^n \cap I_{k,0}^n$, we have
\begin{align*}
& f_{\t\t} \epsilon_{k-1}^\O e(\bi)^\O = 0, &&  f_{\t\t} \epsilon_k^\O e(\bi)^\O = (-1)^{a_k(\bi)} f_{\t\t} + P_k(\t)^{-1} Q_k(\s)^{-1} e_k(\t,\s) f_{\t\s}, \\
& f_{\t\t} \epsilon_{k-1}^\O \epsilon_k^\O e(\bi)^\O = 0, && f_{\t\t} \epsilon_k^\O \epsilon_{k-1}^\O e(\bi)^\O = (-1)^{a_{k-1}(\bi)} P_k(\t)^{-1} Q_k(\s)^{-1} e_k(\t,\s) f_{\t\s},
\end{align*}
by~\autoref{PQ:3},~\autoref{semi:1} and~\autoref{semi:2}. Hence, $f_{\t\t}((-1)^{a_k(\bi)} \epsilon_k^\O + 2(-1)^{a_{k-1}(\bi)} \epsilon_{k-1} - \epsilon_{k-1}^\O\epsilon_k^\O - \epsilon_k^\O \epsilon_{k-1}^\O) e(\bi)^\O$ equals
$$
f_{\t\t} + ((-1)^{a_k(\bi)} - (-1)^{a_{k-1}(\bi)}) P_k(\t)^{-1} Q_k(\s)^{-1} e_k(\t,\s) f_{\t\s}.
$$

By~\autoref{a_ki}, we have $(-1)^{a_k(\bi)} - (-1)^{a_{k-1}(\bi)} = 0$ because $\bi \in I_{k,0}^n$. Hence, we have
\begin{equation} \label{essidem:4:eq1}
f_{\t\t}((-1)^{a_k(\bi)} \epsilon_k^\O + 2(-1)^{a_{k-1}(\bi)} \epsilon_{k-1} - \epsilon_{k-1}^\O\epsilon_k^\O - \epsilon_k^\O \epsilon_{k-1}^\O) e(\bi)^\O = f_{\t\t},
\end{equation}
when $\t(k-1) + \t(k) \neq 0$.

If $\t(k-1) + \t(k) = 0$, then we have $\s(k-1) + \s(k) \neq 0$. As $\bi \in I_{k-1,0}^n \cap I_{k,0}^n$, we have
\begin{align*}
& f_{\t\t} \epsilon_{k-1}^\O e(\bi)^\O = (-1)^{a_{k-1}(\bi)} f_{\t\t}, &&  f_{\t\t} \epsilon_k^\O e(\bi)^\O = (-1)^{a_k(\bi)} f_{\t\t} + P_k(\t)^{-1} Q_k(\s)^{-1} e_k(\t,\s) f_{\t\s}, \\
& f_{\t\t} \epsilon_{k-1}^\O \epsilon_k^\O e(\bi)^\O = (-1)^{a_{k-1}(\bi) + a_k(\bi)} f_{\t\t} + (-1)^{a_{k-1}(\bi)} P_k(\t)^{-1} Q_k(\s)^{-1} e_k(\t,\s) f_{\t\s}, && f_{\t\t} \epsilon_k^\O \epsilon_{k-1}^\O e(\bi)^\O = (-1)^{a_{k-1}(\bi) + a_k(\bi)} f_{\t\t},
\end{align*}
by~\autoref{PQ:3},~\autoref{semi:1} and~\autoref{semi:2}. Hence, $f_{\t\t}((-1)^{a_k(\bi)} \epsilon_k^\O + 2(-1)^{a_{k-1}(\bi)} \epsilon_{k-1} - \epsilon_{k-1}^\O\epsilon_k^\O - \epsilon_k^\O \epsilon_{k-1}^\O) e(\bi)^\O$ equals
$$
((-1)^{2a_k(\bi)} + 2(-1)^{2a_{k-1}(\bi)} - 2(-1)^{a_{k-1}(\bi) + a_k(\bi)}) f_{\t\t} + ((-1)^{a_k(\bi)} - (-1)^{a_{k-1}(\bi)}) P_k(\t)^{-1} Q_k(\s)^{-1} e_k(\t,\s) f_{\t\s}.
$$

By~\autoref{a_ki}, we have $(-1)^{a_k(\bi)} = (-1)^{a_{k-1}(\bi)}$ because $\bi \in I_{k,0}^n$. Hence, we have
\begin{equation} \label{essidem:4:eq2}
f_{\t\t}((-1)^{a_k(\bi)} \epsilon_k^\O + 2(-1)^{a_{k-1}(\bi)} \epsilon_{k-1} - \epsilon_{k-1}^\O\epsilon_k^\O - \epsilon_k^\O \epsilon_{k-1}^\O) e(\bi)^\O = f_{\t\t},
\end{equation}
when $\t(k-1) + \t(k) = 0$.

By (\ref{essidem:4:eq1}) and (\ref{essidem:4:eq2}), we have
$$
f_{\t\t}((-1)^{a_k(\bi)} \epsilon_k^\O + 2(-1)^{a_{k-1}(\bi)} \epsilon_{k-1} - \epsilon_{k-1}^\O\epsilon_k^\O - \epsilon_k^\O \epsilon_{k-1}^\O) e(\bi)^\O = f_{\t\t},
$$
which implies
$$
e(\bi)^\O = (-1)^{a_k(\bi)}e(\bi)^\O \epsilon_k^\O e(\bi)^\O + 2(-1)^{a_{k-1}(\bi)} e(\bi)^\O \epsilon_{k-1}^\O e(\bi)^\O - e(\bi)^\O \epsilon_{k-1}^\O \epsilon_k^\O e(\bi)^\O - e(\bi)^\O \epsilon_k^\O \epsilon_{k-1}^\O e(\bi)^\O \in \BOx.
$$
by~\autoref{idem-semi}. Hence the Lemma holds by lifting the elements into $\B$. \endproof

\begin{Lemma} \label{essidem:5}
Suppose $1 \leq k \leq n-1$ and $\bi \in I_{k,-}^n$ with $i_k + i_{k+1} = 0$. Then we have $e(\bi) (\epsilon_k y_k + y_k \epsilon_k) e(\bi) = (-1)^{a_k(\bi)} e(\bi)$.
\end{Lemma}

\proof Choose an arbitrary $\t \in \Tud_n(\bi)$ and write $\lambda = \t_k$, $\mu = \t_{k+1}$. As $\bi \in I_{k,-}^n$, we have $h_k(\bi) = -2$, and $h_{k+1}(\bi) = -1$ if $i_k = \pm \frac{1}{2}$ and $h_{k+1}(\bi) = 0$ if $i_k \neq \pm \frac{1}{2}$. In both cases, by \eqref{deg:h4:eq2} and \eqref{deg:h4:eq3}, we have $|\mathscr{AR}_\mu(i_{k+1})| = 1$. Let $\alpha = \lambda \ominus \mu$. Because $i_k + i_{k+1} = 0$, by the construction of up-down tableaux, we have $\alpha \in \mathscr{AR}_\mu(i_{k+1})$. Hence, it forces $\t(k) + \t(k+1) = 0$.

Because $h_k(\bi) = -2$, by~\autoref{deg:h4:1}, there exists a unique $\s \in \Tud_n(\bi)$ such that $\s \overset{k}\sim \t$ and $\s \neq \t$, and $c_\t(k) - i_k = -(c_\s(k) - i_k)$. Hence, by~\autoref{PQ:3}, we have
\begin{align*}
f_{\t\t} (\epsilon_k^\O y_k^\O + y_k^\O \epsilon_k^\O) e(\bi)^\O & = 2(c_\t(k) - i_k) f_{\t\t} \epsilon_k^\O \f{\t} + (c_\s(k) - i_k + c_\t(k) - i_k) f_{\t\t} \epsilon_k^\O \f{\s}\\
& = 2(c_\t(k) - i_k) f_{\t\t} \epsilon_k^\O \f{\t} = (-1)^{a_k(\bi)} f_{\t\t}.
\end{align*}

As $\t$ is chosen arbitrary, by~\autoref{idem-semi}, we have
$$
e(\bi)^\O (\epsilon_k^\O y_k^\O + y_k^\O \epsilon_k^\O) e(\bi)^\O = (-1)^{a_k(\bi)} e(\bi)^\O \in \BOx,
$$
which proves the Lemma by lifting the elements into $\B$. \endproof

Next we prove the relation (\ref{rela:5:5}) holds in $\B$. We remind reader that when we prove this relation we do not assume $\bi \in I^n$. In this relation, when $\bi \not\in I^n$, we have $e(\bj) \epsilon_k e(\bi) \epsilon_k e(\bk) = 0$ because $e(\bi) = 0$ by~\autoref{idem-semi}, but it is not easy to see the left hand side of the relation equals $0$, except when $\bi \in P_{k,-}^n$. Here we start by proving the case when $\bi \in P_{k,-}^n$.



\begin{Lemma} \label{essidem:7}
Suppose $1 \leq k \leq n-1$ and $\bi \in P_{k,-}^n$, $\bj, \bk \in I^n$. Then we have $e(\bj) \epsilon_k e(\bi) \epsilon_k e(\bk) = 0$.
\end{Lemma}

\proof When $\bi \not\in I^n$, we have $e(\bi) = 0$ by~\autoref{idem-semi}, which implies $e(\bj) \epsilon_k e(\bi) \epsilon_k e(\bk) = 0$. Hence we assume $\bi \in I^n$ in the rest of the proof. Note that $\bi \in I^n$ and $\bi \in P_{k,-}^n$ is equivalent to $\bi \in I_{k,-}^n$.

Choose arbitrary $\u \in \Tud_n(\bj)$. If $\u(k) + \u(k+1) \neq 0$, we have
\begin{equation} \label{essidem:7:eq1}
f_{\u\u} \epsilon_k^\O e(\bi)^\O \epsilon_k^\O e(\bk)^\O = 0.
\end{equation}

Assume $\u(k) + \u(k+1) = 0$. Because $\bi \in I_{k,-}^n$, we have $h_k(\bi) = -2$. By~\autoref{deg:h4:1}, there exist exactly two up-down tableaux $\s,\t \in \Tud_n(\bi)$ such that $\s \overset{k}\sim \u \overset{k}\sim \t$, and $c_\t(k) - i_k = - (c_\s(k) - i_k)$. Therefore, by~\autoref{semi:B},~\autoref{semi:2} and~\autoref{PQ:3}, we have
\begin{align*}
f_{\u\u} \epsilon_k^\O e(\bi)^\O \epsilon_k^\O e(\bk)^\O & = \sum_{\substack{\v \in \Tud_n(\bk)\\ \v \overset{k}\sim \t}} P_k(\u)^{-1} Q_k(\v)^{-1} e_k(\u,\v) (P_k(\t)^{-1}Q_k(\t)^{-1} e_k(\t,\t) + P_k(\s)^{-1} Q_k(\s)^{-1} e_k(\s,\s))  f_{\u\v}\\
& = \sum_{\substack{\v \in \Tud_n(\bk)\\ \v \overset{k}\sim \t}} P_k(\u)^{-1} Q_k(\v)^{-1} e_k(\u,\v) (\frac{1}{2(c_\t(k) - i_k)} + \frac{1}{2(c_\s(k) - i_k)})  f_{\u\v} = 0.
\end{align*}

By combining the above equality and (\ref{essidem:7:eq1}), we have $f_{\u\u} \epsilon_k^\O e(\bi)^\O \epsilon_k^\O e(\bk)^\O = 0$, and by~\autoref{idem-semi}, we have $e(\bj)^\O \epsilon_k^\O e(\bi)^\O \epsilon_k^\O e(\bk)^\O = 0$. The Lemma follows by lifting the elements into $\B$. \endproof

Next we prove the cases when $\bi \in P_{k,0}^n$. Recall
$$
z_k(\bi) =
\begin{cases}
0, & \text{if $h_k(\bi) < -2$, or $h_k(\bi) \geq 0$ and $i_k \neq 0$,}\\
(-1)^{a_k(\bi)} (1 + \delta_{i_k, -\frac{1}{2}}), & \text{if $-2 \leq h_k(\bi) < 0$,}\\
\frac{1 + (-1)^{a_k(\bi)}}{2}, & \text{if $i_k = 0$.}
\end{cases}
$$

The key point of $z_k(\bi)$ is given in the next Lemma.

\begin{Lemma} \label{essidem:6:2}
Suppose $1 \leq k \leq n-1$ and $\bi \in P_{k,0}^n$, $\bj\in I^n$ with $\bi|_{k-1} = \bj|_{k-1}$. For any $\u \in \Tud_n(\bj)$ with $\u(k) + \u(k+1) = 0$, we have the following properties:
\begin{enumerate}
\item If $h_k(\bi) < -2$, or $h_k(\bi) \geq 0$ and $i_k \neq 0$, then there does not exist $\t \in \Tud_n(\bi)$ such that $\t \overset{k}\sim \u$.

\item If $-2 \leq h_k(\bi) < 0$, then there exists $\t \in \Tud_n(\bi)$ such that $\t \overset{k}\sim \u$.

\item If $i_k = 0$, then
$$
z_k(\bi) =
\begin{cases}
0, & \text{if there does not exist $\t \in \Tud_n(\bi)$ such that $\t \overset{k}\sim \u$,}\\
(-1)^{a_k(\bi)}, & \text{if there exists $\t \in \Tud_n(\bi)$ such that $\t \overset{k}\sim \u$.}
\end{cases}
$$
\end{enumerate}
\end{Lemma}

\proof Write $\lambda = \u_{k-1}$.

(1). By~\autoref{deg:h2}, we have $\bi \not\in I^n$ if $h_k(\bi) < -2$ or $h_k(\bi) > 0$. When $h_k(\bi) = 0$, assume $i_k \neq -\frac{1}{2}$. Then we have $\bi \in P_{k,+}^n$ by the definition and $i_k \neq 0$, which contradicts to the assumptions of the Lemma. Hence we have $i_k = -\frac{1}{2}$, which implies $\bi \not\in I^n$ by~\autoref{deg:h3}. Therefore, in this case we always have $\bi \not\in I^n$, which implies there does not exist $\t \in \Tud_n(\bi)$ such that $\t \overset{k}\sim \u$.

(2). In this case we have $h_k(\bi) = -1$ or $-2$. By~\autoref{deg:h1:1} and~\eqref{deg:h2:ineq} it forces $|\mathscr{AR}_\lambda(i_k)| \geq 1$. Let $\alpha \in \mathscr{AR}_\lambda(i_k)$. Without loss of generality we assume $\alpha \in \mathscr A(\lambda)$. Hence let $\t$ be the up-down tableau such that $\t(k) = \alpha$, $\t(k+1) = -\alpha$ and $\t(\l) = \u(\l)$ for any $\l \neq k,k+1$. Therefore $\t \in \Tud_n(\bi)$ and $\t \overset{k}\sim \u$.

(3). Because $i_k = 0$, we have $h_k(\bi) = 0$ by the definition of $h_k$. By~\autoref{deg:h1:1} and~\eqref{deg:h2:ineq} it forces $|\mathscr{AR}_\lambda(i_k)| = 0$ or $1$. If there exists $\t \in \Tud_n(\bi)$ such that $\t \overset{k}\sim \u$, we require $|\mathscr{AR}_\lambda(0)| = 1$, which implies $z_k(\bi) = 1 = (-1)^{a_k(\bi)}$ by~\autoref{deg:h4:5}; and if there does not exist $\t \in \Tud_n(\bi)$ such that $\t \overset{k}\sim \u$, we have $|\mathscr{AR}_\lambda(0)| = 0$, which implies $z_k(\bi) = 0$ by~\autoref{deg:h4:5}. \endproof

Therefore, under the assumptions of~\autoref{essidem:6:2} we can rewrite $z_k(\bi)$ as
\begin{equation} \label{essidem:6:eq1}
z_k(\bi) =
\begin{cases}
0, & \text{if there does not exist $\t \in \Tud_n(\bi)$ such that $\t \overset{k}\sim \u$,}\\
(-1)^{a_k(\bi)} (1 + \delta_{i_k, -\frac{1}{2}}), & \text{if there exists $\t \in \Tud_n(\bi)$ such that $\t \overset{k}\sim \u$,}
\end{cases}
\end{equation}
which allows us to verify relation~\eqref{rela:5:5} when $\bi \in P_{k,0}^n$ via direct calculations.

\begin{Lemma} \label{essidem:6:1}
Suppose $1 \leq k \leq n-1$ and $\bi \in P_{k,0}^n$, $\bj,\bk \in I^n$. For any $\u \in \Tud_n(\bj)$ with $\u(k) + \u(k+1) = 0$, we have
$$
f_{\u\u} \epsilon_k^\O e(\bi)^\O \epsilon_k^\O e(\bk)^\O = z_k(\bi) f_{\u\u} \epsilon_k^\O e(\bk)^\O.
$$
\end{Lemma}

\proof Suppose there does not exist $\t \in \Tud_n(\bi)$ such that $\t \overset{k}\sim \u$. Then by~\autoref{semi:1} and~\eqref{essidem:6:eq1} we have
$$
f_{\u\u} \epsilon_k^\O e(\bi)^\O \epsilon_k^\O e(\bk)^\O = 0 = z_k(\bi) f_{\u\u} \epsilon_k^\O e(\bk)^\O.
$$

Suppose there exists $\t \in \Tud_n(\bi)$ such that $\t \overset{k}\sim \u$. It requires $\bi \in I_{k,0}^n$. When $h_k(\bi) = -1$, we have $i_k \neq -\frac{1}{2}$, otherwise $\bi \in P_{k,+}^n$ which contradicts to the assumptions of the Lemma. Therefore we have $z_k(\bi) = (-1)^{a_k(\bi)}$.

By~\autoref{deg:h4:1}, there exists a unique $\t \in \Tud_n(\bi)$ such that $\t \overset{k}\sim \u$. Hence, by~\autoref{semi:B},~\autoref{semi:2} and~\autoref{PQ:3}, we have
\begin{align}
f_{\u\u} \epsilon_k^\O e(\bi)^\O \epsilon_k^\O e(\bk)^\O & = \sum_{\substack{\v \in \Tud_n(\bk)\\ \v \overset{k}\sim \t}} P_k(\u)^{-1} e_k(\u,\t) Q_k(\t)^{-1} P_k(\t)^{-1} Q_k(\v)^{-1} e_k(\t,\v) f_{\u\v} \notag \\
& = (-1)^{a_k(\bi)} \sum_{\substack{\v \in \mathscr T^{ud}_n(\bk)\\ \v \overset{k}\sim \u}} P_k(\u)^{-1} e_k(\u,\v) Q_k(\v)^{-1} f_{\u\v} = (-1)^{a_k(\bi)} f_{\u\u} \epsilon_k^\O e(\bk)^\O = z_k(\bi) f_{\u\u} \epsilon_k^\O e(\bk)^\O. \notag
\end{align}

When $h_k(\bi) = 0$, we have $i_k = 0$. Hence following the same argument as when $h_k(\bi) = -1$, we have
$$
f_{\u\u} \epsilon_k^\O e(\bi)^\O \epsilon_k^\O e(\bk)^\O = (-1)^{a_k(\bi)} f_{\u\u} \epsilon_k^\O e(\bk)^\O = z_k(\bi) f_{\u\u} \epsilon_k^\O e(\bk)^\O.
$$

When $h_k(\bi) = -2$, we have $i_k = -\frac{1}{2}$. Hence following the same argument as when $h_k(\bi) = -1$, we have
$$
f_{\u\u} \epsilon_k^\O e(\bi)^\O \epsilon_k^\O e(\bk)^\O = 2(-1)^{a_k(\bi)} f_{\u\u} \epsilon_k^\O e(\bk)^\O = z_k(\bi) f_{\u\u} \epsilon_k^\O e(\bk)^\O,
$$
where the coefficient is doubled because we have $\s,\t \in \Tud_n(\bi)$ such that $\s \overset{k}\sim \u \overset{k}\sim \t$, and when $h_k(\bi) = -1$, there is a unique $\t \in \Tud_n(\bi)$ such that $\t \overset{k}\sim \u$. \endproof

\begin{Lemma} \label{essidem:6}
Suppose $1 \leq k \leq n-1$ and $\bi \in P_{k,0}^n$, $\bj, \bk \in I^n$. Then we have $e(\bj) \epsilon_k e(\bi) \epsilon_k e(\bk) = z_k(\bi) e(\bj) \epsilon_k e(\bk)$.
\end{Lemma}

\proof Choose arbitrary $\u \in \Tud_n(\bj)$. If $\u(k) + \u(k+1) \neq 0$, we have $f_{\u\u} \epsilon_k^\O e(\bi)^\O \epsilon_k^\O e(\bk)^\O = 0 = z_k(\bi) f_{\u\u} \epsilon_k^\O e(\bk)^\O$; and if $\u(k) + \u(k+1) = 0$, we have $f_{\u\u} \epsilon_k^\O e(\bi)^\O \epsilon_k^\O e(\bk)^\O = z_k(\bi) f_{\u\u} \epsilon_k^\O e(\bk)^\O$ by~\autoref{essidem:6:1}. Therefore, we have
$$
e(\bj)^\O \epsilon_k^\O e(\bi)^\O \epsilon_k^\O e(\bk)^\O = z_k(\bi) e(\bj)^\O \epsilon_k^\O e(\bk)^\O \in \BOx,
$$
by~\autoref{idem-semi}, which proves the Lemma by lifting the elements into $\B$. \endproof


To prove (\ref{rela:5:5}) when $\bi \in P_{k,+}^n$, we want to give a result analogue to~\autoref{essidem:6:2}. In more details, we want to prove the next Lemma:

\begin{Lemma} \label{essidem:h5}
Suppose $\bi \in P_{k,+}^n$ and $\bj \in I^n$ with $\bi \overset{k}\sim \bj$ for some $1 \leq k \leq n$. For any $\u \in \Tud_n(\bj)$ with $\u(k) + \u(k+1) = 0$, we have
\begin{eqnarray*}
&& (1 + \delta_{i_k, -\frac{1}{2}}) \left( \sum_{\l \in A_{k,1}^\bi} y_\l^\O - 2 \sum_{\l \in A_{k,2}^\bi} y_\l^\O + \sum_{\l \in A_{k,3}^\bi} y_\l^\O - 2\sum_{\l \in A_{k,4}^\bi} y_\l^\O \right) f_{\u\u}\\
& = &
\begin{cases}
0, & \text{if there does not exist $\t \in \Tud_n(\bi)$ such that $\t \overset{k}\sim \u$,}\\
2(c_\t(k) - i_k) f_{\u\u}, & \text{if there exists $\t \in \Tud_n(\bi)$ such that $\t \overset{k}\sim \u$.}
\end{cases}
\end{eqnarray*}
\end{Lemma}

Under the assumptions of~\autoref{essidem:h5}, write $\lambda = \u_{k-1}$.~\autoref{essidem:h5} explicitly express the connection between $\lambda$ and $\mathscr{AR}_\lambda(i_k)$. Hence we express such connection before we prove~\autoref{essidem:h5}.

For $1 \leq k \leq n$, fix $\bi \in P_{k,+}^n$, $\bj \in I^n$ with $\bi \overset{k}\sim \bj$ and $\u \in \Tud_n(\bj)$. For convenience we write $i = i_k \in P$. Because $\bi \in P_{k,+}^n$ we have $-1 \leq h_k(\bi) \leq 0$.

Let $\u_{k-1} = \lambda$. By~\autoref{deg:h1:1} and~\eqref{deg:h2:ineq}, we have three possible cases:
\begin{enumerate}
\item $h_k(\bi) = 0$ and $|\mathscr{AR}_\lambda(-i)| = |\mathscr{AR}_\lambda(i)| = 0$.

\item $h_k(\bi) = 0$ and $|\mathscr{AR}_\lambda(-i)| = |\mathscr{AR}_\lambda(i)| = 1$.

\item $h_k(\bi) = -1$ and $|\mathscr{AR}_\lambda(-i)| = 0$, $|\mathscr{AR}_\lambda(i)| = 1$.
\end{enumerate}

For $\alpha \in \mathscr{AR}(\lambda)$, define
$$
\cont_\lambda(\alpha) =
\begin{cases}
\cont(\alpha), & \text{if $\alpha \in \mathscr A(\lambda)$,}\\
-\cont(\alpha), & \text{if $\alpha \in \mathscr R(\lambda)$,}
\end{cases}
$$
and denote $\cont_\lambda(i) = \cont_\lambda(\alpha)$ where $\alpha \in \mathscr{AR}_\lambda(i)$ if $|\mathscr{AR}_\lambda(i)| = 1$, and $\cont_\lambda(i) = \frac{x - \delta}{2} + i$ if $|\mathscr{AR}_\lambda(i)| = 0$.

By~\autoref{Naz} and~\autoref{e_k:h1}, set $u$ as an unknown and we have
\begin{equation} \label{essidem:h1:eq1}
\frac{u + c_\u(1)}{u - c_\u(1)} \prod_{\l = 1}^{k-1} \frac{u+c_\u(\l) + 1}{u - c_\u(\l) + 1} \frac{u + c_\u(\l) - 1}{u - c_\u(\l) - 1} \frac{(u - c_\u(\l))^2}{(u + c_\u(\l))^2} = \prod_{\alpha \in \mathscr {AR}(\lambda)} \frac{u + \cont_\lambda(\alpha)}{u - \cont_\lambda(\alpha)}.
\end{equation}

We note that the right hand side of~\eqref{essidem:h1:eq1} contains information of $\mathscr{AR}_\lambda(i)$, which is accessible from the left hand side of~\eqref{essidem:h1:eq1}.

Define subsets of $R[u,x]$ as
\begin{align*}
S_1(\u) & := \{u + c_\u(1)\} \cup \bigcup_{\l=1}^{k-1} \{u + c_\u(\l) + 1, u + c_\u(\l) - 1\}, &
S_2(\u) & := \set{u - c_\u(\l) | 1 \leq \l \leq k-1},\\
S_3(\u) & := \{u - c_\u(1)\} \cup \bigcup_{\l=1}^{k-1} \{u - c_\u(\l) + 1, u - c_\u(\l) - 1\}, &
S_4(\u) & := \set{u + c_\u(\l) | 1 \leq \l \leq k-1},\\
T_1(\u) & := \set{u + \cont_\lambda(\alpha) | \alpha \in \mathscr{AR}(\lambda)}, &
T_2(\u) & := \set{u - \cont_\lambda(\alpha) | \alpha \in \mathscr{AR}(\lambda)},
\end{align*}

By the definitions, for $w[u,x] \in \bigcup_{r = 1}^4 S_{k,r}$, we have $w(i,\delta) \in R$. Define
\begin{align*}
S_1^i(\u) & := \set{w(u,x) \in S_1(\u) | w(i,\delta) = 0 }, &
S_2^i(\u) & := \set{w(u,x) \in S_2(\u) | w(i,\delta) = 0 },\\
S_3^i(\u) & := \set{w(u,x) \in S_3(\u) | w(i,\delta) = 0 }, &
S_4^i(\u) & := \set{w(u,x) \in S_4(\u) | w(i,\delta) = 0 },\\
T_1^i(\u) & := \set{w(u,x) \in T_1(\u) | w(i,\delta) = 0}, &
T_2^i(\u) & := \set{w(u,x) \in T_2(\u) | w(i,\delta) = 0}.
\end{align*}


\begin{Remark} \label{essidem:h1:remark}
By the definitions, we have $c_\u(\l) = \pm \frac{x - \delta}{2} + j$ for some $j \in P$. Therefore, for $w(u,x) \in \bigcup_{r=1}^4 S_r(\u)$, we have $w(u,x) = u \pm \frac{x - \delta}{2} + j$ for some $j \in P$. Hence, $w(i,\delta) = 0$ only if $w(u,x) = u \pm \frac{x - \delta}{2} - i$.

By the definitions, we have $\cont_\lambda(i) = \pm \frac{x - \delta}{2} + i$. Therefore, for $w(u,x) \in \bigcup_{r=1}^4 S^i_r(\u)$, we have $w(\cont_\lambda(i), x) = 2(\cont_\lambda(i) - i)$ or $0$. Similarly, we have the same property for $w(u,x) \in T^i_1(\u) \cup T^i_2(\u)$.
\end{Remark}

\begin{Lemma} \label{essidem:h1}
We have the following equality
$$
\frac{\prod_{w \in S_1^i(\u)} w \prod_{w \in S_2^i(\u)} w^2}{\prod_{w \in S_3^i(\u)} w \prod_{w \in S_4^i(\u)} w^2} = \frac{\prod_{w \in T_1^i(\u)} w}{\prod_{w \in T_2^i(\u)} w}.
$$
\end{Lemma}

\proof For convenience, we define
\begin{align*}
f_1(u,x) & = \frac{\prod_{w \in S_1(\u)} w \prod_{w \in S_2(\u)} w^2}{\prod_{w \in S_3(\u)} w \prod_{w \in S_4(\u)} w^2}, & f_2(u,x) & = \frac{\prod_{w \in T_1(\u)} w}{\prod_{w \in T_2(\u)} w},\\
g_1(u,x) & = \frac{\prod_{w \in S_1^i(\u)} w \prod_{w \in S_2^i(\u)} w^2}{\prod_{w \in S_3^i(\u)} w \prod_{w \in S_4^i(\u)} w^2}, & g_2(u,x) & = \frac{\prod_{w \in T_1^i(\u)} w}{\prod_{w \in T_2^i(\u)} w},
\end{align*}
and $u_1(u,x) = f_1(u,x) / g_1(u,x)$ and $u_2(u,x) = f_2(u,x) / g_2(u,x)$.

By (\ref{essidem:h1:eq1}), we have $f_1(u,x) = f_2(u,x)$. Hence, by~\autoref{essidem:h1:remark}, we can write
$$
f_1(u,\delta) = f_2(u,\delta) = \prod_{j \in P} (u - j)^{a_j},
$$
where $a_j \in \Z$. By the definitions, we have $u_1(i,\delta) \neq 0$ and $u_2(i,\delta) \neq 0$. Therefore, we have
$$
u_1(u,\delta) = \prod_{\substack{j \in P \\ j \neq i}} (u - j)^{a_j} = u_2(u,\delta),
$$
which yields $g_1(u,\delta) = g_2(u,\delta) = (u - i)^{a_i}$.

By~\autoref{essidem:h1:remark}, we can write
$$
g_1(u,x) = (u + \frac{\delta-1}{2} - i)^{k_1} (u - \frac{\delta-1}{2} - i)^{l_1} \qquad \text{and} \qquad g_2(u,x) = (u + \frac{\delta-1}{2} - i)^{k_2} (u - \frac{\delta-1}{2} - i)^{l_2},
$$
where $k_1, l_1, k_2, l_2 \in \Z$ and $k_1 + l_1 = k_2 + l_2 = a_i$. If we have $k_1 > k_2$, then by defining
$$
s_1(u,x) = (u + \frac{\delta-1}{2} - i)^{-k_2}f_1(u,x), \qquad \text{and} \qquad s_2(u,x) = (u + \frac{\delta-1}{2} - i)^{-k_2}f_2(u,x),
$$
we have $s_1(-\frac{\delta-1}{2} + i, x) = 0$ and $s_2(-\frac{\delta-1}{2} + i, x) \neq 0$, which contradicts the fact that $f_1(u,x) = f_2(u,x)$. Hence, we have $k_1 \leq k_2$. Similarly, we have $k_1 \geq k_2$, which implies $k_1 = k_2$. Therefore, we have $l_1 = l_2$, and $g_1(u,x) = g_2(u,x)$. \endproof

\begin{Lemma} \label{essidem:h2}
We have the following equality
$$
\frac{\prod_{w \in S_1^i(\u)} w \prod_{w \in S_2^i(\u)} w^2}{\prod_{w \in S_3^i(\u)} w \prod_{w \in S_4^i(\u)} w^2} =
\begin{cases}
\frac{1}{u - \cont_\lambda(i)}, & \text{if $h_k(\bi) = -1$,}\\
\frac{u + \cont_\lambda(i) - 2i}{u - \cont_\lambda(i)}, & \text{if $h_k(\bi) = 0$ and $|\mathscr{AR}_\lambda(i)| = 1$,}\\
1, & \text{if $h_k(\bi) = 0$ and $|\mathscr{AR}_\lambda(i)| = 0$.}
\end{cases}
$$
\end{Lemma}

\proof By~\autoref{essidem:h1}, it suffices to prove that
$$
\frac{\prod_{w \in T_1^i(\u)} w}{\prod_{w \in T_2^i(\u)} w} =
\begin{cases}
\frac{1}{u - \cont_\lambda(i)}, & \text{if $h_k(\bi) = -1$,}\\
\frac{u + \cont_\lambda(i) - 2i}{u - \cont_\lambda(i)}, & \text{if $h_k(\bi) = 0$ and $|\mathscr{AR}_\lambda(i)| = 1$,}\\
1, & \text{if $h_k(\bi) = 0$ and $|\mathscr{AR}_\lambda(i)| = 0$.}
\end{cases}
$$

For $h_k(\bi) = -1$, by (\ref{deg:h4:eq2}), we have $|\mathscr{AR}_\lambda(i)| = 1$. Hence, there exists a unique $\alpha \in \mathscr{AR}_\lambda(i)$. By the definition, we have $\cont_\lambda(i) = \cont_\lambda(\alpha)$. Hence, we have
$$
\frac{\prod_{w \in T_1^i(\u)} w}{\prod_{w \in T_2^i(\u)} w} = \frac{1}{u - \cont_\lambda(\alpha)} = \frac{1}{u - \cont_\lambda(i)}.
$$

For $h_k(\bi) = 0$ and $|\mathscr{AR}_\lambda(i)| = 1$, by (\ref{deg:h4:eq3}), we have $|\mathscr{AR}_\lambda(i)| = |\mathscr{AR}_\lambda(-i)| = 1$. Therefore, there exist $\alpha \in \mathscr{AR}_\lambda(i)$ and $\beta \in \mathscr{AR}_\lambda(-i)$. By the definition, we have $\cont_\lambda(i) = \cont_\lambda(\alpha)$. If $i = 0$, we have $\alpha = \beta$ because $i = -i$. Hence, we have
$$
\frac{\prod_{w \in T_1^i(\u)} w}{\prod_{w \in T_2^i(\u)} w} = \frac{u + \cont_\lambda(\alpha)}{u - \cont_\lambda(\alpha)} = \frac{u + \cont_\lambda(i) - 2i}{u - \cont_\lambda(i)};
$$
and if $i \neq 0$, by~\autoref{deg:h4}, we have $\alpha, \beta \in \mathscr A(\lambda)$ or $\alpha, \beta \in \mathscr R(\lambda)$. Therefore, we have $\cont_\lambda(\alpha) - i = \cont_\lambda(\beta) + i$, which implies $\cont_\lambda(\beta) = \cont_\lambda(\alpha) - 2i = \cont_\lambda(i) - 2i$. Hence, we have
$$
\frac{\prod_{w \in T_1^i(\u)} w}{\prod_{w \in T_2^i(\u)} w} = \frac{u + \cont_\lambda(\beta)}{u - \cont_\lambda(\alpha)} = \frac{u + \cont_\lambda(i) - 2i}{u - \cont_\lambda(i)}.
$$

For $h_k(\bi) = 0$ and $|\mathscr{AR}_\lambda(i)| = 0$, we have $|\mathscr{AR}_\lambda(i)| = |\mathscr{AR}_\lambda(-i)| = 0$. Hence, we have
$$
\frac{\prod_{w \in T_1^i(\u)} w}{\prod_{w \in T_2^i(\u)} w} = 1,
$$
which completes the proof. \endproof

\autoref{essidem:h2} gives us a method to determine $|\mathscr{AR}_\lambda(i)|$ by giving $\lambda$, but it is not sufficient to determine whether $\alpha \in \mathscr A(\lambda)$ or $\alpha \in \mathscr R(\lambda)$ for $\alpha \in \mathscr{AR}_\lambda(i)$ given $|\mathscr{AR}_\lambda(i)| = 1$.

Recall that by~\autoref{essidem:h1:remark}, for $w(u,x) \in \bigcup_{r = 1}^4 S^i_r(\u)$, we have $w(\cont_\lambda(i),x) = 2(\cont_\lambda(i) - i)$ or $0$. Define
$$
\begin{array}{cc}
a_{1,1} = \#\set{w \in S^i_1(\u) | w(\cont_\lambda(i),x) = 2(\cont_\lambda(i) - i)}, & a_{1,2} = \#\set{w \in S^i_1(\u) | w(\cont_\lambda(i),x) = 0},\\
a_{2,1} = \#\set{w \in S^i_2(\u) | w(\cont_\lambda(i),x) = 2(\cont_\lambda(i) - i)}, & a_{2,2} = \#\set{w \in S^i_2(\u) | w(\cont_\lambda(i),x) = 0},\\
a_{3,1} = \#\set{w \in S^i_3(\u) | w(\cont_\lambda(i),x) = 2(\cont_\lambda(i) - i)}, & a_{3,2} = \#\set{w \in S^i_3(\u) | w(\cont_\lambda(i),x) = 0},\\
a_{4,1} = \#\set{w \in S^i_4(\u) | w(\cont_\lambda(i),x) = 2(\cont_\lambda(i) - i)}, & a_{4,2} = \#\set{w \in S^i_4(\u) | w(\cont_\lambda(i),x) = 0}.
\end{array}
$$

By~\autoref{essidem:h2}, we have
\begin{equation} \label{essidem:eq1}
\begin{array}{rl}
a_{1,1} + 2a_{2,1} - a_{3,1} - 2a_{4,1} & =
\begin{cases}
1, & \text{if $h_k(\bi) = 0$ and $|\mathscr{AR}_\lambda(i)| = 1$,} \\
0, & \text{otherwise;}
\end{cases},\\
a_{1,2} + 2a_{2,2} - a_{3,2} - 2a_{4,2} & =
\begin{cases}
0, & \text{if $h_k(\bi) = 0$ and $|\mathscr{AR}_\lambda(i)| = 0$,}\\
-1, & \text{otherwise.}
\end{cases}
\end{array}
\end{equation}


\begin{Lemma} \label{essidem:h3}
We have the following equalities:
\begin{align*}
\sum_{\l \in A_{k,1}^\bi} y_\l^\O f_{\u\u} & = \left(a_{1,1} (\cont_\lambda(i) - i) - a_{1,2} (\cont_\lambda(i) - i)\right)f_{\u\u},\\
\sum_{\l \in A_{k,2}^\bi} y_\l^\O f_{\u\u} & = \left(a_{2,2} (\cont_\lambda(i) - i) - a_{2,1} (\cont_\lambda(i) - i)\right)f_{\u\u},\\
\sum_{\l \in A_{k,3}^\bi} y_\l^\O f_{\u\u} & = \left(a_{3,2} (\cont_\lambda(i) - i) - a_{3,1} (\cont_\lambda(i) - i)\right)f_{\u\u},\\
\sum_{\l \in A_{k,4}^\bi} y_\l^\O f_{\u\u} & = \left(a_{4,1} (\cont_\lambda(i) - i) - a_{4,2} (\cont_\lambda(i) - i)\right)f_{\u\u}.
\end{align*}
\end{Lemma}

\proof We will only prove the first equality. The rest equalities can be proved following the same argument.

For any polynomial $p(y_1^\O, \ldots, y_{k-1}^\O)$, we write $p(y_1^\O, \ldots, y_{k-1}^\O) f_{\u\u} = \alpha_\u f_{\u\u}$. Recall
$$
A_{k,1}^\bi = \set{1 \leq \l \leq k-1 | i_\l = -i - 1\text{ or } -i + 1, \text{ or $\l = 1$ and $i_\l = -i$}}.
$$

Hence $\sum_{\l \in A_{k,1}^\bi} y_\l^\O$ is a polynomial of $y_1^\O, \ldots, y_{k-1}^\O$. Therefore, it suffices to prove that $\alpha_\u = a_{1,1} (\cont_\lambda(i) - i) - a_{1,2} (\cont_\lambda(i) - i)$.

Choose any $\l \in A_{k,1}^\bi$. If $\l = 1$, there is a unique $w_\l \in \{u + c_\u(1), u + c_\u(1) + 1, u + c_\u(1) - 1\}$ such that $w_\l \in S^i_1(\u)$. Similarly, if $\l > 1$, there is a unique $w_\l \in \{u + c_\u(\l) + 1, u + c_\u(\l) - 1\}$ such that $w_\l \in S^i_1(\u)$. Hence, for any $1 \leq \l \leq k-1$, we have
$$
y_\l^\O f_{\u\u} =
\begin{cases}
(\cont_\lambda(i) - i)f_{\u\u}, & \text{if $w_\l(\cont_\lambda(i),x) = 2(\cont_\lambda(i) - i)$,}\\
-(\cont_\lambda(i) - i)f_{\u\u}, & \text{if $w_\l(\cont_\lambda(i),x) = 0$.}
\end{cases}
$$

Hence by the definitions of $a_{1,1}$ and $a_{1,2}$, we have $\alpha_\u = a_{1,1} (\cont_\lambda(i) - i) - a_{1,2} (\cont_\lambda(i) - i)$. \endproof

Now we are ready to prove~\autoref{essidem:h5}.


\noindent\textit{Proof of~\autoref{essidem:h5}.} By~\autoref{essidem:h3}, we have
\begin{eqnarray*}
&& (\sum_{\l \in A_{k,1}^\bi} y_\l^\O - 2 \sum_{\l \in A_{k,2}^\bi} y_\l^\O + \sum_{\l \in A_{k,3}^\bi} y_\l^\O - 2\sum_{\l \in A_{k,4}^\bi} y_\l^\O) f_{\t\t}\\
& = & (a_{1,1} - a_{1,2} - 2a_{2,2} + 2a_{2,1} + a_{3,2} - a_{3,1} - 2a_{4,1} + 2a_{4,2}) (\cont_\lambda(i) - i) f_{\t\t}\\
& = & \left((a_{1,1} + 2a_{2,1} - a_{3,1} - 2a_{4,1}) - (a_{1,2} + 2a_{2,2} - a_{3,2} - 2a_{4,2})\right)(\cont_\lambda(i) - i) f_{\t\t}.
\end{eqnarray*}

By (\ref{essidem:eq1}), we have
$$
(a_{1,1} + 2a_{2,1} - a_{3,1} - 2a_{4,1}) - (a_{1,2} + 2a_{2,2} - a_{3,2} - 2a_{4,2}) =
\begin{cases}
2, & \text{if $h_k(\bi) = 0$ and $|\mathscr{AR}_\lambda(i)| = 1$,}\\
1, & \text{if $h_k(\bi) = -1$,}\\
0, & \text{if $h_k(\bi) = 0$ and $|\mathscr{AR}_\lambda(i)| = 0$,}
\end{cases}
$$
which implies
\begin{equation} \label{essidem:h5:eq1}
(\sum_{\l \in A_{k,1}^\bi} y_\l^\O - 2 \sum_{\l \in A_{k,2}^\bi} y_\l^\O + \sum_{\l \in A_{k,3}^\bi} y_\l^\O - 2\sum_{\l \in A_{k,4}^\bi} y_\l^\O) f_{\u\u} =
\begin{cases}
2(\cont_\lambda(i) - i) f_{\u\u}, & \text{if $h_k(\bi) = 0$ and $|\mathscr{AR}_\lambda(i)| = 1$,}\\
(\cont_\lambda(i) - i) f_{\u\u}, & \text{if $h_k(\bi) = -1$,}\\
0, & \text{if $h_k(\bi) = 0$ and $|\mathscr{AR}_\lambda(i)| = 0$.}
\end{cases}
\end{equation}

By the construction, there exists $\t \in \Tud_n(\bi)$ such that $\t \overset{k}\sim \u$ if and only if $|\mathscr{AR}_\lambda(i)| = 1$. Moreover, as $\bi \in P_{k,+}^n$, we have $i_k = -\frac{1}{2}$ if $h_k(\bi) = -1$ and $i_k \neq -\frac{1}{2}$ if $h_k(\bi) = 0$. Therefore the Lemma follows by~\eqref{essidem:h5:eq1} and $c_\t(k) = \cont_\lambda(i_k)$. \qed

Finally we prove the relation~\eqref{rela:5:5} when $\bi \in P_{k,+}^n$.


\begin{Lemma} \label{essidem:8}
Suppose $1 \leq k \leq n-1$ and $\bi \in P_{k,+}^n$. Then for any $\bj, \bk \in I^n$, we have
$$
e(\bj) \epsilon_k e(\bi) \epsilon_k e(\bk) =
(-1)^{a_k(\bi)}(1 + \delta_{i_k, -\frac{1}{2}})(\sum_{\l \in A_{k,1}^\bi} y_\l - 2 \sum_{\l \in A_{k,2}^\bi} y_\l + \sum_{\l \in A_{k,3}^\bi} y_\l - 2\sum_{\l \in A_{k,4}^\bi} y_\l) e(\bj) \epsilon_k e(\bk).
$$
\end{Lemma}

\proof Choose arbitrary $\u \in \Tud_n(\bj)$. If $\u(k) + \u(k+1) \neq 0$, by~\autoref{idem:2} we have
\begin{equation} \label{essidem:8:eq1}
f_{\u\u} \epsilon_k^\O e(\bi)^\O \epsilon_k^\O e(\bk)^\O = 0 = (\sum_{\l \in A_{k,1}^\bi} y_\l^\O - 2 \sum_{\l \in A_{k,2}^\bi} y_\l^\O + \sum_{\l \in A_{k,3}^\bi} y_\l^\O - 2\sum_{\l \in A_{k,4}^\bi} y_\l^\O) f_{\u\u} \epsilon_k^\O e(\bk)^\O.
\end{equation}

If $\u(k) + \u(k+1) = 0$, let $\lambda = \u_{k-1}$. If there is no $\t \in \Tud_n(\bi)$ with $\t \overset{k}\sim \u$, we have $|\mathscr{AR}_\lambda(i_k)| = 0$, which forces $h_k(\bi) = 0$ by~\autoref{deg:h1:1}. Hence,
\begin{equation} \label{essidem:8:eq2}
f_{\u\u} \epsilon_k^\O e(\bi)^\O \epsilon_k^\O e(\bk)^\O = 0 = (\sum_{\l \in A_{k,1}^\bi} y_\l^\O - 2 \sum_{\l \in A_{k,2}^\bi} y_\l^\O + \sum_{\l \in A_{k,3}^\bi} y_\l^\O - 2\sum_{\l \in A_{k,4}^\bi} y_\l^\O) f_{\u\u} \epsilon_k^\O e(\bk)^\O,
\end{equation}
by~\autoref{essidem:h5}.

Suppose there exist $\t \in \Tud_n(\bi)$ with $\t \overset{k}\sim \u$. Because $\bi \in P_{k,+}^n$ and $\bi \in I^n$, we have $h_k(\bi) = 0$ or $-1$. Hence, by (\ref{deg:h4:eq1}) and (\ref{deg:h4:eq2}), we always have $|\mathscr{AR}_\lambda(i_k)| = 1$. Therefore, $\t$ is unique, and we have
$$
f_{\u\u} \epsilon_k^\O e(\bi)^\O \epsilon_k^\O e(\bk)^\O  = (-1)^{a_k(\bi)} 2(c_\t(k) - i_k) f_{\u\u} \epsilon_k^\O e(\bk)^\O,
$$
by~\autoref{semi:2} and~\autoref{PQ:3}. So, we have
\begin{equation} \label{essidem:8:eq3}
f_{\u\u} \epsilon_k^\O e(\bi)^\O \epsilon_k^\O e(\bk)^\O = (-1)^{a_k(\bi)}(1 + \delta_{i_k, -\frac{1}{2}})(\sum_{\l \in A_{k,1}^\bi} y_\l^\O - 2 \sum_{\l \in A_{k,2}^\bi} y_\l^\O + \sum_{\l \in A_{k,3}^\bi} y_\l^\O - 2\sum_{\l \in A_{k,4}^\bi} y_\l^\O) f_{\u\u} \epsilon_k^\O e(\bk)^\O,
\end{equation}
by~\autoref{essidem:h5}.

Combining (\ref{essidem:8:eq1}), (\ref{essidem:8:eq2}) and (\ref{essidem:8:eq3}), we have
$$
e(\bj)^\O \epsilon_k^\O e(\bi)^\O \epsilon_k^\O e(\bk)^\O = (-1)^{a_k(\bi)}(1 + \delta_{i_k, -\frac{1}{2}})(\sum_{\l \in A_{k,1}^\bi} y_\l^\O - 2 \sum_{\l \in A_{k,2}^\bi} y_\l^\O + \sum_{\l \in A_{k,3}^\bi} y_\l^\O - 2\sum_{\l \in A_{k,4}^\bi} y_\l^\O) e(\bj)^\O \epsilon_k^\O e(\bk)^\O,
$$
which completes the proof by lifting the elements into $\B$. \endproof

We have proved that (\ref{rela:5:1}) - (\ref{rela:5:4}) hold by~\autoref{esscom:2} and~\autoref{essidem:1} - \ref{essidem:5}; and (\ref{rela:5:5}) holds by~\autoref{essidem:6},
~\autoref{essidem:7} and~\autoref{essidem:8}. Therefore, the essential idempotent relations hold in $\B$.

\begin{Proposition} \label{essidem}
In $\B$, the essential idempotent relations hold.
\end{Proposition}

%
%
%


We close this subsection by proving the untwist relations. We remind the readers that in the rest of this paper we assume $\bi,\bj,\bk \in I^n$.

\begin{Proposition} \label{untwist}
In $\B$, the untwist relation holds.
\end{Proposition}

\proof We need to show that for any $1 \leq k \leq n-1$ and $\bi,\bj \in I^n$, we have
\begin{align*}
e(\bi) \psi_k \epsilon_k e(\bj) & =
\begin{cases}
(-1)^{a_k(\bi)} e(\bi) \epsilon_k e(\bj), & \text{if $\bi \in I_{k,+}^n$ and $i_k \neq -\frac{1}{2}$,}\\
0, & \text{otherwise;}
\end{cases}\\
e(\bj) \epsilon_k \psi_k e(\bi) & = \begin{cases}
(-1)^{a_k(\bi)} e(\bj) \epsilon_k e(\bi), & \text{if $\bi \in I_{k,+}^n$ and $i_k \neq -\frac{1}{2}$,}\\
0, & \text{otherwise.}
\end{cases}
\end{align*}

By~\autoref{idem}, we assume $i_k + i_{k+1} = 0$. Otherwise both sides of the equality will be $0$. We will only prove the first equality, and the second equality follows by the similar argument.

Suppose $\bi \in I_{k,+}^n$ and $i_k \neq 0,-\frac{1}{2}$. Choose arbitrary $\t \in \Tud_n(\bi)$. If $c_\t(k) + c_\t(k+1) \neq 0$, we have
\begin{equation} \label{untwist:eq1}
f_{\t\t} \psi_k^\O \epsilon_k^\O e(\bj)^\O = 0 = (-1)^{a_k(\bi)} f_{\t\t} \epsilon_k^\O e(\bj)^\O,
\end{equation}
by~\autoref{semi:1}. If $c_\t(k) + c_\t(k+1) = 0$, because $i_k \neq 0$ and $\bi \in I_{k,+}^n$, we have $\bi{\cdot}s_k \neq \bi$ and $h_k(\bi) = 0$. Therefore, by~\autoref{deg:h4:1}, there exists a unique $\u \in \Tud_n(\bi{\cdot}s_k)$ such that $\u \overset{k}\sim \t$ and $\u \neq \t$, and $c_\t(k) - i_k = c_\u(k) - i_{k+1}$. We note that $c_\t(k) - i_k = c_\u(k) - i_{k+1}$ implies $c_\t(k) + c_\u(k) = 2(c_\u(k) - i_{k+1})$ because $i_k + i_{k+1} = 0$. Hence, by~\autoref{semi:B},~\autoref{semi:2} and~\autoref{PQ:3}, we have
\begin{align}
f_{\t\t} \psi_k^\O \epsilon_k^\O e(\bj)^\O & = \sum_{\substack{\v \in \Tud_n(\bj) \\ \v \overset{k}\sim \t}} P_k(\t)^{-1} s_k(\t,\u) Q_k(\u)^{-1} P_k(\u)^{-1} e_k(\u,\v) Q_k(\v)^{-1} f_{\t\v} \notag \\
& = \sum_{\substack{\v \in \Tud_n(\bj) \\ \v \overset{k}\sim \t}} \frac{P_k(\u)^{-1}Q_k(\u)^{-1} e_k(\u,\u)}{c_\t(k) + c_\u(k)} P_k(\t)^{-1} e_k(\t,\v) Q_k(\v)^{-1} f_{\t\v} \notag \\
& = \frac{P_k(\u)^{-1}Q_k(\u)^{-1} e_k(\u,\u)}{2(c_\u(k) - i_{k+1})} f_{\t\t} \epsilon_k^\O e(\bj)^\O \notag \\
& = (-1)^{a_k(\bi)} f_{\t\t} \epsilon_k^\O e(\bj)^\O. \label{untwist:eq2}
\end{align}

Because $\t$ is chosen arbitrary, by (\ref{untwist:eq1}), (\ref{untwist:eq2}) and~\autoref{idem-semi}, we have
$$e(\bi)^\O \psi_k^\O \epsilon_k^\O e(\bj)^\O = (-1)^{a_k(\bi)} e(\bi)^\O \epsilon_k^\O e(\bj)^\O \in \BOx,
$$
and therefore $e(\bi) \psi_k \epsilon_k e(\bj) = e(\bi) \epsilon_k e(\bj)$ by lifting the elements into $\B$.

For the rest of the cases, if $h_k(\bi) \neq 0$, we have $\bi{\cdot}s_k \not\in I^n$ by~\autoref{A:9}, which yields $e(\bi) = 0$ by~\autoref{idem:semi:1}. Therefore, we have $e(\bi) \psi_k \epsilon_k e(\bj) = 0$ by~\autoref{com}. If $h_k(\bi) = 0$, we have $i_k = i_{k+1} = 0$, which yields $e(\bi) \psi_k = 0$ by~\autoref{esscom:3}. Therefore, we have $e(\bi) \psi_k \epsilon_k e(\bj) = 0$. \endproof

\subsection{Tangle relations} \label{sec:rela:3}

In this subsection, we will prove the tangle relations hold in $\BOx$ and extend to $\B$ by lifting the elements into $\B$. First we prove (\ref{rela:7:2}).

\begin{Lemma} \label{tangle:2}
Suppose $1 \leq k \leq n-1$ and $\bi, \bj \in I^n$. Then we have $e(\bi)^\O \epsilon_k^\O e(\bj)^\O (y_k^\O + y_{k+1}^\O) = 0$.
\end{Lemma}

\proof By~\autoref{idem:2}, we have $e(\bi)^\O \epsilon_k^\O e(\bj)^\O = 0$ if $i_k + i_{k+1} \neq 0$ or $j_k + j_{k+1} \neq 0$, where the Lemma holds. Suppose $i_k + i_{k+1} = j_k + j_{k+1} = 0$, as $\L_n(\O)$ is a commutative subalgebra of $\BOx$ and $Q_k^\O(\bj)^{-1} e(\bj)^\O \in \L_n(\O)$, we have
$$
e(\bi)^\O \epsilon_k^\O e(\bj)^\O (y_k^\O + y_{k+1}^\O) = e(\bi)^\O P^\O_k(\bi)^{-1} e_k^\O (L_k^\O + L_{k+1}^\O) Q^\O_k(\bj)^{-1} e(\bj)^\O = 0,
$$
which completes the proof. \endproof

\begin{Lemma} \label{tangle:3}
Suppose $1 \leq k \leq n-1$ and $\bi,\bj \in I^n$. Then we have
$$
e(\bi)^\O \epsilon_k^\O \epsilon_{k-1}^\O \epsilon_k^\O e(\bj)^\O = e(\bi)^\O \epsilon_k^\O e(\bj)^\O \in \BOx, \qquad \text{and} \qquad e(\bi)^\O \epsilon_{k-1}^\O \epsilon_k^\O \epsilon_{k-1}^\O e(\bj)^\O = e(\bi)^\O \epsilon_{k-1}^\O e(\bj)^\O \in \BOx.
$$
\end{Lemma}

\proof We only prove the first equality. The second equality follows by the same argument. Choose arbitrary $\t \in \Tud_n(\bi)$. When $\t(k) + \t(k+1) \neq 0$, the equality holds because both sides of the equality are $0$. When $\t(k) + \t(k+1) = 0$, we write $\t = (\alpha_1, \ldots, \alpha_n)$ and define
$$
\s = (\alpha_1, \ldots, \alpha_{k-2}, \alpha_{k-1}, -\alpha_{k-1}, \alpha_{k-1}, \alpha_{k+2}, \alpha_{k+3}, \ldots, \alpha_n).
$$

Because $\alpha_k + \alpha_{k+1} = 0$, $\s$ is an up-down tableau. By~\autoref{PQ:2} and~\autoref{semi:2}, we have
$$
f_{\t\t} \epsilon_k^\O \epsilon_{k-1}^\O \epsilon_k^\O e(\bj)^\O = f_{\t\t} \epsilon_k^\O \f{\s} \epsilon_{k-1}^\O \f{\s} \epsilon_k^\O e(\bj)^\O = f_{\t\t} P_k^\O(\bi)^{-1} e_k^\O e_{k-1}^\O e_k^\O Q_k^\O(\bj)^{-1} e(\bj)^\O =  f_{\t\t} \epsilon_k^\O e(\bj)^\O.
$$

As $\t$ is chosen arbitrary, by~\autoref{idem-semi}, we have $e(\bi)^\O \epsilon_k^\O \epsilon_{k-1}^\O \epsilon_k^\O e(\bj)^\O = e(\bi)^\O \epsilon_k^\O e(\bj)^\O\in \BOx$.
\endproof

%

Now we prove (\ref{rela:7:1}) hold in $\BOx$. We separate the question into several cases based on the values of $i_{k-1}, i_k$ and $i_{k+1}$. In more details, we consider the following three cases:

(6.4.1). When $i_{k-1} + i_{k+1} \neq 0$.

(6.4.2). When $i_{k-1} + i_{k+1} = 0$, $i_k + i_{k-1} \neq 0$ and $i_k + i_{k+1} \neq 0$.

(6.4.3). When $i_{k-1} + i_{k+1} = 0$, and $i_k + i_{k-1} = 0$ or $i_k + i_{k+1} = 0$.

First we consider the case (6.4.1).

\begin{Lemma} \label{tangle:4}
Suppose $1 < k < n$ and $\bi \in I^n$ with $i_{k-1} + i_{k+1} \neq 0$. Then we have
$$
e(\bj)^\O \epsilon_k^\O \epsilon_{k-1}^\O \psi_k^\O e(\bi)^\O = 0 = e(\bj)^\O \epsilon_k^\O \psi_{k-1}^\O e(\bi)^\O, \qquad \text{and} \qquad
e(\bi)^\O \psi_k^\O \epsilon_{k-1}^\O \epsilon_k^\O e(\bj)^\O = 0 = e(\bi)^\O \psi_{k-1}^\O \epsilon_k^\O e(\bj)^\O.
$$
\end{Lemma}

\begin{proof}
By~\autoref{com}, we have
\begin{align*}
e(\bj)^\O \epsilon_k^\O \epsilon_{k-1}^\O \psi_k^\O e(\bi)^\O & = e(\bj)^\O \epsilon_k^\O \epsilon_{k-1}^\O e(\bi{\cdot}s_k)^\O \psi_k^\O e(\bi)^\O, \\
e(\bj)^\O \epsilon_k^\O \psi_{k-1}^\O e(\bi)^\O & = e(\bj)^\O \epsilon_k^\O e(\bi{\cdot}s_{k-1})^\O \psi_{k-1}^\O e(\bi)^\O.
\end{align*}

By~\autoref{idem:2}, both equalities equal to $0$, which proves $e(\bj)^\O \epsilon_k^\O \epsilon_{k-1}^\O \psi_k^\O e(\bi)^\O = 0 = e(\bj)^\O \epsilon_k^\O \psi_{k-1}^\O e(\bi)^\O$. Following the same argument, we can prove $e(\bi)^\O \psi_k^\O \epsilon_{k-1}^\O \epsilon_k^\O e(\bj)^\O = 0 = e(\bi)^\O \psi_{k-1}^\O \epsilon_k^\O e(\bj)^\O$.
\end{proof}

Then we consider the case (6.4.2). In this case, for any $\t \in \Tud_n(\bi)$, we have $\t(k-1) + \t(k) \neq 0$ and $\t(k) + \t(k+1) \neq 0$. Hence, by~\autoref{semi:1}, the actions of $\psi_k^\O$ and $\psi_{k-1}^\O$ on $f_{\t\t}$ are the same as in the KLR algebras. Therefore, for $\bi \in I^n$ with $i_k + i_{k+1} \neq 0$, we have the following Lemma, which is analogue to~\cite[Definition 4.14, Lemma 4.18]{HuMathas:SemiQuiver}.

Recall that when $i_k \neq i_{k+1} + 1$, by~\autoref{L:inv:1}, we have
$$
\frac{1}{1 - L_k^\O + L_{k+1}^\O} e(\bi)^\O = e(\bi)^\O \frac{1}{1 - L_k^\O + L_{k+1}^\O} = \sum_{\t \in \Tud_n(\bi)} \frac{1}{1 - c_\t(k) + c_\t(k+1)} \f{\t} \in \L_n(\O).
$$

\begin{Lemma} \label{tangle:h1}
Suppose $\bi \in I^n$ with $i_k + i_{k+1} \neq 0$. Then
\begin{align*}
e(\bi)^\O \psi_k^\O & =
\begin{cases}
e(\bi)^\O \frac{1}{1 - L_{k+1}^\O + L_k^\O} (s_k^\O - 1), & \text{if $i_k = i_{k+1}$,}\\
e(\bi)^\O ((L_{k+1}^\O - L_k^\O)s_k^\O - 1), & \text{if $i_k = i_{k+1} + 1$,}\\
e(\bi)^\O \frac{1}{1 - L_{k+1}^\O + L_k^\O} ((L_{k+1}^\O - L_k^\O)s_k^\O - 1), & \text{otherwise;}\\
\end{cases} \\
\psi_k^\O e(\bi)^\O & =
\begin{cases}
(s_k^\O + 1) \frac{1}{1 - L_k^\O + L_{k+1}^\O} e(\bi)^\O, & \text{if $i_k = i_{k+1}$,}\\
(s_k^\O (L_k^\O - L_{k+1}^\O) + 1) e(\bi)^\O, & \text{if $i_k = i_{k+1} + 1$,}\\
(s_k^\O (L_k^\O - L_{k+1}^\O) + 1) \frac{1}{1 - L_k^\O + L_{k+1}^\O} e(\bi)^\O, & \text{otherwise.}\\
\end{cases}
\end{align*}
\end{Lemma}

\proof It suffices to prove that for any $\t \in \Tud_n(\bi)$, we have
\begin{align}
f_{\t\t} \psi_k^\O & =
\begin{cases}
f_{\t\t} e(\bi)^\O \frac{1}{1 - L_{k+1}^\O + L_k^\O} (s_k^\O - 1), & \text{if $i_k = i_{k+1}$,}\\
f_{\t\t} ((L_{k+1}^\O - L_k^\O)s_k^\O - 1), & \text{if $i_k = i_{k+1} + 1$,}\\
f_{\t\t} e(\bi)^\O \frac{1}{1 - L_{k+1}^\O + L_k^\O} ((L_{k+1}^\O - L_k^\O)s_k^\O - 1), & \text{otherwise;}\\
\end{cases} \label{tangle:h1:eq1} \\
\psi_k^\O f_{\t\t} & =
\begin{cases}
(s_k^\O + 1) \frac{1}{1 - L_k^\O + L_{k+1}^\O} e(\bi)^\O f_{\t\t}, & \text{if $i_k = i_{k+1}$,}\\
(s_k^\O (L_k^\O - L_{k+1}^\O) + 1) f_{\t\t}, & \text{if $i_k = i_{k+1} + 1$,}\\
(s_k^\O (L_k^\O - L_{k+1}^\O) + 1) \frac{1}{1 - L_k^\O + L_{k+1}^\O} e(\bi)^\O f_{\t\t}, & \text{otherwise.}\\
\end{cases} \label{tangle:h1:eq2}
\end{align}

Because $\t \in \Tud_n(\bi)$ and $i_k + i_{k+1} \neq 0$, we have $\t(k) + \t(k+1) \neq 0$.

If $i_k = i_{k+1}$, because $i_k + i_{k+1} \neq 0$, we have $i_k = i_{k+1} \neq 0$. Hence, by the definition of $h_k$, we have $h_{k+1}(\bi) = h_k(\bi) + 2$ when $i_k \neq \pm\frac{1}{2}$, and $h_{k+1}(\bi) = h_k(\bi) + 3$ when $i_k = \pm \frac{1}{2}$. By~\autoref{deg:h2}, we have $-2 \leq h_k(\bi), h_{k+1}(\bi) \leq 0$, which forces $h_k(\bi) = -2$ and $i_k \neq \pm \frac{1}{2}$. By~\autoref{deg:h4:h}, we have $\t(k) > 0$ and $\t(k+1) < 0$, or $\t(k) < 0$ and $\t(k+1) > 0$. Because $\t(k) + \t(k+1) \neq 0$, by~\autoref{y:h2:8}, we have $\s = \t{\cdot}s_k \in \Tud_n(\bi)$. Then, by~\autoref{semi:B} and~\autoref{semi:1}, we have
$$
f_{\t\t} \psi_k^\O = \frac{1}{c_\t(k+1) - c_\t(k)} f_{\t\t} + \frac{s_k(\t,\s)}{1 - c_\s(k) + c_\s(k+1)} f_{\s\t},
$$
and because $\s = \t{\cdot}s_k$, we have $c_\t(k+1) - c_\t(k) = c_\s(k) - c_\s(k+1)$. Hence,
\begin{align*}
f_{\t\t} e(\bi)^\O \frac{1}{1 - L_{k+1}^\O + L_k^\O} (s_k^\O - 1) & = \frac{1}{1 - c_\t(k+1) + c_\t(k)} (\frac{1}{c_\t(k+1) - c_\t(k)} f_{\t\t} + s_k(\t,\s) f_{\s\t} - 1)\\
& = \frac{1}{c_\t(k+1) - c_\t(k)} f_{\t\t} + \frac{s_k(\t,\s)}{1 - c_\s(k) + c_\s(k+1)} f_{\s\t},
\end{align*}
which proves (\ref{tangle:h1:eq1}) when $i_k = i_{k+1}$.

If $i_k = i_{k+1} + 1$, when $\t{\cdot}s_k$ does not exist, by~\autoref{semi:B} and~\autoref{semi:1}, we have
\begin{align*}
f_{\t\t} ((L_{k+1}^\O - L_k^\O)s_k^\O - 1) & = (\frac{c_\t(k+1) - c_\t(k)}{c_\t(k+1) - c_\t(k)} - 1) f_{\t\t} = 0 = f_{\t\t} \psi_k^\O.
\end{align*}

When $\s = \t{\cdot}s_k$ is an up-down tableau, by~\autoref{semi:B} and~\autoref{semi:1}, we have
$$
f_{\t\t} \psi_k^\O = s_k(\t,\s)(c_\s(k) - c_\s(k+1)) f_{\t\s},
$$
and because $\s = \t{\cdot}s_k$, we have $c_\t(k+1) - c_\t(k) = c_\s(k) - c_\s(k+1)$. Hence,
\begin{align*}
f_{\t\t} ((L_{k+1}^\O - L_k^\O)s_k^\O - 1) & = (\frac{c_\t(k+1) - c_\t(k)}{c_\t(k+1) - c_\t(k)} - 1) f_{\t\t} + (c_\t(k+1) - c_\t(k))s_k(\t,\s) f_{\t\s}\\
& = s_k(\t,\s)(c_\s(k) - c_\s(k+1)) f_{\t\s},
\end{align*}
which proves (\ref{tangle:h1:eq1}) when $i_k = i_{k+1} + 1$.

For the other cases, when $\t{\cdot}s_k$ does not exist, by~\autoref{semi:B} and~\autoref{semi:1}, we have
\begin{align*}
f_{\t\t} e(\bi)^\O \frac{1}{1 - L_{k+1}^\O + L_k^\O} ((L_{k+1}^\O - L_k^\O)s_k^\O - 1) & = \frac{1}{1 - c_\t(k+1) + c_\t(k)}(\frac{c_\t(k+1) - c_\t(k)}{c_\t(k+1) - c_\t(k)} - 1) f_{\t\t} = 0 = f_{\t\t} \psi_k^\O.
\end{align*}

When $\s = \t{\cdot}s_k$ is an up-down tableau, by~\autoref{semi:B} and~\autoref{semi:1}, we have
$$
f_{\t\t} \psi_k^\O = s_k(\t,\s) \frac{c_\s(k) - c_\s(k+1)}{1 - c_\s(k) + c_\s(k+1)} f_{\t\s},
$$
and because $\s = \t{\cdot}s_k$, we have $c_\t(k+1) - c_\t(k) = c_\s(k) - c_\s(k+1)$. Hence,
\begin{eqnarray*}
&&f_{\t\t} e(\bi)^\O \frac{1}{1 - L_{k+1}^\O + L_k^\O} ((L_{k+1}^\O - L_k^\O)s_k^\O - 1) \\
& = & \frac{1}{1 - c_\t(k+1) + c_\t(k)}(\frac{c_\t(k+1) - c_\t(k)}{c_\t(k+1) - c_\t(k)} - 1) f_{\t\t} + \frac{c_\t(k+1) - c_\t(k)}{1 - c_\t(k+1) + c_\t(k)} s_k(\t,\s) f_{\t\s}\\
& = & s_k(\t,\s) \frac{c_\s(k) - c_\s(k+1)}{1 - c_\s(k) + c_\s(k+1)} f_{\t\s},
\end{eqnarray*}
which proves (\ref{tangle:h1:eq1}) for the rest of the cases. Hence, (\ref{tangle:h1:eq1}) holds. Following the similar argument, (\ref{tangle:h1:eq2}) holds. \endproof

The symmetric group $\Sym_n$ acts from left on the rational functions $f \in R(L_1^\O, \ldots, L_n^\O)$ by permuting variables. We denote $s_k{\cdot}f$ by $\leftidx{^{s_k}}f$.

By~\autoref{tangle:h1}, for $\bi \in I^n$ with $i_k + i_{k+1} \neq 0$, we define
\begin{equation} \label{tangle:eq1}
\begin{array}{rlrl}
M_k^\O(\bi) & =
\begin{cases}
\frac{1}{1 - L_k^\O + L_{k+1}^\O} e(\bi)^\O, & \text{if $i_k = i_{k+1}$,}\\
L_k^\O - L_{k+1}^\O, & \text{if $i_k = i_{k+1} + 1$,}\\
\frac{L_k^\O - L_{k+1}^\O}{1 - L_k^\O + L_{k+1}^\O}, & \text{otherwise;}
\end{cases} &
N_k^\O(\bi) & =
\begin{cases}
\frac{1}{1 - L_k^\O + L_{k+1}^\O} e(\bi)^\O, & \text{if $i_k = i_{k+1}$,}\\
1, & \text{if $i_k = i_{k+1} + 1$,}\\
\frac{1}{1 - L_k^\O + L_{k+1}^\O}, & \text{otherwise;}
\end{cases} \\
\widetilde{M}_k^\O(\bi) & =
\begin{cases}
\frac{1}{1 - L_{k+1}^\O + L_k^\O} e(\bi)^\O, & \text{if $i_k = i_{k+1}$,}\\
L_{k+1}^\O - L_k^\O, & \text{if $i_k = i_{k+1} + 1$,}\\
\frac{L_{k+1}^\O - L_k^\O}{1 - L_{k+1}^\O + L_k^\O}, & \text{otherwise;}
\end{cases} &
\widetilde{N}_k^\O(\bi) & =
\begin{cases}
-\frac{1}{1 - L_{k+1}^\O + L_k^\O} e(\bi)^\O, & \text{if $i_k = i_{k+1}$,}\\
-1, & \text{if $i_k = i_{k+1} + 1$,}\\
-\frac{1}{1 - L_{k+1}^\O + L_k^\O}, & \text{otherwise,}
\end{cases}
\end{array}
\end{equation}
such that
\begin{equation} \label{tangle:eq2}
\begin{array}{rl}
\psi_k^\O e(\bi)^\O & = s_k^\O M_k^\O(\bi) e(\bi)^\O + N_k^\O(\bi) e(\bi)^\O,\\
e(\bi)^\O \psi_k^\O & = e(\bi)^\O \widetilde{M}_k^\O(\bi) s_k^\O + e(\bi)^\O \widetilde{N}_k^\O(\bi).
\end{array}
\end{equation}

The following is the technical result we need later.

\begin{Lemma} \label{tangle:h2}
Suppose $\bi \in I^n$ with $i_{k-1} + i_{k+1} = 0$. For any $\t \in \mathscr T^{ud}_n(\bi)$ with $\t(k-1) + \t(k+1) = 0$ and $\t(k) \neq \t(k-1), \t(k+1)$, we have
\begin{align*}
M_k^\O(\bi) Q_{k-1}^\O(\bi{\cdot}s_k)^{-1} f_{\t\t} & = M_{k-1}^\O(\bi) \leftidx{^{s_{k-1}}}Q_k^\O(\bi{\cdot}s_{k-1})^{-1} f_{\t\t},\\
f_{\t\t} P_{k-1}^\O(\bi{\cdot}s_k)^{-1} \widetilde M_k^\O(\bi) & = f_{\t\t} \leftidx{^{s_{k-1}}}P_k^\O(\bi{\cdot}s_{k-1})^{-1} \widetilde M_{k-1}^\O(\bi).
\end{align*}
\end{Lemma}

\proof We only prove the first equality. The second equality follows by the similar method.

Because $\t(k-1) + \t(k+1) = 0$, we have $L_{k-1}^\O f_{\t\t} = - L_{k+1}^\O f_{\t\t}$. Moreover, if $i_{k-1} = i_k = i_{k+1}$, as $i_{k-1} + i_{k+1} = 0$, we have $i_{k-1} = i_k = i_{k+1} = 0$, which forces $\t(k-1) = -\t(k) = \t(k+1)$ by the construction of up-down tableaux. Hence, $i_{k-1} = i_k = i_{k+1}$ is excluded.

By the definition of $Q_k^\O(\bi)$ and $i_{k-1} = - i_{k+1}$, we have
\begin{equation} \label{tangle:h2:eq1}
\leftidx{^{s_{k-1}}}Q_k^\O(\bi{\cdot}s_{k-1})^{-1} f_{\t\t} =
\begin{cases}
\frac{1 - L_{k-1}^\O + L_k^\O}{1 - L_k^\O + L_{k+1}^\O} (L_k^\O - L_{k+1}^\O) Q_{k-1}^\O(\bi{\cdot}s_k)^{-1} f_{\t\t}, & \text{if $i_{k-1} = i_k$,}\\
\frac{1 - L_{k-1}^\O + L_k^\O}{1 - L_k^\O + L_{k+1}^\O} \frac{1}{L_{k-1}^\O - L_k^\O} Q_{k-1}^\O(\bi{\cdot}s_k)^{-1} f_{\t\t}, & \text{if $i_k = i_{k+1}$,}\\
\frac{1}{1 - L_k^\O + L_{k+1}^\O} \frac{L_k^\O - L_{k+1}^\O}{L_{k-1}^\O - L_k^\O} Q_{k-1}^\O(\bi{\cdot}s_k)^{-1} f_{\t\t}, & \text{if $i_{k-1} = i_k + 1$,}\\
(1 - L_{k-1}^\O + L_k^\O) \frac{L_k^\O - L_{k+1}^\O}{L_{k-1}^\O - L_k^\O} Q_{k-1}^\O(\bi{\cdot}s_k)^{-1} f_{\t\t}, & \text{if $i_k = i_{k+1} + 1$,}\\
\frac{L_k^\O - L_{k+1}^\O}{L_{k-1}^\O - L_k^\O} Q_{k-1}^\O(\bi{\cdot}s_k)^{-1} f_{\t\t}, & \text{if $i_{k-1} = i_k + 1$ and $i_k = i_{k+1} + 1$,}\\
\frac{1 - L_{k-1}^\O + L_k^\O}{1 - L_k^\O + L_{k+1}^\O} \frac{L_k^\O - L_{k+1}^\O}{L_{k-1}^\O - L_k^\O} Q_{k-1}^\O(\bi{\cdot}s_k)^{-1} f_{\t\t}, & \text{otherwise.}
\end{cases}
\end{equation}

By (\ref{tangle:eq1}) and (\ref{tangle:h2:eq1}), $M_k^\O(\bi) Q_{k-1}^\O(\bi{\cdot}s_k)^{-1} f_{\t\t} = M_{k-1}^\O(\bi) \leftidx{^{s_{k-1}}}Q_k^\O(\bi{\cdot}s_{k-1})^{-1} f_{\t\t}$ can be verified by direct calculation. \endproof

Now we prove (\ref{rela:7:1}) when $i_{k-1} + i_{k+1} = 0$, $i_k + i_{k-1} \neq 0$ and $i_k + i_{k+1} \neq 0$.

\begin{Lemma} \label{tangle:5:h1}
Suppose $1 < k < n$ and $\bi, \bj \in I^n$. For any $\t \in  \Tud_n(\bi)$ with $\t(k-1) + \t(k+1) \neq 0$, $\t(k-1) + \t(k) \neq 0$ and $\t(k) + \t(k+1) \neq 0$, we have
\begin{align}
e(\bj)^\O \epsilon_k^\O \epsilon_{k-1}^\O \psi_k^\O f_{\t\t} & = 0 = e(\bj)^\O \epsilon_k^\O \psi_{k-1}^\O f_{\t\t}, \label{tangle:5:h1:eq1} \\
f_{\t\t} \psi_k^\O \epsilon_{k-1}^\O \epsilon_k^\O e(\bj)^\O & = 0 = f_{\t\t} \psi_{k-1}^\O \epsilon_k^\O e(\bj)^\O. \label{tangle:5:h1:eq2}
\end{align}
\end{Lemma}

\proof We only prove (\ref{tangle:5:h1:eq1}), and (\ref{tangle:5:h1:eq2}) can be proved following the similar argument.

When $\t{\cdot}s_k$ does not exist, by~\autoref{semi:1}, we have $\psi_k^\O f_{\t\t} = 0$, which implies $e(\bj)^\O \epsilon_k^\O \epsilon_{k-1}^\O \psi_k^\O f_{\t\t} = 0$; and when $\s = \t{\cdot}s_k$ is an up-down tableau, we have $\s(k-1) + \s(k) \neq 0$, which implies $e(\bj)^\O \epsilon_k^\O \epsilon_{k-1}^\O \psi_k^\O f_{\t\t} = 0$ by~\autoref{semi:1}.

When $\t{\cdot}s_{k-1}$ does not exist, by~\autoref{semi:1}, we have $\psi_{k-1}^\O f_{\t\t} = 0$, which implies $e(\bj)^\O \epsilon_k^\O \psi_{k-1}^\O f_{\t\t} = 0$; and when $\s = \t{\cdot}s_{k-1}$ is an up-down tableau, we have $\s(k) + \s(k+1) \neq 0$, which implies $e(\bj)^\O \epsilon_k^\O \psi_{k-1}^\O f_{\t\t} = 0$ by~\autoref{semi:1}. Hence, (\ref{tangle:5:h1:eq1}) follows. \endproof

\begin{Lemma} \label{tangle:5:h2}
Suppose $1 < k < n$ and $\bi \in I^n$ with $i_{k-1} + i_{k+1} = 0$. If $\t \in \Tud_n(\bi)$ with $\t(k-1) + \t(k) \neq 0$ and $\t(k) + \t(k+1) \neq 0$, for any $\bj \in I^n$, we have
\begin{align}
e(\bj)^\O \epsilon_k^\O \epsilon_{k-1}^\O \psi_k^\O f_{\t\t} & = e(\bj)^\O \epsilon_k^\O \psi_{k-1}^\O f_{\t\t}, \label{tangle:5:h2:eq1} \\
f_{\t\t} \psi_k^\O \epsilon_{k-1}^\O \epsilon_k^\O e(\bj)^\O & = f_{\t\t} \psi_{k-1}^\O \epsilon_k^\O e(\bj)^\O. \label{tangle:5:h2:eq2}
\end{align}
\end{Lemma}

\proof We only prove the equality (\ref{tangle:5:h2:eq1}). The equality (\ref{tangle:5:h2:eq2}) can be proved following the similar argument.

When $\t(k-1) + \t(k+1) \neq 0$, because $\t(k-1) + \t(k) \neq 0$ and $\t(k) + \t(k+1) \neq 0$, (\ref{tangle:5:h2:eq1}) holds by~\autoref{tangle:5:h1}.

When $\t(k-1) + \t(k+1) = 0$, because $\t(k-1) + \t(k) \neq 0$ and $\t(k) + \t(k+1) \neq 0$, by~\autoref{com},~\autoref{PQ:2} and (\ref{tangle:eq2}), we have
\begin{align*}
e(\bj)^\O \epsilon_k^\O \epsilon_{k-1}^\O \psi_k^\O f_{\t\t} & = e(\bj)^\O P_k^\O(\bj)^{-1} e_k^\O e_{k-1}^\O Q_{k-1}^\O(\bi{\cdot}s_k)^{-1} \psi_k^\O f_{\t\t} = e(\bj)^\O P_k^\O(\bj)^{-1} e_k^\O e_{k-1}^\O \psi_k^\O Q_{k-1}^\O(\bi{\cdot}s_k)^{-1} f_{\t\t}\\
& = e(\bj)^\O P_k^\O(\bj)^{-1} e_k^\O e_{k-1}^\O s_k^\O M_k^\O(\bi) Q_{k-1}^\O(\bi{\cdot}s_k)^{-1} f_{\t\t} + e(\bj)^\O P_k^\O(\bj)^{-1} e_k^\O e_{k-1}^\O N_k^\O(\bi) Q_{k-1}^\O(\bi{\cdot}s_k)^{-1} f_{\t\t}.
\end{align*}

Because $\t(k-1) + \t(k) \neq 0$, we have $\epsilon_{k-1}^\O f_{\t\t} = 0$ by~\autoref{semi:1}. Hence, as $N_k^\O(\bi) Q_{k-1}^\O(\bi{\cdot}s_k)^{-1} \in \L_n(\O)$, by~\autoref{semi:y}, we have
$$
e(\bj)^\O P_k^\O(\bj)^{-1} e_k^\O e_{k-1}^\O N_k^\O(\bi) Q_{k-1}^\O(\bi{\cdot}s_k)^{-1} f_{\t\t} = 0,
$$
which yields
\begin{equation} \label{tangle:5:h2:eq3}
e(\bj)^\O \epsilon_k^\O \epsilon_{k-1}^\O \psi_k^\O f_{\t\t} = e(\bj)^\O P_k^\O(\bj)^{-1} e_k^\O e_{k-1}^\O s_k^\O M_k^\O(\bi) Q_{k-1}^\O(\bi{\cdot}s_k)^{-1} f_{\t\t}.
\end{equation}

Because $\t(k-1) + \t(k) \neq 0$ and $\t(k) + \t(k+1) \neq 0$, by~\autoref{com} and (\ref{tangle:eq2}), we have
\begin{align*}
e(\bj)^\O \epsilon_k^\O \psi_{k-1}^\O f_{\t\t} & = e(\bj)^\O P_k^\O(\bj)^{-1} e_k^\O \psi_{k-1}^\O \leftidx{^{s_{k-1}}}Q_k^\O(\bi{\cdot}s_{k-1})^{-1} f_{\t\t}\\
& = e(\bj)^\O P_k^\O(\bj)^{-1} e_k^\O s_{k-1}^\O M_{k-1}^\O(\bi) \leftidx{^{s_{k-1}}}Q_k^\O(\bi{\cdot}s_{k-1})^{-1} f_{\t\t} + e(\bj)^\O P_k^\O(\bj)^{-1} e_k^\O N_{k-1}^\O(\bi) \leftidx{^{s_{k-1}}}Q_k^\O(\bi{\cdot}s_{k-1})^{-1} f_{\t\t}.
\end{align*}

Because $\t(k) + \t(k+1) \neq 0$, we have $e_k^\O f_{\t\t} = 0$ by~\autoref{semi:1}. Hence, by~\autoref{esscom}, we have
$$
e(\bj)^\O P_k^\O(\bj)^{-1} e_k^\O N_{k-1}^\O(\bi) \leftidx{^{s_{k-1}}}Q_k^\O(\bi{\cdot}s_{k-1})^{-1} f_{\t\t} = 0,
$$
which yields
\begin{equation} \label{tangle:5:h2:eq4}
e(\bj)^\O \epsilon_k^\O \psi_{k-1}^\O f_{\t\t} = e(\bj)^\O P_k^\O(\bj)^{-1} e_k^\O s_{k-1}^\O M_{k-1}^\O(\bi) \leftidx{^{s_{k-1}}}Q_k^\O(\bi{\cdot}s_{k-1})^{-1} f_{\t\t}.
\end{equation}

By~\autoref{tangle:h2}, we have $M_k^\O(\bi) Q_{k-1}^\O(\bi{\cdot}s_k)^{-1} f_{\t\t} = M_{k-1}^\O(\bi) \leftidx{^{s_{k-1}}}Q_k^\O(\bi{\cdot}s_{k-1})^{-1} f_{\t\t}$. Hence, because $e_k^\O e_{k-1}^\O s_k^\O = e_k^\O s_{k-1}^\O$, the equality (\ref{tangle:5:h2:eq1}) holds by (\ref{tangle:5:h2:eq3}) and (\ref{tangle:5:h2:eq4}). \endproof

\begin{Lemma} \label{tangle:5}
Suppose $1 < k < n$ and $\bi \in I^n$ with $i_{k-1} + i_{k+1} = 0$, $i_k + i_{k-1} \neq 0$ and $i_k + i_{k+1} \neq 0$. Then we have
$$
e(\bj)^\O \epsilon_k^\O \epsilon_{k-1}^\O \psi_k^\O e(\bi)^\O = e(\bj)^\O \epsilon_k^\O \psi_{k-1}^\O e(\bi)^\O \in \BOx, \qquad \text{and} \qquad
e(\bi)^\O \psi_k^\O \epsilon_{k-1}^\O \epsilon_k^\O e(\bj)^\O = e(\bi)^\O \psi_{k-1}^\O \epsilon_k^\O e(\bj)^\O \in \BOx.
$$
\end{Lemma}

\proof Suppose $\t \in \Tud_n(\bi)$. As $i_k + i_{k-1} \neq 0$ and $i_k + i_{k+1} \neq 0$, we have $\t(k) + \t(k-1) \neq 0$ and $\t(k) + \t(k+1) \neq 0$. Hence, by~\autoref{idem-semi} and~\autoref{tangle:5:h2}, as $\t$ is chosen arbitrary, we have
\begin{align*}
e(\bj)^\O \epsilon_k^\O \epsilon_{k-1}^\O \psi_k^\O e(\bi)^\O & = e(\bj)^\O \epsilon_k^\O \psi_{k-1}^\O e(\bi)^\O \in \BOx,\\
e(\bi)^\O \psi_k^\O \epsilon_{k-1}^\O \epsilon_k^\O e(\bj)^\O & = e(\bi)^\O \psi_{k-1}^\O \epsilon_k^\O e(\bj)^\O \in \BOx,
\end{align*}
which proves the Lemma. \endproof

It left us to consider the case (6.4.3).

\begin{Lemma} \label{tangle:6}
Suppose $1 < k < n$ and $\bi \in I^n$ with $i_{k-1} + i_{k+1} = 0$, and $i_k + i_{k-1} = 0$ or $i_k + i_{k+1} = 0$. Then we have $i_k \neq \pm \frac{1}{2}$. Moreover, when $i_k = 0$, for any $\bj \in I^n$, we have
$$
e(\bj)^\O \epsilon_k^\O \epsilon_{k-1}^\O \psi_k^\O e(\bi)^\O = e(\bj)^\O \epsilon_k^\O \psi_{k-1}^\O e(\bi)^\O = 0, \qquad \text{and} \qquad
e(\bi)^\O \psi_k^\O \epsilon_{k-1}^\O \epsilon_k^\O e(\bj)^\O = e(\bi)^\O \psi_{k-1}^\O \epsilon_k^\O e(\bj)^\O = 0.
$$
\end{Lemma}

\proof Suppose $i_k = \pm \frac{1}{2}$. By the assumption of the Lemma, we have $i_k = i_{k-1}$ or $i_k = i_{k+1}$. As $\bi \in I^n$ and $i_k = \pm \frac{1}{2}$, by~\autoref{deg:h3} we have $-2 \leq h_k(\bi) \leq -1$. If $i_{k-1} = i_k$, then $h_{k-1}(\bi) = h_k(\bi) - 3 \leq -4$, which implies $\bi \not\in I^n$ by~\autoref{deg:h2}. If $i_{k+1} = i_k$, $h_{k+1}(\bi) = h_k(\bi) + 3 \geq 1$, which implies $\bi \not\in I^n$ by~\autoref{deg:h2}. Hence we have $\bi \not\in I^n$, which contradicts to the assumption of the Lemma. Hence we have $i_k \neq \pm \frac{1}{2}$.

If $i_k = 0$, then we have $i_{k-1} = i_k = i_{k+1} = 0$. Hence $e(\bi)^\O \psi_k^\O = e(\bi)^\O \psi_{k-1}^\O = 0$ by~\autoref{esscom:3}, which completes the proof. \endproof

\begin{Lemma} \label{tangle:7}
Suppose $1 < k < n$ and $\bi \in I^n$ with $i_{k-1} = i_k = -i_{k+1} \neq 0, \pm \frac{1}{2}$. Then for any $\bj \in I^n$, we have
$$
e(\bj)^\O \epsilon_k^\O \epsilon_{k-1}^\O \psi_k^\O e(\bi)^\O = e(\bj)^\O \epsilon_k^\O \psi_{k-1}^\O e(\bi)^\O \in \BOx, \qquad \text{and} \qquad
e(\bi)^\O \psi_k^\O \epsilon_{k-1}^\O \epsilon_k^\O e(\bj)^\O = e(\bi)^\O \psi_{k-1}^\O \epsilon_k^\O e(\bj)^\O \in \BOx.
$$
\end{Lemma}

\proof We only prove the first equality, and the second equality holds by using similar argument.

In order to prove $e(\bj)^\O \epsilon_k^\O \epsilon_{k-1}^\O \psi_k^\O e(\bi)^\O = e(\bj)^\O \epsilon_k^\O \psi_{k-1}^\O e(\bi)^\O$, it suffices to prove that for any $\t \in \Tud_n(\bi)$, we have
\begin{equation} \label{tangle:7:eq1}
e(\bj)^\O \epsilon_k^\O \epsilon_{k-1}^\O \psi_k^\O f_{\t\t} = e(\bj)^\O \epsilon_k^\O \psi_{k-1}^\O f_{\t\t}.
\end{equation}

Suppose $\t \in \Tud_n(\bi)$. As $i_{k-1} = i_k \neq 0, \pm \frac{1}{2}$, we have $i_{k-1} + i_k \neq 0$, which implies $\t(k-1) + \t(k) \neq 0$. By the construction of up-down tableaux, we have $\t(k-1) \neq \t(k)$, which yields that we will not have $\t(k-1) + \t(k+1) = \t(k) + \t(k+1) = 0$. Hence, we consider the following three cases.


\textbf{Case 1:} $\t(k-1) + \t(k+1) = 0$.

By the construction of up-down tableau, we have $\t(k-1) \neq \t(k)$, which implies $\t(k) + \t(k+1) \neq 0$. Moreover, as $i_{k-1} = i_k \neq 0$, we have $i_{k-1} + i_k \neq 0$, which implies $\t(k-1) + \t(k) \neq 0$. Hence, by~\autoref{tangle:5:h2}, (\ref{tangle:7:eq1}) holds when $\t(k-1) + \t(k+1) = 0$.

\textbf{Case 2:} $\t(k) + \t(k+1) = 0$.

In this case, we can write $\t = (\alpha_1, \ldots, \alpha_n)$ where $\alpha_k = -\alpha_{k+1}$. Define
$$
\s = (\alpha_1, \ldots, \alpha_{k-1}, -\alpha_{k-1}, \alpha_{k-1}, \alpha_{k+2}, \ldots, \alpha_n).
$$

It is easy to see that $\s$ is an up-down tableau and we have $c_\t(k-1) = - c_\s(k)$. As $\s \in \Tud_n(\bi{\cdot}s_k)$, by~\autoref{A:9} we have $h_k(\bi) = 0$, which implies that $\s$ is the unique up-down tableau in $\Tud_n(\bi{\cdot}s_k)$ such that $\s \overset{k}\sim \t$ by~\autoref{deg:h4:1}. By~\autoref{semi:B},~\autoref{semi:1},~\autoref{semi:2} and~\autoref{PQ:2}, we have
\begin{align}
e(\bj)^\O \epsilon_k^\O \epsilon_{k-1}^\O \psi_k^\O f_{\t\t} & = e(\bj)^\O \epsilon_k^\O \epsilon_{k-1}^\O \f{\s} \psi_k^\O f_{\t\t} = e(\bj)^\O P_k^\O(\bj)^{-1} e_k^\O e_{k-1}^\O \f{\s} e_k^\O \frac{Q_k^\O(\bi)^{-1}}{c_\t(k) + c_\s(k)} f_{\t\t} \notag\\
& = e(\bj)^\O P_k^\O(\bj)^{-1} e_k^\O e_{k-1}^\O e_k^\O \frac{Q_k^\O(\bi)^{-1}}{c_\t(k) - c_\t(k-1)} f_{\t\t} = \frac{1}{c_\t(k) - c_\t(k-1)} e(\bj) \epsilon_k^\O f_{\t\t}, \label{tangle:7:eq2}
\end{align}
and by~\autoref{semi:1} and~\autoref{semi:2}, we have
\begin{equation} \label{tangle:7:eq3}
e(\bj)^\O \epsilon_k^\O \psi_{k-1}^\O f_{\t\t} = e(\bj)^\O \epsilon_k^\O \f{\t} \psi_{k-1}^\O f_{\t\t} = \frac{1}{c_\t(k) - c_\t(k-1)} e(\bj) \epsilon_k^\O f_{\t\t}.
\end{equation}

Hence by (\ref{tangle:7:eq2}) and (\ref{tangle:7:eq3}), (\ref{tangle:7:eq1}) holds when $\t(k) + \t(k+1) = 0$.

\textbf{Case 3:} $\t(k-1) + \t(k+1) \neq 0$ and $\t(k) + \t(k+1) \neq 0$.

By~\autoref{tangle:5:h1}, (\ref{tangle:7:eq1}) holds when $\t(k-1) + \t(k+1) \neq 0$ and $\t(k) + \t(k+1) \neq 0$. \endproof

The next Lemma can be proved using the same argument as~\autoref{tangle:7}.

\begin{Lemma} \label{tangle:8}
Suppose $1 < k < n$ and $\bi \in I^n$ with $i_{k-1} = -i_k = -i_{k+1} \neq 0, \pm \frac{1}{2}$. Then for any $\bj \in I^n$, we have
$$
e(\bj)^\O \epsilon_k^\O \epsilon_{k-1}^\O \psi_k^\O e(\bi)^\O = e(\bj)^\O \epsilon_k^\O \psi_{k-1}^\O e(\bi)^\O \in \BOx, \qquad \text{and} \qquad
e(\bi)^\O \psi_k^\O \epsilon_{k-1}^\O \epsilon_k^\O e(\bj)^\O = e(\bi)^\O \psi_{k-1}^\O \epsilon_k^\O e(\bj)^\O \in \BOx.
$$
\end{Lemma}

\begin{Lemma} \label{tangle:O}
In $\BOx$, the tangle relations hold.
\end{Lemma}

\proof The relation (\ref{rela:7:2}) holds by~\autoref{tangle:2} and~\autoref{tangle:3}; and the relation (\ref{rela:7:1}) holds by~\autoref{tangle:4},

\noindent\autoref{tangle:5},~\autoref{tangle:6},~\autoref{tangle:7} and~\autoref{tangle:8}. \endproof

\begin{Corollary} \label{tangle:other}
Suppose $1 \leq k < n$, we have $\epsilon_k^\O \epsilon_{k+1}^\O \psi_k^\O = \epsilon_k^\O \psi_{k+1}^\O$, and $\psi_k^\O \epsilon_{k+1}^\O \epsilon_k^\O = \psi_{k+1}^\O \epsilon_k^\O$.
\end{Corollary}

\proof By~\autoref{tangle:O}, we have $\epsilon_k^\O \epsilon_{k+1}^\O \psi_k^\O = \epsilon_k^\O \epsilon_{k+1}^\O \epsilon_k^\O \psi_{k+1}^\O = \epsilon_k^\O \psi_{k+1}^\O$, and $\psi_k^\O \epsilon_{k+1}^\O \epsilon_k^\O = \psi_{k+1}^\O \epsilon_k^\O \epsilon_{k+1}^\O \epsilon_k^\O = \psi_{k+1}^\O \epsilon_k^\O$, which completes the proof. \endproof

The next Proposition follows by~\autoref{tangle:O}.

\begin{Proposition} \label{tangle}
In $\B$, the tangle relations hold.
\end{Proposition}

\subsection{Braid relations} \label{sec:rela:4}

In this subsection, we prove the braid relations hold in $\B$. The braid relations are determined by the values of $i_{k-1}, i_k$ and $i_{k+1}$. Hence, we separate the question into several cases. In more details, we consider the following cases:

(6.5.1). When $i_{k-1} + i_k \neq 0$, $i_{k-1} + i_{k+1} \neq 0$ and $i_k + i_{k+1} \neq 0$.

(6.5.2). When $i_{k-1} + i_k = 0$ and $i_{k+1} \neq \pm i_{k-1}$, or $i_{k-1} + i_{k+1} = 0$ and $i_k \neq \pm i_{k-1}$, or $i_k + i_{k+1} = 0$ and $i_{k-1} \neq \pm i_k$.

(6.5.3). When $i_{k-1} = i_k = -i_{k+1}$, or $i_{k-1} = -i_k = i_{k+1}$, or $-i_{k-1} = i_k = i_{k+1}$.

First we consider the case (6.5.1).

\begin{Lemma} \label{braid:1}
Suppose $1 < k < n$ and $\bi \in I^n$ satisfies (6.5.1). Then we have
$$
e(\bi) \psi_k \psi_{k-1} \psi_k =
\begin{cases}
e(\bi) \psi_{k-1} \psi_k \psi_{k-1} - e(\bi), & \text{if $i_{k-1} = i_{k+1} = i_k - 1$,}\\
e(\bi) \psi_{k-1} \psi_k \psi_{k-1} + e(\bi), & \text{if $i_{k-1} = i_{k+1} = i_k + 1$,}\\
e(\bi) \psi_{k-1} \psi_k \psi_{k-1}, & \text{otherwise.}
\end{cases}
$$
\end{Lemma}

\proof Suppose $\t \in \mathscr T^{ud}_n(\bi)$. Because $i_{k-1} + i_k \neq 0$, $i_{k-1} + i_{k+1} \neq 0$ and $i_k + i_{k+1} \neq 0$, we have $c_\t(k-1) + c_\t(k) \neq 0$, $c_\t(k-1) + c_\t(k+1) \neq 0$ and $c_\t(k) + c_\t(k+1) \neq 0$. By~\autoref{KLR:1} the Lemma holds. \endproof

Next we consider the case (6.5.2). Note that when $i_k + i_{k+1} = 0$, $i_{k-1} \neq \pm i_k$ is equivalent to $i_{k-1} \neq i_k$ and $i_{k-1} \neq i_{k+1}$. We separate this case further. In more details, we consider the following cases:


(6.5.2.1). When $i_k + i_{k+1} = 0$, and $|i_{k-1} - i_k| > 1$ and $|i_{k-1} - i_{k+1}| > 1$.

(6.5.2.2). When $i_{k-1} + i_k = 0$, and $|i_{k+1} - i_{k-1}| > 1$ and $|i_{k+1} - i_k| > 1$.

(6.5.2.3). When $i_{k-1} + i_{k+1} = 0$, and $|i_k - i_{k-1}| > 1$ and $|i_k - i_{k+1}| > 1$.

(6.5.2.4). When $i_k + i_{k+1} = 0$, and $|i_{k-1} - i_k| = 1$, or $|i_{k-1} - i_{k+1}| = 1$, or $|i_{k-1} - i_k| = |i_{k-1} - i_{k+1}| = 1$.

(6.5.2.5). When $i_{k-1} + i_k = 0$, and $|i_{k+1} - i_{k-1}| = 1$, or $|i_{k+1} - i_k| = 1$, or $|i_{k+1} - i_{k-1}| = |i_{k+1} - i_k| = 1$.

(6.5.2.6). When $i_{k-1} + i_{k+1} = 0$, and $|i_k - i_{k-1}| = 1$, or $|i_k - i_{k+1}| = 1$, or $|i_k - i_{k-1}| = |i_k - i_{k+1}| = 1$.

It is easy to see that if $\bi \in I^n$ satisfies (6.5.2), then $\bi$ satisfies one of (6.5.2.1) - (6.5.2.6). In more details, if $\bi$ satisfy (6.5.2), then we have $r,l \in \{k-1, k, k+1\}$ such that $i_r + i_l = 0$. Let $m \in \{k-1, k, k+1\}$ and $m \neq r, l$. If $|i_m - i_r| = 0$, then $i_m = i_r$, which contradicts to the fact $\bi$ satisfies (6.5.2). Similarly, if $|i_m - i_l| = 0$, it leads to contradiction. Hence, we have $|i_m - i_r| \geq 1$ and $|i_m - i_l| \geq 1$. This shows that $\bi$ satisfies one of (6.5.2.1) - (6.5.2.6).

First we prove (6.5.2.1) - (6.5.2.3). The following Lemmas prove (6.5.2.1).


%
%


\begin{Lemma} \label{braid:2:1:1}
Suppose $1 < k < n$ and $\bi \in I^n$ satisfying (6.5.2.1). Then we have $e(\bi)^\O \psi_k^\O \psi_{k-1}^\O \psi_k^\O = e(\bi)^\O \psi_{k-1}^\O \psi_k^\O \psi_{k-1}^\O$ if $i_k = i_{k+1} = 0$.
\end{Lemma}

\proof Let $\bj = (j_1, \ldots, j_n) = \bi{\cdot}s_k s_{k-1} s_k$. By~\autoref{com}, we have
\begin{equation} \label{braid:2:1:1:eq1}
e(\bi)^\O \psi_k^\O \psi_{k-1}^\O \psi_k^\O = e(\bi)^\O \psi_k^\O \psi_{k-1}^\O \psi_k^\O e(\bj)^\O, \qquad \text{and} \qquad e(\bi)^\O \psi_{k-1}^\O \psi_k^\O \psi_{k-1}^\O = e(\bi)^\O \psi_{k-1}^\O \psi_k^\O \psi_{k-1}^\O e(\bj)^\O.
\end{equation}

When $i_k = i_{k+1} = 0$, we have $j_{k-1} = j_k = 0$. Then by~\autoref{esscom:3}, we have $e(\bi)^\O \psi_k^\O = 0 = \psi_{k-1}^\O e(\bj)^\O$. Therefore, by (\ref{braid:2:1:1:eq1}), we have $e(\bi)^\O \psi_k^\O \psi_{k-1}^\O \psi_k^\O = 0 = e(\bi)^\O \psi_{k-1}^\O \psi_k^\O \psi_{k-1}^\O$, which proves the Lemma. \endproof

\begin{Lemma} \label{braid:2:1:2}
Suppose $1 < k < n$ and $\bi \in I^n$ satisfying (6.5.2.1). For any $\t \in \Tud_n(\bi)$, we have
$$
f_{\t\t} \psi_k^\O \psi_{k-1}^\O \psi_k^\O = f_{\t\t} \psi_{k-1}^\O \psi_k^\O \psi_{k-1}^\O.
$$
\end{Lemma}

\proof Because $|i_{k-1} - i_k| > 1$, $|i_{k-1} - i_{k+1}| > 1$ and $i_k + i_{k+1} = 0$, we have $i_{k-1} + i_k \neq 0$ and $i_{k-1} + i_{k+1} \neq 0$. Therefore, we have $\t(k-1) + \t(k) \neq 0$ and $\t(k-1) + \t(k+1) \neq 0$. It is easy to see that when $\t(k) + \t(k+1) \neq 0$, the Lemma holds by~\autoref{KLR:1}. In the rest of the proof, we consider $\t \in \Tud_n(\bi)$ with $\t(k) + \t(k+1) = 0$.

Set $\bj = \bi{\cdot}s_k s_{k-1} s_k = \bi{\cdot}s_{k-1} s_k s_{k-1}$. By~\autoref{com}, we have
\begin{align}
f_{\t\t} \psi_k^\O \psi_{k-1}^\O \psi_k^\O & = e(\bi)^\O \psi_k^\O e(\bi{\cdot}s_k)^\O \psi_{k-1}^\O e(\bi{\cdot}s_k s_{k-1})^\O \psi_k^\O e(\bi{\cdot}s_k s_{k-1} s_k)^\O, \label{braid:2:1:2:eq1} \\
\qquad \text{and} \qquad f_{\t\t} \psi_{k-1}^\O \psi_k^\O \psi_{k-1}^\O & = e(\bi)^\O \psi_{k-1}^\O e(\bi{\cdot}s_{k-1})^\O \psi_k^\O e(\bi{\cdot}s_{k-1} s_k)^\O \psi_{k-1}^\O e(\bi{\cdot}s_{k-1}s_ks_{k-1})^\O. \label{braid:2:1:2:eq2}
\end{align}

If $\bj \not\in I^n$, (\ref{braid:2:1:2:eq1}) and (\ref{braid:2:1:2:eq2}) both equal to $0$ by~\autoref{idem-semi}. Hence we have
$$
f_{\t\t} \psi_k^\O \psi_{k-1}^\O \psi_k^\O = 0 = f_{\t\t} \psi_{k-1}^\O \psi_k^\O \psi_{k-1}^\O,
$$
where the Lemma holds. Hence, in the rest of the proof, we assume that $\bj \in I^n$.

When $\bj \in I^n$, as $j_{k-1} = i_{k+1} = -i_k$ and $|i_k - i_{k-1}| > 1$, by the definition of $h_k$, we have $h_{k-1}(\bj) = - h_k(\bi)$. Because $\bi, \bj \in I^n$, by~\autoref{deg:h2}, it forces $h_{k-1}(\bj) = h_k(\bi) = 0$. Then by~\autoref{A:9} and~\autoref{A:8}, we have $\bi{\cdot}s_k, \bi{\cdot}s_k s_{k-1}, \bi{\cdot}s_{k-1}, \bi{\cdot}s_{k-1} s_k \in I^n$.

Recall we have $\t(k) + \t(k+1) = 0$. Because we have $h_k(\bi) = 0$, by~\autoref{deg:h4:1}, there exists an unique $\s \in \Tud_n(\bi{\cdot}s_k)$ such that $\t \overset{k}\sim \s$. As $i_k = -i_{k+1} \neq 0$, we have $\s \neq \t$.

Because $|i_{k-1} - i_k| > 1$, $|i_{k-1} - i_{k+1}| > 1$ and $i_k + i_{k+1} = 0$, we have $i_{k-1} + i_k \neq 0$ and $i_{k-1} + i_{k+1} \neq 0$. Hence, by~\autoref{A:8}, there exist $\v = \s{\cdot}s_{k-1} s_k \in \Tud_n(\bj)$. As $h_{k-1}(\bj) = 0$, by~\autoref{deg:h4:1}, there exists a unique $\u \in \Tud_n(\bj{\cdot}s_{k-1})$ such that $\u \overset{k-1}\sim \v$. Because of the uniqueness of $\u$, we have $\u = \t{\cdot}s_{k-1}s_k \in \Tud_n(\bi{\cdot}s_{k-1}s_k)$.

Hence, by~\autoref{semi:1} and~\autoref{semi:2}, we have
\begin{align}
f_{\t\t} \psi_k^\O \psi_{k-1}^\O \psi_k^\O & = f_{\t\t} \psi_k^\O \f{\s} \psi_{k-1}^\O \psi_k^\O \f{\v} = \frac{1}{c_\t(k) + c_\s(k)} f_{\t\t} \epsilon_k^\O \f{\s} \psi_{k-1}^\O \psi_k^\O \f{\v} \notag \\
& = \frac{1}{c_\t(k) + c_\s(k)} f_{\t\t} \epsilon_k^\O \psi_{k-1}^\O \psi_k^\O \f{\v}; \label{braid:2:1:2:eq3}\\
f_{\t\t} \psi_{k-1}^\O \psi_k^\O \psi_{k-1}^\O & = f_{\t\t} \psi_{k-1}^\O \psi_k^\O \f{\u} \psi_{k-1}^\O \f{\v} = \frac{1}{c_\u(k-1) + c_\v(k-1)} f_{\t\t} \psi_{k-1}^\O \psi_k^\O \f{\u} \epsilon_{k-1}^\O \f{\v} \notag \\
& = \frac{1}{c_\u(k-1) + c_\v(k-1)} f_{\t\t} \psi_{k-1}^\O \psi_k^\O \epsilon_{k-1}^\O \f{\v}. \label{braid:2:1:2:eq4}
\end{align}

By~\autoref{tangle:O},~\autoref{tangle:other} and~\autoref{KLR:1}, we have
\begin{align*}
f_{\t\t} \epsilon_k^\O \psi_{k-1}^\O \psi_k^\O \f{\v} & = f_{\t\t} \epsilon_k^\O \epsilon_{k-1}^\O (\psi_k^\O)^2 \f{\v} = f_{\t\t} \epsilon_k^\O \epsilon_{k-1}^\O \f{\v}, \\
f_{\t\t} \psi_{k-1}^\O \psi_k^\O \epsilon_{k-1}^\O \f{\v} & = f_{\t\t} (\psi_{k-1}^\O)^2 \epsilon_k^\O \epsilon_{k-1}^\O \f{\v} = f_{\t\t} \epsilon_k^\O \epsilon_{k-1}^\O \f{\v},
\end{align*}
which yields
\begin{equation} \label{braid:2:1:2:eq5}
f_{\t\t} \psi_k^\O \psi_{k-1}^\O \psi_k^\O = \frac{c_\u(k-1) + c_\v(k-1)}{c_\t(k) + c_\s(k)} f_{\t\t} \psi_{k-1}^\O \psi_k^\O \psi_{k-1}^\O,
\end{equation}
by (\ref{braid:2:1:2:eq3}) and (\ref{braid:2:1:2:eq4}). Because $\u = \t{\cdot}s_{k-1}s_k$, we have $c_\t(k) = c_\u(k-1)$; and because $\v = \s{\cdot}s_{k-1} s_k$, we have $c_\s(k) = c_\v(k-1)$. Therefore, by (\ref{braid:2:1:2:eq5}), we complete the proof. \endproof

\begin{Lemma} \label{braid:2:1}
Suppose $1 < k < n$ and $\bi \in I^n$ satisfying (6.5.2.1). Then we have $e(\bi)^\O \psi_k^\O \psi_{k-1}^\O \psi_k^\O = e(\bi)^\O \psi_{k-1}^\O \psi_k^\O \psi_{k-1}^\O$.
\end{Lemma}

\proof Because $i_k + i_{k+1} = 0$, we have $i_k = -i_{k+1}$. By~\autoref{braid:2:1:1}, the Lemma holds if $i_k = -i_{k+1} = 0$. When $i_k = -i_{k+1} \neq 0$, choose arbitrary $\t \in \Tud_n(\bi)$, we have $f_{\t\t} \psi_k^\O \psi_{k-1}^\O \psi_k^\O = f_{\t\t} \psi_{k-1}^\O \psi_k^\O \psi_{k-1}^\O$ by~\autoref{braid:2:1:2}. As $\t$ is chosen arbitrary, by~\autoref{idem-semi}, the Lemma follows. \endproof

The following Lemmas prove (6.5.2.2) and (6.5.2.3).

\begin{Lemma} \label{braid:2:2}
Suppose $1 < k < n$ and $\bi \in I^n$ satisfying (6.5.2.3). Then we have $e(\bi)^\O \psi_k^\O \psi_{k-1}^\O \psi_k^\O = e(\bi)^\O \psi_{k-1}^\O \psi_k^\O \psi_{k-1}^\O$.
\end{Lemma}

\begin{proof}
Set $\bj = \bi{\cdot}s_k s_{k-1} s_k$, $\bl = \bi{\cdot} s_{k-1}$ and $\bm = \bj{\cdot}s_k$. Then by~\autoref{com} and~\autoref{braid:2:1}, we have
\begin{eqnarray}
&& e(\bl)^\O \psi_{k-1}^\O e(\bi)^\O \psi_k^\O \psi_{k-1}^\O e(\bm)^\O = e(\bl)^\O \psi_{k-1}^\O \psi_k^\O \psi_{k-1}^\O e(\bm)^\O \notag\\
& = & e(\bl)^\O \psi_k^\O \psi_{k-1}^\O \psi_k^\O e(\bm)^\O = e(\bl)^\O \psi_k^\O \psi_{k-1}^\O e(\bj)^\O \psi_k^\O e(\bm)^\O. \label{braid:2:2:eq1}
\end{eqnarray}

Because $|i_k - i_{k-1}| > 1$, $|i_k - i_{k+1}| > 1$ and $i_{k-1} + i_{k+1} = 0$, we have $i_{k-1} + i_k \neq 0$ and $i_k + i_{k+1} \neq 0$. By~\autoref{KLR:1}, we have $(\psi_{k-1}^\O)^2 e(\bi)^\O = e(\bi)^\O$ and $e(\bj)^\O (\psi_k^\O)^2 = e(\bj)^\O$. Hence multiplying $\psi_{k-1}^\O$ from left and $\psi_k^\O$ from right on both sides of (\ref{braid:2:2:eq1}), the Lemma follows.
\end{proof}

The next Lemma can be proved using the same method as~\autoref{braid:2:2}.

\begin{Lemma} \label{braid:2:3}
Suppose $1 < k < n$ and $\bi \in I^n$ satisfying (6.5.2.2). Then we have $e(\bi)^\O \psi_k^\O \psi_{k-1}^\O \psi_k^\O = e(\bi)^\O \psi_{k-1}^\O \psi_k^\O \psi_{k-1}^\O$.
\end{Lemma}

Next we prove (6.5.2.4) - (6.5.2.6). We start by considering some special cases.

\begin{Lemma} \label{braid:3:1}
Suppose $1 < k < n$ and $\bi \in I^n$ satisfying (6.5.2.4). If $i_{k-1} = 0$, or $i_k = i_{k+1} = 0$, or $i_k = -i_{k+1} = \pm \frac{1}{2}$, we have
\begin{align*}
e(\bi) \psi_k \psi_{k-1} \psi_k & = e(\bi) \psi_{k-1} \psi_k \psi_{k-1} = e(\bi) \epsilon_k \epsilon_{k-1} e(\bi{\cdot}s_k s_{k-1} s_k) = 0,\\
\psi_k \psi_{k-1} \psi_k e(\bi) & = \psi_{k-1} \psi_k \psi_{k-1} e(\bi) = e(\bi{\cdot}s_k s_{k-1} s_k) \epsilon_{k-1} \epsilon_k e(\bi) = 0.
\end{align*}
\end{Lemma}

\proof Set $\bj = \bi{\cdot}s_k s_{k-1} s_k.$ First we show that under the assumption of the Lemma, for any $\t \in \Tud_n(\bi)$ and $\s \in \Tud_n(\bj)$, we have $\t|_{k-2} \neq \s|_{k-2}$.

Suppose $i_{k-1} = 0$. By the assumption, we have $i_k = -i_{k+1} = \pm 1$. For any $\t \in \Tud_n(\bi)$, as $i_{k-1} = 0$ and $i_k = -i_{k+1} = \pm 1$, by the construction  of up-down tableaux, we have $\t(k-1), \t(k) > 0$ or $\t(k-1), \t(k) < 0$. If $\t(k-1), \t(k) > 0$ and $i_k = 1$, let $\lambda = \t_{k-2}$ and $(i,j) = \t(k-1)$. Then we have $\t(k) = (i,j+1)$, which implies $(i-1,j) \in [\lambda]$ and $(i-1,j) \not\in \mathscr R(\lambda)$, and $(i+1,j) \not\in [\lambda]$ and $(i+1,j) \not\in \mathscr A(\lambda)$. Therefore, we have $\mathscr{AR}_\lambda(-1) = \emptyset$. Hence, for any $\s \in \Tud_n(\bj)$, we have $\t|_{k-2} \neq \s|_{k-2}$. Following the similar argument, if $\t(k-1), \t(k) < 0$, for any $\s \in \Tud_n(\bj)$, we have $\t|_{k-2} \neq \s|_{k-2}$.

Suppose $i_k = i_{k+1} = 0$. By the assumption, we have $i_{k-1} = \pm 1$. For any $\t \in \Tud_n(\bi)$, as $i_{k-1} = \pm 1$ and $i_k = 0$, by the construction of up-down tableaux, we have $\t(k-1), \t(k) > 0$ or $\t(k-1), \t(k) < 0$. If $\t(k-1), \t(k) > 0$ and $i_{k-1} = 1$, let $\lambda = \t_{k-2}$ and $(i,j) = \t(k-1)$. Then we have $\t(k) = (i+1,j)$, which implies $(i,j-1) \in [\lambda]$ and $(i,j-1) \not\in \mathscr R(\lambda)$, and $(i+1,j) \not\in \mathscr A(\lambda)$. Therefore, we have $\mathscr{AR}_\lambda(0) = \emptyset$. Hence, for any $\s \in \Tud_n(\bj)$, we have $\t|_{k-2} \neq \s|_{k-2}$. Following the similar argument, if $\t(k-1), \t(k) < 0$, for any $\s \in \Tud_n(\bj)$, we have $\t|_{k-2} \neq \s|_{k-2}$.

Suppose $i_k = -i_{k+1} = \pm \frac{1}{2}$. Because $i_{k-1} \neq \pm i_k$ in (6.5.2), we have $i_{k-1} = \pm \frac{3}{2}$. Suppose $\t \in \Tud_n(\bi)$. Following the similar argument as above, by the construction of up-down tableaux, for any $\s \in \Tud_n(\bj)$, we have $\t|_{k-2} \neq \s|_{k-2}$.

Therefore, if $\bi \in I^n$ satisfying the assumptions of the Lemma, for any $\t \in \Tud_n(\bi)$ and $\s \in \Tud_n(\bj)$, we have $\t|_{n-2} \neq \s|_{n-2}$. Hence, we have
\begin{align*}
e(\bi)^\O \psi_k^\O \psi_{k-1}^\O \psi_k^\O & = e(\bi)^\O \psi_{k-1}^\O \psi_k^\O \psi_{k-1}^\O = e(\bi)^\O \epsilon_k^\O \epsilon_{k-1}^\O e(\bi{\cdot}s_k s_{k-1} s_k)^\O = 0 \in \BOx,\\
\psi_k^\O \psi_{k-1}^\O \psi_k e(\bi)^\O & = \psi_{k-1}^\O \psi_k^\O \psi_{k-1}^\O e(\bi)^\O = e(\bi{\cdot}s_k s_{k-1} s_k)^\O \epsilon_{k-1}^\O \epsilon_k^\O e(\bi)^\O = 0 \in \BOx.
\end{align*}

The Lemma follows by lifting the elements into $\B$. \endproof

The next Lemma can be proved by the same argument as~\autoref{braid:3:1}.

\begin{Lemma} \label{braid:3:2}
Suppose $1 < k < n$ and $\bi \in I^n$ satisfying (6.5.2.6). If $i_k = 0$, or $i_{k-1} = i_{k+1} = 0$, or $i_{k-1} = -i_{k+1} = \pm \frac{1}{2}$, we have $e(\bi) \psi_k \psi_{k-1} \psi_k = e(\bi) \psi_{k-1} \psi_k \psi_{k-1} = 0$.
\end{Lemma}


\begin{Remark} \label{braid:remark:1}
Suppose $\bi \in I^n$ satisfying (6.5.2.4) and $\bj = \bi{\cdot}s_k s_{k-1} s_k$. By the definition of $h_k$, as $i_k = -i_{k+1} = -j_{k-1}$, we have
$$
h_k(\bi) = - h_{k-1}(\bj) + \delta_{i_k,-i_{k-1}+1} + \delta_{i_k,-i_{k-1}-1} + 2 \delta_{i_k, i_{k-1}} - \delta_{i_k, i_{k-1} - 1} - \delta_{i_k, i_{k-1} + 1} - 2\delta_{i_k, -i_{k-1}}.
$$

Because $\bi$ satisfies (6.5.2.4), we have $i_k \neq \pm i_{k-1}$, which implies $\delta_{i_k, i_{k-1}} = \delta_{i_k, -i_{k-1}} = 0$; and we have $|i_k - i_{k-1}| = 1$, or $|i_{k+1} - i_{k-1}| = |i_k + i_{k-1}| = 1$, or $|i_k - i_{k-1}| = |i_k + i_{k-1}| = 1$.

Assume $\bi$ is a residue sequence satisfies (6.5.2.4) and does not satisfy the assumption of~\autoref{braid:3:1}. Because $i_k \neq 0$ and $i_{k-1} \neq 0$, we have $|i_k - i_{k-1}| \neq |i_k + i_{k-1}|$.

If $|i_k - i_{k-1}| = 1$, we have $i_k = i_{k-1} \pm 1$, which implies $\delta_{i_k, i_{k-1} - 1} + \delta_{i_k, i_{k-1} + 1} = 1$. It is easy to verify that $\delta_{i_k,-i_{k-1}+1} + \delta_{i_k,-i_{k-1}-1} \neq 0$ only if $i_{k-1} = 0$ or $i_k = 0$. Therefore, by excluding the assumptions of~\autoref{braid:3:1}, we have
$$
h_k(\bi) = - h_{k-1}(\bj) + \delta_{i_k,-i_{k-1}+1} + \delta_{i_k,-i_{k-1}-1} + 2 \delta_{i_k, i_{k-1}} - \delta_{i_k, i_{k-1} - 1} - \delta_{i_k, i_{k-1} + 1} - 2\delta_{i_k, -i_{k-1}} = -h_{k-1}(\bj) - 1,
$$
when $|i_k - i_{k-1}| = 1$.

If $|i_k + i_{k-1}| = 1$, we have $i_k = -i_{k-1} \pm 1$, which implies $\delta_{i_k,-i_{k-1}+1} + \delta_{i_k,-i_{k-1}-1} = 1$. It is easy to verify that $\delta_{i_k, i_{k-1} - 1} + \delta_{i_k, i_{k-1} + 1} \neq 0$ only if $i_{k-1} = 0$ or $i_k = 0$. Therefore, by excluding the assumptions of~\autoref{braid:3:1}, we have
$$
h_k(\bi) = - h_{k-1}(\bj) + \delta_{i_k,-i_{k-1}+1} + \delta_{i_k,-i_{k-1}-1} + 2 \delta_{i_k, i_{k-1}} - \delta_{i_k, i_{k-1} - 1} - \delta_{i_k, i_{k-1} + 1} - 2\delta_{i_k, -i_{k-1}} = -h_{k-1}(\bj) + 1,
$$
when $|i_k + i_{k-1}| = 1$.

In a more concrete form, if $\bi$ is a residue sequence satisfies (6.5.2.4) and does not satisfy the assumption of~\autoref{braid:3:1}, we have $|i_{k-1} - i_k| \neq |i_{k-1} - i_{k+1}|$, and
$$
h_k(\bi) = \begin{cases}
- h_{k-1}(\bj) - 1, & \text{if $|i_{k-1} - i_k| = 1$,}\\
- h_{k-1}(\bj) + 1, & \text{if $|i_{k-1} - i_{k+1}| = 1$.}
\end{cases}
$$

Following the similar argument, if $\bi$ is a residue sequence satisfies (6.5.2.6) and does not satisfy the assumption of~\autoref{braid:3:2}, we have $|i_k - i_{k-1}| \neq |i_k - i_{k+1}|$, and
$$
\begin{cases}
h_{k+1}(\bi) = h_{k-1}(\bj) - 3, & \text{if $|i_k - i_{k+1}| = 1$,}\\
h_{k+1}(\bj) = h_{k-1}(\bi) - 3, & \text{if $|i_k - i_{k-1}| = 1$.}
\end{cases}
$$
\end{Remark}

%

Now we proceed to prove all the other cases.

\begin{Lemma} \label{braid:3:3}
Suppose $1 < k < n$ and $\bi \in I^n$ with $i_k + i_{k+1} = 0$ and $|i_{k-1} - i_{k+1}| = 1$. By excluding the assumptions of~\autoref{braid:3:1}, we have
\begin{align*}
& e(\bi) \psi_k \psi_{k-1} \psi_k = e(\bi) \psi_{k-1} \psi_k \psi_{k-1} = e(\bi) \epsilon_k \epsilon_{k-1} e(\bi{\cdot}s_k s_{k-1} s_k) = 0,\\
& \psi_k \psi_{k-1} \psi_k e(\bi) = \psi_{k-1} \psi_k \psi_{k-1} e(\bi) = e(\bi{\cdot}s_k s_{k-1} s_k) \epsilon_{k-1} \epsilon_k e(\bi) = 0.
\end{align*}
\end{Lemma}

\begin{proof}
Let $\bj = \bi{\cdot}s_k s_{k-1} s_k$. Because the assumptions of~\autoref{braid:3:1} are excluded, by~\autoref{braid:remark:1}, we have $h_k(\bi) = -h_{k-1}(\bj) + 1$. Since $\bi \in I^n$, by~\autoref{deg:h2}, we have $h_k(\bi) \leq 0$, which implies $h_{k-1}(\bj) \geq 1$. Therefore, by~\autoref{deg:h2}, we have $\bj \not\in I^n$. Then by~\autoref{idem-semi}, we have $e(\bj) = 0$. The Lemma follows by~\autoref{com}.
\end{proof}

\begin{Lemma} \label{braid:3:4:h1}
Suppose $1 < k < n$ and $\bi \in I^n$ with $i_k + i_{k+1} = 0$ and $|i_{k-1} - i_k| = 1$. By excluding the assumptions of~\autoref{braid:3:1}, if $\bi{\cdot}s_{k-1} s_k s_{k-1} \in I^n$, then there is exactly one of $\bi{\cdot}s_k$ and $\bi{\cdot}s_{k-1} s_k$ in $I^n$.
\end{Lemma}

\proof Let $\bj = \bi{\cdot} s_{k-1} s_k s_{k-1}$. Because the assumptions of~\autoref{braid:3:1} are excluded, by~\autoref{braid:remark:1}, we have $h_k(\bi) = - h_{k-1}(\bj) - 1$. Moreover, as $\bi{\cdot}s_{k-1} s_k = \bj{\cdot}s_{k-1}$, we have $h_k(\bi) + h_k(\bi{\cdot}s_k) = 0$ and $h_{k-1}(\bj) + h_{k-1}(\bi{\cdot}s_{k-1} s_k) = 0$.

If $\bi{\cdot}s_k \in I^n$, as $\bi \in I^n$, we have $h_k(\bi) = h_k(\bi{\cdot}s_k) = 0$. Hence, we have $h_k(\bj) = -1$, which implies that $h_k(\bi{\cdot}s_{k-1}s_k) = 1 > 0$. Therefore, $\bi{\cdot}s_{k-1} s_k \not\in I^n$.

If $\bi{\cdot}s_k \not\in I^n$, by~\autoref{A:9}, we have $h_k(\bi) \neq 0$. By~\autoref{deg:h2}, we have $-2 \leq h_k(\bi) \leq -1$, which implies $0 \leq h_{k-1}(\bj) \leq 1$. As $\bj \in I^n$, it forces $h_{k-1}(\bj) = 0$ by~\autoref{deg:h2}. Hence, by~\autoref{A:9}, we have $\bi{\cdot}s_{k-1} s_k = \bj{\cdot}s_{k-1} \in I^n$.

Hence, we have $\bi{\cdot}s_k \in I^n$ if and only if $\bi{\cdot}s_{k-1} s_k \not\in I^n$, which proves the Lemma. \endproof

By~\autoref{braid:3:4:h1}, it is equivalent to say that under the assumptions of~\autoref{braid:3:4:h1}, we have either $e(\bi) \psi_k \psi_{k-1} \psi_k = 0$ or $e(\bi) \psi_{k-1} \psi_k \psi_{k-1} = 0$; similarly, we have either $\psi_k \psi_{k-1} \psi_k e(\bi) = 0$ or $\psi_{k-1} \psi_k \psi_{k-1} e(\bi) = 0$.

\begin{Lemma} \label{braid:3:4}
Suppose $1 < k < n$ and $\bi \in I^n$ with $i_k + i_{k+1} = 0$ and $|i_{k-1} - i_k| = 1$. By excluding the assumptions of~\autoref{braid:3:1}, we have
$$
e(\bi) \psi_k \psi_{k-1} \psi_k =
\begin{cases}
e(\bi) \psi_{k-1} \psi_k \psi_{k-1} + e(\bi) \epsilon_k \epsilon_{k-1} e(\bi{\cdot}s_k s_{k-1} s_k), & \text{if $i_{k-1} = i_k - 1$,}\\
e(\bi) \psi_{k-1} \psi_k \psi_{k-1} - e(\bi) \epsilon_k \epsilon_{k-1} e(\bi{\cdot}s_k s_{k-1} s_k), & \text{if $i_{k-1} = i_k + 1$.}
\end{cases}
$$
\end{Lemma}

\proof Define $\bj = \bi{\cdot}s_k s_{k-1} s_k$. By~\autoref{com}, we have $e(\bi) \psi_k \psi_{k-1} \psi_k = e(\bi) \psi_k \psi_{k-1} \psi_k e(\bj)$ and $e(\bi) \psi_{k-1} \psi_k \psi_{k-1} = e(\bi) \psi_{k-1} \psi_k \psi_{k-1} e(\bj)$. We assume $\bj \in I^n$, as otherwise, we have $e(\bi) \psi_k \psi_{k-1} \psi_k = e(\bi) \psi_{k-1} \psi_k \psi_{k-1} = e(\bi) \epsilon_k \epsilon_{k-1} e(\bj) = 0$, and the Lemma follows.

By~\autoref{braid:3:4:h1}, there is exactly one of $\bi{\cdot}s_k$ and $\bi{\cdot}s_{k-1}s_k$ in $I^n$. Suppose $\bi{\cdot}s_k \in I^n$. Then we have $e(\bi) \psi_{k-1} \psi_k \psi_{k-1} = 0$. As $\bi{\cdot}s_k \in I^n$, by~\autoref{A:9}, we have $h_k(\bi) = 0$. Because $|i_{k-1} - i_k| = 1$, by~\autoref{braid:remark:1}, we have $h_{k-1}(\bj) = -h_k(\bi)-1 = -1$. As $j_{k-1} = i_{k+1}$ and the assumptions of~\autoref{braid:3:1} are excluded, we have $j_{k-1} \neq 0, \pm \frac{1}{2}$ as $i_k \neq 0,\pm \frac{1}{2}$. Hence, $\bj \in I_{k-1,0}^n$. Therefore, by (\ref{rela:5:1}),~\autoref{tangle:other},~\autoref{inv} and~\autoref{untwist}, we have
\begin{align*}
e(\bi) \psi_k \psi_{k-1} \psi_k & = (-1)^{a_{k-1}(\bj)} e(\bi) \psi_k \psi_{k-1} \psi_k e(\bj) \epsilon_{k-1} e(\bj) = (-1)^{a_{k-1}(\bj)} e(\bi) \psi_k \psi_{k-1}^2 \epsilon_k e(\bj) \epsilon_{k-1} e(\bj)\\
& = (-1)^{a_{k-1}(\bj)} e(\bi) \psi_k \epsilon_k e(\bj) \epsilon_{k-1} e(\bj) = (-1)^{a_{k-1}(\bj) + a_k(\bi{\cdot}s_k)} e(\bi) \epsilon_k e(\bj) \epsilon_{k-1} e(\bj).
\end{align*}

By direct calculation, we have $(-1)^{a_k(\bi{\cdot}s_k) + a_{k-1}(\bj)} = 1$ when $i_{k-1} = i_k - 1$ and $(-1)^{a_k(\bi{\cdot}s_k) + a_{k-1}(\bj)} = -1$ when $i_{k-1} = i_k + 1$. As $e(\bi) \psi_{k-1} \psi_k \psi_{k-1} = 0$, the Lemma follows when $\bi{\cdot}s_k \in I^n$.

Suppose $\bi{\cdot}s_{k-1}s_k \in I^n$. Then we have $e(\bi) \psi_k \psi_{k-1} \psi_k = 0$. As $\bi{\cdot}s_{k-1}s_k = \bj{\cdot}s_{k-1} \in I^n$, by~\autoref{A:9}, we have $h_{k-1}(\bj) = 0$, which implies $h_k(\bi) = -1$ by~\autoref{braid:remark:1}. Because the assumptions of~\autoref{braid:3:1} are excluded, we have $i_k \neq 0, \pm \frac{1}{2}$. Hence $\bi \in I_{k,0}^n$. Following the similar argument as before, we have
$$
e(\bi) \psi_{k-1} \psi_k \psi_{k-1} = (-1)^{a_k(\bi) + a_{k-1}(\bj{\cdot}s_{k-1})} e(\bi) \epsilon_k \epsilon_{k-1} e(\bj).
$$

By direct calculation, we have $(-1)^{a_k(\bi{\cdot}s_k) + a_{k-1}(\bj)} = -1$ when $i_{k-1} = i_k - 1$ and $(-1)^{a_k(\bi{\cdot}s_k) + a_{k-1}(\bj)} = 1$ when $i_{k-1} = i_k + 1$. As $e(\bi) \psi_k \psi_{k-1} \psi_k = 0$, the Lemma follows when $\bi{\cdot}s_{k-1}s_k \in I^n$. \endproof


The next Lemma follows by almost the same method as~\autoref{braid:3:4}.

\begin{Lemma} \label{braid:3:5}
Suppose $1 < k < n$ and $\bi \in I^n$ with $i_k + i_{k+1} = 0$ and $|i_{k-1} - i_k| = 1$. By excluding the assumptions of~\autoref{braid:3:1}, we have
$$
\psi_k \psi_{k-1} \psi_k e(\bi) =
\begin{cases}
\psi_{k-1} \psi_k \psi_{k-1} e(\bi) + e(\bi{\cdot}s_k s_{k-1} s_k) \epsilon_{k-1} \epsilon_k e(\bi), & \text{if $i_{k-1} = i_k - 1$,}\\
\psi_{k-1} \psi_k \psi_{k-1} e(\bi) - e(\bi{\cdot}s_k s_{k-1} s_k) \epsilon_{k-1} \epsilon_k e(\bi), & \text{if $i_{k-1} = i_k + 1$.}
\end{cases}
$$
\end{Lemma}

By~\autoref{braid:3:1},~\autoref{braid:3:3} and~\autoref{braid:3:4}, (6.5.2.4) has been proved; and by~\autoref{braid:3:1},~\autoref{braid:3:3} and~\autoref{braid:3:5}, (6.5.2.5) has been proved. It left us to prove (6.5.2.6).

\begin{Lemma} \label{braid:3:6}
Suppose $1 < k < n$ and $\bi \in I^n$ satisfying (6.5.2.6). Then we have $e(\bi) \psi_k \psi_{k-1} \psi_k = e(\bi) \psi_{k-1} \psi_k \psi_{k-1} = 0$.
\end{Lemma}

\proof When $\bi$ is under the assumptions of~\autoref{braid:3:2}, the Lemma holds. Hence we consider the cases excluding the assumptions of~\autoref{braid:3:2}.

Let $\bj = \bi{\cdot}s_k s_{k-1} s_k$. If $\bj \not\in I^n$, we have $e(\bj) = 0$ by~\autoref{idem-semi}, and the Lemma follows by~\autoref{com}. Assume $\bj \in I^n$ and write
\begin{align*}
\bi & = (i_1, \ldots, i_{k-2}, i_{k-1}, i_k, i_{k+1}, i_{k+2}, \ldots, i_n),\\
\bj & = (i_1, \ldots, i_{k-2}, i_{k+1}, i_k, i_{k-1}, i_{k+2}, \ldots, i_n).
\end{align*}

As $i_{k-1} + i_{k+1} = 0$, we have $h_{k-1}(\bi) = -h_{k-1}(\bj)$. Because $\bi,\bj \in I^n$, we have $-2 \leq h_{k-1}(\bi) ,h_{k-1}(\bj) \leq 0$ by~\autoref{deg:h2}, which forces $h_{k-1}(\bi) = h_{k-1}(\bj)= 0$. As the assumptions of~\autoref{braid:3:2} are excluded, by~\autoref{braid:remark:1}, we have either $h_{k+1}(\bi) = -3$ or $h_{k+1}(\bj) = -3$. Hence, by~\autoref{deg:h2}, we have $\bi \not\in I^n$ or $\bj \not\in I^n$, which leads to contradiction. Hence, we always have $\bj \not\in I^n$, and the Lemma follows. \endproof

Combining~\autoref{braid:2:1} - \ref{braid:3:6}, we have the following Lemma.

\begin{Lemma} \label{braid:3}
Suppose $1 < k < n$ and $\bi \in I^n$ satisfies (6.5.2). Then we have
\begin{numcases}{e(\bi) \mathcal B_k =}
  e(\bi) \epsilon_k \epsilon_{k-1} e(\bi{\cdot}s_k s_{k-1} s_k), &if $i_k + i_{k+1} = 0$ and $i_{k-1} = \pm (i_k - 1)$,\label{braid:3:r1}\\
- e(\bi) \epsilon_k \epsilon_{k-1} e(\bi{\cdot}s_k s_{k-1} s_k), &if $i_k + i_{k+1} = 0$ and $i_{k-1} = \pm (i_k + 1)$,\label{braid:3:r2}\\
  e(\bi) \epsilon_{k-1} \epsilon_k e(\bi{\cdot}s_k s_{k-1} s_k), &if $i_{k-1} + i_k = 0$ and $i_{k+1} = \pm (i_k - 1)$,\label{braid:3:r3}\\
- e(\bi) \epsilon_{k-1} \epsilon_k e(\bi{\cdot}s_k s_{k-1} s_k), &if $i_{k-1} + i_k = 0$ and $i_{k+1} = \pm (i_k + 1)$,\label{braid:3:r4}\\
0, &otherwise,\label{braid:3:r5}
\end{numcases}
where $\mathcal B_k = \psi_k \psi_{k-1} \psi_k - \psi_{k-1} \psi_k \psi_{k-1}$.
\end{Lemma}

\proof Suppose $i_k + i_{k+1} = 0$, and $i_{k-1} = \pm (i_k - 1)$ or $i_{k-1} = \pm (i_k + 1)$. Then $\bi$ satisfy (6.5.2.4). Hence, by~\autoref{braid:3:1},~\autoref{braid:3:3} and~\autoref{braid:3:4}, (\ref{braid:3:r1}) and (\ref{braid:3:r2}) hold.

Suppose $i_{k-1} + i_k = 0$, and $i_{k+1} = \pm (i_k - 1)$ or $i_{k+1} = \pm (i_k + 1)$. Then $\bi$ satisfy (6.5.2.5). Hence, by~\autoref{braid:3:1},

\noindent\autoref{braid:3:3} and~\autoref{braid:3:5}, (\ref{braid:3:r3}) and (\ref{braid:3:r4}) hold.

For the rest of the cases, when $\bi$ satisfies (6.5.2.1) - (6.5.2.3), by~\autoref{braid:2:1} - \ref{braid:2:3}, we have
$$
e(\bi)^\O \psi_k^\O \psi_{k-1}^\O \psi_k^\O = e(\bi)^\O \psi_{k-1}^\O \psi_k^\O \psi_{k-1}^\O,
$$
which proves $e(\bi) \mathcal B_k = 0$ by lifting the elements into $\B$; and when $\bi$ satisfies (6.5.2.6), by~\autoref{braid:3:6}, we have $e(\bi) \mathcal B_k = 0$. \endproof

Finally, we prove (6.5.3). First we consider some special cases.

\begin{Lemma} \label{braid:4:1}
Suppose $1 < k < n$ and $\bi \in I^n$. If $i_{k-1} = -i_k = i_{k+1} = \pm \frac{1}{2}$, we have $e(\bi) \psi_k \psi_{k-1} \psi_k = e(\bi)\psi_{k-1} \psi_k \psi_{k-1} = 0$.
\end{Lemma}

\begin{proof}
By the definition of $h_k$, as $i_{k-1} = i_{k+1}$, we have $h_k(\bi{\cdot}s_k) = h_{k-1}(\bi{\cdot}s_k) - 3$, which implies $\bi{\cdot}s_k \not\in I^n$ by~\autoref{deg:h2}. Hence $e(\bi{\cdot}s_k) = 0$ by~\autoref{idem-semi} and $e(\bi) \psi_k \psi_{k-1} \psi_k = 0$ by~\autoref{com}. Similarly, we have $h_{k+1}(\bi{\cdot}s_{k-1}) = h_k(\bi{\cdot}s_{k-1}) - 3$. Following the same process we have $e(\bi)\psi_{k-1} \psi_k \psi_{k-1} = 0$.
\end{proof}

\begin{Lemma} \label{braid:4:2}
Suppose $1 < k < n$ and $\bi \in I^n$. If $i_{k-1} = i_k = i_{k+1} = 0$, we have $e(\bi) \psi_k \psi_{k-1} \psi_k = e(\bi)\psi_{k-1} \psi_k \psi_{k-1} = 0$.
\end{Lemma}

\begin{proof}
By~\autoref{esscom:3}, we have $e(\bi) \psi_k = e(\bi) \psi_{k-1} = 0$, which proves the Lemma.
\end{proof}

\begin{Lemma} \label{braid:4:3}
Suppose $1 < k < n$ and $\bi \in I^n$ with $i_{k-1} = i_k = -i_{k+1}$. Then we have
\begin{align*}
e(\bi) \psi_k \psi_{k-1} \psi_k & = e(\bi) \psi_{k-1} \psi_k \psi_{k-1} = 0,\\
\psi_k \psi_{k-1} \psi_k e(\bi) & = \psi_{k-1} \psi_k \psi_{k-1} e(\bi) = 0.
\end{align*}
\end{Lemma}

\begin{proof}
By~\autoref{braid:4:2}, when $i_{k-1} = i_k = 0$ the Lemma follows. When $i_{k-1} = i_k \neq 0$, we have $h_k(\bi) = h_{k-1}(\bi) + 2$ if $i_{k-1} \neq \pm \frac{1}{2}$ and $h_k(\bi) = h_{k-1}(\bi) + 3$ if $i_{k-1} = \pm \frac{1}{2}$. As $\bi \in I^n$, by~\autoref{deg:h2}, we have $-2 \leq h_{k-1}(\bi), h_k(\bi) \leq 0$, which forces $h_{k-1}(\bi) = -2$.

Let $\bj = \bi{\cdot}s_k s_{k-1} s_k$. We have $h_{k-1}(\bj) = -h_{k-1}(\bi) \geq 2$, which implies $\bj \not\in I^n$ by~\autoref{deg:h2}, and hence, $e(\bj) = 0$ by~\autoref{idem-semi}. The Lemma holds by~\autoref{com}.
\end{proof}

\autoref{braid:4:3} proves (6.5.3) when $i_{k-1} = i_k = -i_{k+1}$ and $-i_{k-1} = i_k = i_{k+1}$; and~\autoref{braid:4:1} and~\autoref{braid:4:2} prove (6.5.3) when $i_{k-1} = -i_k = i_{k+1} = 0$ or $\pm\frac{1}{2}$. It only left us to consider when $i_{k-1} = -i_k = i_{k+1} \neq \pm 0, \frac{1}{2}$.

\begin{Lemma} \label{braid:4:4:h1}
Suppose $1 < k < n$ and $\bi \in I^n$ with $i_{k-1} = -i_k = i_{k+1} \neq 0, \pm \frac{1}{2}$. Then at most one of $\bi{\cdot}s_k$ and $\bi{\cdot}s_{k-1}$ is in $I^n$.
\end{Lemma}

\begin{proof}
Because $i_{k-1} = -i_k = i_{k+1} \neq 0, \pm \frac{1}{2}$, we have $h_{k-1}(\bi) = h_{k-1}(\bi{\cdot}s_k) = h_k(\bi{\cdot}s_k) -2$. Assume $\bi{\cdot}s_k \in I^n$. Then by~\autoref{deg:h2}, we have $-2 \leq h_{k-1}(\bi), h_k(\bi{\cdot}s_k) \leq 0$, which forces $h_{k-1}(\bi) = -2$. Hence $h_{k-1}(\bi{\cdot}s_{k-1}) = -h_{k-1}(\bi) = 2$. By~\autoref{deg:h2}, we have $\bi{\cdot}s_{k-1} \not\in I^n$, which completes the proof.
\end{proof}

So we will consider 3 cases: $\bi{\cdot}s_k \in I^n$, $\bi{\cdot}s_{k-1} \in I^n$, and both $\bi{\cdot}s_k, \bi{\cdot}s_{k-1} \not\in I^n$. Note that by~\autoref{A:9}, we have $\bi{\cdot}s_k \in I^n$ if and only if $h_k(\bi) = h_k(\bi{\cdot}s_k) = 0$ and $\bi{\cdot}s_{k-1}$ if and only if $h_{k-1}(\bi) = h_{k-1}(\bi{\cdot}s_{k-1}) = 0$.

\begin{Lemma} \label{braid:4:4:1}
Suppose $1 < k < n$ and $\bi \in I^n$ with $i_{k-1} = - i_k = i_{k+1} \neq 0, \pm \frac{1}{2}$. If $h_k(\bi) = 0$, for $\t \in \Tud_n(\bi)$ with $\t(k-1) = \t(k+1)$, we have
$$
f_{\t\t} \psi_k^\O \psi_{k-1}^\O \psi_k^\O e(\bi)^\O = (-1)^{a_{k-1}(\bi) + 1} f_{\t\t} \epsilon_{k-1}^\O e(\bi)^\O.
$$
\end{Lemma}

\proof Because $h_k(\bi) = 0$ and $i_{k-1} = -i_k$, by~\autoref{deg:h4:2} we have $\t(k) = -\t(k-1)$. Hence, as $\t(k-1) = \t(k+1)$, we have $\t(k-1) = -\t(k) = \t(k+1)$.

Write $\t = (\alpha_1, \ldots, \alpha_n)$ and let $\alpha = \t(k-1)$. Hence, we have $\alpha = \gamma_1$ or $\alpha = -\gamma_2$. Let $\beta$ be a general nodes such that $\beta = -\gamma_2$ if $\alpha = \gamma_1$, and $\beta = \gamma_1$ if $\beta = -\gamma_2$. Define an up-down tableau $\s = (\beta_1, \ldots, \beta_n)$ such that $\s \overset{k-1}\sim \t$, and $\beta_{k-1} = \beta$, $\beta_k = -\beta$. Hence, we can write
\begin{align*}
\t & = (\alpha_1, \ldots, \alpha_{k-2}, \alpha, -\alpha, \alpha, \alpha_{k+2}, \ldots, \alpha_n), \\
\s & = (\alpha_1, \ldots, \alpha_{k-2}, \beta, -\beta, \alpha, \alpha_{k+2}, \ldots, \alpha_n).
\end{align*}

Because $\gamma_1$ and $\gamma_2$ are not adjacent, i.e. $-\beta$ and $\alpha$ are not adjacent, by~\autoref{y:h2:8}, $\v = \s{\cdot}s_k$ is an up-down tableau, and we can write
$$
\v = (\alpha_1, \ldots, \alpha_{k-2}, \beta, \alpha, -\beta, \alpha_{k+2}, \ldots, \alpha_n);
$$
and because $-\beta + \alpha \neq 0$, by~\autoref{y:h2:8}, $\u = \v{\cdot}s_{k-1}$ is an up-down tableau, and we can write
$$
\u = (\alpha_1, \ldots, \alpha_{k-2}, \alpha, \beta, -\beta, \alpha_{k+2}, \ldots, \alpha_n).
$$

Notice that $\u \overset{k}\sim \v$ and $\u \in \Tud_n(\bi{\cdot}s_k)$. Because $h_k(\bi) = 0$, by~\autoref{deg:h4:1}, $\u$ is the unique up-down tableau in $\Tud_n(\bi{\cdot}s_k)$ such that $\u \overset{k}\sim \t$. Hence, by~\autoref{semi:1} and~\autoref{semi:2}, we have
\begin{equation} \label{braid:4:4:eq1}
f_{\t\t} \psi_k^\O \psi_{k-1}^\O \psi_k^\O e(\bi)^\O = f_{\t\t} \psi_k^\O \f{\u} \psi_{k-1}^\O \f{\u} \psi_k^\O \f{\t} + f_{\t\t} \psi_k^\O \f{\u} \psi_{k-1}^\O \f{\v} \psi_k^\O \f{\s}.
\end{equation}

First we work with the second term of (\ref{braid:4:4:eq1}). By~\autoref{semi:2},~\autoref{tangle:5:h2} and~\autoref{inv:1}, we have
\begin{align*}
f_{\t\t} \psi_k^\O \f{\u} \psi_{k-1}^\O \f{\v} \psi_k^\O \f{\s} & = \frac{1}{c_\t(k) + c_\u(k)} f_{\t\t} \epsilon_k^\O \f{\u} \psi_{k-1}^\O \f{\v} \psi_k^\O \f{\s}\\
& = \frac{1}{c_\t(k) + c_\u(k)} f_{\t\t} \epsilon_k^\O \epsilon_{k-1}^\O (\psi_k^\O)^2\f{\s} = \frac{1}{c_\t(k) + c_\u(k)} f_{\t\t} \epsilon_k^\O \epsilon_{k-1}^\O \f{\s}.
\end{align*}

Because $i_k \neq 0, \pm \frac{1}{2}$ and $h_k(\bi) = 0$, we have $\bi \in I_{k,+}^n$. Hence, by~\autoref{PQ:3} and~\autoref{a_ki}, we have
\begin{align}
f_{\t\t} \psi_k^\O \f{\u} \psi_{k-1}^\O \f{\v} \psi_k^\O \f{\s} & = \frac{1}{c_\t(k) + c_\u(k)} f_{\t\t} \epsilon_k^\O \f{\t} \epsilon_{k-1}^\O \f{\s} \notag \\
& = (-1)^{a_k(\bi)} \frac{2(c_\t(k) - i_k)}{c_\t(k) + c_\u(k)} f_{\t\t} \epsilon_{k-1}^\O \f{\s} = (-1)^{a_{k-1}(\bi) + 1} \frac{2(c_\t(k) - i_k)}{c_\t(k) + c_\u(k)} f_{\t\t} \epsilon_{k-1}^\O \f{\s}. \label{braid:4:4:eq2}
\end{align}

Then we work with the first term of (\ref{braid:4:4:eq1}). By~\autoref{semi:2} and~\autoref{inv:1}, we have
$$
f_{\t\t} \psi_k^\O \f{\u} \psi_{k-1}^\O \f{\u} \psi_k^\O \f{\t} = \frac{1}{c_\u(k) - c_\u(k-1)} \f{\t} (\psi_k^\O)^2 \f{\t} = \frac{1}{c_\u(k) - c_\u(k-1)} f_{\t\t}.
$$

Because $i_{k-1} \neq 0, \pm \frac{1}{2}$ and $h_{k-1}(\bi) = -2$, we have $\bi \in I_{k,-}^n$. Hence by~\autoref{PQ:3} and~\autoref{semi:2}, we have
\begin{equation} \label{braid:4:4:eq3}
f_{\t\t} \psi_k^\O \f{\u} \psi_{k-1}^\O \f{\u} \psi_k^\O \f{\t} = \frac{1}{c_\u(k) - c_\u(k-1)} f_{\t\t} = (-1)^{a_{k-1}(\bi)+1} \frac{2(c_\t(k-1) - i_{k-1})}{c_\u(k-1) - c_\u(k)} f_{\t\t} \epsilon_{k-1}^\O \f{\t}.
\end{equation}

By the definitions, we have $\alpha > 0$ if $\beta < 0$ and $\alpha < 0$ if $\beta > 0$. Hence, we have $2(c_\t(k) - i_k) = c_\t(k) + c_\u(k)$ and $2(c_\t(k-1) - i_{k-1}) = c_\u(k-1) - c_\u(k)$. Therefore, we have
$$
\frac{2(c_\t(k) - i_k)}{c_\t(k) + c_\u(k)} = \frac{2(c_\t(k-1) - i_{k-1})}{c_\u(k-1) - c_\u(k)} = 1.
$$

Because $h_{k-1}(\bi) = -2$, by~\autoref{deg:h4:1}, $\s$ is the unique up-down tableau in $\Tud_n(\bi)$ such that $\s \overset{k-1}\sim \t$ and $\s \neq \t$. Hence, we have $f_{\t\t} \epsilon_k^\O e(\bi)^\O = f_{\t\t} \epsilon_k^\O \left( \f{\t} + \f{\s} \right)$. Substituting (\ref{braid:4:4:eq2}) and (\ref{braid:4:4:eq3}) into (\ref{braid:4:4:eq1}), the Lemma follows \endproof

Following the similar argument as above, we have the next Lemma.

\begin{Lemma} \label{braid:4:4:2}
Suppose $1 < k < n$ and $\bi \in I^n$ with $i_{k-1} = - i_k = i_{k+1} \neq 0, \pm \frac{1}{2}$. If $h_k(\bi) = 0$, for $\t \in \Tud_n(\bi)$ with $\t(k-1) \neq \t(k+1)$, we have
$$
f_{\t\t} \psi_k^\O \psi_{k-1}^\O \psi_k^\O e(\bi)^\O = (-1)^{a_{k-1}(\bi) + 1} f_{\t\t} \epsilon_{k-1}^\O e(\bi)^\O.
$$
\end{Lemma}

Combining~\autoref{braid:4:4:1} and~\autoref{braid:4:4:2}, we have the next Lemma.

\begin{Lemma} \label{braid:4:4}
Suppose $1 < k < n$ and $\bi \in I^n$ with $i_{k-1} = - i_k = i_{k+1} \neq 0, \pm \frac{1}{2}$. If $h_k(\bi) = 0$, we have
$$
e(\bi) \psi_k \psi_{k-1} \psi_k = e(\bi) \psi_{k-1} \psi_k \psi_{k-1} - (-1)^{a_{k-1}(\bi)} e(\bi) \epsilon_{k-1} e(\bi).
$$
\end{Lemma}

\proof Choose arbitrary $\t \in \Tud_n(\bi)$. By~\autoref{braid:4:4:1} and~\autoref{braid:4:4:2}, we have
$$
f_{\t\t} \psi_k^\O \psi_{k-1}^\O \psi_k^\O e(\bi)^\O = (-1)^{a_{k-1}(\bi) + 1} f_{\t\t} \epsilon_{k-1}^\O e(\bi)^\O,
$$
which implies
$$
e(\bi)^\O \psi_k^\O \psi_{k-1}^\O \psi_k^\O e(\bi)^\O = (-1)^{a_{k-1}(\bi) + 1} e(\bi)^\O \epsilon_{k-1}^\O e(\bi)^\O.
$$

The Lemma follows by lifting the elements from $\BOx$ to $\B$. \endproof

The next Lemma is proved following the similar argument as~\autoref{braid:4:4}.

\begin{Lemma} \label{braid:4:5}
Suppose $1 < k < n$ and $\bi \in I^n$ with $i_{k-1} = - i_k = i_{k+1} \neq 0, \pm \frac{1}{2}$. If $h_{k-1}(\bi) = 0$, we have
$$
e(\bi) \psi_k \psi_{k-1} \psi_k = e(\bi) \psi_{k-1} \psi_k \psi_{k-1} + (-1)^{a_k(\bi)} e(\bi) \epsilon_k e(\bi).
$$
\end{Lemma}

The only left case is that when both of $\bi{\cdot}s_k$ and $\bi{\cdot}s_{k-1}$ are not in $I^n$.

\begin{Lemma} \label{braid:4:6}
Suppose $1 < k < n$ and $\bi \in I^n$ with $i_{k-1} = - i_k = i_{k+1} \neq 0, \pm \frac{1}{2}$. If $h_{k-1}(\bi), h_k(\bi) \neq 0$, we have
$$
e(\bi) \psi_k \psi_{k-1} \psi_k = e(\bi) \psi_{k-1} \psi_k \psi_{k-1} = 0.
$$
\end{Lemma}

\proof As $\bi{\cdot}s_k, \bi{\cdot}s_{k-1} \not\in I^n$, we have $e(\bi) \psi_k = e(\bi) \psi_{k-1} = 0$ by~\autoref{idem-semi} and~\autoref{com}, which implies $e(\bi) \psi_k \psi_{k-1} \psi_k = e(\bi) \psi_{k-1} \psi_k \psi_{k-1} = 0$. \endproof

\begin{Lemma} \label{braid:4}
Suppose $1 < k < n$ and $\bi \in I^n$ satisfies (6.5.3). Then we have
$$
e(\bi) \mathcal B_k =
\begin{cases}
-(-1)^{a_{k-1}(\bi)} e(\bi) \epsilon_{k-1} e(\bi{\cdot}s_k s_{k-1} s_k), & \text{if $i_{k-1} = -i_k = i_{k+1} \neq 0, \pm \frac{1}{2}$ and $h_k(\bi) = 0$},\\
(-1)^{a_{k-1}(\bi)} e(\bi) \epsilon_k e(\bi{\cdot}s_k s_{k-1} s_k), & \text{if $i_{k-1} = -i_k = i_{k+1} \neq 0, \pm \frac{1}{2}$ and $h_{k-1}(\bi) = 0$},\\
0, & \text{otherwise,}
\end{cases}
$$
where $\mathcal B_k = \psi_k \psi_{k-1} \psi_k - \psi_{k-1} \psi_k \psi_{k-1}$.
\end{Lemma}

\proof If $i_{k-1} = -i_k = i_{k+1} \neq 0, \pm \frac{1}{2}$ and $h_k(\bi) = 0$, the Lemma holds by~\autoref{braid:4:4}; and if $i_{k-1} = -i_k = i_{k+1} \neq 0, \pm \frac{1}{2}$ and $h_{k-1}(\bi) = 0$, the Lemma holds by~\autoref{braid:4:5}; and for the rest of the cases, by~\autoref{braid:4:1} - \ref{braid:4:3} and~\autoref{braid:4:6}, the Lemma holds. \endproof

Therefore, by combining~\autoref{braid:1},~\autoref{braid:3} and~\autoref{braid:4}, we have the following Proposition.

\begin{Proposition} \label{braid}
In $\B$, the braid relations hold.
\end{Proposition}

\subsection{The graded cellular basis of $\B$} \label{sec:rela:main}

Now we are ready to prove our main result of this paper.

\begin{Theorem} \label{main:2}
Suppose $R$ is a field with characteristic $0$ and $\delta \in R$. Then $\B \cong \G{n}$.
\end{Theorem}

\proof We can define a map $\G{n} \longrightarrow \B$ by sending $e(\bi)$ to $e(\bi)$, $y_k$ to $y_k$, $\psi_k$ to $\psi_k$ and $\epsilon_k$ to $\epsilon_k$. By~\autoref{idem},
~\autoref{com},~\autoref{esscom},~\autoref{inv},~\autoref{essidem},~\autoref{untwist},~\autoref{tangle} and
~\autoref{braid}, the map is a homomorphism. By~\autoref{B:span}, the map is surjective. By~\autoref{main:1}, we have $\dim \G{n} \leq (2n-1)!!$. As $\dim\B = (2n-1)!!$, it implies the map is an isomorphism. Hence we have $\B \cong \G{n}$. \endproof

\begin{Theorem} \label{main:3}
Suppose $R$ is a field with characteristic $0$ and $\delta \in R$. Then $\B$ is a graded cellular algebra with a graded cellular basis
$$
B = \set{\psi_{\s\t} | (\lambda,f) \in \widehat B_n, \s,\t \in \Tud_n(\lambda)}.
$$
\end{Theorem}

\begin{proof}
By~\autoref{main:1}, the set $B$ spans $\G{n}$. By~\autoref{main:2}, we have $\dim \G{n} = (2n-1)!!$, which makes $B$ to be a basis of $\G{n}$. The cellularity is proved by~\autoref{basis:main}. The elements $\psi_{\s\t}$'s are homogeneous by the construction and we have a degree function $\deg$ on up-down tableaux such that $\deg \psi_{\s\t} = \deg \s + \deg \t$ by~\autoref{deg:tab}. Finally, as we have a $*$-involution on $\G{n}$ such that $\psi_{\s\t}^* = \psi_{\t\s}$, one can see that $B$ forms a graded cellular basis of $\G{n}$. Finally as $\G{n} \cong \B$ by~\autoref{main:2}, we complete the proof.
\end{proof}

The next Corollary is straightforward by~\autoref{main:3}.

\begin{Corollary} \label{idem:residue}
For any $\bi \in P^n$ and $e(\bi) \in \G{n}$, we have $e(\bi) \neq 0$ if and only if $\bi \in I^n$.
\end{Corollary}

\proof Suppose $\bi$ is the residue sequence of an up-down tableau $\t$. By~\autoref{main:3}, we have $\psi_{\t\t} \neq 0$. Because $\psi_{\t\t} = \psi_{\t\t}e(\bi)$, we have $e(\bi) \neq 0$.

Suppose $\bi \not\in I^n$. For any up-down tableau $\t$, we have $\psi_{\s\t} e(\bi) = 0$. Therefore, we have $\G{n} e(\bi) = 0$, which implies $e(\bi) = 0$. \endproof

Suppose $E_n(\delta)$ is the two-sided ideal of $\G{n}$ generated by $\epsilon_1, \ldots, \epsilon_{n-1}$ and $E'_n(\delta)$ is the two-sided ideal of $\B$ generated by $e_1, \ldots, e_{n-1}$. By~\autoref{remark:generators}, we have $e(\bi) \epsilon_k e(\bj) = e(\bi) P_k(\bi)^{-1} e_k Q_k(\bj)^{-1} e(\bj)$, which implies that $\epsilon_k \in E'(\delta)$ for any $1 \leq k \leq n-1$. Similarly, we have $e(\bi) e_k e(\bj) = e(\bi) P_k(\bi) \epsilon_k Q_k(\bj) e(\bj)$, which implies $e_k \in E(\delta)$ for any $1 \leq k \leq n-1$. Therefore, we have $E_n(\delta) \cong E'_n(\delta)$, with isomorphism compatible with $\G{n} \cong \B$. Hence, we have $\G{n} / E_n(\delta) \cong \B/ E'(\delta)$.

Because $E'_n(\delta)$ is the two-sided ideal of $\B$ generated by $e_1, \ldots, e_{n-1}$, by the definition of $\B$, we have $\B/ E'(\delta) \cong R\Sym_n$. Because $R$ is a field with characteristic $0$, by~\autoref{BK:iso}, we have $R\Sym_n \cong \R(R)$ with $\Lambda = \Lambda_k$, for any $k \in \Z$. The next Theorem connects the cyclotomic Khovanov-Lauda-Rouquier algebra and $\G{n}$.

\begin{Theorem}
Suppose $\Lambda = \Lambda_k$ for some $k \in \Z$ and $R$ is a field with characteristic $0$. Then we have $\G{n} / E_n(\delta) \cong \R(R)$.
\end{Theorem}

Finally, we remark that for Brauer algebras $\B$ over fields $R$ with characteristic $p > 0$, or more precisely, for cyclotomic Nazarov-Wenzl algebras $\mathscr W_{r,n}(\mathbf u)$ over arbitrary field $R$, we should be able to extend the idea of this paper and construct a $\Z$-graded algebra similar to $\G{n}$ isomorphic to $\mathscr W_{r,n}(\mathbf u)$. The algebras are generated with elements
$$
G_n(\delta) = \set{e(\bi)|\bi \in P^n} \cup \set{y_k| 1 \leq k \leq n} \cup \set{\psi_k| 1\leq k \leq n-1} \cup \set{\epsilon_k| 1\leq k \leq n-1},
$$
with degrees similar to $\G{n}$. We are also able to construct a set of homogeneous elements
$$
\set{\psi_{\s\t} | (\lambda,f) \in \widehat B_n, \s,\t \in \Tud_n(\lambda)},
$$
which forms a graded cellular basis of $\mathscr W_{r,n}(\mathbf u)$.

\bibliography{papers}

\end{document}